\documentclass[11pt,a4paper,reqno]{amsart}

\usepackage{geometry}
\usepackage[T1]{fontenc}

\usepackage{amssymb,mathtools,mathrsfs}
\usepackage{enumitem,etoolbox,placeins}
\usepackage{xcolor}
\usepackage[hidelinks,linkcolor=black,anchorcolor=black,citecolor=black]{hyperref}
\usepackage{hyperref}
\numberwithin{equation}{section}
\newtheorem{theorem}{Theorem}[section]
\newtheorem{proposition}[theorem]{Proposition}
\newtheorem{lemma}[theorem]{Lemma}
\newtheorem{corollary}[theorem]{Corollary}
\theoremstyle{remark}
\newtheorem{remark}[theorem]{Remark}
\newcommand{\R}{\mathbb R}
\newcommand{\T}{\mathbb T}
\newcommand{\Z}{\mathbb Z}
\newcommand{\av}[1]{\langle #1\rangle}
\newcommand{\dd}{\,\mathrm d}
\newcommand{\ii}{\mathrm i}
\newcommand{\Rea}{\operatorname{Re}}
\newcommand{\diver}{\operatorname{div}}
\newcommand{\cB}{\mathcal B}

\allowdisplaybreaks[4]

\makeatletter
\@namedef{subjclassname@2020}{\textup{2020} Mathematics Subject Classification}
\makeatother

\title[Stability threshold of MHD Couette flow]{Stability threshold of Couette flow for two-dimensional MHD equations on the infinite channel}

\author[Jiahong Wu]{Jiahong Wu}
\address[J. Wu]{Department of Mathematics, University of Notre Dame, Notre Dame, IN, 46556, USA}
\email{jwu29@nd.edu}  

\author[Tianyuan Yu]{Tianyuan Yu}
\address[T. Yu]{School of Science, Dalian Maritime University, Dalian, 116026, People's Republic of China}
\email{tianyuanyupde@163.com, tianyuanyu@dlmu.edu.cn}

\author[Ning Zhu]{Ning Zhu}
\address[N. Zhu]{School of Mathematics, Shandong University, Jinan, Shandong, 250100, People's Republic of China}
\email{ning.zhu@sdu.edu.cn, mathzhu1@163.com}

\date{}
\begin{document}
\begin{abstract}
In this paper, we study the stability threshold of the two-dimensional incompressible MHD equations near Couette flow and a constant horizontal magnetic field \((\beta,0)\) on \(\mathbb T\times\mathbb R\). For \(N>10\), \(0<\nu,\mu\le1\) and \(|\beta|\ge(\nu+\mu)/(3\sqrt{\nu\mu})\), we prove global stability for initial perturbations of size \(\delta_N\min\{\nu,\mu\}^{1/3}\) in \(H^N\), where \(\delta_N>0\) is independent of \(\nu,\mu,\beta\). We also obtain enhanced dissipation of the nonzero horizontal modes on the time scale \(\min\{\nu,\mu\}^{-1/3}\), together with inviscid damping of the vertical components. The proof combines time-dependent Fourier multipliers with an energy correction constructed using the Alfv\'en oscillations. By combining the nonlinear terms arising from this correction, we obtain the cancellation needed to close the energy estimates at the $1/3$ scale.
\end{abstract}
\maketitle

\section{Introduction}
\label{sec:intro}

We consider the following two-dimensional incompressible MHD equations on the infinite channel $\T\times\R$:
\begin{equation}
\begin{cases}
 \partial_tU+U\cdot\nabla U+\nabla P
       =\nu\Delta U+B\cdot\nabla B,\\
 \partial_tB+U\cdot\nabla B
       =\mu\Delta B+B\cdot\nabla U,\\
 \diver U=\diver B=0,
\end{cases}
\qquad (x,y)\in\T\times\R,
\label{eq:MHD}
\end{equation}
where $\T=\R/(2\pi\Z)$. Here $U$, $B$ and $P$ denote the velocity field, the magnetic field and the scalar pressure, respectively. The constants $\nu>0$ and $\mu>0$ are the fluid viscosity and the magnetic resistivity. The goal of this article is to understand the nonlinear stability for  \eqref{eq:MHD} near the equilibrium 
\[
 U_{sh}(x,y)=(y,0),\qquad B_{bg}=(\beta,0),
 \qquad \beta\in\R\setminus\{0\}.
\]
Writing $U=U_{sh}+u$ and $B=B_{bg}+b$, we obtain the perturbation equations
\begin{equation}
\begin{cases}
 (\partial_t+y\partial_x)u-\nu\Delta u
       +u^2e_1-\beta\partial_xb+\nabla p
       =-u\cdot\nabla u+b\cdot\nabla b,\\
 (\partial_t+y\partial_x)b-\mu\Delta b
       -b^2e_1-\beta\partial_xu
       =-u\cdot\nabla b+b\cdot\nabla u,\\
 \diver u=\diver b=0,\\
 (u,b)|_{t=0}=(u_{\rm in},b_{\rm in}),
\end{cases}
\label{eq:perturbation}
\end{equation}
where $e_1=(1,0)$.

The MHD equations describe the motion of electrically conducting fluids in the presence of a magnetic field. Classical well-posedness results for the viscous and resistive equations were established by Sermange--Temam~\cite{ST83}. The stability of a prescribed flow requires quantitative estimates on the evolution of its perturbations, particularly when the dissipative coefficients are small. For Couette flow, this question is closely related to the transition threshold problem in hydrodynamic stability.

The study of transition from laminar flow to turbulence goes back to Reynolds~\cite{Rey83}. Although plane Couette flow is linearly stable at every Reynolds number~\cite{Rom73}, transient growth can amplify perturbations and affect the nonlinear dynamics. This mechanism was studied by Kelvin~\cite{Kel87} and Orr~\cite{Orr07}, and its role in transition was developed in \cite{Chap02,RSBH98,TTRD93}. For further background, we refer to \cite{BGM19,SH01}. In the Navier--Stokes equations, Bedrossian--Germain--Masmoudi \cite{BGM17} proposed the following {\it transition threshold problem}: {\it given a norm $\|\cdot\|_{X}$, determine $\alpha=\alpha(X)$ such that}
\begin{align*}
&\|u_{\rm in}\|_{X}\lesssim\nu^\alpha\Longrightarrow \text{ stability},\\
&\|u_{\rm in}\|_{X}\gg\nu^\alpha\Longrightarrow \text{ possible instability}.
\end{align*}
The index $\alpha$ is called the {\it transition threshold}. The smaller $\alpha$, the larger the basin of stability. 

For the two-dimensional Euler equations, Bedrossian--Masmoudi~\cite{BM15} proved nonlinear inviscid damping and asymptotic stability near Couette flow for small Gevrey perturbations. The results of Lin--Zeng~\cite{LZ11} on nearby inviscid structures and Deng--Masmoudi~\cite{DM23} on instability in low Gevrey regularity show the importance of the regularity assumption. Linear inviscid damping for a class of monotone shear flows was established by Wei--Zhang--Zhao~\cite{WZZ18}.

For the two-dimensional Navier--Stokes equations, Bedrossian--Masmoudi--Vicol~\cite{BMV16} proved enhanced dissipation and inviscid damping for Gevrey perturbations whose size is independent of $\nu$. In Sobolev spaces, Bedrossian--Vicol--Wang~\cite{BVW18} obtained stability at size $\nu^{1/2}$ for shear flows near Couette. Masmoudi--Zhao~\cite{MZ20} established stability at this scale in a critical space, and subsequently obtained the $\nu^{1/3}$ threshold for Couette flow in a high Sobolev space~\cite{MZ22}. Wei--Zhang~\cite{WZ23} proved nonlinear enhanced dissipation and inviscid damping at the $\nu^{1/3}$ scale by a Fourier multiplier method. Li--Masmoudi--Zhao~\cite{LMZ25} further studied the dependence of the stability threshold on Gevrey regularity. In a finite channel with no-slip boundary conditions, Chen--Li--Wei--Zhang~\cite{CLWZ20} obtained the $\nu^{1/2}$ threshold for $H^2$ velocity perturbations.

For the three-dimensional Navier--Stokes equations, Bedrossian--Germain--Masmoudi proved stability at size $\nu$ in Gevrey spaces~\cite{BGM20} and at size $\nu^{3/2}$ in Sobolev spaces~\cite{BGM17}. They also studied the dynamics above the $\nu$ scale~\cite{BGM22}. Wei--Zhang~\cite{WZ21} obtained the $\nu$ threshold for $H^2$ perturbations on $\T\times\R\times\T$. Chen--Wei--Zhang~\cite{CWZ24} established the same threshold in a finite channel with no-slip boundary conditions.

Shear mixing and enhanced dissipation play a central role in these results. For example, the nonzero horizontal modes of a solution to
\(
 (\partial_t+y\partial_x)f=\nu\Delta f
\)
satisfy
\[
 \|f_{\ne}(t)\|_{L^2}
 \le e^{-\nu t^3/12}\|f_{{\rm in},\ne}\|_{L^2}.
\]
Thus the dissipation acts on the time scale $\nu^{-1/3}$, which is shorter than the diffusive time scale $\nu^{-1}$. We refer to Constantin--Kiselev--Ryzhik--Zlato\v{s}~\cite{CKRZ08} for the connection between mixing and diffusion. In the nonlinear problem, time-dependent Fourier multipliers are used to control interactions near critical times while preserving the enhanced dissipation. Related methods have also been applied to the Boussinesq equations with vertical dissipation~\cite{DWZ21}.

For the MHD equations, the background magnetic field introduces an additional coupling between the velocity and magnetic perturbations. In three dimensions, Liss~\cite{Lis20} proved Sobolev stability at size $\nu$ for equal viscosity and resistivity and a sufficiently strong magnetic field whose direction satisfies a Diophantine condition. Rao--Zhang--Zi~\cite{RZZ25} treated unequal viscosity and resistivity under this condition on the field direction. For rationally aligned magnetic fields, Wang--Xu--Zhang~\cite{WXZ26} obtained the $\nu$ threshold under equal diffusion and a sufficiently strong background field. Knobel~\cite{Kno3D25} studied nonlinear transient growth for a spanwise background magnetic field. These results show that the direction of the magnetic field affects both the growth of perturbations and the stability threshold.

In two dimensions, Dolce~\cite{Dol24} introduced symmetric variables to control the transient growth of vorticity and current, and obtained stability at size $\nu^{2/3}$ in the regime $\mu^3\lesssim\nu\le\mu$. Stability estimates for other regimes of viscosity and resistivity were developed by Jin--Ren--Wei~\cite{JRW25} and Knobel~\cite{KnoNR25}. Wang--Zhang~\cite{WangZhang26} established stability estimates for general viscosity and magnetic resistivity, including the $\nu^{1/2}$ scale when $\nu=\mu$.
Results for partial dissipation include horizontal magnetic dissipation~\cite{KZ25}, vertical dissipation~\cite{LWZ26} and anisotropic dissipation~\cite{WWX}. For echoes and long-time dynamics in other MHD regimes, we refer to \cite{KnoIdeal,KZ23}.

Let $\rho=\min\{\nu,\mu\}$. The purpose of this paper is to establish global nonlinear stability at the scale $\rho^{1/3}$. More precisely, for $N>10$ and
\[
 0<\nu,\mu\le1,\qquad
 |\beta|\ge\frac{\nu+\mu}{3\sqrt{\nu\mu}},
\]
we prove that initial perturbations of size $\delta_N\rho^{1/3}$ in $H^N$ remain globally small in a Sobolev norm following the Couette transport. The constants are uniform in this parameter range. The nonzero horizontal modes of the velocity and magnetic field decay exponentially on the time scale $\rho^{-1/3}$, and their vertical components have an additional $(1+t)^{-1}$ decay. In particular, when $\nu=\mu$, we obtain the same $1/3$ stability exponent as for the two-dimensional Navier--Stokes equations.

We now state our main result. For a function $f=f(x,y)$, let
\[
 f_0(y)=\frac1{2\pi}\int_{\T}f(x,y)\dd x,
 \qquad f_{\ne}=f-f_0,\qquad \Pi_{\ne}f=f_{\ne}.
\]
We denote by $H^s(\T\times\R)$ and $H^s(\R)$ the inhomogeneous Sobolev spaces, and write $H^s$ for $H^s(\T\times\R)$. The horizontal measure is normalized by $\dd x/(2\pi)$, so the two Sobolev norms coincide for an $x$-independent function. For any given norm, we use the convention $\|(f,g)\|^2=\|f\|^2+\|g\|^2$.
Define
\[
 \Lambda_t^2=1-\partial_x^2-(\partial_y+t\partial_x)^2,
 \qquad
 \Lambda_t^a(k,\eta)=(1+k^2+(\eta+kt)^2)^{a/2}.
\]
For every $a\in\R$, the operator $\Lambda_t^a$ commutes with $\partial_t+y\partial_x$.

\begin{theorem}\label{thm:main}
Fix $N>10$. There exist constants $0<\delta_N\le1$ and $C_N>0$, independent of $\nu,\mu,\beta$, with the following property.
Suppose
\begin{equation}
 0<\nu,\mu\le1,\qquad
 |\beta|\ge\frac{\nu+\mu}{3\sqrt{\nu\mu}},\qquad
 \rho=\min\{\nu,\mu\}.
 \label{eq:strongfield}
\end{equation}
If the divergence-free initial data $(u_{\rm in},b_{\rm in})$ satisfy
\begin{equation}
 \varepsilon^2=\|u_{\rm in}\|_{H^N}^2+\|b_{\rm in}\|_{H^N}^2,
 \qquad \varepsilon\le\delta_N\rho^{1/3},
 \label{eq:smallness}
\end{equation}
then the perturbed MHD equations \eqref{eq:perturbation} admit a unique global solution satisfying
\begin{align}
 &\sup_{t\ge0}\Big[
   \|(u_0,b_0)(t)\|_{H^N(\R)}^2
   +e^{\rho^{1/3}t/8}\|\Lambda_t^N(u_{\ne},b_{\ne})(t)\|_{L^2}^2\Big]
       \le C_N\varepsilon^2,\label{eq:mainbound}\\
 &\int_0^\infty\Bigl[
   \nu\|\partial_yu_0\|_{H^N(\R)}^2
   +\mu\|\partial_yb_0\|_{H^N(\R)}^2
   +e^{\rho^{1/3}t/8}\bigl(
      \nu\|\nabla\Lambda_t^N u_{\ne}\|_{L^2}^2
      +\mu\|\nabla\Lambda_t^N b_{\ne}\|_{L^2}^2\bigr)\Bigr]\dd t
       \le C_N\varepsilon^2,\label{eq:maindiff}\\
 &\rho^{1/3}\int_0^\infty e^{\rho^{1/3}t/8}
                  \|\Lambda_t^N(u_{\ne},b_{\ne})\|_{L^2}^2\dd t
       \le C_N\varepsilon^2.\label{eq:mainED}
\end{align}
\end{theorem}

When $\nu=\mu$, the magnetic field condition reduces to $|\beta|\ge2/3$. More generally, if $\mu=a\nu$ for a fixed $a>0$, it becomes
$|\beta|\ge(\sqrt a+a^{-1/2})/3$.
The initial perturbations may have nonzero horizontal means, which are controlled by \eqref{eq:mainbound}--\eqref{eq:maindiff}. Since incompressibility and square integrability imply $u_0^2=b_0^2=0$, we also obtain the following decay estimates.

\begin{corollary}\label{cor:decay}
Under the assumptions of Theorem~\ref{thm:main}, the velocity and magnetic perturbations satisfy, for all $t\ge0$,
\begin{align}
 \|\Lambda_t^N(u_{\ne},b_{\ne})(t)\|_{L^2}
 &\le C_N\varepsilon e^{-\rho^{1/3}t/16},
 \label{eq:fielddecay}\\
 \|\Lambda_t^{N-1}(u^2,b^2)(t)\|_{L^2}
 &\le C_N\varepsilon(1+t)^{-1}e^{-\rho^{1/3}t/16}.
 \label{eq:verticaldecay}
\end{align}
\end{corollary}

Estimate \eqref{eq:fielddecay} gives enhanced dissipation of the nonzero horizontal modes, while \eqref{eq:verticaldecay} gives inviscid damping of the vertical components.

\begin{remark}
The constant $3$ in \eqref{eq:strongfield} is chosen for convenience, and we have not tried to optimize it. The lower bound on $|\beta|$ is used only to obtain the coercivity of the diffusion terms for the pair $(\widetilde u,b)$ and to absorb the linear pressure term into the weight $g$ (see also Steps~1 and 2 in the proof of Proposition~\ref{prop:closedenergy}). 
\end{remark}

We next explain the main difficulties and the ideas of the proof. In the MHD system, the vorticity and current can grow transiently even when the velocity and magnetic field remain bounded. Moreover, the horizontal components do not have the inviscid decay of the velocity in the two-dimensional Navier--Stokes equations. To obtain the $\rho^{1/3}$ stability scale, we use the magnetic coupling in both the quadratic energy estimate and an energy correction.

First, following the adapted velocity introduced in \cite{KnoNR25}, we set
\[
 \widetilde u=u+\beta^{-1}\partial_x^{-1}b^2_{\ne}e_1.
\]
The second component of the magnetic equation cancels the term $u^2e_1$ in the velocity equation, and the magnetic shear term is incorporated into $\beta\partial_x\widetilde u$. We then consider
\[
 E(t)=\|\mathcal A\widetilde u(t)\|_{L^2}^2
       +\|\mathcal Ab(t)\|_{L^2}^2,
\]
where $\mathcal A$ is the time-dependent Fourier multiplier defined in \eqref{eq:Af}. Its bounded factors control the remaining pressure term and generate the weighted space-time estimates used in the nonlinear analysis.

We treat the nonlinear convection reactions together in the Elsasser variables
\[
 Z^s=u+sb,\qquad s=\pm1.
\]
By incompressibility, symmetrization gives a difference
\[
 h(k,\eta)-h(m,\zeta),\qquad
 k=\ell+m,\quad \eta=\xi+\zeta,
\]
where $h$ is the energy weight defined in \eqref{eq:weight}. The Alfv\'en contribution to the time derivative of the corresponding nonlinear product has coefficient
\[
 -s\ii\beta k-s\ii\beta\ell+s\ii\beta m
       =-2s\ii\beta\ell.
\]
For $\ell\ne0$, this identity determines the correction $\mathcal B$ in \eqref{eq:B}. Its time derivative cancels the main nonlinear term in the derivative of $E-\mathcal B$.

Differentiating $\mathcal B$ also produces other nonlinear reactions. The key point is to combine the nonlinear contributions from its first and third factors before estimating them. For $k=j+\ell+p$ and $\eta=\alpha+\xi+\gamma$, this gives the second difference
\[
 h(k,\eta)-h(j+p,\alpha+\gamma)
 -h(j+\ell,\alpha+\xi)+h(j,\alpha),
\]
together with a commutator of the interaction kernels. Using these cancellations and the Fourier multiplier estimates, together with the weighted $L^2$ bounds in time and enhanced dissipation, we bound the remaining nonlinear reactions by
\(
 C_N\rho^{-2/3}\mathcal E(T_1)^2,
\)
where $\mathcal E(T_1)$ is defined in \eqref{eq:controlenergy}. We also estimate the interactions involving the horizontal means and bound both $\mathcal B(0)$ and $\mathcal B(T_1)$ directly. The regularity condition $N>10$ ensures the Fourier $L^1$ bounds used in the nonlinear estimates.

For unequal viscosity and resistivity, we retain the separate coefficients of the physical diffusion terms. Incompressibility and \eqref{eq:strongfield} give coercivity in terms of
\(
 \nu\|\nabla\mathcal Au\|_{L^2}^2
       +\mu\|\nabla\mathcal Ab\|_{L^2}^2
\).
In the diffusion terms arising from $\mathcal B'$, the two derivatives are assigned to distinct amplitudes and estimated with their corresponding viscosity or resistivity. The coefficient is controlled by
\[
 \frac{\nu+\mu}{|\beta|\sqrt{\nu\mu}}\le3.
\]
This yields estimates uniform in the parameters satisfying \eqref{eq:strongfield}.

Combining the nonlinear estimates, we obtain
\begin{equation}
 \mathcal E(T_1)\le C_N\Big[
 \varepsilon^2+\rho^{-1/3}\mathcal E(T_1)^{3/2}
                  +\rho^{-2/3}\mathcal E(T_1)^2\Big].
\label{eq:introclosure}
\end{equation}
Under the bootstrap assumption $\mathcal E(T_1)^{1/2}\le M\varepsilon$, the nonlinear terms are bounded by
\[
 C_N\bigl(M\varepsilon\rho^{-1/3}
       +M^2\varepsilon^2\rho^{-2/3}\bigr)\mathcal E(T_1).
\]
They are absorbed when $\varepsilon\le\delta_N\rho^{1/3}$ and $\delta_N$ is sufficiently small. This proves the stability estimate at the $\rho^{1/3}$ scale.

The rest of the paper is organized as follows. In Section~\ref{sec:equations}, we construct the energy and reduce Theorem~\ref{thm:main} to the estimate \eqref{eq:introclosure}. Section~\ref{sec:identity} derives the corrected energy identity. The Fourier multiplier estimates are proved in Section~\ref{sec:multiplier}. Sections~\ref{sec:cubic} and \ref{sec:quartic} establish estimates for the nonlinear reactions. Finally, Section~\ref{sec:bootstrap} completes the proof of the energy estimate and Corollary~\ref{cor:decay}.

\section{The energy and the bootstrap estimate}
\label{sec:equations}

Throughout Sections~\ref{sec:equations}--\ref{sec:bootstrap}, we assume \eqref{eq:strongfield}. In this section, we first present the notations used throughout the paper. Then we construct the multiplier \(\mathcal A\), the energy \(E\), the correction \(\mathcal B\), and the energy functional \(\mathcal E(t)\).
Proposition~\ref{prop:closedenergy} states the estimate for \(E-\mathcal B\), and the resulting {\it a priori} estimate \eqref{eq:closedenergy}. 
We use this proposition to prove Theorem~\ref{thm:main} at the end of the section. The corresponding energy estimates are established in Sections~\ref{sec:identity}--\ref{sec:bootstrap}.

\subsection{Notations and construction of the Fourier multipliers}

Throughout the paper, constants \(C,C_N\) may change between inequalities and depend only on \(N\), unless stated otherwise.
They are uniform in the parameters satisfying \eqref{eq:strongfield}. Every sum over \(s\) uses \(s=+1,-1\).  We equip $\T$ with measure $\dd x/(2\pi)$ and use $k\in\Z$, $\eta\in\R$ as Fourier variables:
\[
 \widehat f(k,\eta)
 =\frac1{2\pi}\int_\R\int_\T
       f(x,y)e^{-\ii(kx+\eta y)}\dd x\dd y.
\]
For Fourier amplitudes, $\|\cdot\|_p$ denotes the $L^p$ norm with counting measure in $k$ and measure $\dd\eta/(2\pi)$ in $\eta$. The convolution convention is
\begin{equation*}
 (\widehat f*\widehat g)(k,\eta)
 =\sum_{\ell\in\Z}\int_\R
       \widehat f(\ell,\xi)\widehat g(k-\ell,\eta-\xi)
                                      \frac{\dd\xi}{2\pi},
 \qquad \widehat{fg}=\widehat f*\widehat g.
\end{equation*}
We write $|D|=(-\Delta)^{1/2}$. For a Fourier amplitude $F$,
the same notation denotes multiplication by its symbol: \( (|D|F)(k,\eta)=\sqrt{k^2+\eta^2}F(k,\eta) \).
For scalar or vector fields, \(\langle v,w\rangle\) denotes the \(L^2\) inner product, with measure \(\dd x\dd y/(2\pi)\). We also write $\av t=\sqrt{1+t^2}$ and $\av{\rho^{1/3}t}=\sqrt{1+\rho^{2/3}t^2}$. For notational simplicity, we often suppress the time variable and write $f(k,\eta)$ in place of $f(t,k,\eta)$ when no confusion can arise.

Now we construct the Fourier multipliers used throughout the paper. We write
\begin{equation}
 e(k,\eta)=\frac{(-\eta,k)}{\sqrt{k^2+\eta^2}},\qquad
 \kappa(k,\eta)=
 \begin{cases}
 \frac{|k|}{\sqrt{k^2+\eta^2}},&k\ne0,\\
 0,&k=0.
 \end{cases}
 \label{eq:frequencies}
\end{equation}
We also set $e(0,0)=0$.

We define $g$ and $G$ by
\begin{equation}
\begin{aligned}
 g(t,k,\eta)=\sum_{n\in\Z\setminus\{0\}}
   \frac{(1+(k-n)^2)^{-4}}
        {1+(t-(\eta+kt)/n)^2},\quad
 G(t,k,\eta)=\int_0^t g(s,k,\eta+k(t-s))\dd s.
\end{aligned}
\label{eq:g}
\end{equation}
The series is absolutely convergent. Consequently,
\begin{equation}
 0\le G(t,k,\eta)\le\pi\sum_{q\in\Z}(1+q^2)^{-4},\quad
 \sup_{k,\eta}\int_0^\infty g(t,k,\eta-kt)\dd t<\infty,
 \quad \kappa^2(k,\eta)\le g(k,\eta).
\label{eq:gprimitive}
\end{equation}
The last inequality uses the summand $n=k$, whose denominator is $1+(\eta/k)^2$, when $k\ne0$. 
For $0<\sigma\le1$, let
\begin{equation}
 \Phi_\sigma(t,a)=\arctan(\sigma(t-a))+\arctan(\sigma a)
       =\int_0^t\frac{\sigma\dd s}{1+\sigma^2(s-a)^2}.
\label{eq:Phi}
\end{equation}
Set
\begin{equation}
\begin{aligned}
 V(t,k,\eta)&=
 \begin{cases}
 \Phi_{\rho^{1/3}}(t,(\eta+kt)/k),&k\ne0,\\
 0,&k=0,
 \end{cases}\\
 \theta(t,k,\eta)&=\exp\{-2(7G(t,k,\eta)+V(t,k,\eta))\},\\
 h(t,k,\eta)&=\theta(t,k,\eta)\Lambda_t^{2N}(k,\eta).
\end{aligned}
\label{eq:weight}
\end{equation}
The bounded factor satisfies
\begin{equation}
 0<\exp\bigg(-14\pi\sum_{q\in\Z}(1+q^2)^{-4}-2\pi\bigg)
       \le\theta(k,\eta)\le1,\quad
 \theta(0,k,\eta)=1,\quad\theta(-k,-\eta)=\theta(k,\eta).
\label{eq:thetabounds}
\end{equation}
Since $(\partial_t-k\partial_\eta)G=g$, the factor $\exp\{-7G\}$ gives a positive $g$-weighted term in the energy estimate. The inequality $\kappa^2\le g$ allows us to control the linear pressure term by this term.
The function $V$ gives additional damping on the nonzero horizontal modes, which combines with diffusion to yield enhanced dissipation. For $k\ne0$, we have
\begin{equation}
 V(k,\eta)=\arctan\big(\rho^{1/3}(\eta/k+t)\big)
              -\arctan(\rho^{1/3}\eta/k),\quad
(\partial_t-k\partial_\eta)V(k,\eta)=\frac{\rho^{1/3}}{1+\rho^{2/3}\eta^2/k^2}.
 \label{eq:physicalV}
\end{equation}
For $k=0$, both $V(k,\eta)$ and $(\partial_t-k\partial_\eta)V(k,\eta)$ vanish. Also,
\begin{equation}
 G(k,\eta)=\sum_{n\ne0}(1+(k-n)^2)^{-4}
                  \Phi_1(t,(\eta+kt)/n).
 \label{eq:physicalG}
\end{equation}
We now define the operator as follows:
\begin{equation}
 \widehat{\mathcal A(t)w}(k,\eta)
 =e^{\rho^{1/3}t\mathbf 1_{k\ne0}/16}\sqrt{\theta(t,k,\eta)}\Lambda_t^N(k,\eta)\widehat w(k,\eta)=:\mathcal A(t,k,\eta)\widehat w(k,\eta).
 \label{eq:Af}
\end{equation}
We use the same notation $\mathcal A$ for the Fourier multiplier and its symbol.
The uniform positive bounds for $\theta$ imply that $\|\mathcal Aw\|_2$ is equivalent to $\|e^{\rho^{1/3}t\mathbf 1_{k\ne0}/16}\Lambda_t^Nw\|_{2}$. 
Since $(\partial_t-k\partial_\eta)\Lambda_t^N(k,\eta)=0$, we also have 
\begin{equation}
(\partial_t-k\partial_\eta)\mathcal A(t,k,\eta)=\Big(\frac1{16}\rho^{1/3}\mathbf 1_{k\ne0}-7g(k,\eta)-(\partial_t-k\partial_\eta) V(k,\eta)\Big)\mathcal A(t,k,\eta).
 \label{eq:2}
\end{equation}

\subsection{The energy functional and the bootstrap argument}

Incompressibility gives \(u_0^2=b_0^2=0\).
The commutator $[\partial_t+y\partial_x,\diver]v=-\partial_xv^2$ implies $\diver((\partial_t+y\partial_x)u)=\partial_xu^2$.
With the Leray projection
$\mathbb P=I-\nabla\Delta^{-1}\diver$, we obtain
\begin{equation*}
 (\partial_t+y\partial_x)u
   =\nu\Delta u+\beta\partial_xb-u^2e_1
     +2\partial_x\nabla\Delta^{-1}u^2
     +\mathbb P(b\cdot\nabla b-u\cdot\nabla u).
\end{equation*}
The two shear terms in the original system are $-u^2e_1$ in the velocity equation and $+b^2e_1$ in the magnetic equation in \eqref{eq:perturbation}. Testing the velocity shear term directly would give $-2\langle\mathcal Au^1,\mathcal Au^2\rangle$. The pointwise bound for this term contains only one factor $\kappa$ while the multiplier supplies the square $g\ge\kappa^2$ (see \eqref{eq:goodorder} and \eqref{eq:Eprime}). 

The magnetic equation provides a way to remove both shear terms.
Indeed, its second component contains $\beta\partial_xu^2$.
Applying $\beta^{-1}\partial_x^{-1}\Pi_{\ne}$ to that equation therefore produces $+u^2$, which cancels $-u^2$ in the first velocity equation. This is the reason for the corrected velocity ($e_1=(1,0)$)
\begin{equation}
 \widetilde u=u+\beta^{-1}\partial_x^{-1}b^2_{\ne}e_1,
\label{eq:goodunknown}
\end{equation}
which obeys
\begin{align*}
(\partial_t+y\partial_x)\widetilde u={}&\nu\Delta\widetilde u
       +\beta\partial_xb
       +2\partial_x\nabla\Delta^{-1}u^2\\
&+\mathbb P(b\cdot\nabla b-u\cdot\nabla u)
       +\beta^{-1}\Pi_{\ne}(b^1u^2-u^1b^2)e_1+\frac{\mu-\nu}{\beta}
                  \partial_x^{-1}\Delta b^2_{\ne}e_1.	
\end{align*}
Lemma~\ref{lem:goodunknown} derives this equation together with the equation for \(b\). The remaining linear pressure term contains $\kappa^2$ and can be absorbed by $g$ (see also \eqref{eq:pressurebound}). 

We define the energy by
\begin{equation}
 E(t)=\|\mathcal A\widetilde u(t)\|_2^2 + \|\mathcal Ab(t)\|_2^2.
\label{eq:quadraticenergy}
\end{equation}
The transformation \((u,b)\mapsto(\widetilde u,b)\) and its inverse are uniformly bounded under \eqref{eq:strongfield}.
Lemma~\ref{lem:energycomparison} gives the precise comparison used below.
 
However, when applying the energy method to $(\widetilde u,b)$, the weighted energy identity for \(E\), derived in \eqref{eq:Eprime}, contains the two nonlinear inner products
\[
2\langle\mathcal Au,\mathcal A\mathbb P(b\cdot\nabla b-u\cdot\nabla u)\rangle
           +2\langle\mathcal Ab,\mathcal A(b\cdot\nabla u-u\cdot\nabla b)\rangle.
\]
To handle this, we use the Elsasser fields \(Z^+=u+b\) and \(Z^-=u-b\).
For $s\in\{-1,1\}$, write $Z^s=u+s b$. At the Fourier frequency $(k,\eta)$, incompressibility gives $(k,\eta)\cdot(\widehat Z^{s,1},\widehat Z^{s,2})=0$. Thus $\widehat Z^s(k,\eta)$ has a single component along the unit vector 
\[
e(k,\eta)=\frac{(-\eta,k)}{\sqrt{k^2+\eta^2}}.
\]
Hence we can define its scalar amplitude explicitly by
\begin{equation}
\begin{aligned}
 f^s(t,k,\eta)
 :={}&-\ii e(k,\eta)\cdot
       \bigl(\widehat u(t,k,\eta)+s\widehat b(t,k,\eta)\bigr)\\
 ={}&\frac{\ii\eta}{\sqrt{k^2+\eta^2}}
                 \bigl(\widehat u^1(k,\eta)+s\widehat b^1(k,\eta)\bigr)
 -\frac{\ii k}{\sqrt{k^2+\eta^2}}
                 \bigl(\widehat u^2(k,\eta)+s\widehat b^2(k,\eta)\bigr).
\end{aligned}
\label{eq:fexplicit}
\end{equation}
A direct computation gives $f^s(-k,-\eta)=\overline{f^s(k,\eta)}$. The definitions of $f^s(t,k,\eta)$ at $(k,\eta)=(0,0)$ may be chosen arbitrarily. With this scalar amplitude, it holds that
\begin{equation}
 \widehat Z^s(k,\eta)=\ii e(k,\eta)f^s(k,\eta).
\label{eq:elsasser}
\end{equation}
Both original fields are recovered by
\begin{equation}
 \widehat u(k,\eta)=\frac{\ii e(k,\eta)}{2}(f^+(k,\eta)+f^-(k,\eta)),
 \qquad
 \widehat b(k,\eta)=\frac{\ii e(k,\eta)}{2}(f^+(k,\eta)-f^-(k,\eta)).
\label{eq:freconstruction}
\end{equation}
Consequently
\begin{equation}
 |f^s(k,\eta)|=|\widehat{u+s b}(k,\eta)|,
 \qquad
 |f^+(k,\eta)|^2+|f^-(k,\eta)|^2
       =2\bigl(|\widehat u(k,\eta)|^2+|\widehat b(k,\eta)|^2\bigr).
\label{eq:fparseval}
\end{equation}
To relate them to the velocity used in \(E\), we also compute \(\widetilde u\) in terms of \(f^s\). For $k\ne0$, the second component of \eqref{eq:freconstruction} gives
\[
 \frac{\widehat{b^2}(k,\eta)}{\ii\beta k}
 =\frac1{\ii\beta k}\frac{\ii k}{2\sqrt{k^2+\eta^2}}
                              (f^+-f^-)(k,\eta)
 =\frac{(f^+-f^-)(k,\eta)}{2\beta\sqrt{k^2+\eta^2}}.
\]
Substituting this into \eqref{eq:goodunknown} yields
\begin{equation*}
 \widehat{\widetilde u}(k,\eta)
 =\frac{\ii e(k,\eta)}2(f^++f^-)(k,\eta)
  +\frac{(f^+-f^-)(k,\eta)}{2\beta\sqrt{k^2+\eta^2}}e_1,
 \qquad k\ne0.
\end{equation*}
At $k=0$, the correction vanishes and $\widetilde u_0=u_0$.

We then write
\begin{equation}
F^s(k,\eta)=e^{\rho^{1/3}t\mathbf 1_{k\ne0}/16}\Lambda_t^N(k,\eta)f^s(k,\eta),
 \label{eq:weightedamplitudes}
\end{equation}
which has the same Sobolev and exponential factors as $\mathcal A Z^s$, with the bounded factor $\sqrt\theta$ omitted.
We then define the energy functional as follows:
\begin{equation}
\begin{aligned}
 \mathcal E(T_1)={}&\sup_{0\le t\le T_1}\sum_s\|F^s(t)\|_2^2\\
 &+\int_0^{T_1}\sum_s\sum_{k\in\Z}\int_{\R}
       (g(k,\eta)+\rho^{1/3}\mathbf 1_{k\ne0})|F^s(k,\eta)|^2\frac{\dd\eta}{2\pi}\dd t\\
 &+2\int_0^{T_1}\sum_{k\in\Z}\int_{\R} e^{\rho^{1/3}t\mathbf 1_{k\ne0}/8}\Lambda_t^{2N}(k,\eta)(k^2+\eta^2)
          \bigl(\nu|\widehat u(k,\eta)|^2+\mu|\widehat b(k,\eta)|^2\bigr)\frac{\dd\eta}{2\pi}\dd t.
\end{aligned}
\label{eq:controlenergy}
\end{equation}
The inequality
\begin{equation}
\begin{aligned}
\rho\sum_s\||D|F^s(t)\|_2^2
\le 2\sum_{k\in\Z}\int_{\R} e^{\rho^{1/3}t\mathbf 1_{k\ne0}/8}\Lambda_t^{2N}(k,\eta)(k^2+\eta^2)
          \bigl(\nu|\widehat u(k,\eta)|^2+\mu|\widehat b(k,\eta)|^2\bigr)\frac{\dd\eta}{2\pi}.
 \label{eq:minimumcontrol}
\end{aligned}
\end{equation}
shows that the mixed diffusion integral controls $\rho\sum_s\||D|F^s\|_2^2$. 

The following lemma gives the precise comparison.

\begin{lemma}
\label{lem:energycomparison}
Under \eqref{eq:strongfield}, the following pointwise comparison holds:
\begin{equation}
\begin{aligned}
 \frac14\bigl(|\widehat u(k,\eta)|^2+|\widehat b(k,\eta)|^2\bigr)
 \le|\widehat{\widetilde u}(k,\eta)|^2+|\widehat b(k,\eta)|^2
 \le4\bigl(|\widehat u(k,\eta)|^2+|\widehat b(k,\eta)|^2\bigr).
\end{aligned}
\label{eq:physicalenergycomparison}
\end{equation}
Consequently, there exist two universal constants $c,C>0$ such that for $(k,\eta)\neq(0,0)$,
\begin{equation}
 c\sum_s|F^s(k,\eta)|^2
 \le\mathcal A^2(k,\eta)
       \bigl(|\widehat{\widetilde u}(k,\eta)|^2+|\widehat b(k,\eta)|^2\bigr)
 \le C\sum_s|F^s(k,\eta)|^2.
 \label{eq:energyequivalence}
\end{equation}
\end{lemma}
\begin{proof}
The condition \eqref{eq:strongfield} and $\nu+\mu\ge2\sqrt{\nu\mu}$ give
\begin{equation}
 |\beta|\ge\frac23,\qquad
 \frac1{|\beta|}\le\frac32,\qquad
 \frac{\nu+\mu}{|\beta|\sqrt{\nu\mu}}\le3.
 \label{eq:coefficientratios}
\end{equation}
For $k\ne0$, the velocity correction obeys
\[
 |\widehat{\widetilde u-u}(k,\eta)|
 =\frac{|\widehat b^2(k,\eta)|}{|\beta k|}
 \le\frac32|\widehat b(k,\eta)|.
\]
Young's inequality therefore gives
\begin{align*}
 |\widehat{\widetilde u}|^2+|\widehat b|^2
 \le4|\widehat u|^2+\frac43|\widehat{\widetilde u-u}|^2
                                      +|\widehat b|^2
 \le4\bigl(|\widehat u|^2+|\widehat b|^2\bigr).
\end{align*}
Applying the same inequality to $u=\widetilde u-(\widetilde u-u)$ yields
\[
 |\widehat u|^2+|\widehat b|^2
 \le4\bigl(|\widehat{\widetilde u}|^2+|\widehat b|^2\bigr).
\]
At $k=0$ the correction vanishes, so both bounds hold there as well. This proves \eqref{eq:physicalenergycomparison}. Finally,
\[
 \mathcal A^2(k,\eta)
       \bigl(|\widehat u(k,\eta)|^2+|\widehat b(k,\eta)|^2\bigr)
       =\frac{\theta(k,\eta)}2\sum_s|F^s(k,\eta)|^2.
\]
Combining this identity with \eqref{eq:thetabounds} proves \eqref{eq:energyequivalence}.
\end{proof}

Now since \(\mathbb Pu=u\) and \(\mathcal A\) commutes with
\(\mathbb P\), expansion of \(Z^s=u+sb\) gives 
\[
\begin{aligned}
&2\langle\mathcal Au,\mathcal A\mathbb P(b\cdot\nabla b-u\cdot\nabla u)\rangle
 +2\langle\mathcal Ab,\mathcal A(b\cdot\nabla u-u\cdot\nabla b)\rangle
 =-\sum_s\langle\mathcal AZ^s,\mathcal A(Z^{-s}\cdot\nabla Z^s)\rangle\\
={}&\frac12\sum_s\Rea\sum_{\ell,m\in\Z}\int_{\R^{2}}
       (\mathcal A^2(k,\eta)-\mathcal A^2(m,\zeta))a(k,\eta;\ell,\xi)
                 \overline{f^s(k,\eta)}f^{-s}(\ell,\xi)f^s(m,\zeta)\frac{\dd\xi\dd\zeta}{(2\pi)^{2}},
\end{aligned}
\]
where
\begin{equation}
 a(k,\eta;\ell,\xi)
 =\frac{\ell\eta-k\xi}{\sqrt{\ell^2+\xi^2}}\,
  \frac{k(k-\ell)+\eta(\eta-\xi)}
   {\sqrt{k^2+\eta^2}
        \sqrt{(k-\ell)^2+(\eta-\xi)^2}}.
 \label{eq:A}
\end{equation}
This identity is derived in Section~\ref{sec:identity}. The right-hand side of the above identity leads to the following bad term to be cancelled:
\begin{align*}
\frac12\sum_s\Rea\sum_{\substack{\ell,m\in\Z\\\ell\ne0}}
 \int_{\R^2}e^{\rho^{1/3}t\mathbf1_{k\ne0}\mathbf1_{m\ne0}/8}
 (h(k,\eta)-h(m,\zeta))a(k,\eta;\ell,\xi)\overline{f^s(k,\eta)}f^{-s}(\ell,\xi)f^s(m,\zeta)
 \frac{\dd\xi\dd\zeta}{(2\pi)^2}.
\end{align*}
Here \(k=\ell+m\) and \(\eta=\xi+\zeta\).
The magnetic coupling in the equation for \(Z^s\) is \(s\beta\partial_xZ^s\). Its contribution to the time derivative of $\overline{f^s(k,\eta)}, f^{-s}(\ell,\xi), f^s(m,\zeta)$ therefore has coefficient
\[
-s\ii\beta k-s\ii\beta\ell+s\ii\beta m=-2s\ii\beta\ell.
\]
For \(\ell\ne0\), multiplication by \(\ii/(4s\beta\ell)\) turns this coefficient into \(1/2\). This inspires us to define the correction
\begin{equation}
\begin{aligned}
\cB(t)&=\sum_s\Rea\sum_{\substack{\ell,m\in\Z\\\ell\ne0}}\int_{\R^{2}} \frac{\ii}{4s\beta\ell}e^{\frac18\rho^{1/3} t\mathbf 1_{k\ne0}\mathbf 1_{m\ne0}}(h(k,\eta)-h(m,\zeta))a(k,\eta;\ell,\xi)\\
&\qquad\qquad\qquad\qquad\quad\times \overline{f^s(k,\eta)}f^{-s}(\ell,\xi)f^s(m,\zeta)\frac{\dd\xi\dd\zeta}{(2\pi)^{2}}.
\label{eq:B}
\end{aligned}
\end{equation}
Section~\ref{sec:identity} verifies that the Alfv\'en part of \(\mathcal B'\) (see \eqref{eq:Ialfven}) is exactly the preceding bad term, and writes out all remaining terms of that derivative. We consequently estimate \(E-\mathcal B\).

These constructions lead to the following corrected energy estimate and the resulting bootstrap inequality for $\mathcal E$.
\begin{proposition}
\label{prop:closedenergy}
Let \((u,b)\) be a smooth solution of the MHD equations \eqref{eq:perturbation} on \([0,T_1]\), under \eqref{eq:strongfield}, with initial size
\(\varepsilon\le1\) as in \eqref{eq:smallness}. For every
\(0\le t\le T_1\),
\begin{equation}
\begin{aligned}
&E(t)-\mathcal B(t)
 +\frac14\int_0^t
   \bigl(\nu\|\nabla\mathcal A(\tau)u(\tau)\|_2^2
        +\mu\|\nabla\mathcal A(\tau)b(\tau)\|_2^2\bigr)\dd\tau\\
&+\int_0^t\sum_{k\in\Z}\int_\R
 \Big(7g(\tau,k,\eta)+\frac38\rho^{1/3}\mathbf1_{k\ne0}\Big)
 \mathcal A^2(\tau,k,\eta)
 \bigl(|\widehat{\widetilde u}(\tau,k,\eta)|^2
              +|\widehat b(\tau,k,\eta)|^2\bigr)
                    \frac{\dd\eta}{2\pi}\dd\tau\\
\le{} &E(0)-\mathcal B(0)
 +C_N\bigl[\rho^{-1/3}\mathcal E(T_1)^{3/2}
                     +\rho^{-2/3}\mathcal E(T_1)^2\bigr].
\end{aligned}
\label{eq:correctedenergyestimate}
\end{equation}
The correction $\cB(t)$ satisfies
\begin{equation}
\begin{aligned}
\sup_{0\le t\le T_1}|\mathcal B(t)|
 &\le C_N\rho^{-1/3}\mathcal E(T_1)^{3/2},\\
|\mathcal B(0)|&\le C_N\varepsilon^3,\qquad
E(0)\le4\varepsilon^2.
\end{aligned}
\label{eq:correctionboundarycontrol}
\end{equation}
Consequently, there is a constant \(C_0\ge4\), depending only on
\(N\), such that
\begin{equation}
 \mathcal E(T_1)\le C_0\Big[
       \varepsilon^2+\rho^{-1/3}\mathcal E(T_1)^{3/2}
                         +\rho^{-2/3}\mathcal E(T_1)^2\Big].
 \label{eq:closedenergy}
\end{equation}
All constants are independent of \(T_1,\nu,\mu,\beta\) in the
stated range.
\end{proposition}

The proof of Proposition~\ref{prop:closedenergy} is given in Section~\ref{sec:bootstrap}.
We now establish the following bootstrap estimate.

\begin{proposition}[Bootstrap improvement]\label{prop:bootstrap}
Fix $N>10$. There are constants $M\ge8$ and $0<\delta_N\le1$, depending only on $N$, such that the following holds for every $(\nu,\mu,\beta)$ satisfying \eqref{eq:strongfield}.
Let $(u,b)$ be a smooth solution of the MHD equations \eqref{eq:perturbation} on $[0,T_1]$, with initial size $\varepsilon$ as in
\eqref{eq:smallness}. If
\begin{equation}
 \varepsilon\le\delta_N\rho^{1/3},
 \qquad \mathcal E(T_1)^{1/2}\le M\varepsilon,
\label{eq:bootstrapassumption}
\end{equation}
then
\begin{equation}
 \mathcal E(T_1)^{1/2}\le\frac12M\varepsilon.
\label{eq:bootstrapconclusion}
\end{equation}
The constants are independent of $T_1,\nu,\mu,\beta$.
\end{proposition}

\begin{proof}
Fix $C_0\ge4$ from Proposition~\ref{prop:closedenergy}. It depends only on $N$. Choose
\begin{equation}
 M=4\sqrt{C_0},\qquad
 0<\delta_N\le\frac1{16C_0M}\le 1,\qquad
 \varepsilon\le\delta_N\rho^{1/3}.
 \label{eq:bootstrapconstants}
\end{equation}
Under \eqref{eq:bootstrapassumption}, the choice \eqref{eq:bootstrapconstants} gives
\[
 C_0\rho^{-1/3}\mathcal E(T_1)^{1/2}
       +C_0\rho^{-2/3}\mathcal E(T_1)
       \le\frac1{16}+\frac1{256C_0}<\frac14.
\]
Since $\varepsilon\le\delta_N\le1$, Proposition~\ref{prop:closedenergy} gives
\begin{equation*}
 \mathcal E(T_1)\le\frac43C_0\varepsilon^2,
 \qquad
 \mathcal E(T_1)^{1/2}\le\frac{M}{\sqrt{12}}\varepsilon
                              <\frac M2\varepsilon.
\end{equation*}
The proof is thus complete.
\end{proof}

We now deduce Theorem~\ref{thm:main}.

\begin{proof}[Proof of Theorem~\ref{thm:main}]
For fixed positive $\nu$ and $\mu$, the standard Sobolev local theory for viscous, resistive MHD yields a unique maximal $H^N$ perturbation solution on $[0,T_*)$; see \cite{ST83} for the classical local theory.
The same energy construction applies to \eqref{eq:perturbation}: the Couette transport is skew-adjoint in $L^2$, its commutators with derivatives have order at most $N$, and the constant magnetic coupling cancels in the sum of the velocity and magnetic energies.
In particular, a finite $T_*$ can occur only if $\|(u,b)(t)\|_{H^N}$ becomes unbounded as $t\to T_*$.
The estimates below are first justified for smooth data and extend to $H^N$ data by smooth divergence-free approximation and the local theory.

If $\varepsilon=0$, the perturbation is identically zero. Suppose henceforth that $\varepsilon>0$, and choose $M$ and $\delta_N$ as in \eqref{eq:bootstrapconstants}.
Define
\[
 T_2=\sup\bigl\{t\in[0,T_*):
       \mathcal E(\tau)^{1/2}\le M\varepsilon
       \text{ for every }0\le\tau\le t\bigr\}.
\]
The map $t\mapsto\mathcal E(t)$ is continuous on $[0,T_*)$, since its first term is the running supremum of a continuous energy and its remaining terms are locally absolutely continuous.
Since $\mathcal E(0)=2\varepsilon^2$ and $M\ge8$, the bound holds strictly at $t=0$, and hence $T_2>0$.
For every $t<T_2$, Proposition~\ref{prop:bootstrap} gives the improved bound \eqref{eq:bootstrapconclusion}.
If $T_2<T_*$, continuity gives $\mathcal E(T_2)^{1/2}\le M\varepsilon/2$ and extends the original bound to a short interval beyond $T_2$, contradicting its definition. Thus $T_2=T_*$.

It remains to rule out $T_*<\infty$. 
The frequency inequality
\[
 1+k^2+\eta^2
 \le 2(1+t^2)\bigl(1+k^2+(\eta+kt)^2\bigr)
\]
and \eqref{eq:fparseval} imply
\[
 \|(u,b)(t)\|_{H^N}
 \le C_N\av t^N\|\Lambda_t^N(u,b)(t)\|_2
 \le C_N\av{T_*}^N M\varepsilon,
 \qquad 0\le t<T_*.
\]
This contradicts the continuation criterion, so $T_*=\infty$.
Finally, the first, third, and second terms of \eqref{eq:controlenergy}, respectively, give \eqref{eq:mainbound}, \eqref{eq:maindiff}, and \eqref{eq:mainED} by Parseval's identity and \eqref{eq:fparseval}.
\end{proof}	

\section{Derivation of the corrected energy identity}
\label{sec:identity}

In this section, we derive an identity for \(E-\mathcal B\).

\subsection{The corrected equations and the energy identity}

We first establish the equations satisfied by \(\widetilde u\) and \(b\).

\begin{lemma}
\label{lem:goodunknown}
Let \((u,b)\) be a smooth solution of the MHD equations \eqref{eq:perturbation}, and define \(\widetilde u\) by \eqref{eq:goodunknown}. Then
\begin{equation}
\begin{aligned}
 (\partial_t+y\partial_x)\widetilde u={}&\nu\Delta\widetilde u
       +\beta\partial_xb
       +2\partial_x\nabla\Delta^{-1}u^2\\
 &+\mathbb P(b\cdot\nabla b-u\cdot\nabla u)
       +\beta^{-1}\Pi_{\ne}(b^1u^2-u^1b^2)e_1+\frac{\mu-\nu}{\beta}
                  \partial_x^{-1}\Delta b^2_{\ne}e_1,\\
 (\partial_t+y\partial_x)b={}&\mu\Delta b+\beta\partial_x\widetilde u
                  +b\cdot\nabla u-u\cdot\nabla b.
\end{aligned}
\label{eq:good}
\end{equation}
The correction and the pressure satisfy
\begin{equation}
 |\widehat{u^2}(k,\eta)|=\kappa(k,\eta)|\widehat u(k,\eta)|,
 \qquad \sqrt{k^2+\eta^2}|\widehat{\widetilde u-u}(k,\eta)|
                          =|\beta|^{-1}\mathbf1_{k\ne0}|\widehat b(k,\eta)|.
\label{eq:goodorder}
\end{equation}
\begin{equation}
 \big|\widehat{\,2\partial_x\nabla\Delta^{-1}u^2\,}(k,\eta)\big|
       =\frac{2|k|}{\sqrt{k^2+\eta^2}}|\widehat{u^2}(k,\eta)|
       =2\kappa^2(k,\eta)|\widehat u(k,\eta)|.
\label{eq:pressurebound}
\end{equation}
The change $(u,b)\mapsto(\widetilde u,b)$ and its inverse are bounded on $H^a\times H^a$, with norm at most $1+|\beta|^{-1}$.
\end{lemma}
\begin{proof}
Here $\partial_x^{-1}\Pi_{\ne}$ is the Fourier multiplier $(\ii k)^{-1}$ at $k\ne0$ and zero at $k=0$. It has norm at most one on every Sobolev space. The map from $(u,b)$ to $(\widetilde u,b)$ and its inverse are therefore bounded, with operator norm at most $1+|\beta|^{-1}$.

Since $\diver u=\diver b=0$, we have
\begin{align*}
 b\cdot\nabla u^2-u\cdot\nabla b^2
 &=\partial_x(b^1u^2-u^1b^2),\nonumber\\
 \beta^{-1}\partial_x^{-1}\Pi_{\ne}
      (b\cdot\nabla u^2-u\cdot\nabla b^2)e_1
 &=\beta^{-1}\Pi_{\ne}(b^1u^2-u^1b^2)e_1.
\end{align*}
Using this identity and $u_0^2=b_0^2=0$, application of $\partial_t+y\partial_x$ to the correction gives
\[
 (\partial_t+y\partial_x)(\beta^{-1}\partial_x^{-1}b^2_{\ne}e_1)
 =\mu\Delta(\beta^{-1}\partial_x^{-1}b^2_{\ne}e_1)
    +u^2e_1+\beta^{-1}\Pi_{\ne}(b^1u^2-u^1b^2)e_1.
\]
The $u^2e_1$ term cancels the velocity shear term. The diffusion terms combine as
\[
 \nu\Delta u+\frac{\mu}{\beta}\partial_x^{-1}\Delta b^2_{\ne}e_1
 =\nu\Delta\widetilde u
    +\frac{\mu-\nu}{\beta}\partial_x^{-1}\Delta b^2_{\ne}e_1.
\]
In the magnetic equation, $\beta\partial_xu+b^2e_1=\beta\partial_x\widetilde u$.
These two cancellations prove \eqref{eq:good}.

For $k\neq0$, incompressibility gives
\[
 |\widehat{u^2}|=\frac{|k|}{\sqrt{k^2+\eta^2}}|\widehat u|,
 \qquad
 \sqrt{k^2+\eta^2}\frac{|\widehat{b^2}|}{|\beta k|}
       =|\beta|^{-1}|\widehat b|.
\]
For $k=0$, the correction vanishes, proving \eqref{eq:goodorder}.
The symbol of $2\partial_x\nabla\Delta^{-1}$ is $2k(k,\eta)/(k^2+\eta^2)$. Hence combining this with the above identities proves \eqref{eq:pressurebound}.
For $k=0$, the pressure term vanishes.
\end{proof}

We next derive the energy identity for $E(t)$ defined in \eqref{eq:quadraticenergy}.

\begin{proposition}
\label{prop:quadraticidentity}
Let \((u,b)\) be a smooth solution of \eqref{eq:perturbation}. With the notation in \eqref{eq:g}, \eqref{eq:weight}, \eqref{eq:Af}, and \eqref{eq:goodunknown}, we have
\begin{align}
 \frac{\dd}{\dd t}E(t)&+2\sum_{k\in\Z}\int_{\R}\mathcal A^2(k,\eta)(k^2+\eta^2)\big(\nu|\widehat{\widetilde u}(k,\eta)|^2+\mu|\widehat b(k,\eta)|^2\big)\frac{\dd\eta}{2\pi}\nonumber\\
 &+2\sum_{k\in\Z}\int_{\R}(7g(k,\eta)+(\partial_t-k\partial_\eta)V(k,\eta))\mathcal A^2(k,\eta)
           \big(|\widehat{\widetilde u}(k,\eta)|^2+|\widehat b(k,\eta)|^2\big)\frac{\dd\eta}{2\pi}\nonumber\\
 ={}&\frac18\rho^{1/3}\sum_{\substack{k\ne0}}\int_{\R}\mathcal A^2(k,\eta)\big(|\widehat{\widetilde u}(k,\eta)|^2+|\widehat b(k,\eta)|^2\big)\frac{\dd\eta}{2\pi}\nonumber\\
 &+4\langle\mathcal A\widetilde u,
             \mathcal A\partial_x\nabla\Delta^{-1}u^2\rangle
 -\frac{2(\mu-\nu)}{\beta}
   \langle\nabla\mathcal A\widetilde u^1,
                \nabla\mathcal A\partial_x^{-1}b^2_{\ne}\rangle\nonumber\\
 &+2\langle\mathcal Au,\mathcal A\mathbb P(b\cdot\nabla b-u\cdot\nabla u)\rangle
           +2\langle\mathcal Ab,\mathcal A(b\cdot\nabla u-u\cdot\nabla b)\rangle\nonumber\\
 &+2\langle\mathcal A(\widetilde u-u),\mathcal A\mathbb P(b\cdot\nabla b-u\cdot\nabla u)\rangle
           +2\left\langle\mathcal A\widetilde u,\mathcal A\bigl[\beta^{-1}\Pi_{\ne}(b^1u^2-u^1b^2)e_1\bigr]\right\rangle.
\label{eq:Eprime}
\end{align}
\end{proposition}
\begin{proof}
Parseval's identity gives
\begin{align*}
\frac12\frac{\dd}{\dd t}\|\mathcal A\widetilde u\|_2^2
=\sum_{k\in\Z}\int_\R
 \mathcal A(k,\eta)\partial_t\mathcal A(k,\eta)
 |\widehat{\widetilde u}(k,\eta)|^2\frac{\dd\eta}{2\pi}
+\Rea\sum_{k\in\Z}\int_\R
 \mathcal A^2(k,\eta)\overline{\widehat{\widetilde u}(k,\eta)}
 \cdot\partial_t\widehat{\widetilde u}(k,\eta)\frac{\dd\eta}{2\pi}.
\end{align*}
For each fixed \(k\), integration by parts in \(\eta\) gives
\begin{align*}
\Rea\int_\R k\mathcal A^2(k,\eta)
 \overline{\widehat{\widetilde u}(k,\eta)}
 \cdot\partial_\eta\widehat{\widetilde u}(k,\eta)\frac{\dd\eta}{2\pi}
=-\int_\R k\mathcal A(k,\eta)\partial_\eta\mathcal A(k,\eta)
                      |\widehat{\widetilde u}(k,\eta)|^2\frac{\dd\eta}{2\pi}.
\end{align*}
Thus by virtue of \((\partial_t-k\partial_\eta)\widehat w
          =\widehat{(\partial_t+y\partial_x)w}\), we get
\begin{align*}
\frac{\dd}{\dd t}\|\mathcal A\widetilde u\|_2^2
=2\langle\mathcal A\widetilde u,\mathcal A(\partial_t+y\partial_x)\widetilde u\rangle
+2\sum_{k\in\Z}\int_{\R}\frac{(\partial_t-k\partial_\eta)\mathcal A(k,\eta)}{\mathcal A(k,\eta)}|\mathcal A(k,\eta)\widehat{\widetilde u}(k,\eta)|^2\frac{\dd\eta}{2\pi}.
\end{align*}
Similarly, it holds that
\begin{align*}
\frac{\dd}{\dd t}\|\mathcal Ab\|_2^2
=2\langle\mathcal Ab,\mathcal A(\partial_t+y\partial_x)b\rangle
+2\sum_{k\in\Z}\int_{\R}\frac{(\partial_t-k\partial_\eta)\mathcal A(k,\eta)}{\mathcal A(k,\eta)}|\mathcal A(k,\eta)\widehat{b}(k,\eta)|^2\frac{\dd\eta}{2\pi}.
\end{align*}
Integration by parts gives
\[
 2\beta\langle\mathcal A\widetilde u,\partial_x\mathcal Ab\rangle
 +2\beta\langle\mathcal Ab,\partial_x\mathcal A\widetilde u\rangle=0.
\]
The diffusion terms are $-2\nu\|\nabla\mathcal A\widetilde u\|_2^2-2\mu\|\nabla\mathcal Ab\|_2^2$. A direct computation then gives
\begin{equation*}
\frac{2(\mu-\nu)}{\beta}
 \langle\mathcal A\widetilde u^1,
                 \mathcal A\partial_x^{-1}\Delta b^2_{\ne}\rangle=-\frac{2(\mu-\nu)}{\beta}
 \langle\nabla\mathcal A\widetilde u^1,
                 \nabla\mathcal A\partial_x^{-1}b^2_{\ne}\rangle.
\end{equation*}
Substituting the remaining terms of \eqref{eq:good} and using \eqref{eq:2}, we obtain \eqref{eq:Eprime}.
\end{proof}

\subsection{The nonlinear energy term in Elsasser variables}

We next rewrite the two inner products in \eqref{eq:Eprime} that are paired with \(\mathcal Au\) and \(\mathcal Ab\).
The two terms containing the velocity correction will be estimated directly in Section~\ref{sec:cubic}.

\begin{lemma}
\label{lem:scalar}
With the notation in \eqref{eq:fexplicit} and \eqref{eq:elsasser}, for $s=\pm1$, the amplitudes $f^s(k,\eta)$ satisfy
\begin{equation}
\begin{aligned}
(\partial_t-k\partial_\eta)f^s(k,\eta)
={}&
 \Big(s\ii\beta k-\frac{\nu+\mu}{2}(k^2+\eta^2)\Big)f^s(k,\eta)
 +\Big(\frac{k\eta}{(k^2+\eta^2)}
              -\frac{\nu-\mu}{2}(k^2+\eta^2)\Big)f^{-s}(k,\eta)\\
&+\sum_{\ell\in\Z}\int_{\R}a(k,\eta;\ell,\xi)f^{-s}(\ell,\xi)f^s(k-\ell,\eta-\xi)\frac{\dd\xi}{2\pi}.
\label{eq:scalar}
\end{aligned}
\end{equation}
Here $a$ is the explicit kernel \eqref{eq:A}. Recalling \eqref{eq:frequencies}, it also holds that
\begin{equation}
 a(k,\eta;\ell,\xi)
  =\bigl(e(\ell,\xi)\cdot(k-\ell,\eta-\xi)\bigr)
            \bigl(e(k,\eta)\cdot e(k-\ell,\eta-\xi)\bigr).
 \label{eq:kernelgeometric}
\end{equation}
\end{lemma}
\begin{proof}
By \eqref{eq:perturbation}, we can verify that the Elsasser variables $Z^s=u+sb$ satisfy
\[
 (\partial_t+y\partial_x)Z^s=\frac{\nu+\mu}{2}\Delta Z^s+\frac{\nu-\mu}{2}\Delta Z^{-s}+s\beta\partial_xZ^s
          -Z^{-s,2}e_1-\nabla p-Z^{-s}\cdot\nabla Z^s.
\]
Since $e(k,\eta)\cdot(\partial_t-k\partial_\eta)e(k,\eta)=0$ and $\widehat Z^s=\ii e f^s$ by \eqref{eq:elsasser}, taking the Fourier transform and the scalar product with $-\ii e(k,\eta)$ gives
\begin{align*}
(\partial_t-k\partial_\eta)f^s(k,\eta)={}&\ii\widehat{Z^{-s}\cdot\nabla Z^s}(k,\eta)\cdot e(k,\eta)\\
&+\Big(s\ii\beta k-\frac{\nu+\mu}{2}(k^2+\eta^2)\Big)f^s(k,\eta)
 +\Big(\frac{k\eta}{(k^2+\eta^2)}
              -\frac{\nu-\mu}{2}(k^2+\eta^2)\Big)f^{-s}(k,\eta),
\end{align*}
where we have used 
\begin{align*}
&\ii\widehat{\nabla p}(k,\eta)\cdot e(k,\eta)=-(k,\eta)\cdot e(k,\eta)\widehat p(k,\eta)=0,\\
&(-\ii e(k,\eta))\cdot(-\widehat{Z^{-s,2}}(k,\eta)e_1)
       =-e^1(k,\eta)e^2(k,\eta)f^{-s}(k,\eta)
       =\frac{k\eta}{(k^2+\eta^2)}f^{-s}(k,\eta).	
\end{align*}
For the quadratic term, a direct computation gives
\begin{align}
 &\widehat{Z^{-s}\cdot\nabla Z^s}(k,\eta)
 =\sum_{a=1}^2
       (\widehat{Z^{-s,a}}*\widehat{\partial_a Z^s})(k,\eta)\notag\\
=&-\ii\sum_{\ell\in\Z}\int_{\R}(e(\ell,\xi)\cdot (k-\ell,\eta-\xi))e(k-\ell,\eta-\xi) f^{-s}(\ell,\xi)f^s(k-\ell,\eta-\xi)\frac{\dd\xi}{2\pi}. \label{eq:1}
\end{align}
Substituting this expression into the preceding equation gives \eqref{eq:scalar} and \eqref{eq:kernelgeometric}.
\end{proof}

The symmetry needed to obtain a difference of energy weights is the following elementary identity.

\begin{lemma}
\label{lem:kernelidentities}
For \(k=\ell+m\) and \(\eta=\xi+\zeta\), the kernel in
\eqref{eq:A} satisfies
\begin{align}
 a(m,\zeta;-\ell,-\xi)=-a(k,\eta;\ell,\xi),
\label{eq:endpointsym}\\	
 |a(k,\eta;0,\xi)|\le |m|.\label{eq:meanA}
\end{align}
\end{lemma}
\begin{proof}

By $e(k,\eta)=(-\eta,k)/\sqrt{k^2+\eta^2}$ and \eqref{eq:kernelgeometric}, we find that \eqref{eq:endpointsym} follows from $e(-\ell,-\xi)\cdot(m+\ell,\zeta+\xi)=e(-\ell,-\xi)\cdot(m,\zeta)$. 
Moreover,
\[
|a(k,\eta;0,\xi)|\le |e(0,\xi)\cdot(k,\eta-\xi)||e(k,\eta)\cdot e(k,\eta-\xi)|\le |k|=|m|
\]
since $\ell=0$ implies $k=m$.
\end{proof}
Since $\mathbb Pu=u$, expanding $Z^s=u+sb$ gives
\[
\begin{aligned}
&2\langle\mathcal Au,\mathcal A\mathbb P(b\cdot\nabla b-u\cdot\nabla u)\rangle
 +2\langle\mathcal Ab,\mathcal A(b\cdot\nabla u-u\cdot\nabla b)\rangle
 =-\sum_s\langle\mathcal AZ^s,\mathcal A(Z^{-s}\cdot\nabla Z^s)\rangle.
\end{aligned}
\]
Substituting $\widehat Z^s(k,\eta)=\ii e(k,\eta)f^s(k,\eta)$ from \eqref{eq:elsasser} and using \eqref{eq:kernelgeometric} and \eqref{eq:1}, we obtain
\begin{align*}
-\langle\mathcal AZ^s,\mathcal A(Z^{-s}\cdot\nabla Z^s)\rangle &=-\Rea\sum_{k\in\Z}\int_{\R} \mathcal A^2(k,\eta)\overline{\widehat{Z^s}(k,\eta)}\cdot\widehat{(Z^{-s}\cdot\nabla Z^s)}(k,\eta)\frac{\dd\eta}{2\pi}\\
&= \Rea\sum_{k,\ell\in\Z}\int_{\R^{2}}
       \mathcal A^2(k,\eta)a(k,\eta;\ell,\xi)
       \overline{f^s(k,\eta)}f^{-s}(\ell,\xi)f^s(k-\ell,\eta-\xi)
                       \frac{\dd\xi\dd\eta}{(2\pi)^{2}}\\
&= \Rea\sum_{\ell,m\in\Z}\int_{\R^{2}}
       \mathcal A^2(k,\eta)a(k,\eta;\ell,\xi)
       \overline{f^s(k,\eta)}f^{-s}(\ell,\xi)f^s(m,\zeta)
                       \frac{\dd\xi\dd\zeta}{(2\pi)^{2}}.	
\end{align*}
In the last equality, we used the change of variables
\[
m=k-\ell,\qquad \zeta=\eta-\xi,
\qquad\text{so that}\qquad
k=\ell+m,\quad \eta=\xi+\zeta.
\]
Next, we make the change of variables
\[
(\ell,m,\xi,\zeta)\longmapsto
(-\ell,\ell+m,-\xi,\xi+\zeta).
\]
Noting that $\mathcal A$, $a$ and $\langle\mathcal AZ^s,\mathcal A(Z^{-s}\cdot\nabla Z^s)\rangle$ are real and using \eqref{eq:endpointsym} and $f^s(-k,-\eta)=\overline{f^s(k,\eta)}$, we obtain
\[
\begin{aligned}
-\langle\mathcal AZ^s,\mathcal A(Z^{-s}\cdot\nabla Z^s)\rangle=-\Rea\sum_{\ell,m\in\Z}\int_{\R^2}
 \mathcal A^2(m,\zeta)a(k,\eta;\ell,\xi)
 \overline{f^s(k,\eta)}f^{-s}(\ell,\xi)f^s(m,\zeta)
 \frac{\dd\xi\dd\zeta}{(2\pi)^2}.
\end{aligned}
\]
Hence we get
\begin{align}
 &2\langle\mathcal Au,\mathcal A\mathbb P(b\cdot\nabla b-u\cdot\nabla u)\rangle
                  +2\langle\mathcal Ab,\mathcal A(b\cdot\nabla u-u\cdot\nabla b)\rangle\nonumber\\
 ={}&\frac12\sum_s\Rea\sum_{\ell,m\in\Z}\int_{\R^{2}}
       (\mathcal A^2(k,\eta)-\mathcal A^2(m,\zeta))a(k,\eta;\ell,\xi)
                 \overline{f^s(k,\eta)}f^{-s}(\ell,\xi)f^s(m,\zeta)\frac{\dd\xi\dd\zeta}{(2\pi)^{2}}.
\label{eq:Tcommutator}
\end{align}
We split the sum in \eqref{eq:Tcommutator} into $\ell=0$ and $\ell\ne0$. In the latter part, we add and subtract
\[
e^{\rho^{1/3}t\mathbf1_{k\ne0}\mathbf1_{m\ne0}/8}(h(k,\eta)-h(m,\zeta)).
\]
This gives three terms. The first two are
\begin{align}
 &\frac12\sum_s\Rea\sum_{\substack{\ell,m\in\Z\\\ell=0}}\int_{\R^{2}}
       (\mathcal A^2(k,\eta)-\mathcal A^2(m,\zeta))a(k,\eta;\ell,\xi)
                          \overline{f^s(k,\eta)}f^{-s}(\ell,\xi)f^s(m,\zeta)\frac{\dd\xi\dd\zeta}{(2\pi)^{2}},
       \label{eq:T0}\\
 &\frac12\sum_s\Rea\sum_{\substack{\ell,m\in\Z\\\ell\ne0}}\int_{\R^{2}}
       \bigl[\mathcal A^2(k,\eta)-\mathcal A^2(m,\zeta)-e^{\frac18\rho^{1/3} t\mathbf 1_{k\ne0}\mathbf 1_{m\ne0}}(h(k,\eta)-h(m,\zeta))\bigr]
       \nonumber\\
 &\hspace{12mm}\times a(k,\eta;\ell,\xi)\overline{f^s(k,\eta)}f^{-s}(\ell,\xi)f^s(m,\zeta)\frac{\dd\xi\dd\zeta}{(2\pi)^{2}}.
       \label{eq:TJ}
\end{align}
The third term is
\begin{align}
\frac12\sum_s\Rea\sum_{\substack{\ell,m\in\Z\\\ell\ne0}}
 \int_{\R^2}e^{\rho^{1/3}t\mathbf1_{k\ne0}\mathbf1_{m\ne0}/8}
 (h(k,\eta)-h(m,\zeta))a(k,\eta;\ell,\xi)\overline{f^s(k,\eta)}f^{-s}(\ell,\xi)f^s(m,\zeta)
 \frac{\dd\xi\dd\zeta}{(2\pi)^2}.
\label{eq:cancelledcubic}
\end{align}
Thus the right-hand side of \eqref{eq:Tcommutator} is the sum of \eqref{eq:T0}, \eqref{eq:TJ}, and \eqref{eq:cancelledcubic}.
The next subsection shows that \eqref{eq:cancelledcubic} occurs in $\mathcal B'$.

\subsection{The time derivative of $\mathcal B$}

In the double integrals below, \(k=\ell+m\) and \(\eta=\xi+\zeta\).
Differentiating $\mathcal B$ defined in \eqref{eq:B} and integrating by parts in $\xi$ and $\zeta$, we obtain
\begin{equation}\label{ddbt-1}
\begin{split}
\frac{\dd}{\dd t}\mathcal B(t)
={}&\sum_s\Rea\sum_{\substack{\ell,m\in\Z\\\ell\ne0}}
 \int_{\R^2}\frac{\ii}{4s\beta\ell}
 (\partial_t-\ell\partial_\xi-m\partial_\zeta)\Bigl[
 e^{\rho^{1/3}t\mathbf1_{k\ne0}\mathbf1_{m\ne0}/8}
 (h(k,\eta)-h(m,\zeta))\\
&\qquad\times a(k,\eta;\ell,\xi)
 \overline{f^s(k,\eta)}f^{-s}(\ell,\xi)f^s(m,\zeta)\Bigr]
 \frac{\dd\xi\dd\zeta}{(2\pi)^2}.
\end{split}
\end{equation}
The chain rule gives
\begin{align*}
 (\partial_t-\ell\partial_\xi-m\partial_\zeta)
       &\bigl(\overline{f^s(k,\eta)}f^{-s}(\ell,\xi)f^s(m,\zeta)\bigr)
 =\overline{(\partial_t-k\partial_\eta)f^s(k,\eta)}
                   f^{-s}(\ell,\xi)f^s(m,\zeta)\\
& +\overline{f^s(k,\eta)}
          (\partial_t-\ell\partial_\xi)f^{-s}(\ell,\xi)f^s(m,\zeta)
 +\overline{f^s(k,\eta)}f^{-s}(\ell,\xi)
                   (\partial_t-m\partial_\zeta)f^s(m,\zeta).
\end{align*}
Thus substituting \eqref{eq:scalar} of Lemma~\ref{lem:scalar} into \eqref{ddbt-1} gives
\begin{align}\label{eq:Bseventerms}
\frac{\dd}{\dd t}\mathcal B(t)=I_1+I_2+I_3+I_4+I_5+I_6+I_7.	
\end{align}
These seven quantities are defined explicitly as follows. 

\noindent\emph{The coefficient terms \(I_1\).}
\begin{align}
I_1={}&\sum_s\Rea\sum_{\substack{\ell,m\in\Z\\\ell\ne0}}\int_{\R^{2}}\frac{\ii}{4s\beta\ell}
  (\partial_t-\ell\partial_{\xi}-m\partial_{\zeta})
  \bigl[e^{\frac18\rho^{1/3} t\mathbf 1_{k\ne0}\mathbf 1_{m\ne0}}(h(k,\eta)-h(m,\zeta))a(k,\eta;\ell,\xi)\bigr]\nonumber\\
&\qquad\times\overline{f^s(k,\eta)}f^{-s}(\ell,\xi)f^s(m,\zeta)
                       \frac{\dd\xi\dd\zeta}{(2\pi)^{2}}.
\label{eq:R3}
\end{align}

\noindent\emph{The Alfv\'en term \(I_2\).}
\begin{align}
 I_2={}&\sum_s\Rea\sum_{\substack{\ell,m\in\Z\\\ell\ne0}}
 \int_{\R^2}\frac{\ii e^{\rho^{1/3}t\mathbf1_{k\ne0}\mathbf1_{m\ne0}/8}
 (h(k,\eta)-h(m,\zeta))a(k,\eta;\ell,\xi)}{4s\beta\ell}(-s\ii\beta k-s\ii\beta\ell+s\ii\beta m)
 \nonumber\\
 &\qquad\times
 \overline{f^s(k,\eta)}f^{-s}(\ell,\xi)f^s(m,\zeta)
 \frac{\dd\xi\dd\zeta}{(2\pi)^2}\nonumber\\
 ={}&\frac12\sum_s\Rea\sum_{\substack{\ell,m\in\Z\\\ell\ne0}}
 \int_{\R^2}e^{\rho^{1/3}t\mathbf1_{k\ne0}\mathbf1_{m\ne0}/8}
 (h(k,\eta)-h(m,\zeta))a(k,\eta;\ell,\xi)
 \overline{f^s(k,\eta)}f^{-s}(\ell,\xi)f^s(m,\zeta)
 \frac{\dd\xi\dd\zeta}{(2\pi)^2}
 \label{eq:Ialfven}	
\end{align}
where we used $k=\ell+m$. This is exactly \eqref{eq:cancelledcubic}. Thus \(I_2\) cancels that term of \(E'\) when we compute \((E-\mathcal B)'\).

\noindent\emph{The shear terms \(I_3\).}
\begin{equation}
\begin{aligned}
I_3={}&\sum_s\Rea\sum_{\substack{\ell,m\in\Z\\\ell\ne0}}\int_{\R^{2}}\frac{\ii e^{\frac18\rho^{1/3} t\mathbf 1_{k\ne0}\mathbf 1_{m\ne0}}(h(k,\eta)-h(m,\zeta))a(k,\eta;\ell,\xi)}{4s\beta\ell}\\*
&\times\bigg[\frac{k\eta}{(k^2+\eta^2)}\overline{f^{-s}(k,\eta)}f^{-s}(\ell,\xi)f^s(m,\zeta) +\frac{\ell\xi}{(\ell^2+\xi^2)}\overline{f^s(k,\eta)}f^s(\ell,\xi)f^s(m,\zeta)\\*
&\qquad+\frac{m\zeta}{(m^2+\zeta^2)}\overline{f^s(k,\eta)}f^{-s}(\ell,\xi)f^{-s}(m,\zeta)\bigg]\frac{\dd\xi\dd\zeta}{(2\pi)^{2}}.
\label{eq:shearterm}
\end{aligned}
\end{equation}

\noindent\emph{The diffusion term \(I_4\).}
\begin{align}
I_4=&-\sum_s\Rea\sum_{\substack{\ell,m\in\Z\\\ell\ne0}}\int_{\R^{2}}
 \frac{\ii e^{\rho^{1/3}t\mathbf 1_{k\ne0}\mathbf 1_{m\ne0}/8}(h(k,\eta)-h(m,\zeta))a(k,\eta;\ell,\xi)}
      {4s\beta\ell}\nonumber\\
&\times\bigg[\frac{\nu+\mu}{2}((k^2+\eta^2)+(\ell^2+\xi^2)+(m^2+\zeta^2))\overline{f^s(k,\eta)}f^{-s}(\ell,\xi)f^s(m,\zeta)\nonumber\\
&\qquad+\frac{\nu-\mu}{2}\bigl((k^2+\eta^2)\overline{f^{-s}(k,\eta)}f^{-s}(\ell,\xi)f^s(m,\zeta) +(\ell^2+\xi^2)\overline{f^s(k,\eta)}f^s(\ell,\xi)f^s(m,\zeta)\nonumber\\
&\qquad+(m^2+\zeta^2)\overline{f^s(k,\eta)}f^{-s}(\ell,\xi)f^{-s}(m,\zeta)\bigr) \bigg]\frac{\dd\xi\dd\zeta}{(2\pi)^{2}}.
\label{eq:Rnu}
\end{align}
The remaining terms are
\begin{align}
 I_5={}&\sum_s\Rea\sum_{\substack{\ell,m\in\Z\\\ell\ne0}}
 \int_{\R^2}\frac{\ii e^{\rho^{1/3}t\mathbf1_{k\ne0}\mathbf1_{m\ne0}/8}
 (h(k,\eta)-h(m,\zeta))a(k,\eta;\ell,\xi)}{4s\beta\ell}
 \nonumber\\
 &\quad\times\overline{\sum_{p\in\Z}\int_\R
 a(k,\eta;p,\gamma)f^{-s}(p,\gamma)f^s(k-p,\eta-\gamma)
 \frac{\dd\gamma}{2\pi}}f^{-s}(\ell,\xi)f^s(m,\zeta)
 \frac{\dd\xi\dd\zeta}{(2\pi)^2},\nonumber\\
 I_6={}&\sum_s\Rea\sum_{\substack{\ell,m\in\Z\\\ell\ne0}}
 \int_{\R^2}\frac{\ii e^{\rho^{1/3}t\mathbf1_{k\ne0}\mathbf1_{m\ne0}/8}
 (h(k,\eta)-h(m,\zeta))a(k,\eta;\ell,\xi)}{4s\beta\ell}
 \nonumber\\
 &\quad\times\overline{f^s(k,\eta)}
 \bigg(\sum_{p\in\Z}\int_\R a(\ell,\xi;p,\gamma)
 f^s(p,\gamma)f^{-s}(\ell-p,\xi-\gamma)
 \frac{\dd\gamma}{2\pi}\bigg)f^s(m,\zeta)\frac{\dd\xi\dd\zeta}{(2\pi)^2},
 \nonumber\\
 I_7={}&\sum_s\Rea\sum_{\substack{\ell,m\in\Z\\\ell\ne0}}
 \int_{\R^2}\frac{\ii e^{\rho^{1/3}t\mathbf1_{k\ne0}\mathbf1_{m\ne0}/8}
 (h(k,\eta)-h(m,\zeta))a(k,\eta;\ell,\xi)}{4s\beta\ell}
 \nonumber\\
 &\quad\times\overline{f^s(k,\eta)}f^{-s}(\ell,\xi)
 \bigg(\sum_{p\in\Z}\int_\R a(m,\zeta;p,\gamma)
 f^{-s}(p,\gamma)f^s(m-p,\zeta-\gamma)
 \frac{\dd\gamma}{2\pi}\bigg)\frac{\dd\xi\dd\zeta}{(2\pi)^2}.
 \label{eq:R4native}
\end{align}

We first rewrite $I_5$ and $I_7$. For $I_5$, we rename the original \((m,\zeta)\) as \((j,\alpha)\), and replace the original inner frequency \((p,\gamma)\) by \((-p,-\gamma)\). The frequency \((k,\eta)=(\ell+m,\xi+\zeta)\) then becomes \((j+\ell,\alpha+\xi)\). 
The amplitude product is thus converted to
\[
\begin{aligned}
\overline{f^{-s}(-p,-\gamma)}
 \,\overline{f^s(j+\ell+p,\alpha+\xi+\gamma)}
 f^{-s}(\ell,\xi)f^s(j,\alpha)
 =\overline{f^s(j+\ell+p,\alpha+\xi+\gamma)}f^s(j,\alpha)
 f^{-s}(\ell,\xi)f^{-s}(p,\gamma).
\end{aligned}
\]
Here we used \(f^{-s}(-p,-\gamma)=\overline{f^{-s}(p,\gamma)}\). The identity \eqref{eq:endpointsym} in Lemma~\ref{lem:kernelidentities} gives
\[
a(j+\ell,\alpha+\xi;-p,-\gamma)=a(j+\ell+p-p,\alpha+\xi+\gamma-\gamma;-p,-\gamma)
=-a(j+\ell+p,\alpha+\xi+\gamma;p,\gamma).
\]
Hence \(I_5\) becomes
\[
\begin{aligned}
I_5=&\sum_s\Rea\sum_{\substack{j,\ell,p\in\Z\\\ell\ne0}}\int_{\R^{3}} \frac{\ii e^{\frac18\rho^{1/3} t\mathbf 1_{j+\ell\ne0}\mathbf 1_{j\ne0}}(h(j+\ell,\alpha+\xi)-h(j,\alpha))}{4s\beta\ell}a(j+\ell,\alpha+\xi;-p,-\gamma)\\
&\quad \times a(j+\ell,\alpha+\xi;\ell,\xi)\overline{f^s(j+\ell+p,\alpha+\xi+\gamma)}f^s(j,\alpha)f^{-s}(\ell,\xi)f^{-s}(p,\gamma)\frac{\dd\alpha\dd\xi\dd\gamma}{(2\pi)^{3}}.
\end{aligned}
\]
Redefining
\(
k=j+\ell+p,\eta=\alpha+\xi+\gamma,
\)
we get
\[
\begin{aligned}
I_5=&-\sum_s\Rea\sum_{\substack{j,\ell,p\in\Z\\\ell\ne0}}\int_{\R^{3}} \frac{\ii e^{\frac18\rho^{1/3} t\mathbf 1_{j+\ell\ne0}\mathbf 1_{j\ne0}}(h(j+\ell,\alpha+\xi)-h(j,\alpha))}{4s\beta\ell}\\
&\quad \times a(k,\eta;p,\gamma)a(j+\ell,\alpha+\xi;\ell,\xi)\overline {f^s(k,\eta)}f^s(j,\alpha)f^{-s}(\ell,\xi)f^{-s}(p,\gamma)\frac{\dd\alpha\dd\xi\dd\gamma}{(2\pi)^{3}}.
\end{aligned}
\]
For $I_7$, we set
\(
j=m-p,\alpha=\zeta-\gamma.
\)
Thus \(m=j+p\), \(\zeta=\alpha+\gamma\), and 
\[
k=j+\ell+p,\qquad\eta=\alpha+\xi+\gamma.
\]
Consequently,
\[
\begin{aligned}
I_7={}&\sum_s\Rea\sum_{\substack{j,\ell,p\in\Z\\\ell\ne0}}
 \int_{\R^3}\frac{\ii}{4s\beta\ell}
 e^{\rho^{1/3}t\mathbf 1_{k\ne0}\mathbf 1_{j+p\ne0}/8} \bigl(h(k,\eta)-h(j+p,\alpha+\gamma)\bigr)\\
&\quad\times
 a(k,\eta;\ell,\xi)a(j+p,\alpha+\gamma;p,\gamma)
 \overline{f^s(k,\eta)}f^s(j,\alpha)
 f^{-s}(\ell,\xi)f^{-s}(p,\gamma)
                     \frac{\dd\alpha\dd\xi\dd\gamma}{(2\pi)^3}.
\end{aligned}
\]
Adding the displayed expressions for \(I_5\) and \(I_7\) gives
\begin{align}
I_5+I_7&=\sum_s\Rea\sum_{\substack{j,\ell,p\in\Z\\\ell\ne0}}\int_{\R^{3}}\frac{\ii}{4s\beta\ell}\overline {f^s(k,\eta)}f^s(j,\alpha)f^{-s}(\ell,\xi)f^{-s}(p,\gamma)\nonumber\\
 &\quad\times\Bigl[e^{\frac18\rho^{1/3} t\mathbf 1_{k\ne0}\mathbf 1_{j+p\ne0}}(h(k,\eta)-h(j+p,\alpha+\gamma))
 a(k,\eta;\ell,\xi)a(j+p,\alpha+\gamma;p,\gamma)\nonumber\\
 &\qquad-e^{\frac18\rho^{1/3} t\mathbf 1_{j+\ell\ne0}\mathbf 1_{j\ne0}}(h(j+\ell,\alpha+\xi)-h(j,\alpha))
 a(k,\eta;p,\gamma)a(j+\ell,\alpha+\xi;\ell,\xi)\Bigr]\frac{\dd\alpha\dd\xi\dd\gamma}{(2\pi)^{3}}.
\label{eq:pair}
\end{align}
Here \(k=j+\ell+p\), \(\eta=\alpha+\xi+\gamma\) in \eqref{eq:pair}.

We next consider $I_6$. In \eqref{eq:R4native}, set 
\(
q=\ell-p,\quad\delta=\xi-\gamma
\).
Then
\[
\ell=p+q,\quad \xi=\gamma+\delta,\qquad
k=p+q+m,\quad\eta=\gamma+\delta+\zeta.
\]
The constraint \(\ell\ne0\) becomes \(p+q\ne0\). This gives
\begin{equation}
\begin{split}
 I_6={}&\sum_s\Rea\sum_{\substack{p,q,m\in\Z\\p+q\ne0}}\int_{\R^{3}}
 \frac{\ii e^{\frac18\rho^{1/3} t\mathbf 1_{k\ne0}\mathbf 1_{m\ne0}}(h(k,\eta)-h(m,\zeta))}{4s\beta(p+q)}
 a(k,\eta;p+q,\gamma+\delta)\\
 &\qquad\times a(p+q,\gamma+\delta;p,\gamma)
       \overline{f^s(k,\eta)}f^s(p,\gamma)f^{-s}(q,\delta)f^s(m,\zeta)\frac{\dd\gamma\dd\delta\dd\zeta}{(2\pi)^{3}}.
\end{split}
\label{eq:adv}
\end{equation}

\subsection{The identity for $E-\mathcal B$}

Combining the preceding identities, we obtain the following proposition.

\begin{proposition}\label{prop:identity}
With the notation of Section~\ref{sec:equations}, every smooth solution $(u,b)$ of \eqref{eq:perturbation} on $[0,T_1]$ satisfies
\begin{align}
&\frac{\mathrm d}{\mathrm dt}(E-\mathcal B) +2\nu\|\nabla\mathcal A\widetilde u\|_2^2
                         +2\mu\|\nabla\mathcal Ab\|_2^2
+2\sum_k\int_\R(7g+(\partial_t-k\partial_\eta)V)\mathcal A^2
                \big(|\widehat{\widetilde u}|^2+|\widehat b|^2\big)\frac{\dd\eta}{2\pi}\nonumber\\
 &=\frac{\rho^{1/3}}8
       \big(\|\mathcal A\widetilde u_{\ne}\|_2^2+\|\mathcal Ab_{\ne}\|_2^2\big)
   +4\langle\mathcal A\widetilde u,\mathcal A\partial_x\nabla\Delta^{-1}u^2\rangle
-\frac{2(\mu-\nu)}\beta
       \langle\nabla\mathcal A\widetilde u^1,
                    \nabla\mathcal A\partial_x^{-1}b_{\ne}^2\rangle\nonumber\\
 &\quad+2\langle\mathcal A(\widetilde u-u),
                  \mathcal A\mathbb P(b\cdot\nabla b-u\cdot\nabla u)\rangle
       +\frac2\beta\langle\mathcal A\widetilde u^1,
                         \mathcal A\Pi_{\ne}(b^1u^2-u^1b^2)\rangle\nonumber\\
 &\quad+\frac12\sum_s\Rea\sum_m\int_{\R^2}
       (\mathcal A^2(k,\eta)-\mathcal A^2(m,\zeta))a(k,\eta;0,\xi)
       \overline{f^s(k,\eta)}f^{-s}(0,\xi)f^s(m,\zeta)
                                \frac{\dd\xi\dd\zeta}{(2\pi)^2}\nonumber\\
 &\quad+\frac12\sum_s\Rea\sum_{\ell\ne0,m\in\Z}\int_{\R^2}
       [\mathcal A^2(k,\eta)-\mathcal A^2(m,\zeta)
          -e^{\rho^{1/3}t\mathbf 1_{k\ne0}\mathbf 1_{m\ne0}/8}(h(k,\eta)-h(m,\zeta))]\nonumber\\
 &\hspace{35mm}\times a(k,\eta;\ell,\xi)
       \overline{f^s(k,\eta)}f^{-s}(\ell,\xi)f^s(m,\zeta)
                                \frac{\dd\xi\dd\zeta}{(2\pi)^2}\nonumber\\
 &\quad-\sum_s\Rea\sum_{\ell\neq0,m\in\Z}\int_{\R^{2}}\frac{\ii}{4s\beta\ell}
  (\partial_t-\ell\partial_{\xi}-m\partial_{\zeta})
  \big[e^{\frac18\rho^{1/3} t\mathbf 1_{k\ne0}\mathbf 1_{m\ne0}}(h(k,\eta)-h(m,\zeta))a(k,\eta;\ell,\xi)\big]\nonumber\\
&\hspace{35mm}\times\overline{f^s(k,\eta)}f^{-s}(\ell,\xi)f^s(m,\zeta)\frac{\dd\xi\dd\zeta}{(2\pi)^{2}}\nonumber\\
&\quad-\sum_s\Rea\sum_{\ell\neq0,m\in\Z}\int_{\R^{2}}\frac{\ii e^{\frac18\rho^{1/3} t\mathbf 1_{k\ne0}\mathbf 1_{m\ne0}}(h(k,\eta)-h(m,\zeta))a(k,\eta;\ell,\xi)}{4s\beta\ell}\nonumber\\
&\hspace{10mm}\times\bigg[\frac{k\eta}{(k^2+\eta^2)}\overline{f^{-s}(k,\eta)}f^{-s}(\ell,\xi)f^s(m,\zeta) +\frac{\ell\xi}{(\ell^2+\xi^2)}\overline{f^s(k,\eta)}f^s(\ell,\xi)f^s(m,\zeta)\nonumber\\
&\hspace{16mm}+\frac{m\zeta}{(m^2+\zeta^2)}\overline{f^s(k,\eta)}f^{-s}(\ell,\xi)f^{-s}(m,\zeta)\bigg]\frac{\dd\xi\dd\zeta}{(2\pi)^{2}}\nonumber\\
&\quad+\sum_s\Rea\sum_{\ell\neq0,m\in\Z}\int_{\R^{2}}
 \frac{\ii e^{\rho^{1/3}t\mathbf 1_{k\ne0}\mathbf 1_{m\ne0}/8}(h(k,\eta)-h(m,\zeta))a(k,\eta;\ell,\xi)}
      {4s\beta\ell}\nonumber\\
&\hspace{10mm}\times\bigg[\frac{\nu+\mu}{2}((k^2+\eta^2)+(\ell^2+\xi^2)+(m^2+\zeta^2))
                       \overline{f^s(k,\eta)}f^{-s}(\ell,\xi)f^s(m,\zeta)\nonumber\\
&\hspace{16mm}+\frac{\nu-\mu}{2}\bigl( (k^2+\eta^2)\overline{f^{-s}(k,\eta)}f^{-s}(\ell,\xi)f^s(m,\zeta) +(\ell^2+\xi^2)\overline{f^s(k,\eta)}f^s(\ell,\xi)f^s(m,\zeta)\nonumber\\
&\hspace{16mm}+(m^2+\zeta^2)\overline{f^s(k,\eta)}f^{-s}(\ell,\xi)f^{-s}(m,\zeta)\bigr) \bigg]\frac{\dd\xi\dd\zeta}{(2\pi)^{2}}\nonumber\\
&\quad-\sum_s\Rea\sum_{\substack{j,\ell,p\in\Z\\\ell\ne0}}\int_{\R^{3}}\frac{\ii}{4s\beta\ell}\overline{f^s(k,\eta)}f^s(j,\alpha)f^{-s}(\ell,\xi)f^{-s}(p,\gamma) \nonumber\\
&\hspace{10mm}\times \Bigl[e^{\frac18\rho^{1/3} t\mathbf 1_{k\ne0}\mathbf 1_{j+p\ne0}}(h(k,\eta)-h(j+p,\alpha+\gamma))
                   a(k,\eta;\ell,\xi)a(j+p,\alpha+\gamma;p,\gamma)\nonumber\\
 &\hspace{16mm}-e^{\frac18\rho^{1/3} t\mathbf 1_{j+\ell\ne0}\mathbf 1_{j\ne0}}(h(j+\ell,\alpha+\xi)-h(j,\alpha))a(k,\eta;p,\gamma)a(j+\ell,\alpha+\xi;\ell,\xi)\Bigr]\frac{\dd\alpha\dd\xi\dd\gamma}{(2\pi)^{3}}\nonumber\\
&\quad-\sum_s\Rea\sum_{\substack{p,q,m\in\Z\\p+q\ne0}}\int_{\R^{3}}
\frac{\ii e^{\frac18\rho^{1/3} t\mathbf 1_{k\ne0}\mathbf 1_{m\ne0}}(h(k,\eta)-h(m,\zeta))}{4s\beta(p+q)} a(k,\eta;p+q,\gamma+\delta)\nonumber\\
&\hspace{10mm}\times a(p+q,\gamma+\delta;p,\gamma)
       \overline{f^s(k,\eta)}f^s(p,\gamma)f^{-s}(q,\delta)f^s(m,\zeta)\frac{\dd\gamma\dd\delta\dd\zeta}{(2\pi)^{3}}.
\label{eq:identity}
\end{align}
In the double integrals, \(k=\ell+m\) and \(\eta=\xi+\zeta\), with \(\ell=0\) in the first one.
In the triple integral over \((j,\alpha),(\ell,\xi),(p,\gamma)\),
\(k=j+\ell+p\) and \(\eta=\alpha+\xi+\gamma\).
In the last integral, \(k=p+q+m\) and \(\eta=\gamma+\delta+\zeta\).
\end{proposition}

\begin{proof}
By \eqref{eq:Bseventerms}, subtracting \(\mathcal B'\) from
\(E'\) subtracts exactly \(I_1,\ldots,I_7\).
Substituting the decomposition
\eqref{eq:T0}--\eqref{eq:cancelledcubic} into \eqref{eq:Eprime},
we find that its last term equals \(I_2\) and hence cancels.
The terms \(-I_1\), \(-I_3\), \(-I_4\), \(-(I_5+I_7)\), and \(-I_6\)
are the negatives of \eqref{eq:R3}, \eqref{eq:shearterm}, \eqref{eq:Rnu},
\eqref{eq:pair}, and \eqref{eq:adv}, respectively.
All the other terms of \eqref{eq:Eprime} remain as displayed.
\end{proof}

The terms on the right-hand side of \eqref{eq:identity} are estimated in Sections~\ref{sec:cubic}--\ref{sec:bootstrap}.

\section{Fourier multiplier estimates}
\label{sec:multiplier}

In this section, we establish the Fourier multiplier estimates used in the subsequent sections.

\subsection{Estimates for $g$ and $\kappa$}

\begin{lemma}\label{lem:shift}
The multipliers in \eqref{eq:frequencies}, \eqref{eq:g}, and \eqref{eq:weight} satisfy, uniformly in $t\ge0$ and $0<\rho\le1$:
\begin{align}
 g(k+\ell,\eta+\xi)&\le C(1+\ell^2)^4(1+(\xi+\ell t)^2)g(k,\eta),
 \label{eq:shiftgseparated}\\
 g(k+\ell,\eta+\xi)&\le C\Lambda_t^{10}(\ell,\xi)g(k,\eta),\label{eq:shiftg}\\
 \kappa(k+\ell,\eta+\xi)&\le C\Lambda_t^{5}(\ell,\xi)\sqrt{g(k,\eta)},\label{eq:shiftkappa}\\
 \av t\kappa(\ell,\xi)&\le C\Lambda_t(\ell,\xi)\quad(\ell\ne0),\label{eq:tkappa}\\
 \frac{\av t\kappa(\ell,\xi)}{|\ell|}
       (1+\ell^2)^2\sqrt{1+(\xi+\ell t)^2}
 &\le C\Lambda_t^{4}(\ell,\xi)\quad(\ell\ne0),
 \label{eq:tkappaphase}\\
 \av t\kappa(\ell,\xi)(1+\ell^2)^2\sqrt{1+(\xi+\ell t)^2}
 &\le C\Lambda_t^{5}(\ell,\xi)\quad(\ell\ne0).
 \label{eq:tkappatransfer}
\end{align}
\end{lemma}
\begin{proof}
For each $n\ne0$,
\begin{align*}
1+(k-n)^2&\le2(1+(k+\ell-n)^2)(1+\ell^2),\\
1+\Big(t-\frac{\eta+kt}{n}\Big)^2
&\le2\bigg[1+\Big(t-\frac{\eta+kt+\xi+\ell t}{n}\Big)^2\bigg]
            \bigg[1+\Big(\frac{\xi+\ell t}{n}\Big)^2\bigg].
\end{align*}
Since $|n|\ge1$, substitution in the definition of $g$ gives
\begin{align*}
g(k+\ell,\eta+\xi)
=\sum_{n\ne0}\frac{(1+(k+\ell-n)^2)^{-4}}{1+(t-(\eta+kt+\xi+\ell t)/n)^2}\le C(1+\ell^2)^4(1+(\xi+\ell t)^2)
 \sum_{n\ne0}\frac{(1+(k-n)^2)^{-4}}
 {1+(t-(\eta+kt)/n)^2}.
\end{align*}
Hence it holds that
\begin{align*} 
g(k+\ell,\eta+\xi)\le C\Lambda_t^{10}(\ell,\xi)g(k,\eta).
\end{align*}
These are \eqref{eq:shiftgseparated} and \eqref{eq:shiftg}.
For $k+\ell\ne0$, the summand $n=k+\ell$ gives
\[
\kappa^2(k+\ell,\eta+\xi)
\le g(k+\ell,\eta+\xi)
\le C\Lambda_t^{10}(\ell,\xi)g(k,\eta),
\]
which proves \eqref{eq:shiftkappa}. For $k+\ell=0$, its left-hand side is zero.
For the three bounds containing $\kappa(\ell,\xi)$, we can use \eqref{eq:frequencies} to compute
\[
 \av t\kappa(\ell,\xi)
 \le C\sqrt{1+((\xi+\ell t)/\ell)^2}
 \le C\frac{\Lambda_t(\ell,\xi)}{|\ell|}.
\]
It proves \eqref{eq:tkappa}. Since $|\ell|\ge1$,
\[
 \frac{(1+\ell^2)^2}{\ell^2}\le2(1+\ell^2),
 \qquad
 \frac{(1+\ell^2)^2}{|\ell|}
             \le\sqrt2(1+\ell^2)^{3/2}.
\]
Multiplication by $\sqrt{1+(\xi+\ell t)^2}$ proves
\eqref{eq:tkappaphase} and \eqref{eq:tkappatransfer}.
\end{proof}

\subsection{First and second differences of the Fourier multiplier}

We first estimate $G$ and $V$, then pass to $\theta=\exp\{-2(7G+V)\}$ and $h=\theta\Lambda_t^{2N}$.
Let $k=j+\ell+p$, $\eta=\alpha+\xi+\gamma$, and assume
\begin{equation}
8\big[\sqrt{\ell^2+(\xi+\ell t)^2}+\sqrt{p^2+(\gamma+pt)^2}\big]\le \sqrt{j^2+(\alpha+jt)^2}.
\label{eq:rectangle}
\end{equation}
We shall also use the following condition:
\begin{equation}
 t\ge1,\qquad |j|>16(|\ell|+|p|+1),\qquad
                   t/8\le(\alpha+jt)/j\le8t.
\label{eq:criticalregion}
\end{equation}
For $0\le s,r\le1$, the reverse triangle inequality and
\eqref{eq:rectangle} give
\[
\begin{aligned}
&|\Lambda_t(j+s\ell+rp,\alpha+s\xi+r\gamma)
                    -\Lambda_t(j,\alpha)|
\le s\sqrt{\ell^2+(\xi+\ell t)^2}
       +r\sqrt{p^2+(\gamma+pt)^2}
 \le\frac18\Lambda_t(j,\alpha).
\end{aligned}
\]
Consequently,
\begin{equation}
 \frac78\Lambda_t(j,\alpha)
 \le\Lambda_t(j+s\ell+rp,\alpha+s\xi+r\gamma)
 \le\frac98\Lambda_t(j,\alpha),\qquad 0\le s,r\le1.
 \label{eq:referencecomparison}
\end{equation}

\begin{lemma}\label{lem:scalarprimitives}
For $t\ge0$, $a\in\mathbb R$, $0<\sigma\le1$, and $0<\rho\le1$,  the weight $\Phi_\sigma(t,a)$ defined in \eqref{eq:Phi} satisfies
\begin{equation}
 0\le\Phi_\sigma\le\pi,\qquad
 |a|\Phi_\sigma\le C\av t,\qquad
 a^2\Phi_1\le C\av t^2.
\label{eq:prim}
\end{equation}
Moreover, we have
\begin{equation}\label{eq:EDderiv}
\begin{aligned}
 (1+a^2)|\partial_a\Phi_{\rho^{1/3}}|
     &\le C\av t\av{\rho^{1/3}t},\\
 (1+a^2)^{2}|\partial_a^2\Phi_{\rho^{1/3}}|
   +(1+a^2)^{3/2}|\partial_a\Phi_{\rho^{1/3}}|
     &\le C\rho^{-1/3}\av t\av{\rho^{1/3}t}^3,\\
 a^2\Phi_{\rho^{1/3}}(t,a)
       & \le C\rho^{-1/3}t\av{\rho^{1/3}t}.    
\end{aligned}
\end{equation}
\end{lemma}
\begin{proof}
The inequality $0\le\Phi_\sigma\le\pi$ follows from \eqref{eq:Phi}. Then by \eqref{eq:Phi} again,
\begin{equation}\label{eq:3}
\begin{aligned}
 \partial_a\Phi_\sigma(t,a)
 &=\frac{\sigma}{1+\sigma^2a^2}
                   -\frac{\sigma}{1+\sigma^2(t-a)^2}
 =\int_0^t\frac{2\sigma^3(s-a)}{(1+\sigma^2(s-a)^2)^2}\,\mathrm ds,\\
 \partial_a^2\Phi_\sigma(t,a)
 &=-\frac{2\sigma^3a}{(1+\sigma^2a^2)^2}
   -\frac{2\sigma^3(t-a)}{(1+\sigma^2(t-a)^2)^2}
 =\int_0^t\frac{6\sigma^5(s-a)^2-2\sigma^3}
                  {(1+\sigma^2(s-a)^2)^3}\,\mathrm ds.
\end{aligned}
\end{equation}
If $|a|\le2t+2$, then the remaining two bounds in \eqref{eq:prim} are clear. If $|a|>2t+2$, then $|s-a|\ge|a|/2$ for $0\le s\le t$, and
\[
 |a|\Phi_\sigma(t,a)
 \le Ct\frac{\sigma|a|}{1+\sigma^2a^2}\le Ct.
\]
With $\sigma=1$ the same argument gives $a^2\Phi_1\le Ct$ in this region, hence the third bound in \eqref{eq:prim} holds. 

For \eqref{eq:EDderiv}, we use the identity
\begin{equation}\label{scalesharp}
\Phi_{\rho^{1/3}}(t,a)=\Phi_1(\rho^{1/3}t,\rho^{1/3}a).	
\end{equation}
For $0\le s\le t\le1$, due to \( 1+a^2 \le 2(1+(a-s)^2)(1+s^2) \le 4(1+(a-s)^2)\), we have 
\begin{align*}
 (1+a^2)|\partial_a\Phi_1(t,a)|
 &\le
 8\int_0^t\frac{|s-a|}{1+(s-a)^2}\,\mathrm ds
 \le 4t,\\
 (1+a^2)^{3/2}|\partial_a\Phi_1(t,a)|
 &\le
 16\int_0^t\frac{|s-a|}{\sqrt{1+(s-a)^2}}\,\mathrm ds
 \le16t,\\
 (1+a^2)^2|\partial_a^2\Phi_1(t,a)|
 &\le
 16\int_0^t
 \frac{|6(s-a)^2-2|}{1+(s-a)^2}\,\mathrm ds
 \le96t.
\end{align*}
For $t\ge1$ and $|a|\le2t+2$, \eqref{eq:3} gives \(|\partial_a\Phi_1(t,a)|\le1\), and \(|\partial_a^2\Phi_1(t,a)|\le C\). Then it follows that
\begin{equation}\label{eq:bound1}
\begin{split}
 (1+a^2)|\partial_a\Phi_1(t,a)|
 &\le C(1+t^2)
 \le Ct\av t,\\
 (1+a^2)^2|\partial_a^2\Phi_1(t,a)|
 +(1+a^2)^{3/2}|\partial_a\Phi_1(t,a)|
 &\le C\bigl((1+t^2)^2+(1+t^2)^{3/2}\bigr)
 \le Ct\av t^3.
\end{split}
\end{equation}
For $t\ge1$ and $|a|\ge2t+2>1$, we have \(|s-a|\ge |a|-t>\frac{|a|}{2}>1\). Thus 
\begin{align}\label{eq:17}
 |\partial_a\Phi_1(t,a)|
 \le2\int_0^t|s-a|^{-3}\,\mathrm ds
 \le16t|a|^{-3},\quad
 |\partial_a^2\Phi_1(t,a)|
 \le6\int_0^t|s-a|^{-4}\,\mathrm ds
 \le96t|a|^{-4}.
\end{align}
Due to $|a|\ge2t+2>1$, we obtain that \eqref{eq:bound1} also holds true for this case. 
Since $1+a^2\le\rho^{-2/3}(1+\rho^{2/3}a^2)$, using \eqref{scalesharp} and then applying \eqref{eq:bound1} at $(\rho^{1/3}t,\rho^{1/3}a)$ yields
\[
\begin{aligned}
(1+a^2)|\partial_a\Phi_{\rho^{1/3}}(t,a)|
 &\le Ct\av{\rho^{1/3}t},\\
(1+a^2)^2|\partial_a^2\Phi_{\rho^{1/3}}(t,a)|
 +(1+a^2)^{3/2}|\partial_a\Phi_{\rho^{1/3}}(t,a)|
 &\le C\rho^{-1/3}t\av{\rho^{1/3}t}^3.
\end{aligned}
\]
For $t\le1$, $\Phi_1\le t$ and the preceding split in $a$ give $a^2\Phi_1(t,a)\le Ct$. For $t\ge1$, \eqref{eq:prim} gives $a^2\Phi_1(t,a)\le Ct\av t$. Thus the last inequality of \eqref{eq:EDderiv} follows from \eqref{scalesharp}.
\end{proof}

\begin{lemma}\label{lem:primitivesizes}
With the multipliers defined in \eqref{eq:g} and \eqref{eq:weight}, for every $j\in\Z$ and $\alpha\in\R$, it holds that
\begin{equation}
\begin{aligned}
\Lambda_t(j,\alpha)G(j,\alpha)
 &\le C\sqrt{1+j^2}\av t,\\
\Lambda_t^2(j,\alpha)G(j,\alpha)
 &\le C(1+j^2)\av t^2.
\end{aligned}
\label{eq:smallhorizontalG}
\end{equation}
and
\begin{equation}
\begin{aligned}
\Lambda_t(j,\alpha)V(j,\alpha)
 &\le C\sqrt{1+j^2}\av t,\\
\Lambda_t^2(j,\alpha)V(j,\alpha)
 &\le C(1+j^2)\rho^{-1/3}\av t\av{\rho^{1/3}t}.
\end{aligned}
\label{eq:smallhorizontalV}
\end{equation}
Under \eqref{eq:rectangle}, we also have 
\begin{equation}
\begin{aligned}
\Lambda_t(j,\alpha)\bigl(
 |G(j+\ell,\alpha+\xi)-G(j,\alpha)|
 +|V(j+\ell,\alpha+\xi)-V(j,\alpha)|\bigr)
\le C\sqrt{1+(|j|+|j+\ell|)^2}\av t.
\end{aligned}
\label{eq:boundedhorizontalprimitive}
\end{equation}
\end{lemma}
\begin{proof}
For $j\in\mathbb Z$, we first estimate $G(j,\alpha)$. For $r=0,1,2$, due to $1+n^2\le2(1+j^2)(1+(j-n)^2)$, we have
\begin{equation*}
 \sum_{n\ne0}(1+(j-n)^2)^{-4}(1+n^2)^{r/2}
 \le C_r(1+j^2)^{r/2}.
\end{equation*}
By \eqref{eq:prim},
\[
\begin{aligned}
|\alpha+jt|\Phi_1(t,(\alpha+jt)/n)&\le C|n|\av t,\\
(\alpha+jt)^2\Phi_1(t,(\alpha+jt)/n)&\le Cn^2\av t^2.
\end{aligned}
\]
Substitution in \eqref{eq:physicalG} proves \eqref{eq:smallhorizontalG}.

To prove \eqref{eq:smallhorizontalV}, recalling \eqref{eq:weight}, we have $V(j,\alpha)=\Phi_{\rho^{1/3}}(t,a)$. For $j\ne0$, apply \eqref{eq:prim} and \eqref{eq:EDderiv} of Lemma~\ref{lem:scalarprimitives} with $a=(\alpha+jt)/j$. It follows that
\[
 |\alpha+jt|V(j,\alpha)\le C|j|\av t,
 \qquad
 (\alpha+jt)^2V(j,\alpha)
       \le Cj^2\rho^{-1/3}t\av{\rho^{1/3}t}.
\]
Together with $V\le\pi$, we get \eqref{eq:smallhorizontalV}. For $j=0$, the definition $V(0,\alpha)=0$ gives the same conclusion.

To prove \eqref{eq:boundedhorizontalprimitive}, we deduce from $G,V\ge0$, the first inequalities in \eqref{eq:smallhorizontalG} and \eqref{eq:smallhorizontalV}, and \eqref{eq:referencecomparison} that
\begin{align*}
&\Lambda_t(j,\alpha)\bigl(
 |G(j+\ell,\alpha+\xi)-G(j,\alpha)|
 +|V(j+\ell,\alpha+\xi)-V(j,\alpha)|\bigr)\\
\le{}&\Lambda_t(j,\alpha)\bigl[
 G(j+\ell,\alpha+\xi)+G(j,\alpha)
 +V(j+\ell,\alpha+\xi)+V(j,\alpha)\bigr]\\
\le{}& C\av t\bigg[
 \frac{\Lambda_t(j,\alpha)}{\Lambda_t(j+\ell,\alpha+\xi)}
          \sqrt{1+(j+\ell)^2}+\sqrt{1+j^2}\bigg]\\
\le{}& C\av t
       \bigl[\sqrt{1+(j+\ell)^2}+\sqrt{1+j^2}\bigr]\\
\le{}& C\sqrt{1+(|j|+|j+\ell|)^2}\av t,
\end{align*}
which completes the proof.
\end{proof}

\begin{lemma}[First and second differences of $G$]\label{lem:Gdifferences}
With the multiplier $G$ defined in \eqref{eq:g}, let $k=j+\ell+p$ and $\eta=\alpha+\xi+\gamma$, and assume \eqref{eq:rectangle}. Then 
\begin{align}
&\Lambda_t(j,\alpha)|G(j+\ell,\alpha+\xi)-G(j,\alpha)|\nonumber\\
&\le C(1+\Lambda_t(\ell,\xi)+\Lambda_t(p,\gamma))
 \Bigl[\av t+\av t^2
 \mathbf1_{\substack{t\ge1,\ |j|>16(|\ell|+|p|+1)\\t/8\le(\alpha+jt)/j\le8t}}\nonumber\\
&\hspace{35mm}\times\sqrt{g(j,\alpha)+g(j+\ell,\alpha+\xi)
                      +g(j+p,\alpha+\gamma)+g(k,\eta)}\Bigr],
\label{eq:GFD1}\\
&\Lambda_t^2(j,\alpha)|G(k,\eta)-G(j+p,\alpha+\gamma)
                       -G(j+\ell,\alpha+\xi)+G(j,\alpha)|\nonumber\\
&\le C(1+\Lambda_t(\ell,\xi)+\Lambda_t(p,\gamma))^2
 \Bigl[\av t^2+\av t^4
 \mathbf1_{\substack{t\ge1,\ |j|>16(|\ell|+|p|+1)\\t/8\le(\alpha+jt)/j\le8t}}\nonumber\\
&\hspace{35mm}\times\bigl(g(j,\alpha)+g(j+\ell,\alpha+\xi)
                        +g(j+p,\alpha+\gamma)+g(k,\eta)\bigr)\Bigr].
\label{eq:GFD2}
\end{align}
The left-hand side of \eqref{eq:GFD1} can also be replaced by $\Lambda_t(j,\alpha)|(G(j+p,\alpha+\gamma)-G(j,\alpha))|$, $\Lambda_t(j,\alpha)|(G(k,\eta)-G(j+p,\alpha+\gamma))|$ and $\Lambda_t(j,\alpha)|(G(k,\eta)-G(j+\ell,\alpha+\xi))|$.
\end{lemma}
\begin{proof}

\noindent\emph{Case 1. $|j|\le16(|\ell|+|p|+1)$.}
If $|j|\le16(|\ell|+|p|+1)$, then
\begin{equation}
 \max\{|j|,|j+\ell|,|j+p|,|k|\}
       \le17(|\ell|+|p|+1)\le17(1+\Lambda_t(\ell,\xi)+\Lambda_t(p,\gamma)).
\label{eq:fourhorizontalbounds}
\end{equation}
By \eqref{eq:referencecomparison} and \eqref{eq:smallhorizontalG}, we get
\begin{align*}
\Lambda_t(j,\alpha)|G(j+\ell,\alpha+\xi)-G(j,\alpha)|
\le C\sqrt{1+(|j|+|j+\ell|)^2}\av t
\le C\bigl(1+\Lambda_t(\ell,\xi)+\Lambda_t(p,\gamma)\bigr)\av t.
\end{align*}
The other three first differences follow in the same way.
For the second difference, the same calculations give
\begin{align*}
&\Lambda_t^2(j,\alpha)
 |G(k,\eta)-G(j+p,\alpha+\gamma)-G(j+\ell,\alpha+\xi)+G(j,\alpha)|\\
\le{} &C\bigl[(1+k^2)+(1+(j+p)^2)+(1+(j+\ell)^2)+(1+j^2)\bigr]\av t^2\\
\le{} &C\bigl(1+\Lambda_t(\ell,\xi)+\Lambda_t(p,\gamma)\bigr)^2\av t^2.
\end{align*}
\noindent\emph{Case 2. $|j|>16(|\ell|+|p|+1)$.} In this case, reindexing \eqref{eq:physicalG} gives, for $s,r\in\{0,1\}$,
\begin{equation*}
 \begin{split}
 G(j+s\ell+rp,\alpha+s\xi+r\gamma)
 &=\sum_{\substack{q\in\Z\\j+q+s\ell+rp\ne0}}
       (1+q^2)^{-4}\Phi_1(t,a_q(s,r)),\\
 a_q(s,r)&=
 \frac{\alpha+jt+s(\xi+\ell t)+r(\gamma+pt)}{j+q+s\ell+rp}.
 \end{split}
\end{equation*}
\noindent\emph{Step 1. $|q|>|j|/4$.} When $|q|>|j|/4$, both $|j+s\ell+rp|$ and $|j+q+s\ell+rp|$ are at most $C|q|$. Consequently, \eqref{eq:prim}, \eqref{eq:referencecomparison} and $\Lambda_t(j,\alpha)\le C\Lambda_t(j+s\ell+rp,\alpha+s\xi+r\gamma)\le C\av q\av{a_q(s,r)}$ imply
\[
 \begin{aligned}
 \Lambda_t(j,\alpha)\sum_{\substack{|q|>|j|/4\\j+q+s\ell+rp\ne0}}
 (1+q^2)^{-4}\Phi_1(t,a_q(s,r))
 &\le C\av t\sum_q(1+q^2)^{-7/2}\le C\av t,\\
\Lambda_t(j,\alpha)^2\sum_{\substack{|q|>|j|/4\\j+q+s\ell+rp\ne0}}
 (1+q^2)^{-4}\Phi_1(t,a_q(s,r))
 &\le C\av t^2\sum_q(1+q^2)^{-3}\le C\av t^2.
 \end{aligned}
\]
The triangle inequality gives the same bounds for the first and second differences of these sums. It remains to estimate the common range $|q|\le|j|/4$.

\noindent\emph{Step 2. $|q|\le|j|/4$.}
For $|q|\le|j|/4$, extend $a_q$ to $[0,1]^2$. Due to $|j|>16(|\ell|+|p|+1)$, it follows from \eqref{eq:rectangle} that
\[
 \frac{11}{16}\le\frac{j+q+s\ell+rp}{j}\le\frac{21}{16},
 \qquad
 \bigg|\frac{s(\xi+\ell t)+r(\gamma+pt)}j\bigg|
       \le\frac18\av{(\alpha+jt)/j}.
\]
Thus it holds that
\begin{equation}
 \frac23\av{(\alpha+jt)/j}\le\av{a_q(s,r)}\le\frac{18}{11}\av{(\alpha+jt)/j},
 \qquad \Lambda_t(j,\alpha)\approx |j|\av{(\alpha+jt)/j}.
\label{eq:Gquotientcomparison}
\end{equation}
A direct computation then gives
 \begin{align}\label{eq:4}
 \partial_s a_q=\frac{\xi+\ell t-\ell a_q}{j+q+s\ell+rp},\qquad
 \partial_r a_q=\frac{\gamma+pt-pa_q}{j+q+s\ell+rp},\qquad
 \partial_s\partial_r a_q =-\frac{\ell\partial_r a_q+p\partial_s a_q}{j+q+s\ell+rp}.
 \end{align}
Hence we obtain
\begin{equation}
 \begin{aligned}
 \Lambda_t(j,\alpha)|\partial_s a_q|&\le C\Lambda_t(\ell,\xi)\av{(\alpha+jt)/j}^2,
 &\Lambda_t(j,\alpha)^2|\partial_s\partial_r a_q|
 &\le C\Lambda_t(\ell,\xi)\Lambda_t(p,\gamma)\av{(\alpha+jt)/j}^3,\\
 \Lambda_t(j,\alpha)|\partial_r a_q|&\le C\Lambda_t(p,\gamma)\av{(\alpha+jt)/j}^2,
 &\Lambda_t(j,\alpha)^2|\partial_s a_q\partial_r a_q|
 &\le C\Lambda_t(\ell,\xi)\Lambda_t(p,\gamma)\av{(\alpha+jt)/j}^4.
 \end{aligned}
\label{eq:precisechaincoefficients}
\end{equation}
For the first difference, we have
\begin{align}\label{fristdifference}
\Phi_1(t,a_q(1,0))-\Phi_1(t,a_q(0,0))
=\int_0^1\partial_a\Phi_1(t,a_q(s,0))\partial_s a_q(s,0)\dd s.
\end{align}
The other three first differences are integrals of the corresponding first derivatives. For the second difference,
\begin{equation}
 \begin{split}
 &\Phi_1(t,a_q(1,1))-\Phi_1(t,a_q(1,0))-\Phi_1(t,a_q(0,1))+\Phi_1(t,a_q(0,0))\\
 ={}&\int_0^1\!\int_0^1
 \big[\partial_a^2\Phi_1(t,a_q(s,r))\partial_s a_q\partial_r a_q+\partial_a\Phi_1(t,a_q(s,r))\partial_s\partial_r a_q\big]\dd s\dd r.
 \end{split}
\label{eq:rectintegral}
\end{equation}

\emph{Subcase 1. The case $t\le1$, or $t\ge1$ and $(\alpha+jt)/j\notin[t/8,8t]$.}
If $t\le1$, the identities in \eqref{eq:3} give uniform bounds for the weighted derivatives below. 
If $t\ge1$ and $(\alpha+jt)/j\notin[t/8,8t]$, then 
\begin{align*}
\frac{\alpha+jt}{j}<\frac t8
&\quad\Longrightarrow\quad
a_q(s,r) \le\frac{8+\sqrt{65}}{44}\,t<\frac t2,\\
\frac{\alpha+jt}{j}>8t
&\quad\Longrightarrow\quad
a_q(s,r) \ge\frac{128-2\sqrt{65}}{21}\,t>2t.
\end{align*}
Hence $a_q(s,r)\notin[t/2,2t]$.
For such $a$, we have $(1+a^2)/(1+(t-a)^2)\le4$. 
When $|a|\le2t+2$, the identities in \eqref{eq:3} therefore give
\begin{align*}
(1+a^2)|\partial_a\Phi_1(t,a)|\le C,\qquad
(1+a^2)^2|\partial_a^2\Phi_1(t,a)|+(1+a^2)^{3/2}|\partial_a\Phi_1(t,a)|\le Ct.
\end{align*}
When $|a|>2t+2>1$, \eqref{eq:17} gives the same bounds. Combining the two ranges of $a$ with the case $t\le1$, we obtain, at every argument under consideration,
\[
 \begin{aligned}
 (1+a^2)|\partial_a\Phi_1(t,a)|\le C\av t,\qquad 
 (1+a^2)^2|\partial_a^2\Phi_1(t,a)|
 +(1+a^2)^{3/2}|\partial_a\Phi_1(t,a)|&\le C\av t^2.
 \end{aligned}
\]
Combining these bounds with \eqref{eq:Gquotientcomparison}--\eqref{eq:rectintegral}, we get
\begin{align*}
&\Lambda_t(j,\alpha)|\partial_a\Phi_1(t,a_q(s,r))\partial_sa_q(s,r)|\le C\Lambda_t(\ell,\xi)\av t,\\
&\Lambda_t(j,\alpha)|\partial_a\Phi_1(t,a_q(s,r))\partial_ra_q(s,r)| \le C\Lambda_t(p,\gamma)\av t,\\
&\Lambda_t(j,\alpha)^2
 \big|\partial_s\partial_r[\Phi_1(t,a_q(s,r))]\big|
\le C\Lambda_t(\ell,\xi)\Lambda_t(p,\gamma)\av t^2.
\end{align*}
Hence recalling the bounds estimated in Step~1, by virtue of \eqref{fristdifference}, we can sum over $q$ to conclude 
\begin{align*}
&\Lambda_t(j,\alpha)
\big|G(j+\ell,\alpha+\xi)-G(j,\alpha)\big|\\
\le{}&C\av t
+\sum_{|q|\le|j|/4}\frac1{(1+q^2)^4}\Lambda_t(j,\alpha)|\Phi_1(t,a_q(1,0))-\Phi_1(t,a_q(0,0))|\\
\le{}& C\Lambda_t(\ell,\xi)\av t.
\end{align*}
The same argument gives the other three first differences. Similarly, \eqref{eq:rectintegral} gives
\begin{align*}
\Lambda_t^2(j,\alpha)|G(k,\eta)-G(j+p,\alpha+\gamma)
                       -G(j+\ell,\alpha+\xi)+G(j,\alpha)|
\le C\Lambda_t(\ell,\xi)\Lambda_t(p,\gamma)\av t^2.
\end{align*}
This proves \eqref{eq:GFD1} and \eqref{eq:GFD2} in the present case.

\emph{Subcase 2. The case $t\ge1$ and $t/8\le(\alpha+jt)/j\le8t$.} We shall use the expression $\Phi_1(t,a)=\arctan a+\arctan(t-a)$. In this case, \eqref{eq:Gquotientcomparison}, \eqref{eq:4} and the mean value theorem give
\begin{equation}\label{eq:18}
 \begin{split}
 \av{a_q(s,r)}&\approx\av t,\qquad
 \Lambda_t(j,\alpha)\approx |j|\av t,\\
 \sup_{(s,r),(s',r')\in[0,1]^2}|a_q(s,r)-a_q(s',r')|
 &\le D:=\frac{C_*(1+\Lambda_t(\ell,\xi)+\Lambda_t(p,\gamma))\av t}{|j|}, 	
 \end{split}
\end{equation}
where $C_*$ is a fixed sufficiently large constant. 

\noindent\emph{First differences.}
For $\arctan a$, it follows from \eqref{eq:precisechaincoefficients} and \eqref{eq:18} that
\begin{align*}
\Lambda_t(j,\alpha)|\arctan(a_q(1,0))-\arctan(a_q(0,0))|
\le\int_0^1\bigg|\frac{\Lambda_t(j,\alpha)\partial_s a_q(s,0)}
                         {1+a_q(s,0)^2}\bigg|\dd s
\le C\Lambda_t(\ell,\xi).
\end{align*}
Similar estimates hold for the other three differences. For $\arctan(t-a)$, \eqref{eq:18} and the inequality
\begin{equation*}
 |\arctan x-\arctan y|
 \le C|x-y|\bigl(\av x^{-1}+\av y^{-1}\bigr)
\end{equation*}
give, for any $v=(v_1,v_2),w=(w_1,w_2)\in\{0,1\}^2$,
\begin{align*}
&\Lambda_t(j,\alpha)|\arctan(t-a_q(v))-\arctan(t-a_q(w))|\\
\le{} &C(1+\Lambda_t(\ell,\xi)+\Lambda_t(p,\gamma))\av t^2
       \bigl(\av{t-a_q(v)}^{-1}+\av{t-a_q(w)}^{-1}\bigr).
\end{align*}
Then the definition of $g$ further implies
\begin{align}\label{eq:6}
 \sum_{|q|\le|j|/4}
 \frac{(1+q^2)^{-4}}{\av{t-a_q(v)}}
 \le C\bigg(\sum_{|q|\le|j|/4}
       \frac{(1+q^2)^{-4}}{\av{t-a_q(v)}^2}\bigg)^{1/2}
 \le C\sqrt{g(j+v_1\ell+v_2p,\alpha+v_1\xi+v_2\gamma)}.
\end{align}
Hence recalling the bounds estimated in Step~1 and $\Phi_1(t,a)=\arctan a+\arctan(t-a)$, we can sum over $q$ to conclude 
\begin{align*}
&\Lambda_t(j,\alpha)
\big|G(j+\ell,\alpha+\xi)-G(j,\alpha)\big|\\
\le{}&C\av t
+\sum_{|q|\le|j|/4}\frac1{(1+q^2)^4}\Lambda_t(j,\alpha)|\Phi_1(t,a_q(1,0))-\Phi_1(t,a_q(0,0))|\\
\le{}& C(1+\Lambda_t(\ell,\xi)+\Lambda_t(p,\gamma))\\
&\times\Bigl[\av t+\av t^2
 \mathbf1_{\substack{t\ge1,\ |j|>16(|\ell|+|p|+1)\\t/8\le(\alpha+jt)/j\le8t}}\sqrt{g(j,\alpha)+g(j+\ell,\alpha+\xi)
                      +g(j+p,\alpha+\gamma)+g(k,\eta)}\Bigr].
\end{align*}
The same argument gives the other three first differences. 

\noindent\emph{Second difference.}
For $\arctan a$, it follows from \eqref{eq:precisechaincoefficients} and \eqref{eq:18} that (similar to \eqref{eq:rectintegral})
\begin{align}\label{eq:7}
&\Lambda_t^2(j,\alpha)
 \big|\arctan(a_q(1,1))-\arctan(a_q(1,0))
-\arctan(a_q(0,1))+\arctan(a_q(0,0))\big|\nonumber\\
\le{} &C\Lambda_t(\ell,\xi)\Lambda_t(p,\gamma)\av t.
\end{align}
It follows from \eqref{eq:precisechaincoefficients} and \eqref{eq:18} that
\[
\Lambda_t(j,\alpha)^2|\partial_s\partial_r\arctan(t-a_q(s,r))|
       \le C(1+\Lambda_t(\ell,\xi)+\Lambda_t(p,\gamma))^2\av t^4 [1+(t-a_q(s,r))^2]^{-1}.
\]
If the upper bound $D$ in \eqref{eq:18} obeys $D\le1$, then for every $(s,r)\in[0,1]^2$ and $v\in\{0,1\}^2$, we have
\begin{align*}
1+|t-a_q(v)|^2
\le 1+2|t-a_q(s,r)|^2
       +2|a_q(s,r)-a_q(v)|^2
\le 3\bigl(1+|t-a_q(s,r)|^2\bigr).
\end{align*}
If $D>1$ and $|t-a_q(v)|\ge2D$ for every $v=(v_1,v_2)\in\{0,1\}^2$, then
\[
|t-a_q(s,r)|
\ge |t-a_q(v)|-|a_q(s,r)-a_q(v)|
\ge |t-a_q(v)|-D
\ge \frac12|t-a_q(v)|,
\]
and hence $[1+(t-a_q(s,r))^2]^{-1}\le4[1+(t-a_q(v))^2]^{-1}$. 
In either case, the preceding mixed-derivative estimate gives
\begin{align*}
&\Lambda_t(j,\alpha)^2\int_0^1\!\int_0^1
       |\partial_s\partial_r\arctan(t-a_q(s,r))|\dd s\dd r\\
\le{} &C(1+\Lambda_t(\ell,\xi)+\Lambda_t(p,\gamma))^2
       \av t^4\sum_{v\in\{0,1\}^2}[1+(t-a_q(v))^2]^{-1}.
\end{align*}
Integrating the mixed derivative over $[0,1]^2$, we obtain
\begin{equation}
\begin{aligned}
&\Lambda_t(j,\alpha)^2
 \big|\arctan(t-a_q(1,1))-\arctan(t-a_q(1,0))
       -\arctan(t-a_q(0,1))+\arctan(t-a_q(0,0))\big|\\
\le{} &C(1+\Lambda_t(\ell,\xi)+\Lambda_t(p,\gamma))^2
       \av t^4\sum_{v\in\{0,1\}^2}[1+(t-a_q(v))^2]^{-1}.
\end{aligned}
\label{eq:Gatanmixedcorner}
\end{equation}
If $D>1$ and there is a $v_0\in\{0,1\}^2$ such that $|t-a_q(v_0)|<2D$, then the definition of $D$ in \eqref{eq:18} gives
\begin{align*}
(1+\Lambda_t(\ell,\xi)+\Lambda_t(p,\gamma))^2
       \av t^4\sum_{v\in\{0,1\}^2}[1+(t-a_q(v))^2]^{-1}
\ge\frac{(1+\Lambda_t(\ell,\xi)+\Lambda_t(p,\gamma))^2
                \av t^4}{5D^2}
 =\frac{j^2\av t^2}{5C_*^2}.
\end{align*}
Since $|j|\ge1$ and $t/8\le(\alpha+jt)/j\le8t$ in the present case, it holds that
\begin{align*}
\Lambda_t(j,\alpha)^2
=1+j^2+(\alpha+jt)^2
\le Cj^2\av t^2
\le C(1+\Lambda_t(\ell,\xi)+\Lambda_t(p,\gamma))^2
       \av t^4\sum_{v\in\{0,1\}^2}[1+(t-a_q(v))^2]^{-1}.
\end{align*}
Thus, by $|\arctan x|\le\pi/2$, \eqref{eq:Gatanmixedcorner} also holds in this case.
Now recalling $\Phi_1(t,a)=\arctan a+\arctan(t-a)$, by \eqref{eq:7}, \eqref{eq:Gatanmixedcorner} and using the bounds estimated in Step~1, we get
\begin{align*}
&\Lambda_t^2(j,\alpha)
 |G(k,\eta)-G(j+p,\alpha+\gamma)-G(j+\ell,\alpha+\xi)+G(j,\alpha)|\\
\le{} &C\bigl(1+\Lambda_t(\ell,\xi)+\Lambda_t(p,\gamma)\bigr)^2\av t^2\\
&+\sum_{|q|\le|j|/4}\frac1{(1+q^2)^4}\Lambda_t(j,\alpha)^2|\Phi_1(t,a_q(1,1))-\Phi_1(t,a_q(1,0))-\Phi_1(t,a_q(0,1))+\Phi_1(t,a_q(0,0))|\\
\le{} &C\bigl(1+\Lambda_t(\ell,\xi)+\Lambda_t(p,\gamma)\bigr)^2\Bigl[\av t^2+\av t^4
 \mathbf1_{\substack{t\ge1,\ |j|>16(|\ell|+|p|+1)\\t/8\le(\alpha+jt)/j\le8t}}\sum_{|q|\le|j|/4}\frac1{(1+q^2)^4}\sum_{v\in\{0,1\}^2}[1+(t-a_q(v))^2]^{-1}\Bigr].
\end{align*}
Hence the desired estimate follows by \eqref{eq:6}. This proves \eqref{eq:GFD1} and \eqref{eq:GFD2} in the present case.

Combining Cases 1 and 2 completes the proof.
\end{proof}

\begin{lemma}[First and second differences of $V$]\label{lem:Vdifferences}
With the multipliers defined in \eqref{eq:weight}, let $k=j+\ell+p$ and $\eta=\alpha+\xi+\gamma$, and assume \eqref{eq:rectangle}. Then
\begin{equation}
\begin{aligned}
&\Lambda_t(j,\alpha)|V(j+\ell,\alpha+\xi)-V(j,\alpha)|
\le C(1+\Lambda_t(\ell,\xi)+\Lambda_t(p,\gamma))
             \av t\av{\rho^{1/3}t},\\
&\Lambda_t(j,\alpha)^2
 |V(k,\eta)-V(j+p,\alpha+\gamma)-V(j+\ell,\alpha+\xi)+V(j,\alpha)|\\
&\qquad\le C(1+\Lambda_t(\ell,\xi)+\Lambda_t(p,\gamma))^2
             \rho^{-1/3}\av t\av{\rho^{1/3}t}^3.
\end{aligned}
\label{eq:Vdifferencesharp}
\end{equation}
The first inequality also holds, with the same weight
$\Lambda_t(j,\alpha)$, for
\[
\begin{gathered}
 V(j+p,\alpha+\gamma)-V(j,\alpha),\qquad
 V(k,\eta)-V(j+p,\alpha+\gamma),\qquad
 V(k,\eta)-V(j+\ell,\alpha+\xi).
\end{gathered}
\]
\end{lemma}
\begin{proof}
When $|j|>16(|\ell|+|p|+1)$, use \eqref{eq:precisechaincoefficients} with $q=0$ and \eqref{eq:EDderiv}. Integrating the first derivative over $s\in[0,1]$ with $r=0$, and the mixed derivative over $(s,r)\in[0,1]^2$ gives \eqref{eq:Vdifferencesharp}.
When $|j|\le16(|\ell|+|p|+1)$,
\eqref{eq:boundedhorizontalprimitive} and \eqref{eq:fourhorizontalbounds}
give
\[
\begin{aligned}
\Lambda_t(j,\alpha)|V(j+\ell,\alpha+\xi)-V(j,\alpha)|
\le C\sqrt{1+(|j|+|j+\ell|)^2}\av t
\le C\bigl(1+\Lambda_t(\ell,\xi)+\Lambda_t(p,\gamma)\bigr)\av t.
\end{aligned}
\]
For the second difference, use $V\ge0$,
\eqref{eq:referencecomparison}, and \eqref{eq:smallhorizontalV}:
\begin{align*}
&\Lambda_t^2(j,\alpha)
 |V(k,\eta)-V(j+p,\alpha+\gamma)-V(j+\ell,\alpha+\xi)+V(j,\alpha)|\\
\le{} &\Lambda_t^2(j,\alpha)
 \bigl[V(k,\eta)+V(j+p,\alpha+\gamma)+V(j+\ell,\alpha+\xi)+V(j,\alpha)\bigr]\\
\le{} &C\bigl[\Lambda_t^2(k,\eta)V(k,\eta)
 +\Lambda_t^2(j+p,\alpha+\gamma)V(j+p,\alpha+\gamma)\\
&\hspace{15mm}+\Lambda_t^2(j+\ell,\alpha+\xi)V(j+\ell,\alpha+\xi)
 +\Lambda_t^2(j,\alpha)V(j,\alpha)\bigr]\\
\le{} &C\bigl[(1+k^2)+(1+(j+p)^2)+(1+(j+\ell)^2)+(1+j^2)\bigr]
\rho^{-1/3}\av t\av{\rho^{1/3}t}\\
\le{} &C\bigl(1+\Lambda_t(\ell,\xi)+\Lambda_t(p,\gamma)\bigr)^2
       \rho^{-1/3}\av t\av{\rho^{1/3}t}.
\end{align*}
Here the last step is \eqref{eq:fourhorizontalbounds}. Since $\av{\rho^{1/3}t}\ge1$, these bounds imply \eqref{eq:Vdifferencesharp}. For the other three first differences, apply the first bound in \eqref{eq:smallhorizontalV} at the two frequencies and use \eqref{eq:referencecomparison} and \eqref{eq:fourhorizontalbounds}.
\end{proof}

\begin{lemma}[First and second differences of $\theta$]\label{lem:FD}
With the multipliers defined in \eqref{eq:weight}, let $k=j+\ell+p$ and $\eta=\alpha+\xi+\gamma$, and assume \eqref{eq:rectangle}. Then
\begin{equation}
\begin{aligned}
&\Lambda_t(j,\alpha)|\theta(j+\ell,\alpha+\xi)-\theta(j,\alpha)|\\
\le{} &C\big(1+\sqrt{p^2+(\gamma+pt)^2}+\sqrt{\ell^2+(\xi+\ell t)^2}\big)\Bigl[\av t\av{\rho^{1/3}t}\\
&\quad+\av t^2\mathbf1_{\substack{t\ge1,\ |j|>16(|\ell|+|p|+1)\\t/8\le(\alpha+jt)/j\le8t}}\sqrt{g(j,\alpha)+g(j+\ell,\alpha+\xi)+g(j+p,\alpha+\gamma)+g(k,\eta)}\Bigr].
\end{aligned}
\label{eq:FD1}
\end{equation}
The left-hand side of \eqref{eq:FD1} can also be replaced by 
\[
 \Lambda_t(j,\alpha)|\theta(j+p,\alpha+\gamma)-\theta(j,\alpha)|,\quad
 \Lambda_t(j,\alpha)|\theta(k,\eta)-\theta(j+p,\alpha+\gamma)|,\quad
 \Lambda_t(j,\alpha)|\theta(k,\eta)-\theta(j+\ell,\alpha+\xi)|.
\]
Moreover,
\begin{align}
&\Lambda_t(j,\alpha)^2|\theta(k,\eta)-\theta(j+p,\alpha+\gamma)-\theta(j+\ell,\alpha+\xi)+\theta(j,\alpha)|\nonumber\\
\le{} &C\big(1+\sqrt{p^2+(\gamma+pt)^2}+\sqrt{\ell^2+(\xi+\ell t)^2}\big)^2\Big[\rho^{-1/3}\av t\av{\rho^{1/3} t}^3+\av t^2\nonumber\\
&\quad +\mathbf1_{\substack{t\ge1,\ |j|>16(|\ell|+|p|+1)\\t/8\le(\alpha+jt)/j\le8t}}\av t^4(g(j,\alpha)+g(j+\ell,\alpha+\xi)+g(j+p,\alpha+\gamma)+g(k,\eta))\Big].
\label{eq:FD2}
\end{align}
The same estimates hold with each difference of $\theta$ replaced by the corresponding difference of $h$ divided by $\Lambda_t^N(k,\eta)\Lambda_t^N(j,\alpha)$, with constants depending on $N$.
We also have the weaker estimate
\begin{equation}
 \Lambda_t(j,\alpha)|\theta(j+\ell,\alpha+\xi)-\theta(j,\alpha)|
          \le C\big(1+\sqrt{p^2+(\gamma+pt)^2}+\sqrt{\ell^2+(\xi+\ell t)^2}\big)\av t^2.
\label{eq:FDweak}
\end{equation}
We also have the following bound in terms of $j$ and $j+\ell$:
\begin{equation}
 \Lambda_t(j,\alpha)|\theta(j+\ell,\alpha+\xi)-\theta(j,\alpha)|
          \le C\big(1+(|j|+|j+\ell|)^2\big)\big(1+\ell^2+(\xi+\ell t)^2\big)\av t.
\label{eq:FDsmallh}
\end{equation}
Estimates \eqref{eq:FDweak} and \eqref{eq:FDsmallh} also hold with the difference of $\theta$ replaced by the corresponding difference of $h$ divided by $\Lambda_t^N(k,\eta)\Lambda_t^N(j,\alpha)$.
\end{lemma}

\begin{proof}
We shall use $\log\theta=-14G-2V$ and the uniform bounds \eqref{eq:thetabounds}.

\noindent\emph{1. First differences of $\theta$.}
The fundamental theorem of calculus gives
\begin{equation}\label{eq:9}
\begin{split}
\theta(j+\ell,\alpha+\xi)-\theta(j,\alpha)
=\bigl[\log\theta(j+\ell,\alpha+\xi)-\log\theta(j,\alpha)\bigr]
 \int_0^1\theta(j+\ell,\alpha+\xi)^s\theta(j,\alpha)^{1-s}\dd s.
\end{split}
\end{equation}
Since $0<\theta\le1$, it follows from \eqref{eq:9} that
\begin{align*}
&\Lambda_t(j,\alpha)|\theta(j+\ell,\alpha+\xi)-\theta(j,\alpha)|\\
\le{} &\Lambda_t(j,\alpha)
       |\log\theta(j+\ell,\alpha+\xi)-\log\theta(j,\alpha)|\\
\le{} &14\Lambda_t(j,\alpha)|G(j+\ell,\alpha+\xi)-G(j,\alpha)|
    +2\Lambda_t(j,\alpha)|V(j+\ell,\alpha+\xi)-V(j,\alpha)|\\
\le{} &C\bigl(1+\Lambda_t(\ell,\xi)+\Lambda_t(p,\gamma)\bigr)
 \Bigl[\av t\av{\rho^{1/3}t}\\*
&\hspace{15mm}+\av t^2
 \mathbf1_{\substack{t\ge1,\ |j|>16(|\ell|+|p|+1)\\t/8\le(\alpha+jt)/j\le8t}}
 \sqrt{g(j,\alpha)+g(j+\ell,\alpha+\xi)+g(j+p,\alpha+\gamma)+g(k,\eta)}\Bigr].
\end{align*}
Here we used \eqref{eq:GFD1}, \eqref{eq:Vdifferencesharp}, and $\av{\rho^{1/3}t}\ge1$. This proves \eqref{eq:FD1}.
The other three first-difference estimates follow from a similar argument.

\noindent\emph{2. The second difference of $\theta$.}
Direct expansion gives
\[
\begin{aligned}
&\theta(k,\eta)-\theta(j+p,\alpha+\gamma)
 -\theta(j+\ell,\alpha+\xi)+\theta(j,\alpha)\\
={}&\frac{\theta(j+\ell,\alpha+\xi)\theta(j+p,\alpha+\gamma)}{\theta(j,\alpha)}
 \bigg[\frac{\theta(k,\eta)\theta(j,\alpha)}
 {\theta(j+\ell,\alpha+\xi)\theta(j+p,\alpha+\gamma)}-1\bigg]\\
&+\frac{(\theta(j+\ell,\alpha+\xi)-\theta(j,\alpha))
 (\theta(j+p,\alpha+\gamma)-\theta(j,\alpha))}{\theta(j,\alpha)}.
\end{aligned}
\]
Notice that
\[
\begin{aligned}
\frac{\theta(k,\eta)\theta(j,\alpha)}
 {\theta(j+\ell,\alpha+\xi)\theta(j+p,\alpha+\gamma)}
=\exp\bigl(\log\theta(k,\eta)-\log\theta(j+p,\alpha+\gamma)
                  -\log\theta(j+\ell,\alpha+\xi)+\log\theta(j,\alpha)\bigr).
\end{aligned}
\]
By virtue of \eqref{eq:thetabounds}, applying $|e^a-1|\le e^{|a|}|a|$ and \eqref{eq:9}, we obtain
\[
\begin{aligned}
&|\theta(k,\eta)-\theta(j+p,\alpha+\gamma)
 -\theta(j+\ell,\alpha+\xi)+\theta(j,\alpha)|\\
\le {}&C|\log\theta(k,\eta)-\log\theta(j+p,\alpha+\gamma)
 -\log\theta(j+\ell,\alpha+\xi)+\log\theta(j,\alpha)|\\
&+C|\log\theta(j+\ell,\alpha+\xi)-\log\theta(j,\alpha)|
          |\log\theta(j+p,\alpha+\gamma)-\log\theta(j,\alpha)|.
\end{aligned}
\]
For the first term on the right, the identity $\log\theta=-14G-2V$ yields
\begin{align*}
&\Lambda_t^2(j,\alpha)
 |\log\theta(k,\eta)-\log\theta(j+p,\alpha+\gamma)
          -\log\theta(j+\ell,\alpha+\xi)+\log\theta(j,\alpha)|\\
\le{}&14\Lambda_t^2(j,\alpha)
 |G(k,\eta)-G(j+p,\alpha+\gamma)-G(j+\ell,\alpha+\xi)+G(j,\alpha)|\\
&+2\Lambda_t^2(j,\alpha)
 |V(k,\eta)-V(j+p,\alpha+\gamma)-V(j+\ell,\alpha+\xi)+V(j,\alpha)|\\
\le{}&C\bigl(1+\Lambda_t(\ell,\xi)+\Lambda_t(p,\gamma)\bigr)^2
 \Bigl[\av t^2+\rho^{-1/3}\av t\av{\rho^{1/3}t}^3\\*
&\hspace{15mm}+\av t^4
 \mathbf1_{\substack{t\ge1,\ |j|>16(|\ell|+|p|+1)\\t/8\le(\alpha+jt)/j\le8t}}
 \bigl(g(j,\alpha)+g(j+\ell,\alpha+\xi)+g(j+p,\alpha+\gamma)+g(k,\eta)\bigr)\Bigr].
\end{align*}
The last step is \eqref{eq:GFD2} and \eqref{eq:Vdifferencesharp}.
For the second term on the right, we have
\begin{align*}
&\Lambda_t^2(j,\alpha)
 |\log\theta(j+\ell,\alpha+\xi)-\log\theta(j,\alpha)|
 |\log\theta(j+p,\alpha+\gamma)-\log\theta(j,\alpha)|\\
\le{} &C\bigl(1+\Lambda_t(\ell,\xi)+\Lambda_t(p,\gamma)\bigr)^2
 \Bigl[\av t\av{\rho^{1/3}t}\\*
&\hspace{15mm}+\av t^2
 \mathbf1_{\substack{t\ge1,\ |j|>16(|\ell|+|p|+1)\\t/8\le(\alpha+jt)/j\le8t}}
 \sqrt{g(j,\alpha)+g(j+\ell,\alpha+\xi)+g(j+p,\alpha+\gamma)+g(k,\eta)}\Bigr]^2\\
\le{} &C\bigl(1+\Lambda_t(\ell,\xi)+\Lambda_t(p,\gamma)\bigr)^2
 \Bigl[\av t^2\av{\rho^{1/3}t}^2\\
&\hspace{15mm}+\av t^4
 \mathbf1_{\substack{t\ge1,\ |j|>16(|\ell|+|p|+1)\\t/8\le(\alpha+jt)/j\le8t}}
 \bigl(g(j,\alpha)+g(j+\ell,\alpha+\xi)+g(j+p,\alpha+\gamma)+g(k,\eta)\bigr)\Bigr].
\end{align*}
Then notice that $0<\rho=\min\{\nu,\mu\}\le1$ implies
\[
 \av t^2\av{\rho^{1/3}t}^2
 \le\rho^{-1/3}\av t\av{\rho^{1/3}t}^3.
\]
Adding the two displayed bounds proves \eqref{eq:FD2}.

\noindent\emph{3. First and second differences of $h$.}
We first estimate the differences of $\Lambda_t^{2N}$. We have
\[
\begin{aligned}
\partial_s\Lambda_t^{2N}(j+s\ell,\alpha+s\xi)
=2N\bigl[\ell(j+s\ell)+(\xi+\ell t)(\alpha+jt+s(\xi+\ell t))\bigr]\Lambda_t^{2N-2}(j+s\ell,\alpha+s\xi).
\end{aligned}
\]
Hence we get, via \eqref{eq:referencecomparison}, that
\begin{align*}
|\partial_s\Lambda_t^{2N}(j+s\ell,\alpha+s\xi)|
& \le C_N\sqrt{\ell^2+(\xi+\ell t)^2}
       \Lambda_t^{2N-1}(j+s\ell,\alpha+s\xi)\\
& \le C_N\sqrt{\ell^2+(\xi+\ell t)^2}\Lambda_t^{2N-1}(j,\alpha).	
\end{align*}
A direct computation then gives
\[
\begin{aligned}
&\partial_s\partial_r\Lambda_t^{2N}(j+s\ell+rp,\alpha+s\xi+r\gamma)\\
={}&2N\bigl[\ell p+(\xi+\ell t)(\gamma+pt)\bigr]
   \Lambda_t^{2N-2}(j+s\ell+rp,\alpha+s\xi+r\gamma)\\
&+4N(N-1)\bigl[\ell(j+s\ell+rp)
 +(\xi+\ell t)(\alpha+jt+s(\xi+\ell t)+r(\gamma+pt))\bigr]\\
&\quad\times\bigl[p(j+s\ell+rp)
 +(\gamma+pt)(\alpha+jt+s(\xi+\ell t)+r(\gamma+pt))\bigr]
\Lambda_t^{2N-4}(j+s\ell+rp,\alpha+s\xi+r\gamma).
\end{aligned}
\]
Hence we obtain
\begin{align*}
&|\partial_s\partial_r
       \Lambda_t^{2N}(j+s\ell+rp,\alpha+s\xi+r\gamma)|\\
\le{}& C_N\sqrt{\ell^2+(\xi+\ell t)^2}
           \sqrt{p^2+(\gamma+p t)^2}
            \Lambda_t^{2N-2}(j+s\ell+rp,\alpha+s\xi+r\gamma).
\end{align*}
Due to $k=j+\ell+p$, $\eta=\alpha+\xi+\gamma$, integration over $s\in[0,1]$ or $(s,r)\in[0,1]^2$, followed by \eqref{eq:referencecomparison}, therefore gives
\begin{align*}
&\frac{|\Lambda_t^{2N}(j+\ell,\alpha+\xi)-\Lambda_t^{2N}(j,\alpha)|}
       {\Lambda_t^N(k,\eta)\Lambda_t^N(j,\alpha)}
\le\frac{C_N\sqrt{\ell^2+(\xi+\ell t)^2}
                     \Lambda_t^{2N-1}(j,\alpha)}{\Lambda_t^{2N}(j,\alpha)}
 =\frac{C_N\sqrt{\ell^2+(\xi+\ell t)^2}}{\Lambda_t(j,\alpha)},\\
&\frac{|\Lambda_t^{2N}(k,\eta)-\Lambda_t^{2N}(j+p,\alpha+\gamma)
       -\Lambda_t^{2N}(j+\ell,\alpha+\xi)+\Lambda_t^{2N}(j,\alpha)|}
       {\Lambda_t^N(k,\eta)\Lambda_t^N(j,\alpha)}
\le\frac{C_N\sqrt{\ell^2+(\xi+\ell t)^2}
           \sqrt{p^2+(\gamma+p t)^2}}{\Lambda_t^2(j,\alpha)}.
\end{align*}
The first estimate also applies to the other three pairs of frequencies, using \eqref{eq:referencecomparison} throughout.

Since $h=\theta\Lambda_t^{2N}$, the first difference satisfies
\begin{align*}
h(j+\ell,\alpha+\xi)-h(j,\alpha)
&=\theta(j+\ell,\alpha+\xi)
       [\Lambda_t^{2N}(j+\ell,\alpha+\xi)-\Lambda_t^{2N}(j,\alpha)]\\
&\quad+[\theta(j+\ell,\alpha+\xi)-\theta(j,\alpha)]
                                      \Lambda_t^{2N}(j,\alpha).
\end{align*}
Hence by \eqref{eq:thetabounds} and \eqref{eq:referencecomparison},
\begin{equation}\label{eq:10}
\begin{split}
\Lambda_t(j,\alpha)
 \frac{|h(j+\ell,\alpha+\xi)-h(j,\alpha)|}
      {\Lambda_t^N(k,\eta)\Lambda_t^N(j,\alpha)}
\le C_N\sqrt{\ell^2+(\xi+\ell t)^2}
 +C_N\Lambda_t(j,\alpha)|\theta(j+\ell,\alpha+\xi)-\theta(j,\alpha)|.
\end{split}
\end{equation}
Thus \eqref{eq:FD1} bounds the second term. The first term satisfies the same bound, since $\av t\av{\rho^{1/3}t}\ge1$.

For the second difference, use the exact identity
\[
\begin{aligned}
&h(k,\eta)-h(j+p,\alpha+\gamma)-h(j+\ell,\alpha+\xi)+h(j,\alpha)\\
={}&\theta(k,\eta)\bigl[\Lambda_t^{2N}(k,\eta)
 -\Lambda_t^{2N}(j+p,\alpha+\gamma)
 -\Lambda_t^{2N}(j+\ell,\alpha+\xi)+\Lambda_t^{2N}(j,\alpha)\bigr]\\
&+(\theta(k,\eta)-\theta(j+p,\alpha+\gamma))
 (\Lambda_t^{2N}(j+p,\alpha+\gamma)-\Lambda_t^{2N}(j,\alpha))\\
&+(\theta(k,\eta)-\theta(j+\ell,\alpha+\xi))
 (\Lambda_t^{2N}(j+\ell,\alpha+\xi)-\Lambda_t^{2N}(j,\alpha))\\
&+\Lambda_t^{2N}(j,\alpha)
 [\theta(k,\eta)-\theta(j+p,\alpha+\gamma)
 -\theta(j+\ell,\alpha+\xi)+\theta(j,\alpha)].
\end{aligned}
\]
Thanks to \eqref{eq:thetabounds} and the above estimate, the first term on the right is bounded by
\begin{align*}
&\frac{\Lambda_t^2(j,\alpha)|\theta(k,\eta)|}
 {\Lambda_t^N(k,\eta)\Lambda_t^N(j,\alpha)}
 \bigl|\Lambda_t^{2N}(k,\eta)-\Lambda_t^{2N}(j+p,\alpha+\gamma)
       -\Lambda_t^{2N}(j+\ell,\alpha+\xi)+\Lambda_t^{2N}(j,\alpha)\bigr|\\
\le{} &C_N\sqrt{\ell^2+(\xi+\ell t)^2}
           \sqrt{p^2+(\gamma+p t)^2}.
\end{align*}
For the next two terms, the first-difference estimates for
$\Lambda_t^{2N}$ give
\begin{align*}
&\frac{\Lambda_t^2(j,\alpha)}{\Lambda_t^N(k,\eta)\Lambda_t^N(j,\alpha)}
 \Bigl[|\theta(k,\eta)-\theta(j+p,\alpha+\gamma)|
       |\Lambda_t^{2N}(j+p,\alpha+\gamma)-\Lambda_t^{2N}(j,\alpha)|\\*
&\hspace{18mm}+|\theta(k,\eta)-\theta(j+\ell,\alpha+\xi)|
       |\Lambda_t^{2N}(j+\ell,\alpha+\xi)-\Lambda_t^{2N}(j,\alpha)|\Bigr]\\
\le{} &C_N\Lambda_t(j,\alpha)
 \Bigl[\Lambda_t(p,\gamma)|\theta(k,\eta)-\theta(j+p,\alpha+\gamma)|
       +\Lambda_t(\ell,\xi)|\theta(k,\eta)-\theta(j+\ell,\alpha+\xi)|\Bigr]\\
\le{} &C_N\big(1+\sqrt{p^2+(\gamma+pt)^2}+\sqrt{\ell^2+(\xi+\ell t)^2}\big)^2\av t^2.
\end{align*}
In the last step we used \eqref{eq:FD1}, $\av{\rho^{1/3}t}\le\av t$, and
\[
 0\le g(k,\eta)\le\sum_{n\ne0}(1+(k-n)^2)^{-4}
       \le\sum_{q\in\Z}(1+q^2)^{-4}<\infty.
\]
Finally, the fourth term satisfies
\begin{align*}
&\frac{\Lambda_t^2(j,\alpha)\Lambda_t^{2N}(j,\alpha)}
 {\Lambda_t^N(k,\eta)\Lambda_t^N(j,\alpha)}
 |\theta(k,\eta)-\theta(j+p,\alpha+\gamma)
                 -\theta(j+\ell,\alpha+\xi)+\theta(j,\alpha)|\\
\le{} &C_N\Lambda_t^2(j,\alpha)
 |\theta(k,\eta)-\theta(j+p,\alpha+\gamma)
                 -\theta(j+\ell,\alpha+\xi)+\theta(j,\alpha)|.
\end{align*}
The right-hand side of \eqref{eq:FD2} bounds this expression. 
Adding them proves the second-difference estimate for $h$.

\noindent\emph{4. The bounds \eqref{eq:FDweak} and \eqref{eq:FDsmallh}.}
Due to the uniform bound for $g$ and $\av{\rho^{1/3}t}\le\av t$, it follows from \eqref{eq:FD1} that
\begin{align*}
\Lambda_t(j,\alpha)|\theta(j+\ell,\alpha+\xi)-\theta(j,\alpha)|
&\le C\big(1+\sqrt{p^2+(\gamma+pt)^2}+\sqrt{\ell^2+(\xi+\ell t)^2}\big)
       \bigl(\av t\av{\rho^{1/3}t}+\av t^2\bigr)\\
&\le C\big(1+\sqrt{p^2+(\gamma+pt)^2}+\sqrt{\ell^2+(\xi+\ell t)^2}\big)\av t^2.
\end{align*}
This proves \eqref{eq:FDweak}. The corresponding estimate for $h$ follows from \eqref{eq:10}.
For \eqref{eq:FDsmallh}, use \eqref{eq:thetabounds}, \eqref{eq:boundedhorizontalprimitive} and \eqref{eq:9}:
\begin{align*}
&\Lambda_t(j,\alpha)|\theta(j+\ell,\alpha+\xi)-\theta(j,\alpha)|\\
\le{} &14\Lambda_t(j,\alpha)|G(j+\ell,\alpha+\xi)-G(j,\alpha)|
    +2\Lambda_t(j,\alpha)|V(j+\ell,\alpha+\xi)-V(j,\alpha)|\\
\le{} &C\sqrt{1+(|j|+|j+\ell|)^2}\av t.
\end{align*}
For $h$, \eqref{eq:10} consequently gives
\begin{align*}
\Lambda_t(j,\alpha)
 \frac{|h(j+\ell,\alpha+\xi)-h(j,\alpha)|}
      {\Lambda_t^N(k,\eta)\Lambda_t^N(j,\alpha)}
&\le C_N\Bigl[\sqrt{\ell^2+(\xi+\ell t)^2}
       +\sqrt{1+(|j|+|j+\ell|)^2}\av t\Bigr]\\
&\le C_N\bigl(1+(|j|+|j+\ell|)^2\bigr)
            \bigl(1+\ell^2+(\xi+\ell t)^2\bigr)\av t.
\end{align*}
These inequalities prove \eqref{eq:FDsmallh} for both weights.
The proof is thus complete.
\end{proof}

The operator $\partial_t-\ell\partial_{\xi}-m\partial_{\zeta}$ in \eqref{eq:R3} also acts on $h$.
The next lemma bounds the resulting differences of $g$ and $(\partial_t-k\partial_\eta)V$.

\begin{lemma}\label{lem:timediff}
With the multipliers defined in \eqref{eq:g} and \eqref{eq:weight}, if $k=m+\ell,\eta=\zeta+\xi$ and $8\sqrt{\ell^2+(\xi+\ell t)^2}\le\sqrt{m^2+(\zeta+mt)^2}$, then
\begin{align}
& \Lambda_t(m,\zeta)|g(k,\eta)-g(m,\zeta)|\le C\Lambda_t(\ell,\xi)
          \big[1+\av t^2(g(k,\eta)+g(m,\zeta))\big],\label{eq:gdiff}\\
& \Lambda_t(m,\zeta)|(\partial_t-k\partial_{\eta})V(k,\eta)-(\partial_t-m\partial_{\zeta})V(m,\zeta)|\le C\Lambda_t(\ell,\xi)\av{\rho^{1/3} t}^2.
\label{eq:vdiff}
\end{align}
In addition, for $k\ne0$,
\begin{equation}
 |k||\partial_\eta\theta(k,\eta)|\le C.
\label{eq:verticalLip}
\end{equation}
\end{lemma}

\begin{proof}
Similar to \eqref{eq:referencecomparison}, due to $k=m+\ell,\eta=\zeta+\xi$ and $8\sqrt{\ell^2+(\xi+\ell t)^2}\le\sqrt{m^2+(\zeta+mt)^2}$, we have 
\begin{align}\label{eq:11}
\frac78\Lambda_t(m,\zeta)
\le\Lambda_t(m+s\ell,\zeta+s\xi)
\le\frac98\Lambda_t(m,\zeta),\qquad 0\le s\le1.
\end{align}

We first prove \eqref{eq:gdiff}. First suppose that $|m|\le16(|\ell|+1)$. Since
\[
\Lambda_t(m,\zeta)\le\sqrt{1+m^2}+|\zeta+mt|,
\qquad
\frac{|a|}{1+(t-a)^2}
\le\frac{t}{1+(t-a)^2}+\frac12,
\]
the definition of $g$ in \eqref{eq:g} and its uniform bound give
\begin{align*}
\Lambda_t(m,\zeta)g(m,\zeta)
&\le \sqrt{1+m^2}g(m,\zeta)
 +\frac12\sum_{n\ne0}(1+(m-n)^2)^{-4}|n|+t\sum_{n\ne0}
 \frac{(1+(m-n)^2)^{-4}|n|}
      {1+(t-(\zeta+mt)/n)^2}\\
&\le C\sqrt{1+m^2}+Ct\bigg[\sum_{n\ne0}(1+(m-n)^2)^{-4}n^2\bigg]^{1/2}\sqrt{g(m,\zeta)}\\
&\le C\Lambda_t(\ell,\xi)\bigl(1+t^2g(m,\zeta)\bigr).
\end{align*}
Here we also used $\sqrt{1+m^2}\le C\sqrt{1+\ell^2}\le C\Lambda_t(\ell,\xi)$.
Since $|k|\le17(|\ell|+1)$, the same calculation and \eqref{eq:11} give
\begin{align*}
\Lambda_t(m,\zeta)g(k,\eta)
\le\frac87\Lambda_t(k,\eta)g(k,\eta)\le C\Lambda_t(\ell,\xi)\bigl(1+t^2g(k,\eta)\bigr).
\end{align*}
Adding the last two bounds yields \eqref{eq:gdiff}.

Now suppose that $|m|>16(|\ell|+1)>0$. Reindexing the two series gives
\begin{align*}
g(m,\zeta)
=\sum_{\substack{q\in\Z\\m+q\ne0}}
 \frac{(1+q^2)^{-4}}{1+(t-(\zeta+mt)/(m+q))^2},\quad
g(k,\eta)
=\sum_{\substack{q\in\Z\\k+q\ne0}}
 \frac{(1+q^2)^{-4}}{1+(t-(\eta+kt)/(k+q))^2}.
\end{align*}
We estimate separately the sums over $|q|\le|m|/4$ and $|q|>|m|/4$.
For the case $|q|\le|m|/4$, one has 
\[
\frac34|m|\le|m+q|\le\frac54|m|,
\qquad
\frac{11}{16}|m|\le|k+q|\le\frac{21}{16}|m|.
\]
For brevity, we write
\[
a=\frac{\zeta+mt}{m+q},\qquad
b=\frac{\eta+kt}{k+q}.
\]
Then we have
\begin{align*}
\Lambda_t(m,\zeta)\le C|m|\sqrt{1+a^2},\qquad
\Lambda_t(m,\zeta)\le\frac87\Lambda_t(k,\eta)
 \le C|m|\sqrt{1+b^2}.
\end{align*}
Using $\eta+kt=\zeta+mt+\xi+\ell t$, we obtain
\begin{align*}
b-a
=\frac{(\xi+\ell t)(m+q)-\ell(\zeta+mt)}{(m+q)(k+q)}
=\frac{\xi+\ell t-\ell a}{k+q}
 =\frac{\xi+\ell t-\ell b}{m+q}.
\end{align*}
Thus it holds that
\begin{align*}
\Lambda_t(m,\zeta)|b-a|
&\le C\sqrt{1+a^2}|\xi+\ell t-\ell a|
\le C\sqrt{\ell^2+(\xi+\ell t)^2}(1+a^2),\\
\Lambda_t(m,\zeta)|b-a|
&\le C\sqrt{1+b^2}|\xi+\ell t-\ell b|
\le C\sqrt{\ell^2+(\xi+\ell t)^2}(1+b^2).
\end{align*}
On the other hand, direct subtraction yields
\begin{align*}
\bigg|\frac1{1+(t-a)^2}-\frac1{1+(t-b)^2}\bigg|
&\le\frac{|a-b|}{2}
 \bigg(\frac1{1+(t-a)^2}+\frac1{1+(t-b)^2}\bigg).
\end{align*}
Hence it follows that
\begin{align*}
\Lambda_t(m,\zeta)
 \bigg|\frac1{1+(t-a)^2}-\frac1{1+(t-b)^2}\bigg|
&\le C\Lambda_t(\ell,\xi)
 \bigg(\frac{1+a^2}{1+(t-a)^2}
       +\frac{1+b^2}{1+(t-b)^2}\bigg)\\
&\le C\Lambda_t(\ell,\xi)
 \bigg[1+t^2\bigg(\frac1{1+(t-a)^2}
                 +\frac1{1+(t-b)^2}\bigg)\bigg].
\end{align*}
Since $|k|=|\ell+m|\ge|m|-|\ell|\ge15|m|/16>|m|/4$ and $|q|\le|m|/4$, we have $k+q\neq0$ and $m+q\neq0$. Hence multiplying by $(1+q^2)^{-4}$ and summing gives
\begin{align*}
&\Lambda_t(m,\zeta)
 \bigg|\sum_{|q|\le|m|/4}(1+q^2)^{-4}
 \bigg(\frac1{1+(t-(\eta+kt)/(k+q))^2}
       -\frac1{1+(t-(\zeta+mt)/(m+q))^2}\bigg)\bigg|\\
\le{} &C\Lambda_t(\ell,\xi)
 \bigl[1+t^2\bigl(g(k,\eta)+g(m,\zeta)\bigr)\bigr].
\end{align*}
For the case $|q|>|m|/4$, the inequalities $|m|<4|q|$ and $|k|\le17|m|/16$ imply
\begin{align*}
\sum_{|q|>|m|/4}(1+q^2)^{-4}
 \Bigl[\sqrt{1+m^2}+\sqrt{1+k^2}
       +|m+q|+|k+q|+(m+q)^2+(k+q)^2\Bigr]\le C.
\end{align*}
Thus we further conclude that 
\begin{align*}
\Lambda_t(m,\zeta)
 \sum_{\substack{|q|>|m|/4\\m+q\ne0}}
 \frac{(1+q^2)^{-4}}{1+(t-(\zeta+mt)/(m+q))^2}
&\le C+t\sum_{\substack{|q|>|m|/4\\m+q\ne0}}
 \frac{(1+q^2)^{-4}|m+q|}{1+(t-(\zeta+mt)/(m+q))^2}\\
&\le C\bigl(1+t\sqrt{g(m,\zeta)}\bigr)
 \le C\bigl(1+t^2g(m,\zeta)\bigr).
\end{align*}
Likewise, using \eqref{eq:11}, we get
\begin{align*}
\Lambda_t(m,\zeta)
 \sum_{\substack{|q|>|m|/4\\k+q\ne0}}
 \frac{(1+q^2)^{-4}}{1+(t-(\eta+kt)/(k+q))^2}
\le C\bigl(1+t^2g(k,\eta)\bigr).
\end{align*}
Since $\Lambda_t(\ell,\xi)\ge1$, adding these estimates proves \eqref{eq:gdiff}.

We next prove \eqref{eq:vdiff}.
First suppose that $|m|\le16(|\ell|+1)$. For $m\ne0$, by \eqref{eq:physicalV}
\begin{align*}
\Lambda_t(m,\zeta)(\partial_t-m\partial_\zeta)V(m,\zeta)
&\le \sqrt{1+m^2}
 \frac{\rho^{1/3}\sqrt{1+(t+\zeta/m)^2}}
      {1+\rho^{2/3}\zeta^2/m^2}\\
&\le \sqrt{1+m^2}
 \bigg[\rho^{1/3}(1+t)
       +\frac{\rho^{1/3}|\zeta/m|}{1+\rho^{2/3}\zeta^2/m^2}\bigg]\\
&\le C\sqrt{1+m^2}\av{\rho^{1/3}t}.
\end{align*}
This bound also holds for $m=0$, since the left-hand side is zero.
Applying it at $(k,\eta)$ as well, we obtain
\begin{align*}
\Lambda_t(m,\zeta)
 \bigl|(\partial_t-k\partial_\eta)V(k,\eta)
       -(\partial_t-m\partial_\zeta)V(m,\zeta)\bigr|
\le C\bigl[\sqrt{1+m^2}+\sqrt{1+k^2}\bigr]\av{\rho^{1/3}t}
\le C\Lambda_t(\ell,\xi)\av{\rho^{1/3}t}.
\end{align*}

Now suppose that $|m|>16(|\ell|+1)$. For $0\le s\le1$, we have $|m+s\ell|\approx|m|$. Then write
\[
a(s)=\frac{\zeta+mt+s(\xi+\ell t)}{m+s\ell}.
\]
Recalling \eqref{eq:physicalV} and $k=m+\ell,\eta=\zeta+\xi$, we have 
\[
(\partial_t-k\partial_\eta)V(k,\eta)
       -(\partial_t-m\partial_\zeta)V(m,\zeta)=\frac{\rho^{1/3}}{1+\rho^{2/3}(t-a(1))^2}-\frac{\rho^{1/3}}{1+\rho^{2/3}(t-a(0))^2}.
\]
Direct differentiation gives
\[
\frac{\dd}{\dd s}
 \frac{\rho^{1/3}}{1+\rho^{2/3}(t-a(s))^2}
=\frac{2\rho(t-a(s))a'(s)}
 {[1+\rho^{2/3}(t-a(s))^2]^2}.
\]
Due to
\begin{align*}
a'(s)
=\frac{m(\xi+\ell t)-\ell(\zeta+mt)}{(m+s\ell)^2}
=\frac{\xi+\ell t-\ell a(s)}{m+s\ell},
\end{align*}
we infer from \eqref{eq:11} and $\Lambda_t(m+s\ell,\zeta+s\xi) \le C|m+s\ell|[1+a(s)^2]^{1/2}$ that
\begin{align*}
\Lambda_t(m,\zeta)|a'(s)|
\le C[1+a(s)^2]^{1/2}|\xi+\ell t-\ell a(s)|
\le C\Lambda_t(\ell,\xi)[1+a(s)^2].
\end{align*}
Then for every real $a$, using $0<\rho\le1$, we have
\begin{align*}
\rho^{2/3}(1+a^2)
\le\rho^{2/3}+2\rho^{2/3}t^2+2\rho^{2/3}(t-a)^2
\le2\av{\rho^{1/3}t}^2[1+\rho^{2/3}(t-a)^2],
\end{align*}
and therefore
\begin{align*}
(1+a^2)\frac{2\rho|t-a|}{[1+\rho^{2/3}(t-a)^2]^2}
\le4\av{\rho^{1/3}t}^2
 \frac{\rho^{1/3}|t-a|}{1+\rho^{2/3}(t-a)^2}
\le2\av{\rho^{1/3}t}^2.
\end{align*}
Hence we obtain
\begin{align*}
\Lambda_t(m,\zeta)
 \bigl|(\partial_t-k\partial_\eta)V(k,\eta)
       -(\partial_t-m\partial_\zeta)V(m,\zeta)\bigr|
&\le\int_0^1
 \frac{2\rho|t-a(s)|\Lambda_t(m,\zeta)|a'(s)|}
 {[1+\rho^{2/3}(t-a(s))^2]^2}\dd s\\
&\le C\Lambda_t(\ell,\xi)\int_0^1
 \frac{2\rho|t-a(s)|[1+a(s)^2]}
 {[1+\rho^{2/3}(t-a(s))^2]^2}\dd s\\
&\le C\Lambda_t(\ell,\xi)\av{\rho^{1/3}t}^2.
\end{align*}
This proves \eqref{eq:vdiff}.

Finally, we show \eqref{eq:verticalLip}. Due to $\theta=e^{-14G-2V}$ and $0<\theta\le1$, we have
\[
|k||\partial_\eta\theta(k,\eta)|
\le14|k|\,|\partial_\eta G(k,\eta)|
   +2|k|\,|\partial_\eta V(k,\eta)|.
\]
By \eqref{eq:physicalG}, we have
\begin{align*}
\partial_\eta G(k,\eta)
&= \sum_{n\ne0}\frac{(1+(k-n)^2)^{-4}}{n}
 \bigg[\frac1{1+((\eta+kt)/n)^2}
       -\frac1{1+(t-(\eta+kt)/n)^2}\bigg],\\
|\partial_\eta G(k,\eta)|&\le2\sum_{n\ne0}\frac{(1+(k-n)^2)^{-4}}{|n|}.
\end{align*}
For $k\ne0$, split the last sum as follows:
\begin{align*}
\sum_{n\ne0}\frac{(1+(k-n)^2)^{-4}}{|n|}
&=\sum_{\substack{n\ne0\\|n-k|\le|k|/2}}
 \frac{(1+(k-n)^2)^{-4}}{|n|}
 +\sum_{\substack{n\ne0\\|n-k|>|k|/2}}
 \frac{(1+(k-n)^2)^{-4}}{|n|}\\
&\le\frac2{|k|}\sum_{q\in\Z}(1+q^2)^{-4}
 +\frac2{|k|}
 \sum_{q\in\Z}(1+q^2)^{-7/2}\\
& \le\frac C{|k|}.
\end{align*}
Similarly, \eqref{eq:physicalV} gives
\begin{align*}
\partial_\eta V(k,\eta)
=\frac{\rho^{1/3}}{k}
 \bigg[\frac1{1+\rho^{2/3}(\eta/k+t)^2}
       -\frac1{1+\rho^{2/3}\eta^2/k^2}\bigg],\qquad
|\partial_\eta V(k,\eta)|\le\frac{2\rho^{1/3}}{|k|}.
\end{align*}
Adding these two bounds proves \eqref{eq:verticalLip} and completes the proof.
\end{proof}

\subsection{Comparison of $g$ at related frequencies}

We next compare the values of $g$ at different frequencies.
\begin{lemma}\label{lem:kerneltransfer}
With the multipliers defined in \eqref{eq:g} and \eqref{eq:weight}, if $k=m+\ell$ and $\eta=\zeta+\xi$, then
\begin{equation}
\begin{aligned}
&g(k,\eta)+g(m,\zeta)\le C\Lambda_t^{5}(\ell,\xi)\sqrt{g(k,\eta)g(m,\zeta)},\\
&\sqrt{g(k,\eta)+g(m,\zeta)}\,\kappa(a,\omega)
\le C\Lambda_t^{5}(\ell,\xi)\sqrt{g(k,\eta)g(m,\zeta)},\qquad (a,\omega)\in\{(k,\eta),(m,\zeta)\}.
 \label{eq:twotransfer}
\end{aligned}
\end{equation}
If $k=j+\ell+p$ and $\eta=\alpha+\xi+\gamma$, then
\begin{equation}
\begin{aligned}
&g(j,\alpha)+g(j+\ell,\alpha+\xi)
                +g(j+p,\alpha+\gamma)+g(k,\eta)\\
&\le C(1+\ell^2)^2(1+p^2)^2
       \sqrt{1+(\xi+\ell t)^2}\sqrt{1+(\gamma+pt)^2}
       \sqrt{g(k,\eta)g(j,\alpha)}.
 \label{eq:doubletransferseparated}
\end{aligned}
\end{equation}
\end{lemma}
\begin{proof}
First let $k=m+\ell$ and $\eta=\zeta+\xi$.
Applying \eqref{eq:shiftg} gives
\[
 g(k,\eta)\le C\Lambda_t^{10}(\ell,\xi)g(m,\zeta),
 \qquad g(m,\zeta)\le C\Lambda_t^{10}(\ell,\xi)g(k,\eta).
\]
Multiplying the first inequality by $g(k,\eta)$ and the
second by $g(m,\zeta)$, and then taking square roots, yields
\begin{align*}
 g(k,\eta)^2
 &\le C\Lambda_t^{10}(\ell,\xi)g(k,\eta)g(m,\zeta),\\
 g(m,\zeta)^2
 &\le C\Lambda_t^{10}(\ell,\xi)g(k,\eta)g(m,\zeta),\\
 g(k,\eta)+g(m,\zeta)
 &\le C\Lambda_t^5(\ell,\xi)\sqrt{g(k,\eta)g(m,\zeta)}.
\end{align*}
This is the first inequality in \eqref{eq:twotransfer}.
For the second one, use $\kappa^2\le g$ from
\eqref{eq:gprimitive}:
\begin{align*}
 \bigl(g(k,\eta)+g(m,\zeta)\bigr)\kappa(k,\eta)^2
 &\le\bigl(g(k,\eta)+g(m,\zeta)\bigr)g(k,\eta)\\
 &\le\bigl(C\Lambda_t^{10}(\ell,\xi)+1\bigr)g(k,\eta)g(m,\zeta)
 \le C\Lambda_t^{10}(\ell,\xi)g(k,\eta)g(m,\zeta),\\
 \bigl(g(k,\eta)+g(m,\zeta)\bigr)\kappa(m,\zeta)^2
 &\le\bigl(g(k,\eta)+g(m,\zeta)\bigr)g(m,\zeta)\\
 &\le\bigl(1+C\Lambda_t^{10}(\ell,\xi)\bigr)g(k,\eta)g(m,\zeta)
 \le C\Lambda_t^{10}(\ell,\xi)g(k,\eta)g(m,\zeta).
\end{align*}
This proves the required bound for both choices of $(a,\omega)$.

Next let $k=j+\ell+p$ and $\eta=\alpha+\xi+\gamma$.
We prove \eqref{eq:doubletransferseparated} by estimating its four terms separately. For $(j+\ell,\alpha+\xi)$,
\eqref{eq:shiftgseparated} gives
\begin{align*}
 g(j+\ell,\alpha+\xi)
 &\le C(1+\ell^2)^4(1+(\xi+\ell t)^2)g(j,\alpha),\\
 g(j+\ell,\alpha+\xi)= g(k-p,\eta-\gamma)
 &\le C(1+p^2)^4(1+(\gamma+pt)^2)g(k,\eta).
\end{align*}
Multiplying the two inequalities and taking the
square root, we obtain
\begin{align*}
 g(j+\ell,\alpha+\xi)
 \le C(1+\ell^2)^2(1+p^2)^2
       \sqrt{1+(\xi+\ell t)^2}\sqrt{1+(\gamma+pt)^2}
       \sqrt{g(k,\eta)g(j,\alpha)}.
\end{align*}
For $(j+p,\alpha+\gamma)$, the corresponding two estimates are
\begin{align*}
 g(j+p,\alpha+\gamma)
 &\le C(1+p^2)^4(1+(\gamma+pt)^2)g(j,\alpha),\\
 g(j+p,\alpha+\gamma) = g(k-\ell,\eta-\xi)
 &\le C(1+\ell^2)^4(1+(\xi+\ell t)^2)g(k,\eta).
\end{align*}
Their product gives
\begin{align*}
 g(j+p,\alpha+\gamma)
 \le C(1+\ell^2)^2(1+p^2)^2
       \sqrt{1+(\xi+\ell t)^2}\sqrt{1+(\gamma+pt)^2}
       \sqrt{g(k,\eta)g(j,\alpha)}.
\end{align*}
For the remaining two frequencies, use $k=j+\ell+p$ and $\eta=\alpha+\xi+\gamma$:
\begin{align*}
 g(k,\eta)= g(j+\ell+p,\alpha+\xi+\gamma)
 &\le C(1+(\ell+p)^2)^4
       (1+(\xi+\ell t+\gamma+pt)^2)g(j,\alpha),\\
 g(j,\alpha)= g(k-(\ell+p),\eta-(\xi+\gamma))
 &\le C(1+(\ell+p)^2)^4
       (1+(\xi+\ell t+\gamma+pt)^2)g(k,\eta).
\end{align*}
As in the proof of \eqref{eq:twotransfer}, these imply
\begin{align*}
g(k,\eta)+g(j,\alpha)
&\le C(1+(\ell+p)^2)^2
       \sqrt{1+(\xi+\ell t+\gamma+pt)^2}
       \sqrt{g(k,\eta)g(j,\alpha)}\\
&\le C(1+\ell^2)^2(1+p^2)^2
       \sqrt{1+(\xi+\ell t)^2}\sqrt{1+(\gamma+pt)^2}
       \sqrt{g(k,\eta)g(j,\alpha)}.
\end{align*}
Adding the bounds for $g(j+\ell,\alpha+\xi)$, $g(j+p,\alpha+\gamma)$, and $g(k,\eta)+g(j,\alpha)$ proves \eqref{eq:doubletransferseparated}.
\end{proof}

\begin{remark}
Since
\begin{align*}
 (1+\ell^2)^2\sqrt{1+(\xi+\ell t)^2}
 &\le(1+\ell^2+(\xi+\ell t)^2)^{5/2}=\Lambda_t^5(\ell,\xi),\\
 (1+p^2)^2\sqrt{1+(\gamma+pt)^2}
 &\le(1+p^2+(\gamma+pt)^2)^{5/2}=\Lambda_t^5(p,\gamma),
\end{align*}
we also obtain
\begin{equation}
\begin{aligned}
g(j,\alpha)+g(j+\ell,\alpha+\xi)
                  +g(j+p,\alpha+\gamma)+g(k,\eta)
\le C\Lambda_t^{5}(\ell,\xi)
       \Lambda_t^{5}(p,\gamma)
       \sqrt{g(k,\eta)g(j,\alpha)}.
 \label{eq:doubletransfer}
\end{aligned}
\end{equation}
For the three frequencies $(p,\gamma)$, $(p+j,\gamma+\alpha)$, and $(k,\eta)$, it follows from \eqref{eq:doubletransfer} that
\begin{equation}
\begin{split}
g(k,\eta)+g(j+p,\alpha+\gamma)+g(p,\gamma)
&\le g(p,\gamma)+g(p+j,\gamma+\alpha)
            +g(p+\ell,\gamma+\xi)+g(k,\eta)\\
&\le C\Lambda_t^5(j,\alpha)\Lambda_t^5(\ell,\xi) \sqrt{g(k,\eta)g(p,\gamma)}.	
\end{split}
 \label{eq:triangletransfer}
\end{equation}	
\end{remark}

\section{Bounds for the nonlinear reactions}
\label{sec:cubic}

In this section, we establish estimates for \(\mathcal B\) and the nonlinear reactions.

\subsection{The symbol estimates}
We first compute the derivative of $a(k,\eta;\ell,\xi)$ and establish the following pointwise estimates.
\begin{lemma}
\label{lem:kerneldifferentiation}
With the notation in \eqref{eq:kernelgeometric}, for \(k=\ell+m\) and \(\eta=\xi+\zeta\), one has
\begin{align}
(\partial_t-\ell\partial_{\xi}-m\partial_{\zeta})a(k,\eta;\ell,\xi)={}&(\kappa^2(m,\zeta)-\kappa^2(k,\eta))\frac{(\ell\zeta-m\xi)^2}{\sqrt{\ell^2+\xi^2}\sqrt{k^2+\eta^2}\sqrt{m^2+\zeta^2}}\nonumber\\
&+\frac{\ell\xi}{(\ell^2+\xi^2)}a(k,\eta;\ell,\xi).
\label{eq:Adot}
\end{align}
\end{lemma}
\begin{proof}
By \eqref{eq:kernelgeometric} and the relations $k=\ell+m$ and $\eta=\xi+\zeta$, we have
\[
a(k,\eta;\ell,\xi)=\bigl(e(\ell,\xi)\cdot(m,\zeta)\bigr)
            \bigl(e(k,\eta)\cdot e(m,\zeta)\bigr),
\]
and therefore
\begin{align*}
(\partial_t-\ell\partial_{\xi}-m\partial_{\zeta})a(k,\eta;\ell,\xi)={}&\frac{\ell\zeta-m\xi}{\sqrt{\ell^2+\xi^2}}\bigg[(-\ell\partial_{\xi}-m\partial_{\zeta}) \Big(\frac{(k,\eta)}{\sqrt{k^2+\eta^2}}\cdot\frac{(m,\zeta)}{\sqrt{m^2+\zeta^2}}\Big)\bigg]\\
&+\bigg[(-\ell\partial_{\xi}-m\partial_{\zeta})\Big(\frac{-m\xi+\ell\zeta}{\sqrt{\ell^2+\xi^2}}\Big)\bigg]\bigl(e(k,\eta)\cdot e(m,\zeta)\bigr).
\end{align*}
For the first term, use $k=\ell+m$, $\eta=\xi+\zeta$, and \eqref{eq:frequencies} to obtain 
\[
\begin{aligned}
&(-\ell\partial_{\xi}-m\partial_{\zeta})\Big(\frac{(k,\eta)}{\sqrt{k^2+\eta^2}}\cdot\frac{(m,\zeta)}{\sqrt{m^2+\zeta^2}}\Big)\\
= &-k\partial_{\eta}\Big(\frac{(k,\eta)}{\sqrt{k^2+\eta^2}}\Big)\cdot\frac{(m,\zeta)}{\sqrt{m^2+\zeta^2}}+\frac{(k,\eta)}{\sqrt{k^2+\eta^2}}\cdot(-m\partial_{\zeta})\Big(\frac{(m,\zeta)}{\sqrt{m^2+\zeta^2}}\Big)\\
= &-\kappa^2(k,\eta)e(k,\eta)\cdot\frac{(m,\zeta)}{\sqrt{m^2+\zeta^2}}
              -\kappa^2(m,\zeta)\frac{(k,\eta)}{\sqrt{k^2+\eta^2}}\cdot e(m,\zeta)\\
={} &(\kappa^2(m,\zeta)-\kappa^2(k,\eta))\frac{\ell\zeta-m\xi}{\sqrt{k^2+\eta^2}\sqrt{m^2+\zeta^2}},
\end{aligned}
\]
where we used
\[
\begin{aligned}
&(-k\partial_\eta)\Big(\frac{(k,\eta)}{\sqrt{k^2+\eta^2}}\Big)=-\kappa^2(k,\eta)e(k,\eta).
\end{aligned}
\]
For the second term, direct computation gives 
\begin{align*}
\bigg[(-\ell\partial_{\xi}-m\partial_{\zeta})\Big(\frac{-m\xi+\ell\zeta}{\sqrt{\ell^2+\xi^2}}\Big)\bigg]\bigl(e(k,\eta)\cdot e(m,\zeta)\bigr)&=\frac{\ell\xi}{(\ell^2+\xi^2)}\bigl(e(\ell,\xi)\cdot(m,\zeta)\bigr)\bigl(e(k,\eta)\cdot e(m,\zeta)\bigr)\\
&=\frac{\ell\xi}{(\ell^2+\xi^2)}a(k,\eta;\ell,\xi).
\end{align*}
This proves the identity.
\end{proof}

\begin{lemma}
\label{lem:singlesymbol}
For integers $c\ne0$ and $d$, and $\omega,\sigma\in\R$,
\begin{equation}
\begin{aligned}
 \frac{|c\sigma-d\omega|}{\sqrt{c^2+\omega^2}}
 &\le C\Lambda_t(d,\sigma)(1+\av t\kappa(c,\omega)),\\
 \frac{\Lambda_t(c,\omega)}{\sqrt{c^2+\omega^2}}
 &\le C(1+\av t\kappa(c,\omega)).
\end{aligned}
\label{eq:elementarysymbols}
\end{equation}
Let $k=\ell+m$, $\eta=\xi+\zeta$, and $\ell\ne0$. Then
\begin{equation}
\begin{aligned}
 \frac{|a(k,\eta;\ell,\xi)|}{|\ell|}
 &\le C\Lambda_t(\ell,\xi)\Lambda_t(m,\zeta)\kappa(\ell,\xi),\\
 \frac{|(\partial_t-\ell\partial_\xi-m\partial_\zeta)a(k,\eta;\ell,\xi)|}{|\ell|}
 &\le C\Lambda_t(\ell,\xi)\Lambda_t(m,\zeta)\kappa(\ell,\xi)
 [\kappa(k,\eta)+\kappa(\ell,\xi)+\kappa(m,\zeta)].
\end{aligned}
 \label{eq:lowA}
\end{equation}
If $8\sqrt{\ell^2+(\xi+\ell t)^2}>\sqrt{m^2+(\zeta+mt)^2}$, then 
\begin{align}
 &\frac{|h(k,\eta)-h(m,\zeta)|\,|a(k,\eta;\ell,\xi)|}
 {|\ell|\Lambda_t^N(k,\eta)\Lambda_t^N(\ell,\xi)\Lambda_t^N(m,\zeta)}
 \le C_N[\Lambda_t^{1-N}(k,\eta)+\Lambda_t^{1-N}(m,\zeta)]
 (1+\av t\kappa(\ell,\xi)),\label{eq:reactioncubic}\\
 &\frac{|h(k,\eta)-h(m,\zeta)|\,
 |(\partial_t-\ell\partial_\xi-m\partial_\zeta)a(k,\eta;\ell,\xi)|}
 {|\ell|\Lambda_t^N(k,\eta)\Lambda_t^N(\ell,\xi)\Lambda_t^N(m,\zeta)}\nonumber\\
 &\le C_N[\Lambda_t^{1-N}(k,\eta)+\Lambda_t^{1-N}(m,\zeta)]
 (1+\av t\kappa(\ell,\xi))
 [\kappa(k,\eta)+\kappa(\ell,\xi)+\kappa(m,\zeta)].
 \label{eq:reactionkerneldot}
\end{align}
\end{lemma}
\begin{proof}
The identities $c\sigma-d\omega=c(\sigma+dt)-d(\omega+ct)$ and $|c|\ge1$ give
\begin{align*}
 \frac{|c\sigma-d\omega|}{\sqrt{c^2+\omega^2}}
 &\le |\sigma+dt|\kappa(c,\omega)
       +|d|\frac{|\omega|+t|c|}{\sqrt{c^2+\omega^2}}
 \le C\Lambda_t(d,\sigma)(1+\av t\kappa(c,\omega)),\\
 \frac{\Lambda_t(c,\omega)}{\sqrt{c^2+\omega^2}}
 &\le\frac{\sqrt{1+c^2+\omega^2}+t|c|}{\sqrt{c^2+\omega^2}}
 \le\sqrt2+t\kappa(c,\omega).
\end{align*}
This proves \eqref{eq:elementarysymbols}. Next,
\begin{align*}
 |\ell\zeta-m\xi|
 &=|\ell(\zeta+mt)-m(\xi+\ell t)|
 \le\Lambda_t(\ell,\xi)\Lambda_t(m,\zeta),\\
 \frac{|a(k,\eta;\ell,\xi)|}{|\ell|}
 &\le\frac{|\ell\zeta-m\xi|}{|\ell|\sqrt{\ell^2+\xi^2}}
 \le\Lambda_t(\ell,\xi)\Lambda_t(m,\zeta)\kappa(\ell,\xi).
\end{align*}
Due to $k=\ell+m$ and $\eta=\xi+\zeta$, we have $\ell\zeta-m\xi=k\zeta-m\eta$ and thus
\[
 \frac{|\ell\xi|}{\ell^2+\xi^2}\le\kappa(\ell,\xi),\qquad
 \frac{|\ell\zeta-m\xi|}{\sqrt{k^2+\eta^2}\sqrt{m^2+\zeta^2}}\le1,
 \qquad \kappa^2\le\kappa.
\]
It follows from \eqref{eq:Adot} that
\begin{align*}
 \frac{|(\partial_t-\ell\partial_\xi-m\partial_\zeta)a(k,\eta;\ell,\xi)|}{|\ell|}
 \le\frac{|\ell\zeta-m\xi|}{|\ell|\sqrt{\ell^2+\xi^2}}
 [\kappa(k,\eta)+\kappa(\ell,\xi)+\kappa(m,\zeta)],
\end{align*}
which proves \eqref{eq:lowA}.
If $8\sqrt{\ell^2+(\xi+\ell t)^2}>\sqrt{m^2+(\zeta+mt)^2}$, then we have
\[
 \Lambda_t(m,\zeta)\le8\Lambda_t(\ell,\xi),\qquad
 \Lambda_t(k,\eta)\le9\Lambda_t(\ell,\xi).
\]
Moreover, $\ell\zeta-m\xi=\ell\eta-k\xi$, so it follows from \eqref{eq:elementarysymbols} that
\begin{align*}
 \frac{|\ell\zeta-m\xi|}{|\ell|\sqrt{\ell^2+\xi^2}}
 \le C\Lambda_t(m,\zeta)(1+\av t\kappa(\ell,\xi)),\qquad
 \frac{|\ell\zeta-m\xi|}{|\ell|\sqrt{\ell^2+\xi^2}}
 \le C\Lambda_t(k,\eta)(1+\av t\kappa(\ell,\xi)).
\end{align*}
Since $0<\theta\le1$ and $h=\theta\Lambda_t^{2N}$, these yield
\begin{align*}
 &\frac{h(k,\eta)|a(k,\eta;\ell,\xi)|}
 {|\ell|\Lambda_t^N(k,\eta)\Lambda_t^N(\ell,\xi)\Lambda_t^N(m,\zeta)}
 \le C\Big(\frac{\Lambda_t(k,\eta)}{\Lambda_t(\ell,\xi)}\Big)^N
 \Lambda_t^{1-N}(m,\zeta)(1+\av t\kappa(\ell,\xi))
 \le C_N\Lambda_t^{1-N}(m,\zeta)(1+\av t\kappa(\ell,\xi)),\\
 &\frac{h(m,\zeta)|a(k,\eta;\ell,\xi)|}
 {|\ell|\Lambda_t^N(k,\eta)\Lambda_t^N(\ell,\xi)\Lambda_t^N(m,\zeta)}
 \le C\Big(\frac{\Lambda_t(m,\zeta)}{\Lambda_t(\ell,\xi)}\Big)^N
 \Lambda_t^{1-N}(k,\eta)(1+\av t\kappa(\ell,\xi))
 \le C_N\Lambda_t^{1-N}(k,\eta)(1+\av t\kappa(\ell,\xi)).
\end{align*}
Adding the above estimates proves \eqref{eq:reactioncubic}. Applying a similar argument to $(\partial_t-\ell\partial_\xi-m\partial_\zeta)a$ proves \eqref{eq:reactionkerneldot}. 
\end{proof}

\subsection{Estimates for $\mathcal B$, $I_1$, and $I_3$}

\begin{proposition}
\label{prop:cubic}
Let \(k=\ell+m\) and \(\eta=\xi+\zeta\). With the notation of Section~\ref{sec:equations}, let $(u,b)$ be a smooth solution of the MHD equations \eqref{eq:perturbation} on $[0,T_1]$. Then we have  
\begin{align}
 |\cB(t)|&\le \frac{C_N}{|\beta|}\av t e^{-\frac1{16}\rho^{1/3} t} \Big(\sum_s\|F^s(t)\|_2^2\Big)^{3/2}, \qquad 0\le t\le T_1,\label{eq:Bbound}\\
 |\cB(0)|&\le \frac{C_N}{|\beta|}\varepsilon^3.\label{eq:Binitial}
\end{align}
Moreover, the nonlinear reactions \eqref{eq:R3} and \eqref{eq:shearterm} satisfy
\begin{align}
&\int_0^{T_1}\bigg|\sum_s\Rea\sum_{\substack{\ell,m\in\Z\\\ell\ne0}}\int_{\R^{2}}\frac{\ii}{4s\beta\ell}
  (\partial_t-\ell\partial_{\xi}-m\partial_{\zeta})
  \bigl[e^{\frac18\rho^{1/3} t\mathbf 1_{k\ne0}\mathbf 1_{m\ne0}}(h(k,\eta)-h(m,\zeta))a(k,\eta;\ell,\xi)\bigr]\nonumber\\
&\qquad\times\overline{f^s(k,\eta)}f^{-s}(\ell,\xi)f^s(m,\zeta)
 \frac{\dd\xi\dd\zeta}{(2\pi)^{2}}\bigg|\dd t\nonumber\\
&+\int_0^{T_1}\bigg|\sum_s\Rea\sum_{\substack{\ell,m\in\Z\\\ell\ne0}}\int_{\R^{2}}\frac{\ii e^{\frac18\rho^{1/3} t\mathbf 1_{k\ne0}\mathbf 1_{m\ne0}}(h(k,\eta)-h(m,\zeta))a(k,\eta;\ell,\xi)}{4s\beta\ell}\nonumber\\
&\qquad\times\bigg[\frac{k\eta}{(k^2+\eta^2)}\overline{f^{-s}(k,\eta)}f^{-s}(\ell,\xi)f^s(m,\zeta) +\frac{\ell\xi}{(\ell^2+\xi^2)}\overline{f^s(k,\eta)}f^s(\ell,\xi)f^s(m,\zeta)\nonumber\\
&\qquad\quad +\frac{m\zeta}{(m^2+\zeta^2)}\overline{f^s(k,\eta)}f^{-s}(\ell,\xi)f^{-s}(m,\zeta)\bigg]\frac{\dd\xi\dd\zeta}{(2\pi)^{2}}\bigg|\dd t
\nonumber\\
\le{}&\frac{C_N}{|\beta|}\rho^{-1/3}\mathcal E(T_1)^{3/2}.
\label{eq:R3bound}
\end{align}
\end{proposition}

\begin{proof}[Proof of Proposition~\ref{prop:cubic}]
Recalling \eqref{eq:B} with $k=\ell+m$ and $\eta=\zeta+\xi$, since $\ell\ne0$,
\begin{align}\label{eq:14}
\frac{e^{\rho^{1/3}t\mathbf 1_{k\ne0}\mathbf 1_{m\ne0}/8}}
 {e^{\rho^{1/3}t(\mathbf 1_{k\ne0}+\mathbf 1_{\ell\ne0}+\mathbf 1_{m\ne0})/16}}
 =\begin{cases}
 e^{-\rho^{1/3}t/16},&km\ne0,\\
 e^{-\rho^{1/3}t/8},&km=0.
 \end{cases}
\end{align}
Hence expressing the three factors $f^s$ in terms of $F^s$ leaves at most $e^{-\rho^{1/3}t/16}$ in \eqref{eq:B}, \eqref{eq:R3}, \eqref{eq:shearterm}, and \eqref{eq:Rnu}.

\emph{Estimate on $\cB(t)$.} First suppose that $8\sqrt{\ell^2+(\xi+\ell t)^2}\le\sqrt{m^2+(\zeta+mt)^2}$.
We deduce from \eqref{eq:tkappa}, \eqref{eq:FDweak} (with $j=m$, $\alpha=\zeta$, and $p=\gamma=0$), and \eqref{eq:lowA} that
\begin{align}\label{eq:12}
 \frac{|h(k,\eta)-h(m,\zeta)|\,|a(k,\eta;\ell,\xi)|}
 {|\ell|\Lambda_t^N(k,\eta)\Lambda_t^N(\ell,\xi)\Lambda_t^N(m,\zeta)}
 \le C_N\Lambda_t^{2-N}(\ell,\xi)\av t^2\kappa(\ell,\xi)
 \le C_N\Lambda_t^{3-N}(\ell,\xi)\av t.
\end{align}
Next suppose that $8\sqrt{\ell^2+(\xi+\ell t)^2}>\sqrt{m^2+(\zeta+mt)^2}$. It follows from \eqref{eq:reactioncubic} that
\begin{align}\label{eq:13}
\frac{|h(k,\eta)-h(m,\zeta)|\,|a(k,\eta;\ell,\xi)|}
 {|\ell|\Lambda_t^N(k,\eta)\Lambda_t^N(\ell,\xi)\Lambda_t^N(m,\zeta)}
 \le C_N[\Lambda_t^{1-N}(k,\eta)+\Lambda_t^{1-N}(m,\zeta)]
 (1+\av t\kappa(\ell,\xi)).
\end{align}
Using the definitions in \eqref{eq:fexplicit}, \eqref{eq:B}, and \eqref{eq:weightedamplitudes}, we obtain
\begin{align*}
|\cB(t)|\le &\frac{C_N}{|\beta|}\av t e^{-\rho^{1/3}t/16}
 \sum_s\sum_k\int_\R |F^s(k,\eta)|
 [|F^s|*(\Lambda_t^{3-N}|F^{-s}|)](k,\eta)\frac{\dd\eta}{2\pi}\\
 &+\frac{C_N}{|\beta|}\av t e^{-\rho^{1/3}t/16}\sum_s\sum_k\int_\R |F^s(k,\eta)|
 [|F^{-s}|*(\Lambda_t^{1-N}|F^s|)](k,\eta)\frac{\dd\eta}{2\pi}\\
 &+\frac{C_N}{|\beta|}\av t e^{-\rho^{1/3}t/16}\sum_s\sum_k\int_\R\Lambda_t^{1-N}(k,\eta)|F^s(k,\eta)|
 [|F^{-s}|*|F^s|](k,\eta)\frac{\dd\eta}{2\pi}\\
 \le &\frac{C_N}{|\beta|}\av t e^{-\rho^{1/3}t/16}
 \sum_s\Big(\|F^s\|_2^2\|\Lambda_t^{3-N}F^{-s}\|_1+\|F^s\|_2\|F^{-s}\|_2\|\Lambda_t^{1-N}F^s\|_1\Big)\\
 \le &\frac{C_N}{|\beta|}\av t e^{-\rho^{1/3}t/16}
                   \Bigl(\sum_s\|F^s\|_2^2\Bigr)^{3/2}.	
\end{align*}
Here we used for $r\le9$, $N>10$ and $s=\pm1$,
\begin{align}\label{eq:15}
 \|\Lambda_t^{r-N}F^s\|_1
 \le C_N\|F^s\|_2,\quad
 \|\Lambda_t^{r-N}\sqrt gF^s\|_1
 \le C_N\|\sqrt gF^s\|_2.
\end{align}
This proves \eqref{eq:Bbound}.
At $t=0$, $\sum_s\|F^s(0)\|_2^2=2\varepsilon^2$, so
\eqref{eq:Binitial} follows as well.

\emph{Estimates on $I_1$ and $I_3$.}
A direct computation gives 
\begin{align*}
\partial_t e^{\rho^{1/3}t\mathbf1_{k\ne0}\mathbf1_{m\ne0}/8}
 &=\frac{\rho^{1/3}}8\mathbf1_{k\ne0}\mathbf1_{m\ne0}
 e^{\rho^{1/3}t\mathbf1_{k\ne0}\mathbf1_{m\ne0}/8}.
\end{align*}
Then recalling \eqref{eq:frequencies},
\[
 \frac{|k\eta|}{k^2+\eta^2}\le\kappa(k,\eta),\qquad
 \frac{|\ell\xi|}{\ell^2+\xi^2}\le\kappa(\ell,\xi),\qquad
 \frac{|m\zeta|}{m^2+\zeta^2}\le\kappa(m,\zeta).
\]
Hence it follows from \eqref{eq:R3} and \eqref{eq:shearterm} that
\[
 \int_0^{T_1}\bigl(|I_1(t)|+|I_3(t)|\bigr)\dd t
 \le J_1+J_2+J_3+J_4,
\]
where 
\begin{align*}
 J_1&=\int_0^{T_1}\sum_{\sigma_1,\sigma_2,\sigma_3=\pm1}
 \sum_{\substack{\ell,m\in\Z\\\ell\ne0}}\int_{\R^2}
 \frac{e^{\rho^{1/3}t\mathbf1_{k\ne0}\mathbf1_{m\ne0}/8}}{4|\beta\ell|}|h(k,\eta)-h(m,\zeta)|\,|(\partial_t-\ell\partial_\xi-m\partial_\zeta)a(k,\eta;\ell,\xi)|\\*
 &\qquad\times |f^{\sigma_1}(k,\eta)||f^{\sigma_2}(\ell,\xi)||f^{\sigma_3}(m,\zeta)|\frac{\dd\xi\dd\zeta}{(2\pi)^2}\dd t,\\
 J_2&=\int_0^{T_1}\sum_{\sigma_1,\sigma_2,\sigma_3=\pm1}
 \sum_{\substack{\ell,m\in\Z\\\ell\ne0}}\int_{\R^2}
 \frac{e^{\rho^{1/3}t\mathbf1_{k\ne0}\mathbf1_{m\ne0}/8}}{4|\beta\ell|}|h(k,\eta)-h(m,\zeta)||a(k,\eta;\ell,\xi)|\\*
 &\qquad\times[\kappa(k,\eta)+\kappa(\ell,\xi)+\kappa(m,\zeta)]
 |f^{\sigma_1}(k,\eta)||f^{\sigma_2}(\ell,\xi)||f^{\sigma_3}(m,\zeta)|\frac{\dd\xi\dd\zeta}{(2\pi)^2}\dd t,\\
 J_3&=\int_0^{T_1}\sum_{\sigma_1,\sigma_2,\sigma_3=\pm1}
 \sum_{\substack{\ell,m\in\Z\\\ell\ne0}}\int_{\R^2}
 \frac{e^{\rho^{1/3}t\mathbf1_{k\ne0}\mathbf1_{m\ne0}/8}}{4|\beta\ell|}|(\partial_t-k\partial_\eta)h(k,\eta)-(\partial_t-m\partial_\zeta)h(m,\zeta)|\\*
 &\qquad\times |a(k,\eta;\ell,\xi)||f^{\sigma_1}(k,\eta)||f^{\sigma_2}(\ell,\xi)||f^{\sigma_3}(m,\zeta)|\frac{\dd\xi\dd\zeta}{(2\pi)^2}\dd t,\\
 J_4&=\int_0^{T_1}\sum_{\sigma_1,\sigma_2,\sigma_3=\pm1}
 \sum_{\substack{\ell,m\in\Z\\\ell\ne0}}\int_{\R^2}
 \frac{e^{\rho^{1/3}t\mathbf1_{k\ne0}\mathbf1_{m\ne0}/8}}{4|\beta\ell|}\frac{\rho^{1/3}}8\mathbf1_{k\ne0}\mathbf1_{m\ne0}|h(k,\eta)-h(m,\zeta)||a(k,\eta;\ell,\xi)|\\*
 &\qquad\times |f^{\sigma_1}(k,\eta)||f^{\sigma_2}(\ell,\xi)||f^{\sigma_3}(m,\zeta)|\frac{\dd\xi\dd\zeta}{(2\pi)^2}\dd t.
\end{align*}
For $J_1$ and $J_2$, first suppose that $8\sqrt{\ell^2+(\xi+\ell t)^2}\le\sqrt{m^2+(\zeta+mt)^2}$.
In this case, we use \eqref{eq:tkappa}, \eqref{eq:FD1} (with $j=m$, $\alpha=\zeta$, and $p=\gamma=0$), and \eqref{eq:lowA} to obtain
\begin{equation*}
\begin{aligned}
 &e^{-\rho^{1/3}t/16}
 \frac{|h(k,\eta)-h(m,\zeta)|\,
 |(\partial_t-\ell\partial_\xi-m\partial_\zeta)a(k,\eta;\ell,\xi)|}
 {|\ell|\Lambda_t^N(k,\eta)\Lambda_t^N(\ell,\xi)\Lambda_t^N(m,\zeta)}\\
 &+e^{-\rho^{1/3}t/16}\frac{|h(k,\eta)-h(m,\zeta)|\,|a(k,\eta;\ell,\xi)|}
 {|\ell|\Lambda_t^N(k,\eta)\Lambda_t^N(\ell,\xi)\Lambda_t^N(m,\zeta)}[\kappa(k,\eta)+\kappa(\ell,\xi)+\kappa(m,\zeta)]\\
 \le{} &C_N\Lambda_t^{2-N}(\ell,\xi)e^{-\rho^{1/3}t/16}
 [\av t\av{\rho^{1/3}t}+\av t^2\sqrt{g(k,\eta)+g(m,\zeta)}]
 \kappa(\ell,\xi)
 [\kappa(k,\eta)+\kappa(\ell,\xi)+\kappa(m,\zeta)]\\
\le{} &C_N\Lambda_t^{3-N}(\ell,\xi)e^{-\rho^{1/3}t/16}
 [\av{\rho^{1/3}t}+\av t\sqrt{g(k,\eta)+g(m,\zeta)}]
 [\kappa(k,\eta)+\kappa(\ell,\xi)+\kappa(m,\zeta)]\\
\le{} &C_N\av{\rho^{1/3}t} e^{-\rho^{1/3}t/16}\Lambda_t^{3-N}(\ell,\xi)\big[\sqrt{g(k,\eta)}+\sqrt{g(\ell,\xi)}+\sqrt{g(m,\zeta)}\big]\\
&+C_N \av te^{-\rho^{1/3}t/16}\Big[\Lambda_t^{3-N}(\ell,\xi)\bigl[\sqrt{g(k,\eta)g(\ell,\xi)}+\sqrt{g(\ell,\xi)g(m,\zeta)}\bigr]+\Lambda_t^{8-N}(\ell,\xi)\sqrt{g(k,\eta)g(m,\zeta)}\Big],
\end{aligned}
\end{equation*}
where we also used \eqref{eq:twotransfer} in the last inequality to conclude
\begin{align*}
 &\sqrt{g(k,\eta)+g(m,\zeta)}[\kappa(k,\eta)+\kappa(m,\zeta)]
 \le C\Lambda_t^5(\ell,\xi)\sqrt{g(k,\eta)g(m,\zeta)}.
 \end{align*}
Then suppose that $8\sqrt{\ell^2+(\xi+\ell t)^2}>\sqrt{m^2+(\zeta+mt)^2}$. In this case, applying \eqref{eq:reactionkerneldot} gives
\begin{align*}
&e^{-\rho^{1/3}t/16}
 \frac{|h(k,\eta)-h(m,\zeta)|\,
 |(\partial_t-\ell\partial_\xi-m\partial_\zeta)a(k,\eta;\ell,\xi)|}
 {|\ell|\Lambda_t^N(k,\eta)\Lambda_t^N(\ell,\xi)\Lambda_t^N(m,\zeta)}\\
 &+e^{-\rho^{1/3}t/16}\frac{|h(k,\eta)-h(m,\zeta)|\,|a(k,\eta;\ell,\xi)|}
 {|\ell|\Lambda_t^N(k,\eta)\Lambda_t^N(\ell,\xi)\Lambda_t^N(m,\zeta)}[\kappa(k,\eta)+\kappa(\ell,\xi)+\kappa(m,\zeta)]\\
 \le{} &C_N e^{-\rho^{1/3}t/16}
 [\Lambda_t^{1-N}(m,\zeta)+\Lambda_t^{1-N}(k,\eta)]\\
 &\quad\times\bigl[\kappa(k,\eta)+\kappa(\ell,\xi)+\kappa(m,\zeta)+\av t\kappa(\ell,\xi)\kappa(k,\eta)
 +\av t\kappa(\ell,\xi)^2
 +\av t\kappa(\ell,\xi)\kappa(m,\zeta)\bigr]\\
 \le{} &C_N e^{-\rho^{1/3}t/16}
 \big[\Lambda_t^{1-N}(m,\zeta)+\Lambda_t^{1-N}(k,\eta)\big]\big[\sqrt{g(k,\eta)}+\sqrt{g(\ell,\xi)}+\sqrt{g(m,\zeta)}\big]\\
 &+C_N \av te^{-\rho^{1/3}t/16}\big[\Lambda_t^{1-N}(m,\zeta)+\Lambda_t^{1-N}(k,\eta)\big]\bigl[\sqrt{g(k,\eta)g(\ell,\xi)}
 +\sqrt{g(\ell,\xi)g(m,\zeta)}\bigr]\\
 &+C_N\av t e^{-\rho^{1/3}t/16}\big[\Lambda_t^{6-N}(m,\zeta)\sqrt{g(\ell,\xi)g(k,\eta)}+\Lambda_t^{6-N}(k,\eta)\sqrt{g(\ell,\xi)g(m,\zeta)}\big],
\end{align*}
where we used \eqref{eq:gprimitive} and 
\begin{align*}
 \Lambda_t^{1-N}(m,\zeta)\kappa(\ell,\xi)^2
 &\le C\Lambda_t^{6-N}(m,\zeta)\sqrt{g(\ell,\xi)g(k,\eta)},\\
 \Lambda_t^{1-N}(k,\eta)\kappa(\ell,\xi)^2
 &\le C\Lambda_t^{6-N}(k,\eta)\sqrt{g(\ell,\xi)g(m,\zeta)}.
\end{align*}
due to \eqref{eq:twotransfer}. In the second line we used $\ell=-m+k$, $\xi=-\zeta+\eta$, and $g(-m,-\zeta)=g(m,\zeta)$.
Expressing $f^s$ in terms of $F^s$ by \eqref{eq:weightedamplitudes} and using \eqref{eq:gprimitive} give
\begin{align*}
|\beta|(J_1+J_2)&\le \sum_{\sigma_1,\sigma_2,\sigma_3}C_N\int_0^{T_1} e^{-\rho^{1/3}t/16}\av{\rho^{1/3}t}
 \Bigl[\|\sqrt gF^{\sigma_1}\|_2\|\Lambda_t^{3-N}F^{\sigma_2}\|_1\|F^{\sigma_3}\|_2\\
 &\qquad+\|F^{\sigma_1}\|_2\|\Lambda_t^{3-N}\sqrt gF^{\sigma_2}\|_1\|F^{\sigma_3}\|_2
 +\|F^{\sigma_1}\|_2\|\Lambda_t^{3-N}F^{\sigma_2}\|_1\|\sqrt gF^{\sigma_3}\|_2\Bigr]\dd t\\
 &+\sum_{\sigma_1,\sigma_2,\sigma_3}C_N\int_0^{T_1}\av t e^{-\rho^{1/3}t/16}
 \Bigl[\|\sqrt gF^{\sigma_1}\|_2\|F^{\sigma_3}\|_2
       +\|F^{\sigma_1}\|_2\|\sqrt gF^{\sigma_3}\|_2\Bigr]
           \|\Lambda_t^{3-N}\sqrt gF^{\sigma_2}\|_1\dd t\\
 &+ \sum_{\sigma_1,\sigma_2,\sigma_3}C_N\int_0^{T_1}\av t e^{-\rho^{1/3}t/16}
 \|\sqrt gF^{\sigma_1}\|_2\|\Lambda_t^{8-N}F^{\sigma_2}\|_1
  \|\sqrt gF^{\sigma_3}\|_2\dd t\\
 &+ \sum_{\sigma_1,\sigma_2,\sigma_3}C_N\int_0^{T_1} e^{-\rho^{1/3}t/16}
 \Bigl[\|\sqrt gF^{\sigma_1}\|_2\|F^{\sigma_2}\|_2\|F^{\sigma_3}\|_2
       +\|F^{\sigma_1}\|_2\|\sqrt gF^{\sigma_2}\|_2\|F^{\sigma_3}\|_2\\
 &\hspace{48mm}+\|F^{\sigma_1}\|_2\|F^{\sigma_2}\|_2\|\sqrt gF^{\sigma_3}\|_2\Bigr]\dd t\\
 &+ \sum_{\sigma_1,\sigma_2,\sigma_3}C_N\int_0^{T_1}\av t e^{-\rho^{1/3}t/16}
 \|\sqrt gF^{\sigma_2}\|_2
 \Bigl[\|\sqrt gF^{\sigma_1}\|_2\|F^{\sigma_3}\|_2
       +\|F^{\sigma_1}\|_2\|\sqrt gF^{\sigma_3}\|_2\Bigr]\dd t.
\end{align*}
Here $\sigma_i\in\{-,+\}$. For $a,b\ge0$, we have
\[
  \|\av t^a\av{\rho^{1/3} t}^b e^{-\frac1{16}\rho^{1/3} t}\|_{L^p(0,\infty)}
       \le C_{a,b,p}\rho^{-a/3-1/(3p)},\qquad 1\le p\le\infty.
\]
Hence the term containing $\Lambda_t^{8-N}$ is bounded by
\begin{align*}
 &\int_0^{T_1}\av t e^{-\rho^{1/3}t/16}
   \|\sqrt gF^{\sigma_1}\|_2\|\Lambda_t^{8-N}F^{\sigma_2}\|_1
   \|\sqrt gF^{\sigma_3}\|_2\dd t\\
 \le{} & C_N\|\av t e^{-\rho^{1/3}t/16}\|_{L^\infty(0,\infty)}
   \|\sqrt gF^{\sigma_1}\|_{L^2_tL^2}
   \|F^{\sigma_2}\|_{L^\infty_tL^2}
   \|\sqrt gF^{\sigma_3}\|_{L^2_tL^2}\\
 \le{} & C_N\rho^{-1/3}\mathcal E(T_1)^{3/2}.
\end{align*}
Applying H\"older's inequality to the remaining terms, we obtain
\begin{align*}
J_1+J_2&\le C_N\frac{1}{|\beta|}\|e^{-\rho^{1/3}t/16}\av{\rho^{1/3}t}\|_{L^2_t}
                  \mathcal E(T_1)^{3/2}+ C_N\frac{1}{|\beta|}\|\av t e^{-\rho^{1/3}t/16}\|_{L^\infty_t}\mathcal E(T_1)^{3/2}\\
&\le C_N\frac{1}{|\beta|}\rho^{-1/6}\mathcal E(T_1)^{3/2}+C_N\frac{1}{|\beta|}\rho^{-1/3}\mathcal E(T_1)^{3/2}\\
&\le C_N\frac{1}{|\beta|}\rho^{-1/3}\mathcal E(T_1)^{3/2}.
\end{align*}
For $J_3$, since $(\partial_t-k\partial_\eta)\Lambda_t=0$, the definition of $h$ in \eqref{eq:weight} gives
\begin{equation}
(\partial_t-k\partial_\eta)h(k,\eta)
 =-2\bigl(7g(k,\eta)+(\partial_t-k\partial_\eta)V(k,\eta)\bigr)h(k,\eta).
\label{eq:coefficientderivatives}
\end{equation}
This gives the exact identity
\begin{equation}
\begin{aligned}
&(\partial_t-k\partial_{\eta})h(k,\eta)-(\partial_t-m\partial_{\zeta})h(m,\zeta)\\
=&-2(7g(k,\eta)+(\partial_t-k\partial_{\eta})V(k,\eta))(h(k,\eta)-h(m,\zeta))\\*
&-2\{7(g(k,\eta)-g(m,\zeta))+(\partial_t-k\partial_{\eta})V(k,\eta)-(\partial_t-m\partial_{\zeta})V(m,\zeta)\}h(m,\zeta).
\label{eq:hdotdifference}
\end{aligned}
\end{equation}
First suppose that $8\sqrt{\ell^2+(\xi+\ell t)^2}\le\sqrt{m^2+(\zeta+mt)^2}$. In this case, $\Lambda_t(k,\eta)\approx\Lambda_t(m,\zeta)$ due to $k=\ell+m$ and $\eta=\xi+\zeta$. For the first term on the right of \eqref{eq:hdotdifference}, it follows from \eqref{eq:physicalV}, \eqref{eq:FDweak} (with $j=m$, $\alpha=\zeta$, and $p=\gamma=0$), and \eqref{eq:lowA} that
\begin{align*}
\big|2(7g(k,\eta)+(\partial_t-k\partial_{\eta})V(k,\eta))\big|\frac{|a(k,\eta;\ell,\xi)||h(k,\eta)-h(m,\zeta)|}{|\ell|\Lambda_t^N(k,\eta)\Lambda_t^N(\ell,\xi)\Lambda_t^N(m,\zeta)}
\le C_N\Lambda_t^{2-N}(\ell,\xi)\av t^2\kappa(\ell,\xi)\big[g(k,\eta)+\rho^{1/3}\big].	
\end{align*}
For the second term on the right of \eqref{eq:hdotdifference}, using $h(m,\zeta)\le\Lambda_t^{2N}(m,\zeta)$,
\eqref{eq:gdiff}, \eqref{eq:vdiff} and \eqref{eq:lowA}, we obtain
\begin{align*}
&\big|2\{7(g(k,\eta)-g(m,\zeta))+(\partial_t-k\partial_{\eta})V(k,\eta)-(\partial_t-m\partial_{\zeta})V(m,\zeta)\}h(m,\zeta)\big|
\frac{|a(k,\eta;\ell,\xi)|}{|\ell|\Lambda_t^N(k,\eta)\Lambda_t^N(\ell,\xi)\Lambda_t^N(m,\zeta)}\\
\le{} &C_N\Lambda_t^{1-N}(\ell,\xi)\Lambda_t(m,\zeta)\kappa(\ell,\xi)
 \bigl[|g(k,\eta)-g(m,\zeta)|
 +|(\partial_t-k\partial_\eta)V(k,\eta)
             -(\partial_t-m\partial_\zeta)V(m,\zeta)|\bigr]\\
 \le{} &C_N\Lambda_t^{2-N}(\ell,\xi)\kappa(\ell,\xi)
 \bigl[1+\av t^2(g(k,\eta)+g(m,\zeta))+\av{\rho^{1/3}t}^2\bigr].
\end{align*}
Thus, by \eqref{eq:tkappa} and \eqref{eq:twotransfer}, we get
\begin{equation*}
\begin{aligned}
&\frac{\big|(\partial_t-k\partial_{\eta})h(k,\eta)-(\partial_t-m\partial_{\zeta})h(m,\zeta)\big|\,|a(k,\eta;\ell,\xi)|}{|\ell|\Lambda_t^N(k,\eta)\Lambda_t^N(\ell,\xi)\Lambda_t^N(m,\zeta)}\\
\le{} &C_N\Lambda_t^{3-N}(\ell,\xi)
 \bigl[\av t(g(k,\eta)+g(m,\zeta))+\rho^{1/3}\av t
       +\kappa(\ell,\xi)(1+\av{\rho^{1/3}t}^2)\bigr]\\
\le{} &C_N\Lambda_t^{3-N}(\ell,\xi)
 \bigl[\av t\Lambda_t^5(\ell,\xi)\sqrt{g(k,\eta)g(m,\zeta)}+\rho^{1/3}\av t
       +\kappa(\ell,\xi)(1+\av{\rho^{1/3}t}^2)\bigr].
\end{aligned}
\end{equation*}
Then suppose that $8\sqrt{\ell^2+(\xi+\ell t)^2}>\sqrt{m^2+(\zeta+mt)^2}$. 
We use \eqref{eq:coefficientderivatives} to estimate the two terms in \eqref{eq:hdotdifference} separately. The two computations in the proof of \eqref{eq:reactioncubic} give
\begin{align*}
 &\frac{|(\partial_t-k\partial_\eta)h(k,\eta)
 -(\partial_t-m\partial_\zeta)h(m,\zeta)|\,|a(k,\eta;\ell,\xi)|}
 {|\ell|\Lambda_t^N(k,\eta)\Lambda_t^N(\ell,\xi)\Lambda_t^N(m,\zeta)}\\
 \le{} &C_N\bigl[
 \Lambda_t^{1-N}(m,\zeta)(g(k,\eta)+\rho^{1/3})
 +\Lambda_t^{1-N}(k,\eta)(g(m,\zeta)+\rho^{1/3})\bigr]
 (1+\av t\kappa(\ell,\xi)).
\end{align*}
Here we also used $0\le(\partial_t-k\partial_\eta)V\le\rho^{1/3}$ due to \eqref{eq:physicalV}.
Since $g\le C$ and $\kappa\le\min\{\sqrt g,1\}$, recalling the definition of $f^s$ and $F^s$ in \eqref{eq:fexplicit} and \eqref{eq:weightedamplitudes}, we conclude
\begin{align*}
|\beta|J_3&\le \sum_{\sigma_1,\sigma_2,\sigma_3}C_N\int_0^{T_1} e^{-\rho^{1/3}t/16}\av t
 \|\sqrt gF^{\sigma_1}\|_2\|\Lambda_t^{8-N}F^{\sigma_2}\|_1\|\sqrt gF^{\sigma_3}\|_2\dd t\\
 &+\sum_{\sigma_1,\sigma_2,\sigma_3}C_N\int_0^{T_1} e^{-\rho^{1/3}t/16}
 [\rho^{1/3}\av t+1+\av{\rho^{1/3}t}^2]
 \|F^{\sigma_1}\|_2\|\Lambda_t^{3-N}F^{\sigma_2}\|_1\|F^{\sigma_3}\|_2\dd t\\
 &+\sum_{\sigma_1,\sigma_2,\sigma_3}C_N\rho^{1/3}\int_0^{T_1} e^{-\rho^{1/3}t/16}(1+\av t)
 \|F^{\sigma_2}\|_2\Bigl[\|F^{\sigma_1}\|_2\|\Lambda_t^{1-N}F^{\sigma_3}\|_1
             +\|\Lambda_t^{1-N}F^{\sigma_1}\|_1\|F^{\sigma_3}\|_2\Bigr]\dd t\\
 &+\sum_{\sigma_1,\sigma_2,\sigma_3}C_N\int_0^{T_1} e^{-\rho^{1/3}t/16}\|F^{\sigma_2}\|_2
 \Bigl[\|\sqrt gF^{\sigma_1}\|_2\|\Lambda_t^{1-N}F^{\sigma_3}\|_1
       +\|\Lambda_t^{1-N}F^{\sigma_1}\|_1\|\sqrt gF^{\sigma_3}\|_2\Bigr]\dd t\\
 &+\sum_{\sigma_1,\sigma_2,\sigma_3}C_N\int_0^{T_1}\av t e^{-\rho^{1/3}t/16}\|\sqrt gF^{\sigma_2}\|_2
 \Bigl[\|\sqrt gF^{\sigma_1}\|_2\|\Lambda_t^{1-N}F^{\sigma_3}\|_1
       +\|\Lambda_t^{1-N}F^{\sigma_1}\|_1\|\sqrt gF^{\sigma_3}\|_2\Bigr]\dd t.
\end{align*}
By virtue of \eqref{eq:15}, we further deduce
\begin{align*}
J_3\le{}&\frac{C_N}{|\beta|}\Bigl[
 \|\av t e^{-\rho^{1/3}t/16}\|_{L^\infty(0,\infty)}
 +\|e^{-\rho^{1/3}t/16}\|_{L^2(0,\infty)}\\
 &\qquad+\rho^{1/3}\|(1+\av t)e^{-\rho^{1/3}t/16}\|_{L^1(0,\infty)}
 +\|(1+\av{\rho^{1/3}t}^2)e^{-\rho^{1/3}t/16}\|_{L^1(0,\infty)}
 \Bigr]\mathcal E(T_1)^{3/2}\\
\le{}&\frac{C_N}{|\beta|}
 [\rho^{-1/3}+\rho^{-1/6}+\rho^{1/3}\rho^{-2/3}+\rho^{-1/3}]
 \mathcal E(T_1)^{3/2}\le\frac{C_N}{|\beta|}\rho^{-1/3}\mathcal E(T_1)^{3/2}.
\end{align*}
For $J_4$, arguing as in the proof of \eqref{eq:Bbound}, we obtain 
\[
 J_4\le\frac{C_N}{|\beta|}\rho^{1/3}
    \int_0^{T_1}\av t e^{-\rho^{1/3}t/16}\dd t\,
                         \mathcal E(T_1)^{3/2}
 \le\frac{C_N}{|\beta|}\rho^{-1/3}\mathcal E(T_1)^{3/2}.
\]
Finally, adding the four estimates yields
\[
 \int_0^{T_1}\big(|I_1(t)|+|I_3(t)|\big)\dd t
 \le J_1+J_2+J_3+J_4
 \le\frac{C_N}{|\beta|}\rho^{-1/3}\mathcal E(T_1)^{3/2},
\]
which proves \eqref{eq:R3bound}.
\end{proof}

\subsection{\texorpdfstring{Estimate for the diffusion term $I_4$}{Estimate for the diffusion term I4}}

\begin{proposition}
\label{lem:mixeddiffusion}
Let \(k=\ell+m\) and \(\eta=\xi+\zeta\). With the notation of Section~\ref{sec:equations}, let $(u,b)$ be a smooth solution of the MHD equations \eqref{eq:perturbation} on $[0,T_1]$. Then the nonlinear reaction \eqref{eq:Rnu} satisfies 
\begin{align}
&\int_0^{T_1}\bigg|-\sum_s\Rea\sum_{\substack{\ell,m\in\Z\\\ell\ne0}}\int_{\R^{2}}
 \frac{\ii e^{\rho^{1/3}t\mathbf 1_{k\ne0}\mathbf 1_{m\ne0}/8}(h(k,\eta)-h(m,\zeta))a(k,\eta;\ell,\xi)}
      {4s\beta\ell}\nonumber\\
&\hspace{10mm}\times\biggl[\frac{\nu+\mu}{2}((k^2+\eta^2)+(\ell^2+\xi^2)+(m^2+\zeta^2))\overline{f^s(k,\eta)}f^{-s}(\ell,\xi)f^s(m,\zeta)\nonumber\\
&\hspace{16mm}+\frac{\nu-\mu}{2}\bigl( (k^2+\eta^2)\overline{f^{-s}(k,\eta)}f^{-s}(\ell,\xi)f^s(m,\zeta) +(\ell^2+\xi^2)\overline{f^s(k,\eta)}f^s(\ell,\xi)f^s(m,\zeta)\nonumber\\
&\hspace{16mm}+(m^2+\zeta^2)\overline{f^s(k,\eta)}f^{-s}(\ell,\xi)f^{-s}(m,\zeta)\bigr) \bigg]\frac{\dd\xi\dd\zeta}{(2\pi)^{2}}\bigg|\dd t
\nonumber\\
\le{} &C_N\frac{\nu+\mu}{|\beta|\sqrt{\nu\mu}}\rho^{-1/3}\mathcal E(T_1)^{3/2}.
\label{eq:mixeddiffusionbound}
\end{align}
\end{proposition}
\begin{proof}
Using \eqref{eq:fexplicit}, we can rewrite \eqref{eq:Rnu} as follows:
\begin{align}
I_4=&-\Rea\sum_{\ell\ne0,m}\int_{\R^2}
 \frac{e^{\rho^{1/3}t\mathbf1_{k\ne0}\mathbf1_{m\ne0}/8}
 (h(k,\eta)-h(m,\zeta))a(k,\eta;\ell,\xi)}{2\beta\ell}\nonumber\\
&\times\Bigl[
 (\mu(k^2+\eta^2)+\nu(\ell^2+\xi^2)+\nu(m^2+\zeta^2))
\bigl(\overline{e(k,\eta)\cdot\widehat b(k,\eta)}\bigr)\bigl(e(\ell,\xi)\cdot\widehat u(\ell,\xi)\bigr)\bigl(e(m,\zeta)\cdot\widehat u(m,\zeta)\bigr)\nonumber\\
&\quad-(\nu(k^2+\eta^2)+\mu(\ell^2+\xi^2)+\nu(m^2+\zeta^2))
\bigl(\overline{e(k,\eta)\cdot\widehat u(k,\eta)}\bigr)\bigl(e(\ell,\xi)\cdot\widehat b(\ell,\xi)\bigr)\bigl(e(m,\zeta)\cdot\widehat u(m,\zeta)\bigr)\nonumber\\
&\quad+(\nu(k^2+\eta^2)+\nu(\ell^2+\xi^2)+\mu(m^2+\zeta^2))
\bigl(\overline{e(k,\eta)\cdot\widehat u(k,\eta)}\bigr)\bigl(e(\ell,\xi)\cdot\widehat u(\ell,\xi)\bigr)\bigl(e(m,\zeta)\cdot\widehat b(m,\zeta)\bigr)\nonumber\\
&\quad-\mu(k^2+\eta^2+\ell^2+\xi^2+m^2+\zeta^2)
\bigl(\overline{e(k,\eta)\cdot\widehat b(k,\eta)}\bigr)\bigl(e(\ell,\xi)\cdot\widehat b(\ell,\xi)\bigr)\bigl(e(m,\zeta)\cdot\widehat b(m,\zeta)\bigr)\Bigr]
 \frac{\dd\xi\dd\zeta}{(2\pi)^2}.
\label{eq:mixedphysicalexpansion}
\end{align}
The conditions \(k=\ell+m\), \(\eta=\xi+\zeta\) then give
\begin{align*}
\sqrt{k^2+\eta^2}&\le \sqrt{\ell^2+\xi^2}+\sqrt{m^2+\zeta^2},\\
\sqrt{\ell^2+\xi^2}&\le \sqrt{k^2+\eta^2}+\sqrt{m^2+\zeta^2},\\
\sqrt{m^2+\zeta^2}&\le \sqrt{k^2+\eta^2}+\sqrt{\ell^2+\xi^2}.
\end{align*}
Multiplying them by \(\mu \sqrt{k^2+\eta^2},\mu \sqrt{\ell^2+\xi^2},\mu \sqrt{m^2+\zeta^2}\), respectively, and adding, we obtain
\[
\mu(k^2+\eta^2+\ell^2+\xi^2+m^2+\zeta^2)
 \le 2\mu\bigl[\sqrt{k^2+\eta^2}\sqrt{\ell^2+\xi^2}
 +\sqrt{k^2+\eta^2}\sqrt{m^2+\zeta^2}
 +\sqrt{\ell^2+\xi^2}\sqrt{m^2+\zeta^2}\bigr].
\]
On the other hand, multiplying them by $\mu\sqrt{k^2+\eta^2}$, $\nu\sqrt{\ell^2+\xi^2}$, and $\nu\sqrt{m^2+\zeta^2}$, respectively, gives
\begin{align*}
&\mu(k^2+\eta^2)+\nu(\ell^2+\xi^2)+\nu(m^2+\zeta^2)\\
\le{} &(\nu+\mu)\sqrt{k^2+\eta^2}\sqrt{\ell^2+\xi^2}
 +(\nu+\mu)\sqrt{k^2+\eta^2}\sqrt{m^2+\zeta^2}
 +2\nu\sqrt{\ell^2+\xi^2}\sqrt{m^2+\zeta^2}\\
\le{} &\frac{\nu+\mu}{\sqrt{\nu\mu}}\Bigl[
 \sqrt\mu\sqrt{k^2+\eta^2}\sqrt\nu\sqrt{\ell^2+\xi^2}+\sqrt\mu\sqrt{k^2+\eta^2}\sqrt\nu\sqrt{m^2+\zeta^2}
 +\nu\sqrt{\ell^2+\xi^2}\sqrt{m^2+\zeta^2}\Bigr].
\end{align*}
Here we also used $(\nu+\mu)/\sqrt{\nu\mu}\ge2$. Similarly,
\begin{align*}
 &\nu(k^2+\eta^2)+\mu(\ell^2+\xi^2)+\nu(m^2+\zeta^2)\\
 \le{}&\frac{\nu+\mu}{\sqrt{\nu\mu}}\bigl[
 \sqrt\nu\sqrt{k^2+\eta^2}\sqrt\mu\sqrt{\ell^2+\xi^2}
 +\nu\sqrt{k^2+\eta^2}\sqrt{m^2+\zeta^2}
 +\sqrt\mu\sqrt{\ell^2+\xi^2}\sqrt\nu\sqrt{m^2+\zeta^2}\bigr],\\
 &\nu(k^2+\eta^2)+\nu(\ell^2+\xi^2)+\mu(m^2+\zeta^2)\\
 \le{}&\frac{\nu+\mu}{\sqrt{\nu\mu}}\bigl[
 \nu\sqrt{k^2+\eta^2}\sqrt{\ell^2+\xi^2}
 +\sqrt\nu\sqrt{k^2+\eta^2}\sqrt\mu\sqrt{m^2+\zeta^2}
 +\sqrt\nu\sqrt{\ell^2+\xi^2}\sqrt\mu\sqrt{m^2+\zeta^2}\bigr].
\end{align*}
First suppose that $8\sqrt{\ell^2+(\xi+\ell t)^2}\le\sqrt{m^2+(\zeta+mt)^2}$. It follows from \eqref{eq:14} and \eqref{eq:12} that
\begin{equation*}
\frac{e^{\rho^{1/3}t\mathbf 1_{k\ne0}\mathbf 1_{m\ne0}/8}}
 {e^{\rho^{1/3}t(\mathbf 1_{k\ne0}+\mathbf 1_{\ell\ne0}+\mathbf 1_{m\ne0})/16}}\frac{|h(k,\eta)-h(m,\zeta)|\,|a(k,\eta;\ell,\xi)|}
 {|\ell|\Lambda_t^N(k,\eta)\Lambda_t^N(\ell,\xi)\Lambda_t^N(m,\zeta)}
 \le C_Ne^{-\rho^{1/3}t/16}\av t\Lambda_t^{3-N}(\ell,\xi).
\end{equation*}
Then suppose that $8\sqrt{\ell^2+(\xi+\ell t)^2}>\sqrt{m^2+(\zeta+mt)^2}$. We deduce from \eqref{eq:14}, \eqref{eq:13}, and $\kappa\le1$ that
\begin{align*}
 \frac{e^{\rho^{1/3}t\mathbf 1_{k\ne0}\mathbf 1_{m\ne0}/8}}
 {e^{\rho^{1/3}t(\mathbf 1_{k\ne0}+\mathbf 1_{\ell\ne0}+\mathbf 1_{m\ne0})/16}}\frac{|h(k,\eta)-h(m,\zeta)|\,|a(k,\eta;\ell,\xi)|}
 {|\ell|\Lambda_t^N(k,\eta)\Lambda_t^N(\ell,\xi)\Lambda_t^N(m,\zeta)}
 \le C_Ne^{-\rho^{1/3}t/16}\av t\,[\Lambda_t^{1-N}(m,\zeta)+\Lambda_t^{1-N}(k,\eta)].
\end{align*}
To use the norms of $\mathcal Au$ and $\mathcal Ab$, we apply \eqref{eq:Af} and \eqref{eq:thetabounds} to conclude
\[
\frac{e^{\rho^{1/3}t(\mathbf1_{k\ne0}+\mathbf1_{\ell\ne0}+\mathbf1_{m\ne0})/16}
 \Lambda_t^N(k,\eta)\Lambda_t^N(\ell,\xi)\Lambda_t^N(m,\zeta)}
 {\mathcal A(k,\eta)\mathcal A(\ell,\xi)\mathcal A(m,\zeta)}
=\frac1{\sqrt{\theta(k,\eta)\theta(\ell,\xi)\theta(m,\zeta)}}\le C.
\]
The constant is independent of $t,\nu,\mu,\beta$.
Now using $|e|=1$, recalling \eqref{eq:Af} and applying Young's inequality, we obtain from \eqref{eq:mixedphysicalexpansion} that
\begin{align*}
\int_0^{T_1}|I_4(t)|\dd t
&\le\frac{C_N(\nu+\mu)}{|\beta|\sqrt{\nu\mu}}
 \int_0^{T_1}\av t e^{-\rho^{1/3}t/16}
 \Bigl[\sqrt{\nu\mu}\|\nabla\mathcal Ab\|_2
                    \|\nabla\mathcal Au\|_2\|\mathcal Au\|_2
       +\nu\|\nabla\mathcal Au\|_2^2\|\mathcal Ab\|_2\Bigr]\dd t\\
&\quad+\frac{C_N}{|\beta|}\int_0^{T_1}\av t e^{-\rho^{1/3}t/16}
       \mu\|\nabla\mathcal Ab\|_2^2\|\mathcal Ab\|_2\dd t.
\end{align*}
Here we used $\Lambda_t^{r-N}\in L^2$ for $N-r>1$. 
By \eqref{eq:thetabounds}, \eqref{eq:fparseval} and \eqref{eq:controlenergy},
\begin{equation*}
\sup_{t\le T_1}(\|\mathcal Au\|_2+\|\mathcal Ab\|_2)
 +\sqrt\nu\|\nabla\mathcal Au\|_{L^2_tL^2}
 +\sqrt\mu\|\nabla\mathcal Ab\|_{L^2_tL^2}
 \le C\mathcal E(T_1)^{1/2}.
\end{equation*}
Then it follows from $(\nu+\mu)/\sqrt{\nu\mu}\ge2$ and $\|\av t e^{-\rho^{1/3}t/16}\|_{L^\infty_t}\le C\rho^{-1/3}$ that
\begin{align*}
 \int_0^{T_1}|I_4(t)|\dd t&\le\frac{C_N(\nu+\mu)}{|\beta|\sqrt{\nu\mu}}\rho^{-1/3}\mathcal E(T_1)^{3/2}+\frac{C_N}{|\beta|}\rho^{-1/3}\mathcal E(T_1)^{3/2}\le\frac{C_N(\nu+\mu)}{|\beta|\sqrt{\nu\mu}}\rho^{-1/3}\mathcal E(T_1)^{3/2},
\end{align*}
which is \eqref{eq:mixeddiffusionbound}.
\end{proof}

\subsection{Estimates for the nonlinear reactions involving the velocity correction}

\begin{proposition}
\label{lem:additionalcubic}
With the notation of Section~\ref{sec:equations}, let $(u,b)$ be a smooth solution of the MHD equations \eqref{eq:perturbation} on $[0,T_1]$. Then we have
\begin{align}
 &\big|2\langle\mathcal A(\widetilde u-u),
        \mathcal A\mathbb P(b\cdot\nabla b-u\cdot\nabla u)\rangle\big|
 +\big|2\langle\mathcal A\widetilde u,
         \mathcal A[\beta^{-1}\Pi_{\ne}(b^1u^2-u^1b^2)e_1]\rangle\big|
        \nonumber\\
 \le{}&\frac{C_N}{|\beta|}\Big(\sum_s\|F^s(t)\|_2^2\Big)^{1/2}
                              \sum_s\|F^s_{\ne}(t)\|_2^2.
\label{eq:Npointwise}
\end{align}
Consequently,
\begin{equation}
 \int_0^{T_1}\Big(
 2\big|\langle\mathcal A(\widetilde u-u),
              \mathcal A\mathbb P(b\cdot\nabla b-u\cdot\nabla u)\rangle\big|+\frac2{|\beta|}
       \big|\langle\mathcal A\widetilde u^1,
                         \mathcal A\Pi_{\ne}(b^1u^2-u^1b^2)\rangle\big|
                         \Big)\dd t
 \le \frac{C_N}{|\beta|}\rho^{-1/3}\mathcal E(T_1)^{3/2}.
 \label{eq:Nbound}
\end{equation}
\end{proposition}
\begin{proof}
We write
\begin{align*}
 J_5=2\big|\langle\mathcal A(\widetilde u-u),
                  \mathcal A\mathbb P(b\cdot\nabla b-u\cdot\nabla u)\rangle\big|,\quad
 J_6=\frac2{|\beta|}\big|\langle\mathcal A\widetilde u^1,
                  \mathcal A\Pi_{\ne}(b^1u^2-u^1b^2)\rangle\big|.
\end{align*}
The left-hand side of \eqref{eq:Npointwise} is $J_5+J_6$.
For $k=\ell+m$ and $\eta=\xi+\zeta$, recalling \eqref{eq:thetabounds} and \eqref{eq:Af}, we have
\begin{align*}
 \Lambda_t^N(k,\eta)&\le C_N[\Lambda_t^N(\ell,\xi)+\Lambda_t^N(m,\zeta)],\\
 \mathbf1_{k\ne0}&\le\mathbf1_{\ell\ne0}+\mathbf1_{m\ne0},\\
 \frac{\mathcal A(k,\eta)}{\mathcal A(\ell,\xi)\mathcal A(m,\zeta)}
 &\le C_N[\Lambda_t^{-N}(\ell,\xi)+\Lambda_t^{-N}(m,\zeta)].
\end{align*}
Hence, using $\Lambda_t^{-N}\mathcal A\widehat f,\Lambda_t^{-N}\mathcal A\widehat g\in L^1$ and Young's inequality for convolution, we obtain
\begin{equation}
 \|\mathcal A\Pi_{\ne}(fg)\|_2
 \le C_N\bigl(\|\mathcal Af\|_2\|\mathcal Ag_{\ne}\|_2
             +\|\mathcal Af_{\ne}\|_2\|\mathcal Ag\|_2\bigr).
 \label{eq:weightedproduct}
\end{equation}
By incompressibility and \eqref{eq:goodorder},
\[
 \|\nabla\mathcal A\mathbb P(\widetilde u-u)\|_2
 \le\frac C{|\beta|}\|\mathcal Ab_{\ne}\|_2.
\]
Integration by parts and $\widetilde u=u+\beta^{-1}\partial_x^{-1}b^2_{\ne}e_1$ therefore yield
\begin{align*}
 J_5=2\big|\langle\nabla\mathcal A\mathbb P(\widetilde u-u),
                \mathcal A\Pi_{\ne}(b\otimes b-u\otimes u)\rangle\big|
 \le\frac{C_N}{|\beta|}\|\mathcal Ab_{\ne}\|_2
   \big(\|\mathcal Au\|_2+\|\mathcal Ab\|_2\big)
   \big(\|\mathcal Au_{\ne}\|_2+\|\mathcal Ab_{\ne}\|_2\big).
\end{align*}
Then by virtue of \eqref{eq:weightedproduct}, we have
\begin{align*}
 J_6\le\frac{C_N}{|\beta|}\|\mathcal A\widetilde u_{\ne}\|_2
   \big(\|\mathcal Au\|_2+\|\mathcal Ab\|_2\big)
   \big(\|\mathcal Au_{\ne}\|_2+\|\mathcal Ab_{\ne}\|_2\big).
\end{align*}
Using \eqref{eq:thetabounds}, \eqref{eq:Af}, \eqref{eq:fparseval} and \eqref{eq:physicalenergycomparison}, we obtain
\[
 J_5+J_6\le  \frac{C_N}{|\beta|}
 \Bigl(\sum_s\|F^s(t)\|_2^2\Bigr)^{1/2}
                   \sum_s\|F^s_{\ne}(t)\|_2^2.
\]
This proves \eqref{eq:Npointwise}. Integrating the two terms gives
\begin{align*}
 \int_0^{T_1}(J_5(t)+J_6(t))\dd t\le\frac{C_N}{|\beta|}
 \sup_{t\le T_1}\Bigl(\sum_s\|F^s(t)\|_2^2\Bigr)^{1/2}
                \int_0^{T_1}\sum_s\|F^s_{\ne}(t)\|_2^2\dd t
\le\frac{C_N}{|\beta|}\rho^{-1/3}\mathcal E(T_1)^{3/2},
\end{align*}
which is \eqref{eq:Nbound}.
\end{proof}

\subsection{\texorpdfstring{Estimate for the interactions with the horizontal means}{Estimate for the interactions with the horizontal means}}

\begin{proposition}
\label{lem:meanadvection}
Let \(k=\ell+m\) and \(\eta=\xi+\zeta\). With the notation of Section~\ref{sec:equations}, let $(u,b)$ be a smooth solution of the MHD equations \eqref{eq:perturbation} on $[0,T_1]$. Then the nonlinear reaction \eqref{eq:T0} satisfies
\begin{align}
&\int_0^{T_1}\bigg|\frac12\sum_s\Rea\sum_{\substack{\ell,m\in\Z\\\ell=0}}\int_{\R^{2}} (\mathcal A^2(k,\eta)-\mathcal A^2(m,\zeta))a(k,\eta;\ell,\xi)\overline{f^s(k,\eta)}f^{-s}(\ell,\xi)f^s(m,\zeta)\frac{\dd\xi\dd\zeta}{(2\pi)^{2}}\bigg|\dd t \nonumber\\
\le{} &C\rho^{-1/3}\mathcal E(T_1)^{3/2}.
\label{eq:T0bound}
\end{align}
\end{proposition}
\begin{proof}
Since $\ell=0$ in \eqref{eq:T0}, we have \(k=\ell+m=m\) and \(\eta=\xi+\zeta\). Recalling \eqref{eq:weight} and \eqref{eq:Af}, we have
\[
(\mathcal A^2(k,\eta)-\mathcal A^2(m,\zeta))a(k,\eta;\ell,\xi)=e^{\rho^{1/3}t\mathbf 1_{k\ne0}/8}(h(k,\eta)-h(k,\zeta))a(k,\eta;0,\xi).
\]
For $k=0$, the kernel $a(0,\eta;0,\xi)$ vanishes. For $k\ne0$, $|a(k,\eta;0,\xi)|\le|k|$ by \eqref{eq:meanA}.
The definition \eqref{eq:weightedamplitudes} gives the exact cancellation
\[
e^{\rho^{1/3}t/8}
 |f^s(k,\eta)f^{-s}(0,\xi)f^s(k,\zeta)|
=\frac{|F^s(k,\eta)F^{-s}(0,\xi)F^s(k,\zeta)|}
 {\Lambda_t^N(k,\eta)\Lambda_t^N(0,\xi)\Lambda_t^N(k,\zeta)}.
\]
First suppose that $8|\xi|\le\sqrt{k^2+(\zeta+kt)^2}$.
Then $\Lambda_t(k,\zeta+r\xi)\approx\Lambda_t(k,\zeta)$ for $0\le r\le1$, and in particular $\Lambda_t(k,\eta)\approx\Lambda_t(k,\zeta)$. Recalling \eqref{eq:weight}, we write
\[
\begin{aligned}
h(k,\eta)-h(m,\zeta)
=(\theta(k,\eta)-\theta(m,\zeta))\Lambda_t^{2N}(m,\zeta)
+\theta(k,\eta)(\Lambda_t^{2N}(k,\eta)-\Lambda_t^{2N}(m,\zeta)).
\end{aligned}
\]
Since $\eta-\zeta=\xi$ and $k=m\ne0$, it follows from \eqref{eq:verticalLip} that
\begin{align*}
|k||\theta(k,\zeta+\xi)-\theta(k,\zeta)|&\le C|\xi|,\\
\Lambda_t^{2N}(k,\zeta+\xi)-\Lambda_t^{2N}(k,\zeta)
 &=2N\xi\int_0^1(\zeta+kt+r\xi)
       (1+k^2+(\zeta+kt+r\xi)^2)^{N-1}\dd r,\\
\frac{|k|\,|\Lambda_t^{2N}(k,\eta)-\Lambda_t^{2N}(k,\zeta)|}
 {\Lambda_t^N(k,\eta)\Lambda_t^N(k,\zeta)}
 &\le C_N\frac{|k||\xi|}{\Lambda_t(k,\zeta)}\le C_N|\xi|.
\end{align*}
The bounds \eqref{eq:thetabounds} imply
\begin{align*}
\frac{|(h(k,\eta)-h(m,\zeta))a(k,\eta;0,\xi)|}{\Lambda_t^N(k,\eta)\Lambda_t^N(0,\xi)\Lambda_t^N(m,\zeta)}\le  C_N|\xi|\Lambda_t^{-N}(0,\xi) \le C_N\Lambda_t^{1-N}(0,\xi).	
\end{align*}
Then suppose that $8|\xi|>\sqrt{k^2+(\zeta+kt)^2}$. One has
\[
 \Lambda_t(k,\zeta)\le8\Lambda_t(0,\xi),\qquad
 \Lambda_t(k,\eta)\le9\Lambda_t(0,\xi).
\]
Thus it holds that
\begin{align*}
 \frac{|h(k,\eta)-h(k,\zeta)|\,|a(k,\eta;0,\xi)|}
 {\Lambda_t^N(k,\eta)\Lambda_t^N(0,\xi)\Lambda_t^N(k,\zeta)}
 \le C_N|k|\,[\Lambda_t^{-N}(k,\zeta)+\Lambda_t^{-N}(k,\eta)]
 \le C_N[\Lambda_t^{1-N}(k,\zeta)+\Lambda_t^{1-N}(k,\eta)].
\end{align*}
Consequently, expressing $f^s$ in terms of $F^s$ by \eqref{eq:weightedamplitudes}, we deduce from Young's inequality, H\"older's inequality and \eqref{eq:controlenergy} that
\begin{align*}
&\int_0^{T_1}\bigg|\frac12\sum_s\Rea\sum_{\substack{\ell,m\in\Z\\\ell=0}}\int_{\R^{2}} (\mathcal A^2(k,\eta)-\mathcal A^2(m,\zeta))a(k,\eta;\ell,\xi)
\overline{f^s(k,\eta)}f^{-s}(\ell,\xi)f^s(m,\zeta)\frac{\dd\xi\dd\zeta}{(2\pi)^{2}}\bigg|\dd t \\
\le{} &C_N\sum_s\int_0^{T_1}\|F^s_{\ne}\|_2^2
             \|\Lambda_t^{1-N}(0,\xi)F^{-s}_0(0,\xi)\|_1\dd t\\
 &+C_N\sum_s\int_0^{T_1}\big[\|F^s_{\ne}\|_2\|F^{-s}_0\|_2\|\Lambda_t^{1-N}F^s_{\ne}\|_1+\|\Lambda_t^{1-N}F^s_{\ne}\|_1\|F^{-s}_0\|_2\|F^s_{\ne}\|_2\big]\dd t\\
\le{} &C_N\sum_s\|F^{-s}_0\|_{L^\infty_tL^2}
                  \|F^s_{\ne}\|_{L^2_tL^2}^2+C_N\sum_s\|F^{-s}_0\|_{L^\infty_tL^2}
                         \|F^s_{\ne}\|_{L^2_tL^2}^2\\
\le{} &C_N\rho^{-1/3}\mathcal E(T_1)^{3/2},
\end{align*}
which proves \eqref{eq:T0bound}.
\end{proof}

\subsection{\texorpdfstring{Estimate for the remainder from the exponential weights}{Estimate for the remainder from the exponential weights}}

\begin{proposition}
\label{lem:exponentialjump}
Let \(k=\ell+m\) and \(\eta=\xi+\zeta\). With the notation of Section~\ref{sec:equations}, let $(u,b)$ be a smooth solution of the MHD equations \eqref{eq:perturbation} on $[0,T_1]$. Then the nonlinear reaction \eqref{eq:TJ} satisfies
\begin{equation}
\begin{aligned}
 &\int_0^{T_1}\bigg|\frac12\sum_s\Rea\sum_{\ell\ne0,m}\int_{\R^2}
 [\mathcal A^2(k,\eta)-\mathcal A^2(m,\zeta)
       -e^{\rho^{1/3}t\mathbf 1_{k\ne0}\mathbf 1_{m\ne0}/8}(h(k,\eta)-h(m,\zeta))]\\
 &\hspace{12mm}\times a(k,\eta;\ell,\xi)
       \overline{f^s(k,\eta)}f^{-s}(\ell,\xi)f^s(m,\zeta)
                   \frac{\dd\xi\dd\zeta}{(2\pi)^2}\bigg|\dd t
 \le C_N\rho^{-1/3}\mathcal E(T_1)^{3/2}.
 \label{eq:TJbound}
\end{aligned}
\end{equation}
\end{proposition}
\begin{proof}
A direct computation gives
\[
\begin{aligned}
\mathcal A^2(k,\eta)-\mathcal A^2(m,\zeta)
 -e^{\frac18\rho^{1/3}t\mathbf1_{k\ne0}\mathbf1_{m\ne0}}
       (h(k,\eta)-h(m,\zeta))=\begin{cases}
 0,&km\ne0,\\
 (e^{\rho^{1/3}t/8}-1)h(k,\eta),&m=0,\\
 -(e^{\rho^{1/3}t/8}-1)h(m,\zeta),&k=0.
 \end{cases}
\end{aligned}
\]
Hence it suffices to consider the cases $m=0$ and $k=0$.
For $m=0$, substitute $k=\ell$ and $\eta=\xi+\zeta$ in \eqref{eq:TJ}.
For $k=0$, substitute $m=-\ell$ and change variables from $(\xi,\zeta)$ to $(\xi,\eta)$ in \eqref{eq:TJ}, where $\eta=\xi+\zeta$.
Thus expressing $f^s$ in terms of $F^s$ by \eqref{eq:weightedamplitudes}, we obtain 
\begin{align*}
&\int_0^{T_1}\bigg|\frac12\sum_s\Rea\sum_{\ell\ne0,m}\int_{\R^2}
 [\mathcal A^2(k,\eta)-\mathcal A^2(m,\zeta)
       -e^{\rho^{1/3}t\mathbf 1_{k\ne0}\mathbf 1_{m\ne0}/8}(h(k,\eta)-h(m,\zeta))]\\
 &\hspace{12mm}\times a(k,\eta;\ell,\xi)
       \overline{f^s(k,\eta)}f^{-s}(\ell,\xi)f^s(m,\zeta)
                   \frac{\dd\xi\dd\zeta}{(2\pi)^2}\bigg|\dd t\\
\le{} &\frac12\sum_s\int_0^{T_1}\sum_{\ell\ne0}\int_{\R^2}
  \frac{\big|(1-e^{-\rho^{1/3}t/8})a(\ell,\xi+\zeta;\ell,\xi)h(\ell,\xi+\zeta)\big|}
       {\Lambda_t^{N}(\ell,\xi+\zeta)\Lambda_t^{N}(\ell,\xi)\Lambda_t^{N}(0,\zeta)}
 |F^s(\ell,\xi+\zeta)||F^{-s}(\ell,\xi)||F^s(0,\zeta)|
 \frac{\dd\xi\dd\zeta}{(2\pi)^2}\dd t\\
&+\frac12\sum_s\int_0^{T_1}\sum_{\ell\ne0}\int_{\R^2}
 \frac{\big|(1-e^{-\rho^{1/3}t/8})a(0,\eta;\ell,\xi)h(-\ell,\eta-\xi)\big|}
      {\Lambda_t^{N}(0,\eta)\Lambda_t^{N}(\ell,\xi)\Lambda_t^{N}(-\ell,\eta-\xi)}
 |F^s(0,\eta)||F^{-s}(\ell,\xi)||F^s(-\ell,\eta-\xi)|
 \frac{\dd\xi\dd\eta}{(2\pi)^2}\dd t.
\end{align*}

If $m=0$, then $k=\ell\ne0$ and $|a(k,\eta;\ell,\xi)|\le|\zeta|\kappa(\ell,\xi)$ by \eqref{eq:kernelgeometric} and \eqref{eq:frequencies}.
Using $h=\theta\Lambda_t^{2N}$ and $0<\theta\le1$, we obtain
\begin{equation}
\begin{aligned}
\frac{\big|(1-e^{-\rho^{1/3}t/8})a(k,\eta;\ell,\xi)h(k,\eta)\big|}
 {\Lambda_t^{N}(k,\eta)\Lambda_t^{N}(\ell,\xi)\Lambda_t^{N}(0,\zeta)}
\le C(1-e^{-\frac18\rho^{1/3}t})
       \frac{\Lambda_t^{N}(k,\eta)}{\Lambda_t^{N}(\ell,\xi)\Lambda_t^{N}(0,\zeta)}|\zeta|\kappa(\ell,\xi).
\end{aligned}
\label{eq:jumpsymbol}
\end{equation}
First suppose that $8\sqrt{\ell^2+(\xi+\ell t)^2}\le|\zeta|$. Then we have $\Lambda_t(k,\eta)\approx\Lambda_t(0,\zeta)$ because $k=\ell$ and $\eta=\xi+\zeta$. Moreover, \eqref{eq:tkappa} implies
\[
 (1-e^{-\frac18\rho^{1/3} t})\kappa(\ell,\xi)
       \le C\rho^{1/3} t\kappa(\ell,\xi)\le C\rho^{1/3}\Lambda_t(\ell,\xi).
\]
Thus the left-hand side of \eqref{eq:jumpsymbol} can be further bounded by
\begin{equation*}
\frac{\big|(1-e^{-\rho^{1/3}t/8})a(k,\eta;\ell,\xi)h(k,\eta)\big|}
 {\Lambda_t^{N}(k,\eta)\Lambda_t^{N}(\ell,\xi)\Lambda_t^{N}(0,\zeta)}
\le C_N\rho^{1/3}|\zeta|\Lambda_t^{1-N}(\ell,\xi).
\end{equation*}
Then suppose that  $8\sqrt{\ell^2+(\xi+\ell t)^2}>|\zeta|$. We have $\Lambda_t(k,\eta)\le9\Lambda_t(\ell,\xi)$.
Since $0\le1-e^{-\rho^{1/3}t/8}\le1$ and $\kappa(\ell,\xi)\le1$, the left-hand side of \eqref{eq:jumpsymbol} can be bounded by
\[
\frac{\big|(1-e^{-\rho^{1/3}t/8})a(k,\eta;\ell,\xi)h(k,\eta)\big|}
 {\Lambda_t^{N}(k,\eta)\Lambda_t^{N}(\ell,\xi)\Lambda_t^{N}(0,\zeta)}\le C_N|\zeta|\Lambda_t^{-N}(0,\zeta)
\le C_N\Lambda_t^{1-N}(0,\zeta).
\]
By \eqref{eq:controlenergy} and \eqref{eq:minimumcontrol}, for each $s$,
\[
\|F^s\|_{L^\infty_tL^2}
 +\rho^{1/6}\|F^s_{\ne}\|_{L^2_tL^2}
 +\rho^{1/2}\||\eta|F^s_0(0,\eta)\|_{L^2_tL^2}
\le C\mathcal E(T_1)^{1/2}.
\]
Thus we deduce from Young's inequality, and H\"older's inequality that
\begin{align*}
&\frac12\sum_s\int_0^{T_1}\sum_{\ell\ne0}\int_{\R^2}
  \frac{\big|(1-e^{-\rho^{1/3}t/8})a(\ell,\xi+\zeta;\ell,\xi)h(\ell,\xi+\zeta)\big|}
       {\Lambda_t^{N}(\ell,\xi+\zeta)\Lambda_t^{N}(\ell,\xi)\Lambda_t^{N}(0,\zeta)}
 |F^s(\ell,\xi+\zeta)||F^{-s}(\ell,\xi)||F^s(0,\zeta)|
 \frac{\dd\xi\dd\zeta}{(2\pi)^2}\dd t\\
\le{} &C_N\rho^{1/3}\sum_s\int_0^{T_1}
 \|F^s_{\ne}\|_2\||\zeta|F_0^s(0,\zeta)\|_2
                    \|\Lambda_t^{1-N}(\ell,\xi)F_{\ne}^{-s}\|_1\dd t
+C_N\sum_s\int_0^{T_1}\|F^s_{\ne}\|_2\|F^{-s}_{\ne}\|_2
                      \|\Lambda_t^{1-N}F^s_0\|_1\dd t\\
\le{} &C_N\rho^{1/3}\sum_s\|F^s_{\ne}\|_{L^\infty_tL^2}
            \||\zeta|F_0^s(0,\zeta)\|_{L^2_tL^2}\|F_{\ne}^{-s}\|_{L^2_tL^2}
+C_N\sum_s\|F^s_0\|_{L^\infty_tL^2}
       \|F^s_{\ne}\|_{L^2_tL^2}\|F^{-s}_{\ne}\|_{L^2_tL^2}\\
\le{} &C_N\rho^{1/3}\rho^{-1/2}\rho^{-1/6}\mathcal E(T_1)^{3/2}
       +C_N\rho^{-1/3}\mathcal E(T_1)^{3/2}
 \le C_N\rho^{-1/3}\mathcal E(T_1)^{3/2}.
\end{align*}

It remains to treat $k=0$. Now $m=-\ell\ne0$, $\zeta=\eta-\xi$, and \eqref{eq:kernelgeometric} gives
\[
|a(0,\eta;\ell,\xi)|
\le|e(\ell,\xi)\cdot(-\ell,\eta-\xi)|
=\frac{|\ell\eta|}{\sqrt{\ell^2+\xi^2}}
=|\eta|\kappa(\ell,\xi).
\]
Hence
\[
\begin{aligned}
\frac{\big|(1-e^{-\rho^{1/3}t/8})a(0,\eta;\ell,\xi)h(m,\zeta)\big|}
 {\Lambda_t^{N}(0,\eta)\Lambda_t^{N}(\ell,\xi)\Lambda_t^{N}(m,\zeta)}
\le C_N(1-e^{-\rho^{1/3}t/8})
 \frac{\Lambda_t^N(m,\zeta)}{\Lambda_t^N(\ell,\xi)\Lambda_t^N(0,\eta)}
 |\eta|\kappa(\ell,\xi).
\end{aligned}
\]
If $8\sqrt{\ell^2+(\xi+\ell t)^2}\le|\eta|$, then $\Lambda_t(m,\zeta)\le C\Lambda_t(0,\eta)$, so this is bounded by $C_N\rho^{1/3}|\eta|\Lambda_t^{1-N}(\ell,\xi)$ by \eqref{eq:tkappa}.
If $8\sqrt{\ell^2+(\xi+\ell t)^2}>|\eta|$, then $\Lambda_t(m,\zeta)\le9\Lambda_t(\ell,\xi)$, and this is bounded by $C_N\Lambda_t^{1-N}(0,\eta)$.
Thus
\begin{align*}
&\frac12\sum_s\int_0^{T_1}\sum_{\ell\ne0}\int_{\R^2}
 \frac{\big|(1-e^{-\rho^{1/3}t/8})a(0,\eta;\ell,\xi)h(-\ell,\eta-\xi)\big|}
      {\Lambda_t^{N}(0,\eta)\Lambda_t^{N}(\ell,\xi)\Lambda_t^{N}(-\ell,\eta-\xi)}
 |F^s(0,\eta)||F^{-s}(\ell,\xi)||F^s(-\ell,\eta-\xi)|
 \frac{\dd\xi\dd\eta}{(2\pi)^2}\dd t\\
\le{} &C_N\rho^{1/3}\sum_s\int_0^{T_1}
 \||\eta|F^s_0(0,\eta)\|_2
 \|\Lambda_t^{1-N}F^{-s}_{\ne}\|_1\|F^s_{\ne}\|_2\dd t
+C_N\sum_s\int_0^{T_1}\|\Lambda_t^{1-N}F^s_0\|_1
              \bigl\||F^{-s}_{\ne}|*|F^s_{\ne}|\bigr\|_\infty\dd t\\
\le{} &C_N\rho^{1/3}\sum_s
 \||\eta|F^s_0(0,\eta)\|_{L^2_tL^2}
 \|F^{-s}_{\ne}\|_{L^2_tL^2}\|F^s_{\ne}\|_{L^\infty_tL^2}
+C_N\sum_s\|F^s_0\|_{L^\infty_tL^2}
       \|F^{-s}_{\ne}\|_{L^2_tL^2}\|F^s_{\ne}\|_{L^2_tL^2}\\
\le{} &C_N\rho^{1/3}\rho^{-1/2}\rho^{-1/6}\mathcal E(T_1)^{3/2}
       +C_N\rho^{-1/3}\mathcal E(T_1)^{3/2}\\
\le{} &C_N\rho^{-1/3}\mathcal E(T_1)^{3/2}.
\end{align*}
Adding the two estimates proves \eqref{eq:TJbound}.
\end{proof}

\section{Estimates for the remaining nonlinear reactions}
\label{sec:quartic}

In this section, we establish the following estimates for $I_5+I_7$ and $I_6$ in \eqref{eq:pair} and \eqref{eq:adv}, respectively.

\begin{proposition}
\label{prop:quartic}
Let $N>10$ and let $(u,b)$ be a smooth solution of the MHD equations \eqref{eq:perturbation} on $[0,T_1]$. With the notation of Section~\ref{sec:equations}, the terms \eqref{eq:pair} and \eqref{eq:adv} satisfy
\begin{align}
&\int_0^{T_1}\bigg|\sum_s\Rea\sum_{\substack{j,\ell,p\in\Z\\\ell\ne0}}\int_{\R^{3}}\frac{\ii}{4s\beta\ell}\nonumber\\
&\qquad\times        \Bigl[e^{\frac18\rho^{1/3} t\mathbf 1_{k\ne0}\mathbf 1_{j+p\ne0}}(h(k,\eta)-h(j+p,\alpha+\gamma))
                   a(k,\eta;\ell,\xi)a(j+p,\alpha+\gamma;p,\gamma)\nonumber\\
 &\hspace{28mm}-e^{\frac18\rho^{1/3} t\mathbf 1_{j+\ell\ne0}\mathbf 1_{j\ne0}}(h(j+\ell,\alpha+\xi)-h(j,\alpha))
                   a(k,\eta;p,\gamma)a(j+\ell,\alpha+\xi;\ell,\xi)\Bigr]\nonumber\\
 &\hspace{12mm}\times\overline{f^s(k,\eta)}f^s(j,\alpha)f^{-s}(\ell,\xi)f^{-s}(p,\gamma)
\frac{\dd\alpha\dd\xi\dd\gamma}{(2\pi)^{3}}\bigg|\dd t\nonumber\\
&+\int_0^{T_1}\bigg|\sum_s\Rea\sum_{\substack{p,q,m\in\Z\\p+q\ne0}}\int_{\R^{3}}\frac{\ii e^{\frac18\rho^{1/3} t\mathbf1_{k\ne0}\mathbf 1_{m\ne0}}(h(k,\eta)-h(m,\zeta))}{4s\beta(p+q)}
 a(k,\eta;p+q,\gamma+\delta)\nonumber\\
 &\hspace{8mm}\times a(p+q,\gamma+\delta;p,\gamma)
       \overline{f^s(k,\eta)}f^s(p,\gamma)f^{-s}(q,\delta)f^s(m,\zeta)\frac{\dd\gamma\dd\delta\dd\zeta}{(2\pi)^{3}}\bigg|\dd t\nonumber\\
\le{}&\frac{C_N}{|\beta|}\rho^{-2/3}\mathcal E(T_1)^2.
\label{eq:quarticbound}
\end{align}
In the first integral, $k=j+\ell+p$ and
$\eta=\alpha+\xi+\gamma$; in the second, $k=p+q+m$ and
$\eta=\gamma+\delta+\zeta$.
\end{proposition}

\subsection{Estimates for $I_5+I_7$}
\label{sec:quarticfirst}

For \eqref{eq:pair}, $k=j+\ell+p$, $\eta=\alpha+\xi+\gamma$,
and $\ell\ne0$. Hence
\[
 (k,\eta+kt)=(j,\alpha+jt)+(\ell,\xi+\ell t)+(p,\gamma+pt).
\]
At each fixed time, we divide the admissible frequencies into
\emph{seven disjoint regions}. Regions 1--6 are defined, in order,
by the following six inequalities:
\begin{equation}
\begin{gathered}
8\bigl(\sqrt{\ell^2+(\xi+\ell t)^2}+\sqrt{p^2+(\gamma+pt)^2}\bigr)\le\min\big\{\sqrt{k^2+(\eta+kt)^2},\sqrt{j^2+(\alpha+jt)^2}\big\},\\
8\bigl(\sqrt{k^2+(\eta+kt)^2}+\sqrt{j^2+(\alpha+jt)^2}\bigr)\le\min\big\{\sqrt{\ell^2+(\xi+\ell t)^2},\sqrt{p^2+(\gamma+pt)^2}\big\},\\
8\bigl(\sqrt{j^2+(\alpha+jt)^2}+\sqrt{p^2+(\gamma+pt)^2}\bigr)\le\min\bigl\{\sqrt{k^2+(\eta+kt)^2},\sqrt{\ell^2+(\xi+\ell t)^2}\bigr\},\\
8\bigl(\sqrt{k^2+(\eta+kt)^2}+\sqrt{p^2+(\gamma+pt)^2}\bigr)\le\min\bigl\{\sqrt{j^2+(\alpha+jt)^2},\sqrt{\ell^2+(\xi+\ell t)^2}\bigr\},\\
8\bigl(\sqrt{j^2+(\alpha+jt)^2}+\sqrt{\ell^2+(\xi+\ell t)^2}\bigr)\le\min\bigl\{\sqrt{k^2+(\eta+kt)^2},\sqrt{p^2+(\gamma+pt)^2}\bigr\},\\
8\bigl(\sqrt{k^2+(\eta+kt)^2}+\sqrt{\ell^2+(\xi+\ell t)^2}\bigr)\le\min\bigl\{\sqrt{j^2+(\alpha+jt)^2},\sqrt{p^2+(\gamma+pt)^2}\bigr\}.
\end{gathered}
\label{eq:pairpartition}
\end{equation}
Region 7 is the complement of the union of regions 1--6 within the admissible frequencies. Equivalently, all six inequalities in \eqref{eq:pairpartition} fail there.
The relation above and $\ell\ne0$ imply that at least two of the four lengths are positive. In any of regions 1--6, the two selected larger lengths are therefore positive.
If two different selections held simultaneously, a length selected only in the first would be at least eight times a length selected only in the second, and conversely, a contradiction.
This proves disjointness. Adding region 7 gives a partition.
For the exponential weights, the frequency relation gives
\begin{equation}
\begin{aligned}
 \frac{e^{\rho^{1/3}t\mathbf1_{k\ne0}\mathbf1_{j+p\ne0}/8}}
 {e^{\rho^{1/3}t(\mathbf1_{k\ne0}+\mathbf1_{j\ne0}
                    +\mathbf1_{\ell\ne0}+\mathbf1_{p\ne0})/16}}
 &\le e^{-\rho^{1/3}t/16},\qquad
 \frac{e^{\rho^{1/3}t\mathbf1_{j+\ell\ne0}\mathbf1_{j\ne0}/8}}
 {e^{\rho^{1/3}t(\mathbf1_{k\ne0}+\mathbf1_{j\ne0}
                    +\mathbf1_{\ell\ne0}+\mathbf1_{p\ne0})/16}}
 \le e^{-\rho^{1/3}t/16}.
\end{aligned}
\label{eq:quarticexpdecay}
\end{equation}
Indeed, a numerator equal to $e^{\rho^{1/3}t/8}$ forces at least three of $k,j,\ell,p$ to be nonzero. Otherwise, at least two are nonzero, since $\ell\ne0$ and $k=j+\ell+p$.

The following lemmas give pointwise kernel bounds on these seven regions. In region 1, we have
\begin{equation}
8\bigl(\sqrt{\ell^2+(\xi+\ell t)^2}+\sqrt{p^2+(\gamma+pt)^2}\bigr)\le\min\big\{\sqrt{k^2+(\eta+kt)^2},\sqrt{j^2+(\alpha+jt)^2}\big\}.
\label{eq:twohighset}
\end{equation}
The coefficient in \eqref{eq:pair} has the exact decomposition
\begin{align}
&e^{\frac18\rho^{1/3} t\mathbf 1_{k\ne0}\mathbf 1_{j+p\ne0}}(h(k,\eta)-h(j+p,\alpha+\gamma))a(k,\eta;\ell,\xi)a(j+p,\alpha+\gamma;p,\gamma)\nonumber\\
&-e^{\frac18\rho^{1/3} t\mathbf 1_{j+\ell\ne0}\mathbf 1_{j\ne0}}(h(j+\ell,\alpha+\xi)-h(j,\alpha))a(k,\eta;p,\gamma)a(j+\ell,\alpha+\xi;\ell,\xi)\nonumber\\
={}&e^{\frac18\rho^{1/3} t\mathbf 1_{k\ne0}\mathbf 1_{j+p\ne0}}(h(k,\eta)-h(j+p,\alpha+\gamma))
       \Bigl[a(k,\eta;\ell,\xi)a(j+p,\alpha+\gamma;p,\gamma)\nonumber\\
&\qquad -a(k,\eta;p,\gamma)a(j+\ell,\alpha+\xi;\ell,\xi)\Bigr]\nonumber\\
&+e^{\frac18\rho^{1/3} t\mathbf 1_{k\ne0}\mathbf 1_{j+p\ne0}}(h(k,\eta)-h(j+p,\alpha+\gamma)-h(j+\ell,\alpha+\xi)+h(j,\alpha))\nonumber\\*
&\qquad \times a(k,\eta;p,\gamma)a(j+\ell,\alpha+\xi;\ell,\xi)\nonumber\\
&+(e^{\frac18\rho^{1/3} t\mathbf 1_{k\ne0}\mathbf 1_{j+p\ne0}}-e^{\frac18\rho^{1/3} t\mathbf 1_{j+\ell\ne0}\mathbf 1_{j\ne0}})(h(j+\ell,\alpha+\xi)-h(j,\alpha))\nonumber\\
&\qquad \times a(k,\eta;p,\gamma)a(j+\ell,\alpha+\xi;\ell,\xi).
\label{eq:Ssplit}
\end{align}
The first two terms are treated in Lemmas~\ref{lem:lowpair} and \ref{lem:meanpair}, according to whether $p\ne0$ or $p=0$.
Lemma~\ref{lem:quarticinterface} treats the last term.

\begin{remark}\label{rem:1}
Notice that if $p=0$, then $k=j+\ell$, and $\eta=\alpha+\xi+\gamma$.
Since $j+p=j$ and $j+\ell=k$,
\[
 e^{\rho^{1/3}t\mathbf1_{k\ne0}\mathbf1_{j+p\ne0}/8}
 =e^{\rho^{1/3}t\mathbf1_{j+\ell\ne0}\mathbf1_{j\ne0}/8}.
\]
For $p=0$, the coefficient in \eqref{eq:pair} is the sum of the first two terms on the right-hand side of \eqref{eq:Ssplit}. Thus it suffices to consider $p\neq0$ in Lemma~\ref{lem:quarticinterface}.
\end{remark}

\begin{lemma}\label{lem:lowpair}
Let $k=j+\ell+p$, $\eta=\alpha+\xi+\gamma$, and $\ell,p\ne0$.
On the region \eqref{eq:twohighset}, the weights and kernels defined in \eqref{eq:weight} and \eqref{eq:kernelgeometric} satisfy
\begin{align}
&\frac{e^{\rho^{1/3}t\mathbf1_{k\ne0}\mathbf1_{j+p\ne0}/8}}{e^{\rho^{1/3}t(\mathbf1_{k\ne0}+\mathbf1_{j\ne0}+\mathbf1_{\ell\ne0}+\mathbf1_{p\ne0})/16}}
\frac{1}{|\ell|\Lambda_t^N(k,\eta)\Lambda_t^N(j,\alpha)\Lambda_t^N(\ell,\xi)\Lambda_t^N(p,\gamma)}\nonumber\\
&\times\Bigl[
 |h(k,\eta)-h(j+p,\alpha+\gamma)|
 \bigl|a(k,\eta;\ell,\xi)a(j+p,\alpha+\gamma;p,\gamma)
       -a(k,\eta;p,\gamma)a(j+\ell,\alpha+\xi;\ell,\xi)\bigr|\nonumber\\
&\qquad+
 |h(k,\eta)-h(j+p,\alpha+\gamma)-h(j+\ell,\alpha+\xi)+h(j,\alpha)|
 |a(k,\eta;p,\gamma)a(j+\ell,\alpha+\xi;\ell,\xi)|
 \Bigr]\nonumber\\
\le{} &C_Ne^{-\rho^{1/3}t/16}
 \Lambda_t^{9-N}(\ell,\xi)\Lambda_t^{9-N}(p,\gamma)\Bigl[
 \big[\av t^2+\rho^{-1/3}\av t\av{\rho^{1/3}t}^4\big]
       \sqrt{g(\ell,\xi)g(p,\gamma)}
 +\av t^2\sqrt{g(k,\eta)g(j,\alpha)}\Bigr].
\label{eq:lowpairtarget}
\end{align}
\end{lemma}
\begin{proof}
\emph{Step 1. Estimates for the two products.} 
By \eqref{eq:A}, we have
\begin{align}
&a(k,\eta;\ell,\xi)a(j+p,\alpha+\gamma;p,\gamma)\nonumber\\
={} &\frac{(\ell(\alpha+\gamma)-(j+p)\xi)(p\alpha-j\gamma)}
 {\sqrt{\ell^2+\xi^2}\sqrt{p^2+\gamma^2}}\,
 [e(k,\eta)\cdot e(j+p,\alpha+\gamma)]
 [e(j+p,\alpha+\gamma)\cdot e(j,\alpha)],\nonumber\\
&a(k,\eta;p,\gamma)a(j+\ell,\alpha+\xi;\ell,\xi)\nonumber\\
={} &\frac{(p(\alpha+\xi)-(j+\ell)\gamma)(\ell\alpha-j\xi)}
 {\sqrt{\ell^2+\xi^2}\sqrt{p^2+\gamma^2}}\,
 [e(k,\eta)\cdot e(j+\ell,\alpha+\xi)]
 [e(j+\ell,\alpha+\xi)\cdot e(j,\alpha)].
\label{eq:quarticproductidentities}
\end{align}
Since $\ell,p\ne0$, we have
\begin{align}
|\ell\alpha-j\xi|
 &=|\ell(\alpha+jt)-j(\xi+\ell t)|
 \le\Lambda_t(\ell,\xi)\Lambda_t(j,\alpha),\nonumber\\
|p(\alpha+\xi)-(j+\ell)\gamma|
 &=|p(\alpha+jt+\xi+\ell t)-(j+\ell)(\gamma+pt)|
 \le C\Lambda_t(p,\gamma)\Lambda_t(j,\alpha),\nonumber\\
\frac1{\sqrt{\ell^2+\xi^2}\sqrt{p^2+\gamma^2}}
 &=\frac{\kappa(\ell,\xi)\kappa(p,\gamma)}{|\ell p|}
 \le\kappa(\ell,\xi)\kappa(p,\gamma).
\label{eq:quarticdetbounds}
\end{align}
Here the second inequality uses \eqref{eq:twohighset}. The second identity in \eqref{eq:quarticproductidentities},
\eqref{eq:quarticdetbounds}, and $|e|=1$ yield
\begin{align}
|a(k,\eta;p,\gamma)a(j+\ell,\alpha+\xi;\ell,\xi)|
\le C\Lambda_t(\ell,\xi)\Lambda_t(p,\gamma)
       \Lambda_t(j,\alpha)^2\kappa(\ell,\xi)\kappa(p,\gamma).
\label{eq:quarticproductbound}
\end{align}
For the difference of the two products, first expand
\begin{align}
&(\ell(\alpha+\gamma)-(j+p)\xi)(p\alpha-j\gamma)
 -(p(\alpha+\xi)-(j+\ell)\gamma)(\ell\alpha-j\xi)\nonumber\\
={}&[(\ell\alpha-j\xi)+(\ell\gamma-p\xi)](p\alpha-j\gamma)
 -[(p\alpha-j\gamma)-(\ell\gamma-p\xi)](\ell\alpha-j\xi)\nonumber\\
={}&(\ell\gamma-p\xi)((\ell+p)\alpha-j(\xi+\gamma)).
\label{eq:detcancel}
\end{align}
Subtracting the two identities in \eqref{eq:quarticproductidentities} and using \eqref{eq:detcancel}, we obtain
\begin{align}
&a(k,\eta;\ell,\xi)a(j+p,\alpha+\gamma;p,\gamma)
 -a(k,\eta;p,\gamma)a(j+\ell,\alpha+\xi;\ell,\xi)\nonumber\\
={}&\frac{(\ell\gamma-p\xi)((\ell+p)\alpha-j(\xi+\gamma))}
 {\sqrt{\ell^2+\xi^2}\sqrt{p^2+\gamma^2}}
 [e(k,\eta)\cdot e(j+p,\alpha+\gamma)]
 [e(j+p,\alpha+\gamma)\cdot e(j,\alpha)]\nonumber\\
&+\frac{(p(\alpha+\xi)-(j+\ell)\gamma)(\ell\alpha-j\xi)}
 {\sqrt{\ell^2+\xi^2}\sqrt{p^2+\gamma^2}}
 \Bigl[
 [e(k,\eta)\cdot e(j+p,\alpha+\gamma)]
 [e(j+p,\alpha+\gamma)\cdot e(j,\alpha)]\nonumber\\
&\hspace{35mm}
 -[e(k,\eta)\cdot e(j+\ell,\alpha+\xi)]
 [e(j+\ell,\alpha+\xi)\cdot e(j,\alpha)]
 \Bigr].
\label{eq:quarticproductdifference}
\end{align}
The coefficient of the first term satisfies
\begin{align}
\frac{|(\ell\gamma-p\xi)((\ell+p)\alpha-j(\xi+\gamma))|}
 {\sqrt{\ell^2+\xi^2}\sqrt{p^2+\gamma^2}}
&=\frac{|\ell(\gamma+pt)-p(\xi+\ell t)|}
 {\sqrt{\ell^2+\xi^2}\sqrt{p^2+\gamma^2}}\,
 |(\ell+p)(\alpha+jt)-j(\xi+\ell t+\gamma+pt)|\nonumber\\
&\le C\Lambda_t(\ell,\xi)\Lambda_t(p,\gamma)
 [\Lambda_t(\ell,\xi)+\Lambda_t(p,\gamma)]
 \Lambda_t(j,\alpha)\kappa(\ell,\xi)\kappa(p,\gamma).
\label{eq:quarticdetdifferencebound}
\end{align}
For the difference of the dot products,
$e(a,\omega)=(-\omega,a)/\sqrt{a^2+\omega^2}$ gives
\[
 |e(a,\omega)-e(b,\sigma)|
 \le\frac{2\sqrt{(a-b)^2+(\omega-\sigma)^2}}{\sqrt{a^2+\omega^2}}.
\]
Hence
\begin{equation}
\begin{split}
&\big|[e(k,\eta)\cdot e(a,\omega)][e(a,\omega)\cdot e(j,\alpha)]
 -[e(k,\eta)\cdot e(b,\sigma)][e(b,\sigma)\cdot e(j,\alpha)]\big|\\
\le{}& 2|e(a,\omega)-e(b,\sigma)| \le C\frac{\sqrt{(a-b)^2+(\omega-\sigma)^2}}{\sqrt{a^2+\omega^2}}.
\end{split}
\label{eq:projectordiff}
\end{equation}
Due to $k=j+\ell+p$, $\eta=\alpha+\xi+\gamma$, we have
\begin{align*}
&|\Lambda_t(j+p,\alpha+\gamma)-\Lambda_t(j,\alpha)|\le \sqrt{p^2+(pt+\gamma)^2}, \quad |\Lambda_t(j+\ell,\alpha+\xi)-\Lambda_t(j,\alpha)|\le \sqrt{\ell^2+(\ell t+\xi)^2},\\
&|\Lambda_t(k,\eta)-\Lambda_t(j,\alpha)|\le \sqrt{(\ell+p)^2+(\ell t+p t+\xi+\gamma)^2} \le \sqrt{p^2+(pt+\gamma)^2}+\sqrt{\ell^2+(\ell t+\xi)^2}.
\end{align*}
Applying \eqref{eq:twohighset} yields
\begin{align}
&\frac78\Lambda_t(j,\alpha)
 \le\Lambda_t(j+p,\alpha+\gamma),\
    \Lambda_t(j+\ell,\alpha+\xi),\ \Lambda_t(k,\eta)
 \le\frac98\Lambda_t(j,\alpha).
\label{eq:quartichighcomparison}
\end{align}
Notice that $\ell,p\ne0$ and \eqref{eq:twohighset} imply $\sqrt{j^2+(\alpha+jt)^2}\ge16$.  Thus
\begin{align}
\Lambda_t(j,\alpha)
 &\le C\sqrt{(j+p)^2+(\alpha+\gamma+(j+p)t)^2}\nonumber\\
 &\le C\bigl[\sqrt{(j+p)^2+(\alpha+\gamma)^2}+t|j+p|\bigr],\nonumber\\
\frac{\Lambda_t(j,\alpha)}{\sqrt{(j+p)^2+(\alpha+\gamma)^2}}
 &\le C[1+\av t\kappa(j+p,\alpha+\gamma)].
\label{eq:quartichighphysical}
\end{align}
This also holds for $j+p=0$ due to \eqref{eq:quartichighcomparison}.
Moreover,
\[
 \sqrt{(p-\ell)^2+(\gamma-\xi)^2}
 \le C\av t[\Lambda_t(\ell,\xi)+\Lambda_t(p,\gamma)].
\]
Applying \eqref{eq:projectordiff}, \eqref{eq:quartichighphysical},
and $\kappa^2\le g$ from \eqref{eq:gprimitive}, we obtain
\begin{align}
&\big|[e(k,\eta)\cdot e(j+p,\alpha+\gamma)]
 [e(j+p,\alpha+\gamma)\cdot e(j,\alpha)]\nonumber\\
&\qquad-[e(k,\eta)\cdot e(j+\ell,\alpha+\xi)]
 [e(j+\ell,\alpha+\xi)\cdot e(j,\alpha)]\big|\nonumber\\
\le{} & C\frac{\Lambda_t(\ell,\xi)+\Lambda_t(p,\gamma)}
              {\Lambda_t(j,\alpha)}
       \big[\av t+\av t^2\sqrt{g(j+p,\alpha+\gamma)}\big].
\label{eq:quarticprojectionbound}
\end{align}
Now it follows from \eqref{eq:quarticproductdifference}, \eqref{eq:quarticdetdifferencebound}, \eqref{eq:quarticdetbounds}, and \eqref{eq:quarticprojectionbound} that
\begin{align}
&|a(k,\eta;\ell,\xi)a(j+p,\alpha+\gamma;p,\gamma)
 -a(k,\eta;p,\gamma)a(j+\ell,\alpha+\xi;\ell,\xi)|\nonumber\\
\le{}&C\Lambda_t(\ell,\xi)\Lambda_t(p,\gamma)
 \big[\Lambda_t(\ell,\xi)+\Lambda_t(p,\gamma)]
 \Lambda_t(j,\alpha)\kappa(\ell,\xi)\kappa(p,\gamma)
 [\av t+\av t^2\sqrt{g(j+p,\alpha+\gamma)}\big].
\label{eq:quarticcommutatorbound}
\end{align}
Substitution in \eqref{eq:quarticproductbound} and
\eqref{eq:quarticcommutatorbound}, followed by positivity of $g$,
gives
\begin{equation}
\begin{aligned}
&|a(k,\eta;p,\gamma)a(j+\ell,\alpha+\xi;\ell,\xi)|\\
&\le C(1+\sqrt{\ell^2+(\xi+\ell t)^2}+\sqrt{p^2+(\gamma+pt)^2})^2\Lambda_t(j,\alpha)^2\kappa(\ell,\xi)\kappa(p,\gamma),\\
&|a(k,\eta;\ell,\xi)a(j+p,\alpha+\gamma;p,\gamma)-a(k,\eta;p,\gamma)a(j+\ell,\alpha+\xi;\ell,\xi)|\\
&\le C(1+\sqrt{\ell^2+(\xi+\ell t)^2}+\sqrt{p^2+(\gamma+pt)^2})^4\Lambda_t(j,\alpha)\kappa(\ell,\xi)\kappa(p,\gamma)\\
&\qquad \times \big(\av t+\av t^2\sqrt{g(j,\alpha)+g(j+\ell,\alpha+\xi)+g(j+p,\alpha+\gamma)+g(k,\eta)}\big).
\label{eq:UWnonzero}
\end{aligned}
\end{equation}
\noindent\emph{Step 2. The first and second differences.}
For the first difference, due to \eqref{eq:twohighset}, we have 
\[
\begin{aligned}
8\sqrt{\ell^2+(\xi+\ell t)^2}
 &\le \sqrt{j^2+(\alpha+jt)^2}
                    -8\sqrt{p^2+(\gamma+pt)^2}
 \le \sqrt{(j+p)^2+(\alpha+jt+\gamma+pt)^2}.
\end{aligned}
\]
Applying \eqref{eq:FD1} at $(j+p,\alpha+\gamma)$ and $(k,\eta)$, and using \eqref{eq:quartichighcomparison}, we obtain 
\begin{equation}
\begin{aligned}
\Lambda_t(j,\alpha)\frac{|h(k,\eta)-h(j+p,\alpha+\gamma)|}
 {\Lambda_t^{N}(k,\eta)\Lambda_t^{N}(j,\alpha)} \le{} &C_N\Lambda_t(j+p,\alpha+\gamma)\frac{|h(k,\eta)-h(j+p,\alpha+\gamma)|}
 {\Lambda_t^{N}(k,\eta)\Lambda_t^{N}(j+p,\alpha+\gamma)}\\
\le{} &C_N\Lambda_t(\ell,\xi)
 \Bigl[\av t\av{\rho^{1/3}t}
       +\av t^2\sqrt{g(k,\eta)+g(j+p,\alpha+\gamma)}\Bigr].
\label{eq:lowfirstdifference}
\end{aligned}
\end{equation}
Multiplying \eqref{eq:lowfirstdifference} by the commutator bound in \eqref{eq:UWnonzero}, and cancelling $\Lambda_t(j,\alpha)^{-1}\Lambda_t(j,\alpha)$, gives
\begin{equation}
\begin{aligned}
&\frac{|h(k,\eta)-h(j+p,\alpha+\gamma)|}
 {\Lambda_t^{N}(k,\eta)\Lambda_t^{N}(j,\alpha)}
 \frac{|a(k,\eta;\ell,\xi)a(j+p,\alpha+\gamma;p,\gamma)
       -a(k,\eta;p,\gamma)a(j+\ell,\alpha+\xi;\ell,\xi)|}{|\ell|}\\
\le{} &\frac{C_N\Lambda_t(\ell,\xi)}{|\ell|}
 \Bigl[\Lambda_t^{4}(\ell,\xi)
                   +\Lambda_t^{4}(p,\gamma)\Bigr]
             \kappa(\ell,\xi)\kappa(p,\gamma)\\
&\quad\times\Bigl[
  \av t^2\av{\rho^{1/3}t}^2
  +\av t^4\bigl(g(j,\alpha)+g(j+\ell,\alpha+\xi)
                       +g(j+p,\alpha+\gamma)+g(k,\eta)\bigr)\Bigr].
\label{eq:lowfirstpolynomial}
\end{aligned}
\end{equation}
For the second term, \eqref{eq:FD2} and the first line of \eqref{eq:UWnonzero} cancel $\Lambda_t(j,\alpha)^{-2}\Lambda_t(j,\alpha)^2$. These give
\begin{equation}
\begin{aligned}
&\frac{|h(k,\eta)-h(j+p,\alpha+\gamma)
                    -h(j+\ell,\alpha+\xi)+h(j,\alpha)|}
 {\Lambda_t^{N}(k,\eta)\Lambda_t^{N}(j,\alpha)}
 \frac{|a(k,\eta;p,\gamma)a(j+\ell,\alpha+\xi;\ell,\xi)|}{|\ell|}\\
\le{} &\frac{C_N\Lambda_t(\ell,\xi)}{|\ell|}
 \Bigl[\Lambda_t^{4}(\ell,\xi)
                    +\Lambda_t^{4}(p,\gamma)\Bigr]
              \kappa(\ell,\xi)\kappa(p,\gamma)\Bigl[\av t^2+\rho^{-1/3}\av t\av{\rho^{1/3}t}^4\\
&\qquad+\av t^4
\bigl(g(j,\alpha)+g(j+\ell,\alpha+\xi)
                         +g(j+p,\alpha+\gamma)+g(k,\eta)\bigr)\Bigr].
\label{eq:lowsecondpolynomial}
\end{aligned}
\end{equation}
Since $0<\rho\le1$, we have $\av t\le\rho^{-1/3}\av{\rho^{1/3}t}$ and $\av t^2\av{\rho^{1/3}t}^2
 \le\rho^{-1/3}\av t\av{\rho^{1/3}t}^3$. Thus the left-hand sides of \eqref{eq:lowfirstpolynomial} and \eqref{eq:lowsecondpolynomial} are bounded by
\begin{equation}
\begin{aligned}
&\frac{C_N\Lambda_t(\ell,\xi)}{|\ell|}
 \Bigl[\Lambda_t^{4}(\ell,\xi)
                   +\Lambda_t^{4}(p,\gamma)\Bigr]
                  \kappa(\ell,\xi)\kappa(p,\gamma)\\
&\times\Bigl[\av t^2
 +\av t^4\bigl(g(j,\alpha)+g(j+\ell,\alpha+\xi)
                         +g(j+p,\alpha+\gamma)+g(k,\eta)\bigr)
 +\rho^{-1/3}\av t\av{\rho^{1/3}t}^4\Bigr].
\label{eq:twolowmajorant}
\end{aligned}
\end{equation}
It remains to divide by the weights at $(\ell,\xi),(p,\gamma)$.

\noindent\emph{Step 3. Completion of the kernel estimate.}
For the first and last terms in \eqref{eq:twolowmajorant}, \eqref{eq:gprimitive} and $\ell\neq0$ give
\begin{equation}
\begin{aligned}
&\frac{\Lambda_t(\ell,\xi)[\Lambda_t^4(\ell,\xi)+\Lambda_t^4(p,\gamma)]}
 {|\ell|\Lambda_t^N(\ell,\xi)\Lambda_t^N(p,\gamma)}
 \kappa(\ell,\xi)\kappa(p,\gamma)\big[\av t^2+\rho^{-1/3}\av t\av{\rho^{1/3}t}^4\big]\\
\le{} &2\Lambda_t^{9-N}(\ell,\xi)\Lambda_t^{9-N}(p,\gamma)
       \big[\av t^2+\rho^{-1/3}\av t\av{\rho^{1/3}t}^4\big]\sqrt{g(\ell,\xi)g(p,\gamma)}.
\label{eq:8}
\end{aligned}
\end{equation}
For the term containing $g$ in \eqref{eq:twolowmajorant}, substituting \eqref{eq:doubletransferseparated} and then applying
\eqref{eq:tkappaphase} and \eqref{eq:tkappatransfer}, we obtain
\begin{equation*}
\begin{aligned}
&\frac{\av t^4}{|\ell|}\kappa(\ell,\xi)\kappa(p,\gamma)
 \bigl(g(j,\alpha)+g(j+\ell,\alpha+\xi)
                     +g(j+p,\alpha+\gamma)+g(k,\eta)\bigr)\\
\le{} &C\av t^2
 \bigg[\frac{\av t\kappa(\ell,\xi)}{|\ell|}
       (1+\ell^2)^2\sqrt{1+(\xi+\ell t)^2}\bigg]
 \big[\av t\kappa(p,\gamma)
       (1+p^2)^2\sqrt{1+(\gamma+pt)^2}\big]
 \sqrt{g(k,\eta)g(j,\alpha)}\\
\le{} &C\av t^2\Lambda_t^{4}(\ell,\xi)
                    \Lambda_t^{5}(p,\gamma)
                    \sqrt{g(k,\eta)g(j,\alpha)}.
\end{aligned}
\end{equation*}
Thus
\begin{equation}
\begin{aligned}
&\frac{\Lambda_t(\ell,\xi)[\Lambda_t^4(\ell,\xi)+\Lambda_t^4(p,\gamma)]}
 {|\ell|\Lambda_t^N(\ell,\xi)\Lambda_t^N(p,\gamma)}
               \av t^4\kappa(\ell,\xi)\kappa(p,\gamma)
 \bigl(g(j,\alpha)+g(j+\ell,\alpha+\xi)
                     +g(j+p,\alpha+\gamma)+g(k,\eta)\bigr)\\
\le{} &C\av t^2
 \Bigl[\Lambda_t^{9-N}(\ell,\xi)
                         \Lambda_t^{5-N}(p,\gamma)
+\Lambda_t^{5-N}(\ell,\xi)
                         \Lambda_t^{9-N}(p,\gamma)\Bigr]
                    \sqrt{g(k,\eta)g(j,\alpha)}\\
\le{} &C\av t^2 \Lambda_t^{9-N}(\ell,\xi)\Lambda_t^{9-N}(p,\gamma)\sqrt{g(k,\eta)g(j,\alpha)}.
\label{eq:lowpairtransfer}
\end{aligned}
\end{equation}
Combining \eqref{eq:lowpairtransfer} with \eqref{eq:lowfirstpolynomial}, \eqref{eq:lowsecondpolynomial}, \eqref{eq:twolowmajorant}, \eqref{eq:8}, and the first inequality of \eqref{eq:quarticexpdecay} proves \eqref{eq:lowpairtarget}.
\end{proof}

\begin{lemma}\label{lem:meanpair}
Let $k=j+\ell+p$, $\eta=\alpha+\xi+\gamma$, and $\ell\ne0$, $p=0$.
On the region \eqref{eq:twohighset}, the weights and kernels defined in \eqref{eq:weight} and \eqref{eq:kernelgeometric} satisfy
\begin{align}
&\frac{e^{\rho^{1/3}t\mathbf1_{k\ne0}\mathbf1_{j+p\ne0}/8}}{e^{\rho^{1/3}t(\mathbf1_{k\ne0}+\mathbf1_{j\ne0}+\mathbf1_{\ell\ne0}+\mathbf1_{p\ne0})/16}}
\frac{1}
 {|\ell|
 \Lambda_t^N(k,\eta)\Lambda_t^N(j,\alpha)
 \Lambda_t^N(\ell,\xi)\Lambda_t^N(0,\gamma)}\nonumber\\
&\quad\times\bigl|
 (h(k,\eta)-h(j,\alpha+\gamma))
       a(k,\eta;\ell,\xi)a(j,\alpha+\gamma;0,\gamma)\nonumber\\
&\qquad-
 (h(j+\ell,\alpha+\xi)-h(j,\alpha))
       a(k,\eta;0,\gamma)a(j+\ell,\alpha+\xi;\ell,\xi)\bigr|\nonumber\\
\le{} &C_Ne^{-\rho^{1/3}t/16}
 \Lambda_t^{9-N}(\ell,\xi)\Lambda_t^{9-N}(0,\gamma)
 \Bigl[\av{\rho^{1/3}t}^2+\av t+\rho^{-1/3}\av{\rho^{1/3}t}^4
       +\av t^2\sqrt{g(k,\eta)g(j,\alpha)}\Bigr].
\label{eq:meanpairtarget}
\end{align}
\end{lemma}

\begin{proof}
\emph{Step 1. Estimates for the two products.}
When $p=0$ and $\gamma\ne0$,
\[
 \frac{p(\alpha+\xi)-(j+\ell)\gamma}{\sqrt{p^2+\gamma^2}}
 =-\operatorname{sgn}(\gamma)(j+\ell)=-\operatorname{sgn}(\gamma)k.
\]
It follows from \eqref{eq:quarticproductidentities} and \eqref{eq:quarticdetbounds} that
\begin{equation}
\begin{aligned}
|a(k,\eta;0,\gamma)a(j+\ell,\alpha+\xi;\ell,\xi)|
&\le |k|\frac{|\ell\alpha-j\xi|}{\sqrt{\ell^2+\xi^2}}
 \le\Lambda_t(\ell,\xi)\Lambda_t(j,\alpha)|k|\kappa(\ell,\xi)\\
&\le C\big[1+\sqrt{\ell^2+(\xi+\ell t)^2}+|\gamma|\big]^2
       \Lambda_t(j,\alpha)|k|\kappa(\ell,\xi).
\end{aligned}
\label{eq:Wmean}
\end{equation}
The case $\gamma=0$ follows from $a(k,\eta;0,0)=0$.
To estimate the difference of the two kernel products, use \eqref{eq:quarticproductdifference}. Both products vanish when $\gamma=0$, so assume $\gamma\ne0$. A direct computation gives
\begin{align}
\frac{|\ell\gamma|\,|\ell\alpha-j(\xi+\gamma)|}
 {\sqrt{\ell^2+\xi^2}|\gamma|}
=\kappa(\ell,\xi)
 |\ell(\alpha+jt)-j(\xi+\ell t+\gamma)|
\le[\Lambda_t(\ell,\xi)+|\gamma|]\Lambda_t(j,\alpha)
           \kappa(\ell,\xi).
\label{eq:meandeterminantdifference}
\end{align}
Then by \eqref{eq:quarticdetbounds}, we obtain
\[
\frac{|-(j+\ell)\gamma(\ell\alpha-j\xi)|}
 {\sqrt{\ell^2+\xi^2}\sqrt{\gamma^2}}=\frac{|k(\ell(\alpha+jt)-j(\xi+\ell t))|}
 {\sqrt{\ell^2+\xi^2}}\le C|k|\Lambda_t(\ell,\xi)\Lambda_t(j,\alpha)\kappa(\ell,\xi).
\]
Since $|e|=1$ at these frequencies, we write the angles of $e(k,\eta),e(j,\alpha)$ as $\varphi_1,\varphi_2$, and those of $e(j+p,\alpha+\gamma),e(j+\ell,\alpha+\xi)$ as $\psi_1,\psi_2$. Hence
\begin{align*}
&\big|[e(k,\eta)\cdot e(j+p,\alpha+\gamma)]
 [e(j+p,\alpha+\gamma)\cdot e(j,\alpha)]\\
&\qquad-[e(k,\eta)\cdot e(j+\ell,\alpha+\xi)]
 [e(j+\ell,\alpha+\xi)\cdot e(j,\alpha)]\big|\\
={}&|\sin(\varphi_1+\varphi_2-\psi_1-\psi_2)
           \sin(\psi_1-\psi_2)|\\
\le{}&|\sin(\psi_1-\psi_2)|
 =\frac{|j(\alpha+\xi)-(j+\ell)(\alpha+\gamma)|}
 {\sqrt{j^2+(\alpha+\gamma)^2}\sqrt{(j+\ell)^2+(\alpha+\xi)^2}}.
\end{align*}
By \eqref{eq:quartichighcomparison}, we have
\begin{align}
|j(\alpha+\xi)-(j+\ell)(\alpha+\gamma)|
 =|j(\xi+\ell t)-\ell(\alpha+jt)-(j+\ell)\gamma|
 \le C[\Lambda_t(\ell,\xi)+|\gamma|]\Lambda_t(j,\alpha).
\label{eq:meanangulardeterminant}
\end{align}
Since \eqref{eq:twohighset} and $\ell\neq0$ imply $\sqrt{j^2+(\alpha+jt)^2}\ge8$, the same triangle inequalities used in \eqref{eq:quartichighphysical} and $k=j+\ell+p=j+\ell$ give
\begin{align}
\frac{\Lambda_t(j,\alpha)}{\sqrt{j^2+(\alpha+\gamma)^2}}
 \le C[1+\av t\kappa(j,\alpha+\gamma)],\qquad
\frac{|k|}{\sqrt{(j+\ell)^2+(\alpha+\xi)^2}}
 \le1.
\label{eq:meanphysicalquotients}
\end{align}
Thus it follows from \eqref{eq:quarticproductdifference}, \eqref{eq:meandeterminantdifference}, \eqref{eq:meanangulardeterminant}, and \eqref{eq:meanphysicalquotients} that
\begin{equation}
\begin{aligned}
&|a(k,\eta;\ell,\xi)a(j,\alpha+\gamma;0,\gamma)-a(k,\eta;0,\gamma)a(j+\ell,\alpha+\xi;\ell,\xi)|\\
\le{} &C(1+\sqrt{\ell^2+(\xi+\ell t)^2}+|\gamma|)^4\Lambda_t(j,\alpha)\kappa(\ell,\xi)\\
& \times \big(1+\av t\sqrt{g(j,\alpha)+g(j+\ell,\alpha+\xi)+g(j+p,\alpha+\gamma)+g(k,\eta)}\big).
\label{eq:UWmean}
\end{aligned}
\end{equation}
\noindent\emph{Step 2. The first difference.}
The derivation of \eqref{eq:lowfirstdifference} uses only
\eqref{eq:twohighset} and \eqref{eq:FD1}, so it also applies when $p=0$. Hence we deduce from $\av{\rho^{1/3}t}\ge1$, \eqref{eq:lowfirstdifference} with $p=0$ and \eqref{eq:UWmean} that
\begin{equation}
\begin{aligned}
&\frac{|h(k,\eta)-h(j,\alpha+\gamma)|}
 {\Lambda_t^{N}(k,\eta)\Lambda_t^{N}(j,\alpha)}
 \frac{|a(k,\eta;\ell,\xi)a(j,\alpha+\gamma;0,\gamma)
       -a(k,\eta;0,\gamma)a(j+\ell,\alpha+\xi;\ell,\xi)|}{|\ell|}\\
\le{}&\frac{C_N\Lambda_t(\ell,\xi)}{|\ell|}
 \Bigl[\Lambda_t^{4}(\ell,\xi)+(1+\gamma^2)^2\Bigr] \av t\kappa(\ell,\xi)\\
&\times\Bigl[
 \av{\rho^{1/3}t}^2
 +\av t^2\bigl(g(j,\alpha)+g(j+\ell,\alpha+\xi)
               +g(j,\alpha+\gamma)+g(k,\eta)\bigr)\Bigr].
\label{eq:meanfirstmajorant}
\end{aligned}
\end{equation}
For $p=0$, \eqref{eq:doubletransferseparated} reads
\[
\begin{aligned}
g(j,\alpha)+g(j+\ell,\alpha+\xi)+g(j,\alpha+\gamma)+g(k,\eta)
\le C(1+\ell^2)^2\sqrt{1+(\xi+\ell t)^2}
          \sqrt{1+\gamma^2}\sqrt{g(k,\eta)g(j,\alpha)}.
\end{aligned}
\]
Multiplying by $\av t\kappa(\ell,\xi)/|\ell|$ and applying
\eqref{eq:tkappaphase} gives
\begin{equation}
\begin{aligned}
\frac{\av t\kappa(\ell,\xi)}{|\ell|}
 \bigl(g(j,\alpha)+g(j+\ell,\alpha+\xi)
                         +g(j,\alpha+\gamma)+g(k,\eta)\bigr)
\le C\Lambda_t^{4}(\ell,\xi)
             \sqrt{1+\gamma^2}\sqrt{g(k,\eta)g(j,\alpha)}.
\label{eq:meanjointtransfer}
\end{aligned}
\end{equation}
Therefore we get, via \eqref{eq:meanfirstmajorant},
\eqref{eq:tkappa} and $|\ell|\ge1$ that 
\begin{equation}
\begin{aligned}
&\frac{|h(k,\eta)-h(j,\alpha+\gamma)|}
 {\Lambda_t^{N}(k,\eta)\Lambda_t^{N}(j,\alpha)}
 \frac{|a(k,\eta;\ell,\xi)a(j,\alpha+\gamma;0,\gamma)
       -a(k,\eta;0,\gamma)a(j+\ell,\alpha+\xi;\ell,\xi)|}{|\ell|}\\
\le{}&C_N\av t^2\Bigl[
 \Lambda_t^{9}(\ell,\xi)\sqrt{1+\gamma^2}+\Lambda_t^{5}(\ell,\xi)(1+\gamma^2)^{5/2}\Bigr]
                    \sqrt{g(k,\eta)g(j,\alpha)}\\
&+C_N\Bigl[\Lambda_t^{6}(\ell,\xi)
       +\Lambda_t^{2}(\ell,\xi)(1+\gamma^2)^2\Bigr]
                  \av{\rho^{1/3}t}^2\\
\le{}&C_N \Lambda_t^9(\ell,\xi)\Lambda_t^5(0,\gamma)\big[\av{\rho^{1/3}t}^2+\av t^2\sqrt{g(k,\eta)g(j,\alpha)}\big].
\label{eq:meanfirstpolynomial}
\end{aligned}
\end{equation}
\noindent\emph{Step 3. The second difference.}
Applying \eqref{eq:FD2} with $p=0$, and then multiplying by
\eqref{eq:Wmean}, gives
\begin{equation}
\begin{aligned}
&\frac{|h(k,\eta)-h(j,\alpha+\gamma)
                   -h(j+\ell,\alpha+\xi)+h(j,\alpha)|}
 {\Lambda_t^{N}(k,\eta)\Lambda_t^{N}(j,\alpha)}
 \frac{|a(k,\eta;0,\gamma)a(j+\ell,\alpha+\xi;\ell,\xi)|}{|\ell|}\\
\le{} &\frac{C_N}{|\ell|}
 \Bigl[\Lambda_t^{4}(\ell,\xi)+(1+\gamma^2)^2\Bigr]
                      \frac{|k|}{\Lambda_t(j,\alpha)}\kappa(\ell,\xi) 
\Bigl[\av t^2+\rho^{-1/3}\av t\av{\rho^{1/3}t}^4\\
&\qquad+\mathbf1_{\substack{t\ge1,\ |j|>16(|\ell|+|p|+1)\\t/8\le(\alpha+jt)/j\le8t}}\av t^4\bigl(g(j,\alpha)+g(j+\ell,\alpha+\xi)+g(j,\alpha+\gamma)+g(k,\eta)\bigr)\Bigr].
\label{eq:meansecond}
\end{aligned}
\end{equation}
For the scalar terms in \eqref{eq:meansecond}, $\ell\neq0$, \eqref{eq:quartichighcomparison} and \eqref{eq:tkappa} give
\[
 \frac{|k|}{\Lambda_t(j,\alpha)}
 \frac{\av t\kappa(\ell,\xi)}{|\ell|}
 \le\frac{\Lambda_t(k,\eta)}{\Lambda_t(j,\alpha)}
 \frac{\av t\kappa(\ell,\xi)}{|\ell|}\le C\Lambda_t(\ell,\xi).
\]
For the term containing $g$, the conditions in
\eqref{eq:criticalregion}, with $p=0$, imply
\begin{align*}
 |\alpha+jt|\ge |j|t/8,\quad 
 |k|=|j+\ell|\le\frac{17}{16}|j|, \quad t\ge1,\quad
 \frac{|k|}{\Lambda_t(j,\alpha)}
 \le\frac{17|j|/16}{|j|t/8}
 =\frac{17}{2t}\le\frac{17\sqrt2}{2\av t}.
\end{align*}
Combining this ratio bound with \eqref{eq:meanjointtransfer}, we obtain
\begin{equation*}
\begin{aligned}
&\frac{\mathbf1_{\substack{t\ge1,\ |j|>16(|\ell|+|p|+1)\\t/8\le(\alpha+jt)/j\le8t}}|k|}{\Lambda_t(j,\alpha)|\ell|}\av t^4\kappa(\ell,\xi)
 \bigl(g(j,\alpha)+g(j+\ell,\alpha+\xi)
                         +g(j,\alpha+\gamma)+g(k,\eta)\bigr)\\
\le{} &C\av t^2\Lambda_t^{4}(\ell,\xi)
             \sqrt{1+\gamma^2}\sqrt{g(k,\eta)g(j,\alpha)}.
\end{aligned}
\end{equation*}
If one of $j,j+\ell,j+p,k$ is zero, the characteristic function in \eqref{eq:meansecond} vanishes, and the remaining terms satisfy the preceding bounds.
Consequently, we get
\begin{equation}
\begin{aligned}
&\frac{|h(k,\eta)-h(j,\alpha+\gamma)
                   -h(j+\ell,\alpha+\xi)+h(j,\alpha)|}
 {\Lambda_t^{N}(k,\eta)\Lambda_t^{N}(j,\alpha)}
 \frac{|a(k,\eta;0,\gamma)a(j+\ell,\alpha+\xi;\ell,\xi)|}{|\ell|}\\
\le{} &C_N\Bigl[\Lambda_t^{5}(\ell,\xi)
 +\Lambda_t(\ell,\xi)(1+\gamma^2)^2\Bigr]
       \bigl[\av t+\rho^{-1/3}\av{\rho^{1/3}t}^4\bigr]\\
&+C_N\av t^2\Bigl[
 \Lambda_t^{8}(\ell,\xi)\sqrt{1+\gamma^2} +\Lambda_t^{4}(\ell,\xi)(1+\gamma^2)^{5/2}\Bigr]
                    \sqrt{g(k,\eta)g(j,\alpha)}\\
\le{}&C_N \Lambda_t^9(\ell,\xi)\Lambda_t^5(0,\gamma)\big[\av t+\rho^{-1/3}\av{\rho^{1/3}t}^4+\av t^2\sqrt{g(k,\eta)g(j,\alpha)}\big].
\label{eq:meansecondscalarpolynomial}
\end{aligned}
\end{equation}
Dividing \eqref{eq:meanfirstpolynomial} and \eqref{eq:meansecondscalarpolynomial} by $\Lambda_t^N(\ell,\xi)\Lambda_t^N(0,\gamma)$ and using \eqref{eq:quarticexpdecay}, \eqref{eq:Ssplit}, and Remark~\ref{rem:1}, we obtain \eqref{eq:meanpairtarget}.
\end{proof}

\begin{lemma}\label{lem:quarticinterface}
Let $k=j+\ell+p$, $\eta=\alpha+\xi+\gamma$, and $\ell,p\ne0$.
On the region \eqref{eq:twohighset}, the weights and kernels defined in \eqref{eq:weight} and \eqref{eq:kernelgeometric} satisfy
\begin{align}
&\frac{\bigl|e^{\rho^{1/3}t\mathbf1_{k\ne0}\mathbf1_{j+p\ne0}/8}
             -e^{\rho^{1/3}t\mathbf1_{j+\ell\ne0}\mathbf1_{j\ne0}/8}\bigr|}{e^{\rho^{1/3}t(\mathbf1_{k\ne0}+\mathbf1_{j\ne0}
+\mathbf1_{\ell\ne0}+\mathbf1_{p\ne0})/16}}\frac{|h(j+\ell,\alpha+\xi)-h(j,\alpha)|}
 {|\ell|
\Lambda_t^N(k,\eta)\Lambda_t^N(j,\alpha)
\Lambda_t^N(\ell,\xi)\Lambda_t^N(p,\gamma)}\nonumber\\
&\times\big|a(k,\eta;p,\gamma)a(j+\ell,\alpha+\xi;\ell,\xi)\big|\nonumber\\
\le{} &C_Ne^{-\rho^{1/3}t/16}
 \Lambda_t^{9-N}(\ell,\xi)\Lambda_t^{9-N}(p,\gamma)
 \Bigl[
 \sqrt{j^2+\alpha^2}\sqrt{g(\ell,\xi)}
       +\av t^2\sqrt{g(\ell,\xi)g(p,\gamma)}\Bigr].
\label{eq:interfacetarget}
\end{align}
\end{lemma}
\begin{proof}
If $jk(j+p)(j+\ell)\neq0$, then the left-hand side of \eqref{eq:interfacetarget} is zero.
It remains to consider the case in which at least one of $j,k,j+p,j+\ell$ is zero. In this case, the relation $k=j+\ell+p$ gives
\begin{equation}
 |j|\le|\ell|+|p|,\qquad
 |j|+|j+\ell|\le3(|\ell|+|p|).
\label{eq:interfacehorizontal}
\end{equation}
Therefore
\begin{align}
 \big[1+(|j|+|j+\ell|)^2\big]\Lambda_t^2(\ell,\xi)
 &\le C(1+\ell^2)(1+p^2)\Lambda_t^2(\ell,\xi)
 \le C\Lambda_t^4(\ell,\xi)\Lambda_t^2(p,\gamma).
\label{eq:interfaceweightpowers}
\end{align}
Applying \eqref{eq:FDsmallh} to the difference at $(j+\ell,\alpha+\xi)$ and $(j,\alpha)$, and then using  \eqref{eq:interfaceweightpowers}, gives
\begin{equation}
 \frac{|h(j+\ell,\alpha+\xi)-h(j,\alpha)|}
 {\Lambda_t^N(k,\eta)\Lambda_t^N(j,\alpha)}
 \le \frac{C_N\Lambda_t^4(\ell,\xi)\Lambda_t^2(p,\gamma)\av t}
 {\Lambda_t(j,\alpha)}.
\label{eq:interfacefirstdifference}
\end{equation}
Here \eqref{eq:twohighset} verifies the frequency condition of \eqref{eq:FDsmallh}. 
Multiplying \eqref{eq:interfacefirstdifference} by the first inequality of \eqref{eq:UWnonzero}, and using $|\ell|\ge1$, we find
\begin{equation}
\begin{aligned}
 &\frac{|h(j+\ell,\alpha+\xi)-h(j,\alpha)|
 |a(k,\eta;p,\gamma)a(j+\ell,\alpha+\xi;\ell,\xi)|}
 {|\ell|\Lambda_t^N(k,\eta)\Lambda_t^N(j,\alpha)}\\
 \le{} &C_N\Lambda_t^6(\ell,\xi)\Lambda_t^4(p,\gamma)
       \av t\Lambda_t(j,\alpha)\kappa(\ell,\xi)\kappa(p,\gamma).
\end{aligned}
\label{eq:interfacemajorant}
\end{equation}
To estimate the remaining factor $\Lambda_t(j,\alpha)$, \eqref{eq:interfacehorizontal} gives
\begin{align*}
\Lambda_t(j,\alpha)
 \le1+\sqrt{j^2+\alpha^2}+t|j|
 \le\sqrt{j^2+\alpha^2}
        +C\av t\Lambda_t(\ell,\xi)\Lambda_t(p,\gamma).
\end{align*}
Substituting this into \eqref{eq:interfacemajorant} and applying \eqref{eq:tkappa} to $\av t\kappa(p,\gamma)$, we obtain
\begin{align}
&C_N\Lambda_t^6(\ell,\xi)\Lambda_t^4(p,\gamma)
 \av t\Lambda_t(j,\alpha)\kappa(\ell,\xi)\kappa(p,\gamma)\nonumber\\
\le{} &C_N\Lambda_t^6(\ell,\xi)\Lambda_t^4(p,\gamma)
 \av t\sqrt{j^2+\alpha^2}\kappa(\ell,\xi)\kappa(p,\gamma)
 +C_N\Lambda_t^7(\ell,\xi)\Lambda_t^5(p,\gamma)
 \av t^2\kappa(\ell,\xi)\kappa(p,\gamma)\nonumber\\
\le{} &C_N\Lambda_t^6(\ell,\xi)\Lambda_t^5(p,\gamma)
 \sqrt{j^2+\alpha^2}\kappa(\ell,\xi)
 +C_N\Lambda_t^7(\ell,\xi)\Lambda_t^5(p,\gamma)
 \av t^2\kappa(\ell,\xi)\kappa(p,\gamma).
\label{eq:interfacepolynomial}
\end{align}
Hence by \eqref{eq:gprimitive}, division of \eqref{eq:interfacepolynomial} by $\Lambda_t^N(\ell,\xi)\Lambda_t^N(p,\gamma)$ gives
\begin{align}
&\frac{|h(j+\ell,\alpha+\xi)-h(j,\alpha)|
 |a(k,\eta;p,\gamma)a(j+\ell,\alpha+\xi;\ell,\xi)|}
 {|\ell|\Lambda_t^N(k,\eta)\Lambda_t^N(j,\alpha)
          \Lambda_t^N(\ell,\xi)\Lambda_t^N(p,\gamma)}\nonumber\\
\le{} &C_N\Lambda_t^{9-N}(\ell,\xi)\Lambda_t^{9-N}(p,\gamma)
 \big[\sqrt{j^2+\alpha^2}\sqrt{g(\ell,\xi)}
       +\av t^2\sqrt{g(\ell,\xi)g(p,\gamma)}\big].
\label{eq:interfacenormalized}
\end{align}
Finally, if the exponential difference vanishes there is nothing to prove. Otherwise exactly one numerator exponential is
$e^{\rho^{1/3}t/8}$. Together with $\ell,p\ne0$, this forces at least three of $k,j,\ell,p$ to be nonzero. Therefore
\[
\begin{aligned}
\frac{|e^{\rho^{1/3}t\mathbf1_{k\ne0}\mathbf1_{j+p\ne0}/8}
       -e^{\rho^{1/3}t\mathbf1_{j+\ell\ne0}\mathbf1_{j\ne0}/8}|}
 {e^{\rho^{1/3}t(\mathbf1_{k\ne0}+\mathbf1_{j\ne0}
                     +\mathbf1_{\ell\ne0}+\mathbf1_{p\ne0})/16}}
\le2e^{\rho^{1/3}t/8-3\rho^{1/3}t/16}
 =2e^{-\rho^{1/3}t/16}.
\end{aligned}
\]
Combining this inequality with \eqref{eq:interfacenormalized} proves \eqref{eq:interfacetarget}.
\end{proof}

\begin{lemma}\label{lem:otherpairs}
Let $k=j+\ell+p$, $\eta=\alpha+\xi+\gamma$, and $\ell\ne0$.
On the last five regions in \eqref{eq:pairpartition}, the weights and kernels defined in \eqref{eq:weight} and \eqref{eq:kernelgeometric} satisfy
\begin{align}
&\frac{1}{|\ell|e^{\rho^{1/3}t(\mathbf1_{k\ne0}+\mathbf1_{j\ne0}
+\mathbf1_{\ell\ne0}+\mathbf1_{p\ne0})/16}
\Lambda_t^N(k,\eta)\Lambda_t^N(j,\alpha)
\Lambda_t^N(\ell,\xi)\Lambda_t^N(p,\gamma)}\nonumber\\
&\times\Bigl|
 e^{\rho^{1/3}t\mathbf1_{k\ne0}\mathbf1_{j+p\ne0}/8}
 (h(k,\eta)-h(j+p,\alpha+\gamma))
 a(k,\eta;\ell,\xi)a(j+p,\alpha+\gamma;p,\gamma)\nonumber\\
&\quad-
 e^{\rho^{1/3}t\mathbf1_{j+\ell\ne0}\mathbf1_{j\ne0}/8}
 (h(j+\ell,\alpha+\xi)-h(j,\alpha))
 a(k,\eta;p,\gamma)a(j+\ell,\alpha+\xi;\ell,\xi)\Bigr|\nonumber\\
\le{} &C_Ne^{-\rho^{1/3}t/16}
\begin{cases}
\begin{aligned}
 &\Lambda_t^{9-N}(k,\eta)\Lambda_t^{9-N}(j,\alpha)\\
 &\quad\times\big[1+\av t(\sqrt{g(\ell,\xi)}+\sqrt{g(p,\gamma)})
                    +\av t^2\sqrt{g(\ell,\xi)g(p,\gamma)}\big],
\end{aligned}&\text{region }2,\\[3pt]
\begin{aligned}
 &\Lambda_t^{9-N}(j,\alpha)\Lambda_t^{9-N}(p,\gamma)\\
 &\quad\times\big[1+\av t\sqrt{g(\ell,\xi)}
       +\sqrt{\ell^2+\xi^2}\sqrt{g(p,\gamma)}
       +\av t^2\sqrt{g(\ell,\xi)g(p,\gamma)}\big],
\end{aligned}&\text{region }3,\\[3pt]
\begin{aligned}
 &\Lambda_t^{9-N}(k,\eta)\Lambda_t^{9-N}(p,\gamma)\\
 &\quad\times\big[1+\av t\sqrt{g(\ell,\xi)}
       +\sqrt{\ell^2+\xi^2}\sqrt{g(p,\gamma)}
       +\av t^2\sqrt{g(\ell,\xi)g(p,\gamma)}\big],
\end{aligned}&\text{region }4,\\[3pt]
\begin{aligned}
 &\Lambda_t^{9-N}(j,\alpha)\Lambda_t^{9-N}(\ell,\xi)\\
 &\quad\times\big[\av{\rho^{1/3}t}^2+\av t\sqrt{g(p,\gamma)}
                     +\av t^2\sqrt{g(k,\eta)g(p,\gamma)}\big],
\end{aligned}&\text{region }5,\\[3pt]
\begin{aligned}
 &\Lambda_t^{9-N}(k,\eta)\Lambda_t^{9-N}(\ell,\xi)\\
 &\quad\times\big[\av{\rho^{1/3}t}^2+\av t\sqrt{g(p,\gamma)}
                     +\av t^2\sqrt{g(j,\alpha)g(p,\gamma)}\big],
\end{aligned}&\text{region }6.
\end{cases}
\label{eq:otherpairstarget}
\end{align}
The region numbers refer to the corresponding inequalities in \eqref{eq:pairpartition}.
\end{lemma}

\begin{proof}
The triangle inequality and the six conditions in
\eqref{eq:pairpartition} give
\begin{equation}
\begin{aligned}
&\big|\sqrt{\ell^2+(\xi+\ell t)^2}-\sqrt{p^2+(\gamma+pt)^2}\big|
 \le\sqrt{k^2+(\eta+kt)^2}+\sqrt{j^2+(\alpha+jt)^2}
 ,
 &&\text{region }2,\\
&\big|\sqrt{k^2+(\eta+kt)^2}-\sqrt{\ell^2+(\xi+\ell t)^2}\big|
 \le\sqrt{j^2+(\alpha+jt)^2}+\sqrt{p^2+(\gamma+pt)^2}
 ,
 &&\text{region }3,\\
&\big|\sqrt{j^2+(\alpha+jt)^2}-\sqrt{\ell^2+(\xi+\ell t)^2}\big|
 \le\sqrt{k^2+(\eta+kt)^2}+\sqrt{p^2+(\gamma+pt)^2}
 ,
 &&\text{region }4,\\
&\big|\sqrt{k^2+(\eta+kt)^2}-\sqrt{p^2+(\gamma+pt)^2}\big|
 \le\sqrt{j^2+(\alpha+jt)^2}+\sqrt{\ell^2+(\xi+\ell t)^2}
 ,
 &&\text{region }5,\\
&\big|\sqrt{j^2+(\alpha+jt)^2}-\sqrt{p^2+(\gamma+pt)^2}\big|
 \le\sqrt{k^2+(\eta+kt)^2}+\sqrt{\ell^2+(\xi+\ell t)^2}
 ,
 &&\text{region }6.
\end{aligned}
\label{eq:otherpairsizecomparisons}
\end{equation}
In each line, the right-hand side is at most $1/8$ times the smaller of the two lengths on the left, by the corresponding condition in \eqref{eq:pairpartition}. 
Thus the two larger Sobolev weights in each region are comparable. 
Since $\Lambda_t(j+p,\alpha+\gamma)\le\Lambda_t(j,\alpha)+\Lambda_t(p,\gamma)$, $\Lambda_t(j+\ell,\alpha+\xi)\le\Lambda_t(j,\alpha)+\Lambda_t(\ell,\xi)$, we obtain from \eqref{eq:weight}, \eqref{eq:thetabounds} and \eqref{eq:otherpairsizecomparisons} that
\begin{equation}
\begin{aligned}
h(k,\eta)+h(j+p,\alpha+\gamma)+h(j+\ell,\alpha+\xi)+h(j,\alpha)
\le C_N
\begin{cases}
\Lambda_t^N(\ell,\xi)\Lambda_t^N(p,\gamma),&\text{region }2,\\
\Lambda_t^N(k,\eta)\Lambda_t^N(\ell,\xi),&\text{region }3,\\
\Lambda_t^N(j,\alpha)\Lambda_t^N(\ell,\xi),&\text{region }4,\\
\Lambda_t^N(k,\eta)\Lambda_t^N(p,\gamma),&\text{region }5,\\
\Lambda_t^N(j,\alpha)\Lambda_t^N(p,\gamma),&\text{region }6.
\end{cases}
\end{aligned}
\label{eq:otherpairenergyweights}
\end{equation}

\noindent\emph{Estimate on Region 2.}
By virtue of \eqref{eq:A}, we can use \eqref{eq:elementarysymbols} to conclude that 
\begin{align*}
 |a(k,\eta;\ell,\xi)|&\le C\Lambda_t(k,\eta)(1+\av t\kappa(\ell,\xi)),\\
 |a(j+p,\alpha+\gamma;p,\gamma)|&\le C\Lambda_t(j,\alpha)(1+\av t\kappa(p,\gamma)),\\
 |a(k,\eta;p,\gamma)|&\le C\Lambda_t(k,\eta)(1+\av t\kappa(p,\gamma)),\\
 |a(j+\ell,\alpha+\xi;\ell,\xi)|&\le C\Lambda_t(j,\alpha)(1+\av t\kappa(\ell,\xi)).
\end{align*}
For $p=0$, \eqref{eq:meanA} gives the same bound. Thus we can further obtain
\[
\begin{aligned}
&|a(k,\eta;\ell,\xi)a(j+p,\alpha+\gamma;p,\gamma)|
\le C\Lambda_t(k,\eta)\Lambda_t(j,\alpha)(1+\av t\kappa(\ell,\xi))(1+\av t\kappa(p,\gamma)),\\
&|a(k,\eta;p,\gamma)a(j+\ell,\alpha+\xi;\ell,\xi)|
\le C\Lambda_t(k,\eta)\Lambda_t(j,\alpha)(1+\av t\kappa(p,\gamma))(1+\av t\kappa(\ell,\xi)).
\end{aligned}
\]
Combining these bounds with region 2 of \eqref{eq:otherpairenergyweights} and $|\ell|\ge1$, we obtain
\begin{equation}
\begin{aligned}
 &\frac{|h(k,\eta)-h(j+p,\alpha+\gamma)|
       |a(k,\eta;\ell,\xi)a(j+p,\alpha+\gamma;p,\gamma)|}
 {|\ell|\Lambda_t^N(\ell,\xi)\Lambda_t^N(p,\gamma)}\\
 &+\frac{|h(j+\ell,\alpha+\xi)-h(j,\alpha)|
       |a(k,\eta;p,\gamma)a(j+\ell,\alpha+\xi;\ell,\xi)|}
 {|\ell|\Lambda_t^N(\ell,\xi)\Lambda_t^N(p,\gamma)}\\
 \le{} &C_N\Lambda_t(k,\eta)\Lambda_t(j,\alpha)
 \big[1+\av t\kappa(\ell,\xi)+\av t\kappa(p,\gamma)
                        +\av t^2\kappa(\ell,\xi)\kappa(p,\gamma)\big].
\end{aligned}
\label{eq:pairLP}
\end{equation}
Dividing \eqref{eq:pairLP} by $\Lambda_t^N(k,\eta)\Lambda_t^N(j,\alpha)$ and applying \eqref{eq:gprimitive} and \eqref{eq:quarticexpdecay}, we bound the left-hand side of \eqref{eq:otherpairstarget} in region 2 by
\begin{align*}
C_Ne^{-\rho^{1/3}t/16}\Lambda_t^{1-N}(k,\eta)\Lambda_t^{1-N}(j,\alpha)
\big[1+\av t(\sqrt{g(\ell,\xi)}+\sqrt{g(p,\gamma)})
                 +\av t^2\sqrt{g(\ell,\xi)g(p,\gamma)}\big].
\end{align*}
Since $\Lambda_t^{1-N}\le\Lambda_t^{9-N}$, this proves region 2 of \eqref{eq:otherpairstarget}.

\noindent\emph{Estimate on Region 3.} 
For $p\ne0$, the region condition gives
\[
 |\Lambda_t(j+\ell,\alpha+\xi)-\Lambda_t(\ell,\xi)|
 \le\sqrt{j^2+(\alpha+jt)^2}\le\frac18\Lambda_t(\ell,\xi).
\]
In particular, $\Lambda_t(j+\ell,\alpha+\xi)\approx\Lambda_t(\ell,\xi)$.
Since $(k,\eta)=(j+\ell,\alpha+\xi)+(p,\gamma)$,
\eqref{eq:A}, \eqref{eq:elementarysymbols}, and \eqref{eq:lowA} imply
\begin{align*}
 |a(k,\eta;\ell,\xi)|&\le C\Lambda_t(j+p,\alpha+\gamma)(1+\av t\kappa(\ell,\xi)),\\
 |a(j+p,\alpha+\gamma;p,\gamma)|& \le C\Lambda_t(j,\alpha)\Lambda_t(p,\gamma),\\
 |a(k,\eta;p,\gamma)|&\le C\Lambda_t(p,\gamma) \Lambda_t(j+\ell,\alpha+\xi)\kappa(p,\gamma)\le C\Lambda_t(p,\gamma) \Lambda_t(\ell,\xi)\kappa(p,\gamma),\\
 |a(j+\ell,\alpha+\xi;\ell,\xi)|&\le C\frac{\Lambda_t(j,\alpha) \Lambda_t(\ell,\xi)}{\sqrt{\ell^2+\xi^2}}.
\end{align*}
Since $\ell\ne0$, we have
\begin{align}\label{eq:16}
\Lambda_t(\ell,\xi)^2
 =1+\ell^2+(\xi+\ell t)^2
 \le C(\ell^2+\xi^2+t^2\ell^2),	
\end{align}
which implies
\begin{equation}
\begin{aligned}
\frac{\Lambda_t(\ell,\xi)^2}{|\ell|\sqrt{\ell^2+\xi^2}}\kappa(p,\gamma)
&\le C\bigg[\frac{\sqrt{\ell^2+\xi^2}}{|\ell|}
           +t^2\frac{|\ell|}{\sqrt{\ell^2+\xi^2}}\bigg]\kappa(p,\gamma)\\
&\le C\big[\sqrt{\ell^2+\xi^2}\kappa(p,\gamma)
             +\av t^2\kappa(\ell,\xi)\kappa(p,\gamma)\big].
\end{aligned}
\label{eq:remainingderivative}
\end{equation}
For $p=0$, \eqref{eq:meanA} implies
\[
|a(j,\alpha+\gamma;0,\gamma)|\le|j|,\qquad |a(k,\eta;0,\gamma)|\le|j+\ell|,\qquad
\frac{|j+\ell|}{|\ell|}\le1+|j|.
\]
Multiplying the first two kernel bounds and using $|\ell|\ge1$ gives 
\[
\begin{aligned}
\frac{|a(k,\eta;\ell,\xi)a(j+p,\alpha+\gamma;p,\gamma)|}{|\ell|}
&\le C\Lambda_t(j+p,\alpha+\gamma)\Lambda_t(p,\gamma)
\Lambda_t(j,\alpha)(1+\av t\kappa(\ell,\xi))\\
&\le C\Lambda_t^2(p,\gamma)\Lambda_t^2(j,\alpha)(1+\av t\kappa(\ell,\xi)).
\end{aligned}
\]
This inequality also holds for $p=0$. For the second product with $p\ne0$, we obtain
\[
\begin{aligned}
\frac{|a(k,\eta;p,\gamma)a(j+\ell,\alpha+\xi;\ell,\xi)|}{|\ell|}
\le C\Lambda_t(j,\alpha)\Lambda_t(p,\gamma)
\big[\sqrt{\ell^2+\xi^2}\kappa(p,\gamma)
             +\av t^2\kappa(\ell,\xi)\kappa(p,\gamma)\big].
\end{aligned}
\]
When $p=0$, its bound instead is $C\Lambda_t^2(j,\alpha)(1+\av t\kappa(\ell,\xi))$, since $|j+\ell|/|\ell|\le1+|j|$ and \eqref{eq:16}. Combining these bounds with region 3 of \eqref{eq:otherpairenergyweights} and \eqref{eq:gprimitive}, we obtain 
\begin{align*}
 &\frac{|h(k,\eta)-h(j+p,\alpha+\gamma)|
       |a(k,\eta;\ell,\xi)a(j+p,\alpha+\gamma;p,\gamma)|}
 {|\ell|\Lambda_t^N(k,\eta)\Lambda_t^N(j,\alpha)\Lambda_t^N(\ell,\xi)\Lambda_t^N(p,\gamma)}\\
 &+\frac{|h(j+\ell,\alpha+\xi)-h(j,\alpha)|
       |a(k,\eta;p,\gamma)a(j+\ell,\alpha+\xi;\ell,\xi)|}
 {|\ell|\Lambda_t^N(k,\eta)\Lambda_t^N(j,\alpha)\Lambda_t^N(\ell,\xi)\Lambda_t^N(p,\gamma)}\\
 \le{} &C_N\Lambda_t^{2-N}(j,\alpha)\Lambda_t^{2-N}(p,\gamma)
\big[1+\av t\sqrt{g(\ell,\xi)}
 +\sqrt{\ell^2+\xi^2}\sqrt{g(p,\gamma)}
 +\av t^2\sqrt{g(\ell,\xi)g(p,\gamma)}\big].
\end{align*}
Using $\Lambda_t^{2-N}\le\Lambda_t^{9-N}$ and \eqref{eq:quarticexpdecay}, we obtain the bound for region 3 in \eqref{eq:otherpairstarget}.

\smallskip\noindent\emph{Estimate on Region 4.}
By \eqref{eq:A} and \eqref{eq:meanA}, we have $|a(k,\eta;p,\gamma)|\le C\Lambda_t(k,\eta)\Lambda_t(p,\gamma)$.
By \eqref{eq:endpointsym}, \eqref{eq:elementarysymbols} and $j+\ell=k-p$, $\alpha+\xi=\eta-\gamma$,
\[
 |a(j+\ell,\alpha+\xi;\ell,\xi)|
 \le C\Lambda_t(j+\ell,\alpha+\xi)(1+\av t\kappa(\ell,\xi))
 \le C\Lambda_t(k,\eta)\Lambda_t(p,\gamma)(1+\av t\kappa(\ell,\xi)).
\]
Multiplication and $|\ell|\ge1$ give
\[
\begin{aligned}
\frac{|a(k,\eta;p,\gamma)a(j+\ell,\alpha+\xi;\ell,\xi)|}{|\ell|}
&\le C\Lambda_t^2(k,\eta)\Lambda_t^2(p,\gamma)
                       (1+\av t\kappa(\ell,\xi)).
\end{aligned}
\]
For the other product, when $p\ne0$, $j+p=k-\ell$, $\alpha+\gamma=\eta-\xi$, \eqref{eq:A}, \eqref{eq:endpointsym}, \eqref{eq:elementarysymbols} and \eqref{eq:pairpartition} give
\begin{align*}
 |a(k,\eta;\ell,\xi)|=|a(j+p,\alpha+\gamma;-\ell,-\xi)|
       &\le\frac{|k(\xi+\ell t)-\ell(\eta+kt)|}{\sqrt{\ell^2+\xi^2}}\le C\frac{\Lambda_t(k,\eta)\Lambda_t(\ell,\xi)}{\sqrt{\ell^2+\xi^2}},\\
 |a(j+p,\alpha+\gamma;p,\gamma)|&\le C\Lambda_t(p,\gamma)\Lambda_t(k-\ell,\eta-\xi)\kappa(p,\gamma)\le C\Lambda_t(p,\gamma)\Lambda_t(\ell,\xi)\kappa(p,\gamma).
\end{align*}
Dividing by $|\ell|$ and using \eqref{eq:remainingderivative},
we get, for $p\ne0$,
\[
\begin{aligned}
 \frac{|a(k,\eta;\ell,\xi)a(j+p,\alpha+\gamma;p,\gamma)|}{|\ell|}
 \le C\Lambda_t(k,\eta)\Lambda_t(p,\gamma)
 \bigl[\sqrt{\ell^2+\xi^2}\kappa(p,\gamma)
             +\av t^2\kappa(\ell,\xi)\kappa(p,\gamma)\bigr].
\end{aligned}
\]
For $p=0$, the corresponding bound is $C\Lambda_t^2(k,\eta)(1+\av t\kappa(\ell,\xi))$, by \eqref{eq:meanA}, $|k-\ell|/|\ell|\le1+|k|$, and \eqref{eq:16}. Combining these bounds with region 4 of \eqref{eq:otherpairenergyweights} and \eqref{eq:gprimitive}, we obtain
\begin{align*}
 &\frac{|h(k,\eta)-h(j+p,\alpha+\gamma)|
       |a(k,\eta;\ell,\xi)a(j+p,\alpha+\gamma;p,\gamma)|}
 {|\ell|\Lambda_t^N(k,\eta)\Lambda_t^N(j,\alpha)\Lambda_t^N(\ell,\xi)\Lambda_t^N(p,\gamma)}\\
 &+\frac{|h(j+\ell,\alpha+\xi)-h(j,\alpha)|
       |a(k,\eta;p,\gamma)a(j+\ell,\alpha+\xi;\ell,\xi)|}
 {|\ell|\Lambda_t^N(k,\eta)\Lambda_t^N(j,\alpha)\Lambda_t^N(\ell,\xi)\Lambda_t^N(p,\gamma)}\\
 \le{} &C_N\Lambda_t^{2-N}(k,\eta)\Lambda_t^{2-N}(p,\gamma)
 \big[1+\av t\sqrt{g(\ell,\xi)}
 +\sqrt{\ell^2+\xi^2}\sqrt{g(p,\gamma)}
 +\av t^2\sqrt{g(\ell,\xi)g(p,\gamma)}\big].
\end{align*}
Using $\Lambda_t^{2-N}\le\Lambda_t^{9-N}$ and \eqref{eq:quarticexpdecay}, we obtain the bound for region 4 in \eqref{eq:otherpairstarget}.

\smallskip\noindent\emph{Estimate on Region 5.}
The fifth condition in \eqref{eq:pairpartition} gives
\[
 8\sqrt{\ell^2+(\xi+\ell t)^2}
 \le\sqrt{p^2+(\gamma+pt)^2}-8\sqrt{j^2+(\alpha+jt)^2}
 \le\sqrt{(j+p)^2+(\alpha+\gamma+(j+p)t)^2}.
\]
The region condition also gives
\[
 |\Lambda_t(j+p,\alpha+\gamma)-\Lambda_t(p,\gamma)|
 \le\sqrt{j^2+(\alpha+jt)^2}\le\frac18\Lambda_t(p,\gamma).
\]
Thus $\Lambda_t(j+p,\alpha+\gamma)\approx\Lambda_t(p,\gamma)$.
Applying \eqref{eq:FD1} at $(j+p,\alpha+\gamma)$ and $(k,\eta)$, we obtain
\begin{align}
\frac{|h(k,\eta)-h(j+p,\alpha+\gamma)|}
 {\Lambda_t^N(k,\eta)\Lambda_t^N(p,\gamma)}
\le\frac{C_N\Lambda_t(\ell,\xi)}
              {\Lambda_t(p,\gamma)}
 \big[\av t\av{\rho^{1/3}t}+\av t^2\sqrt{g(k,\eta)+g(j+p,\alpha+\gamma)}\big].
\label{eq:pairKPfirstdifference}
\end{align}
For the two kernels, \eqref{eq:lowA}, \eqref{eq:elementarysymbols}, and \eqref{eq:meanA} give
\begin{align*}
 \frac{|a(k,\eta;\ell,\xi)|}{|\ell|}
 &\le C\Lambda_t(\ell,\xi)\Lambda_t(j+p,\alpha+\gamma)\kappa(\ell,\xi)\le C\Lambda_t(\ell,\xi)\Lambda_t(p,\gamma)\kappa(\ell,\xi),\\
 |a(j+p,\alpha+\gamma;p,\gamma)|
 &\le C\Lambda_t(j,\alpha)(1+\av t\kappa(p,\gamma)).
\end{align*}
Multiplying these bounds by \eqref{eq:pairKPfirstdifference} gives
\begin{align}
&\frac{|h(k,\eta)-h(j+p,\alpha+\gamma)|
 |a(k,\eta;\ell,\xi)a(j+p,\alpha+\gamma;p,\gamma)|}
 {|\ell|\Lambda_t^N(k,\eta)\Lambda_t^N(p,\gamma)}\nonumber\\
\le{} &C_N\Lambda_t(j,\alpha)\Lambda_t^2(\ell,\xi)
 \big[\av t\av{\rho^{1/3}t}+\av t^2
 \sqrt{g(k,\eta)+g(j+p,\alpha+\gamma)}\big]\,
 \kappa(\ell,\xi)(1+\av t\kappa(p,\gamma)).
\label{eq:pairKPraw}
\end{align}
By \eqref{eq:tkappa}, \eqref{eq:gprimitive},
\begin{align}
&\big[\av t\av{\rho^{1/3}t}+\av t^2
 \sqrt{g(k,\eta)+g(j+p,\alpha+\gamma)}\big]\,
 \kappa(\ell,\xi)(1+\av t\kappa(p,\gamma))\nonumber\\
\le{} &C\Lambda_t(\ell,\xi)
 \big[\av{\rho^{1/3}t}+\av t\sqrt{g(k,\eta)+g(j+p,\alpha+\gamma)}\big]
 \big[1+\av t\sqrt{g(p,\gamma)}\big]\nonumber\\
\le{} &C\Lambda_t(\ell,\xi)
 \big[\av{\rho^{1/3}t}^2+\av t^2
       (g(k,\eta)+g(j+p,\alpha+\gamma)+g(p,\gamma))\big].
\label{eq:pairKP}
\end{align}
Next, \eqref{eq:triangletransfer} gives
\[
 g(k,\eta)+g(j+p,\alpha+\gamma)+g(p,\gamma)
 \le C\Lambda_t^5(j,\alpha)\Lambda_t^5(\ell,\xi)
                 \sqrt{g(k,\eta)g(p,\gamma)}.
\]
Substituting this and \eqref{eq:pairKP} into \eqref{eq:pairKPraw},
we obtain
\begin{align}
&\frac{|h(k,\eta)-h(j+p,\alpha+\gamma)|
 |a(k,\eta;\ell,\xi)a(j+p,\alpha+\gamma;p,\gamma)|}
 {|\ell|\Lambda_t^N(k,\eta)\Lambda_t^N(p,\gamma)}\nonumber\\
\le{} &C_N\Lambda_t(j,\alpha)\Lambda_t^3(\ell,\xi)\av{\rho^{1/3}t}^2
 +C_N\Lambda_t^6(j,\alpha)\Lambda_t^8(\ell,\xi)
           \av t^2\sqrt{g(k,\eta)g(p,\gamma)}.
\label{eq:pairKPpolynomial}
\end{align}
For the other product, \eqref{eq:elementarysymbols}, \eqref{eq:meanA}, and \eqref{eq:lowA} give
\begin{align*}
 |a(k,\eta;p,\gamma)|
 &\le C\Lambda_t(j+\ell,\alpha+\xi)(1+\av t\kappa(p,\gamma)),\\
 \frac{|a(j+\ell,\alpha+\xi;\ell,\xi)|}{|\ell|}
 &\le C\Lambda_t(\ell,\xi)\Lambda_t(j,\alpha).
\end{align*}
Since $\Lambda_t(j+\ell,\alpha+\xi)\le
2\Lambda_t(j,\alpha)\Lambda_t(\ell,\xi)$, region 5 of
\eqref{eq:otherpairenergyweights} now implies
\begin{align}
\frac{|h(j+\ell,\alpha+\xi)-h(j,\alpha)|
 |a(k,\eta;p,\gamma)a(j+\ell,\alpha+\xi;\ell,\xi)|}
 {|\ell|\Lambda_t^N(k,\eta)\Lambda_t^N(p,\gamma)}
\le C_N\Lambda_t^2(j,\alpha)\Lambda_t^2(\ell,\xi)
        (1+\av t\kappa(p,\gamma)).
\label{eq:pairKPsecondproduct}
\end{align}
All powers are at most $6$ at $(j,\alpha)$ and $8$ at $(\ell,\xi)$. Adding \eqref{eq:pairKPpolynomial} and \eqref{eq:pairKPsecondproduct}, dividing by $\Lambda_t^N(j,\alpha)\Lambda_t^N(\ell,\xi)$, and using \eqref{eq:gprimitive} and \eqref{eq:quarticexpdecay}, we bound the left-hand side of \eqref{eq:otherpairstarget} in region 5 by
\begin{align*}
&C_Ne^{-\rho^{1/3}t/16}\Lambda_t^{9-N}(j,\alpha)\Lambda_t^{9-N}(\ell,\xi)
\big[\av{\rho^{1/3}t}^2+\av t\sqrt{g(p,\gamma)}+\av t^2\sqrt{g(k,\eta)g(p,\gamma)}\big],
\end{align*}
which proves region 5 of \eqref{eq:otherpairstarget}.

\smallskip\noindent\emph{Estimate on Region 6.}
The sixth condition of \eqref{eq:pairpartition} ensures
\[
 8\sqrt{\ell^2+(\xi+\ell t)^2}\le\sqrt{j^2+(\alpha+jt)^2}.
\]
The last inequality of \eqref{eq:otherpairsizecomparisons} and \eqref{eq:pairpartition} yield $|\Lambda_t(j+\ell,\alpha+\xi)-\Lambda_t(j,\alpha)|\le\frac18\Lambda_t(j,\alpha)$ and $|\Lambda_t(j,\alpha)-\Lambda_t(p,\gamma)|\le\frac18\Lambda_t(p,\gamma)$, thus $\Lambda_t(j+\ell,\alpha+\xi)\approx\Lambda_t(j,\alpha)\approx\Lambda_t(p,\gamma)$. 
It follows from \eqref{eq:FD1} and  \eqref{eq:otherpairsizecomparisons} that
\begin{align}
\frac{|h(j+\ell,\alpha+\xi)-h(j,\alpha)|}
 {\Lambda_t^N(j,\alpha)\Lambda_t^N(p,\gamma)}
\le\frac{C_N\Lambda_t(k,\eta)\Lambda_t(\ell,\xi)}
              {\Lambda_t(p,\gamma)}
 \big[\av t\av{\rho^{1/3}t}+\av t^2
                  \sqrt{g(j+\ell,\alpha+\xi)+g(j,\alpha)}\big].
\label{eq:pairJPfirstdifference}
\end{align}
By \eqref{eq:endpointsym}, \eqref{eq:elementarysymbols},
\eqref{eq:meanA}, and \eqref{eq:lowA},
\begin{align*}
 |a(k,\eta;p,\gamma)|
 &=|a(j+\ell,\alpha+\xi;-p,-\gamma)|\le C\Lambda_t(k,\eta)(1+\av t\kappa(p,\gamma)),\\
 \frac{|a(j+\ell,\alpha+\xi;\ell,\xi)|}{|\ell|}
 &\le C\Lambda_t(\ell,\xi)\Lambda_t(j,\alpha)\kappa(\ell,\xi)\le C\Lambda_t(\ell,\xi)\Lambda_t(p,\gamma)\kappa(\ell,\xi).
\end{align*}
Multiplying these bounds and \eqref{eq:pairJPfirstdifference}, and then using \eqref{eq:tkappa} and \eqref{eq:gprimitive}, yields
\begin{align}
&\frac{|h(j+\ell,\alpha+\xi)-h(j,\alpha)|
 |a(k,\eta;p,\gamma)a(j+\ell,\alpha+\xi;\ell,\xi)|}
 {|\ell|\Lambda_t^N(j,\alpha)\Lambda_t^N(p,\gamma)}\nonumber\\
\le{} &C_N\Lambda_t^2(k,\eta)\Lambda_t^3(\ell,\xi)
 \big[\av{\rho^{1/3}t}+\av t\sqrt{g(j+\ell,\alpha+\xi)+g(j,\alpha)}\big]
 \big[1+\av t\sqrt{g(p,\gamma)}\big]\nonumber\\
\le{} &C_N\Lambda_t^2(k,\eta)\Lambda_t^3(\ell,\xi)
 \big[\av{\rho^{1/3}t}^2+\av t^2
       (g(j+\ell,\alpha+\xi)+g(j,\alpha)+g(p,\gamma))\big].
\label{eq:pairJPbeforetransfer}
\end{align}
Changing the summation index $n$ to $-n$ in \eqref{eq:g} gives
$g(-p,-\gamma)=g(p,\gamma)$.
Since $j=k-\ell-p$, $\alpha=\eta-\xi-\gamma$, $j+\ell=k-p$, and $\alpha+\xi=\eta-\gamma$, applying \eqref{eq:doubletransfer} yields
\[
 g(j+\ell,\alpha+\xi)+g(j,\alpha)+g(p,\gamma)
 \le C\Lambda_t^5(k,\eta)\Lambda_t^5(\ell,\xi)
                    \sqrt{g(j,\alpha)g(p,\gamma)}.
\]
Substitution in \eqref{eq:pairJPbeforetransfer} gives
\begin{align}
&\frac{|h(j+\ell,\alpha+\xi)-h(j,\alpha)|
 |a(k,\eta;p,\gamma)a(j+\ell,\alpha+\xi;\ell,\xi)|}
 {|\ell|\Lambda_t^N(j,\alpha)\Lambda_t^N(p,\gamma)}\nonumber\\
\le{} &C_N\Lambda_t^2(k,\eta)\Lambda_t^3(\ell,\xi)\av{\rho^{1/3}t}^2
 +C_N\Lambda_t^7(k,\eta)\Lambda_t^8(\ell,\xi)
             \av t^2\sqrt{g(j,\alpha)g(p,\gamma)}.
\label{eq:pairJPpolynomial}
\end{align}
For the first product, $\ell\ne0$ and \eqref{eq:A} give
\[
 \frac{|a(k,\eta;\ell,\xi)|}{|\ell|}
 \le \frac{|\ell\eta-k\xi|}{|\ell|\sqrt{\ell^2+\xi^2}}
 \le C\Lambda_t(k,\eta)\Lambda_t(\ell,\xi).
\]
For the second kernel, the identity $p\alpha-j\gamma=p(\alpha+\gamma)-(j+p)\gamma$
and \eqref{eq:elementarysymbols} give, for $p\ne0$,
\[
 |a(j+p,\alpha+\gamma;p,\gamma)|
 \le C\Lambda_t(j+p,\alpha+\gamma)
                      (1+\av t\kappa(p,\gamma)).
\]
Recalling $k=j+p+\ell$ and $\eta=\alpha+\gamma+\xi$, we obtain
\[
|a(j+p,\alpha+\gamma;p,\gamma)|
 \le C\Lambda_t(k-\ell,\eta-\xi)
                      (1+\av t\kappa(p,\gamma))\le C\Lambda_t(k,\eta)\Lambda_t(\ell,\xi)(1+\av t\kappa(p,\gamma)).
\]
At $p=0$, the same conclusion follows from \eqref{eq:meanA}, since $|j|=|k-\ell|\le\Lambda_t(k-\ell,\eta-\xi)$.
Combining them with region 6 of \eqref{eq:otherpairenergyweights}
gives
\begin{align}
\frac{|h(k,\eta)-h(j+p,\alpha+\gamma)|
 |a(k,\eta;\ell,\xi)a(j+p,\alpha+\gamma;p,\gamma)|}
 {|\ell|\Lambda_t^N(j,\alpha)\Lambda_t^N(p,\gamma)}
\le C_N\Lambda_t^2(k,\eta)\Lambda_t^2(\ell,\xi)
                         (1+\av t\kappa(p,\gamma)).
\label{eq:pairJPfirstproduct}
\end{align}
Adding \eqref{eq:pairJPpolynomial} and \eqref{eq:pairJPfirstproduct}, dividing by $\Lambda_t^N(k,\eta)\Lambda_t^N(\ell,\xi)$, and using \eqref{eq:gprimitive} and \eqref{eq:quarticexpdecay}, the left-hand side of \eqref{eq:otherpairstarget} in region 6 is at most
\begin{align*}
C_Ne^{-\rho^{1/3}t/16}\Lambda_t^{9-N}(k,\eta)\Lambda_t^{9-N}(\ell,\xi)
\big[\av{\rho^{1/3}t}^2+\av t\sqrt{g(p,\gamma)}
                           +\av t^2\sqrt{g(j,\alpha)g(p,\gamma)}\big].
\end{align*}
This is region 6 of \eqref{eq:otherpairstarget}.
\end{proof}

\begin{lemma}\label{lem:highpartition}
Let $k=j+\ell+p$, $\eta=\alpha+\xi+\gamma$, and $\ell\ne0$.
If all six inequalities in \eqref{eq:pairpartition} fail, then  the weights and kernels defined in \eqref{eq:weight} and \eqref{eq:kernelgeometric} satisfy
\begin{align}
&\frac{1}{|\ell|e^{\rho^{1/3}t(\mathbf1_{k\ne0}+\mathbf1_{j\ne0}
+\mathbf1_{\ell\ne0}+\mathbf1_{p\ne0})/16}
\Lambda_t^N(k,\eta)\Lambda_t^N(j,\alpha)
\Lambda_t^N(\ell,\xi)\Lambda_t^N(p,\gamma)}\nonumber\\
&\times\Bigl|
 e^{\rho^{1/3}t\mathbf1_{k\ne0}\mathbf1_{j+p\ne0}/8}
 (h(k,\eta)-h(j+p,\alpha+\gamma))
 a(k,\eta;\ell,\xi)a(j+p,\alpha+\gamma;p,\gamma)\nonumber\\
&\quad-
 e^{\rho^{1/3}t\mathbf1_{j+\ell\ne0}\mathbf1_{j\ne0}/8}
 (h(j+\ell,\alpha+\xi)-h(j,\alpha))
 a(k,\eta;p,\gamma)a(j+\ell,\alpha+\xi;\ell,\xi)\Bigr|\nonumber\\
\le{} &C_Ne^{-\rho^{1/3}t/16}\Bigl[
\Lambda_t^{-2}(k,\eta)\Lambda_t^{-2}(j,\alpha)
 +\Lambda_t^{-2}(k,\eta)\Lambda_t^{-2}(\ell,\xi)
 +\Lambda_t^{-2}(k,\eta)\Lambda_t^{-2}(p,\gamma)\nonumber\\
&\qquad+\Lambda_t^{-2}(j,\alpha)\Lambda_t^{-2}(\ell,\xi)
 +\Lambda_t^{-2}(j,\alpha)\Lambda_t^{-2}(p,\gamma)
 +\Lambda_t^{-2}(\ell,\xi)\Lambda_t^{-2}(p,\gamma)\Bigr].
\label{eq:threehightarget}
\end{align}
\end{lemma}
\begin{proof}
For $p\ne0$, \eqref{eq:A}, \eqref{eq:endpointsym} and $|\ell|,|p|\ge1$ give
\begin{equation}
\begin{split}
|a(k,\eta;\ell,\xi)|
 &\le\Lambda_t(\ell,\xi)\Lambda_t(j+p,\alpha+\gamma),\\
|a(j+p,\alpha+\gamma;p,\gamma)|
 &\le\Lambda_t(p,\gamma)\Lambda_t(j,\alpha),\\
|a(k,\eta;p,\gamma)|
 &\le\Lambda_t(p,\gamma)\Lambda_t(j+\ell,\alpha+\xi),\\
|a(j+\ell,\alpha+\xi;\ell,\xi)|
 &\le\Lambda_t(\ell,\xi)\Lambda_t(j,\alpha).
\label{eq:roughpairkernels}
\end{split}
\end{equation}
For $p=0$, \eqref{eq:meanA} replaces the second and third bounds by $|j|$ and $|k|$, respectively. The other two are unchanged.
Since \(\Lambda_t(j+p,\alpha+\gamma)\le\Lambda_t(j,\alpha)+\Lambda_t(p,\gamma)\) and \(\Lambda_t(j+\ell,\alpha+\xi)\le\Lambda_t(j,\alpha)+\Lambda_t(\ell,\xi)\),
each product in \eqref{eq:roughpairkernels} is bounded by
\begin{equation}
 C\max\{\Lambda_t(k,\eta),\Lambda_t(j,\alpha),
                   \Lambda_t(\ell,\xi),\Lambda_t(p,\gamma)\}^{4}.
\label{eq:roughpairproduct}
\end{equation}
Also, $h=\theta\Lambda_t^{2N}$ in \eqref{eq:weight} and
$0<\theta\le1$ in \eqref{eq:thetabounds} imply
\begin{align}
&h(k,\eta)+h(j+p,\alpha+\gamma)
       +h(j+\ell,\alpha+\xi)+h(j,\alpha)\nonumber\\
\le{} &C_N\max\{\Lambda_t(k,\eta),\Lambda_t(j,\alpha),
                   \Lambda_t(\ell,\xi),\Lambda_t(p,\gamma)\}^{2N}.
\label{eq:roughpairweights}
\end{align}
Now consider the frequencies satisfying
\[
 \sqrt{k^2+(\eta+kt)^2}\ge\sqrt{j^2+(\alpha+jt)^2}
 \ge\sqrt{\ell^2+(\xi+\ell t)^2}\ge\sqrt{p^2+(\gamma+pt)^2}.
\]
The relation $k=j+\ell+p$, $\eta=\alpha+\xi+\gamma$ gives
\[
 \sqrt{k^2+(\eta+kt)^2}
 \le\sqrt{j^2+(\alpha+jt)^2}
    +\sqrt{\ell^2+(\xi+\ell t)^2}
    +\sqrt{p^2+(\gamma+pt)^2}
 \le3\sqrt{j^2+(\alpha+jt)^2}.
\]
Since the first inequality of \eqref{eq:pairpartition} fails,
\begin{align*}
 \sqrt{j^2+(\alpha+jt)^2}
 <8\big[\sqrt{\ell^2+(\xi+\ell t)^2}
                  +\sqrt{p^2+(\gamma+pt)^2}\big]
 \le16\sqrt{\ell^2+(\xi+\ell t)^2}.
\end{align*}
Consequently,
\begin{equation}
 \sqrt{\ell^2+(\xi+\ell t)^2}
 >\frac1{16}\sqrt{j^2+(\alpha+jt)^2}
 \ge\frac1{48}\sqrt{k^2+(\eta+kt)^2}.
\label{eq:threehighsizes}
\end{equation}
Thus $\Lambda_t(k,\eta)\le 3\Lambda_t(j,\alpha)\le48\Lambda_t(\ell,\xi)$.
Combining \eqref{eq:roughpairproduct}, \eqref{eq:roughpairweights}, and $|\ell|\ge1$, we deduce that 
\begin{equation}
\begin{aligned}
 &\frac{|h(k,\eta)-h(j+p,\alpha+\gamma)|
       |a(k,\eta;\ell,\xi)a(j+p,\alpha+\gamma;p,\gamma)|}
 {|\ell|\Lambda_t^N(k,\eta)\Lambda_t^N(j,\alpha)\Lambda_t^N(\ell,\xi)\Lambda_t^N(p,\gamma)}\\
 &+\frac{|h(j+\ell,\alpha+\xi)-h(j,\alpha)|
       |a(k,\eta;p,\gamma)a(j+\ell,\alpha+\xi;\ell,\xi)|}
 {|\ell|\Lambda_t^N(k,\eta)\Lambda_t^N(j,\alpha)\Lambda_t^N(\ell,\xi)\Lambda_t^N(p,\gamma)}\\
 \le{} &C_N\Lambda_t^{4-N}(\ell,\xi)\Lambda_t^{-N}(p,\gamma)
 \le C_N\Lambda_t^{-2}(\ell,\xi)\Lambda_t^{-2}(p,\gamma).
\end{aligned}
\label{eq:threehigh}
\end{equation}
The last inequality uses $N\ge6$ and $\Lambda_t\ge1$.

For any other ordering, the frequency relation expresses the largest vector as a signed sum of the other three. Its length is therefore at most three times the second largest length.
The inequality in \eqref{eq:pairpartition} selecting those two largest lengths fails in region 7. Hence the second largest length is less than eight times the sum of the remaining two, and therefore less than sixteen times the third largest.
These are precisely the two comparisons used in \eqref{eq:threehighsizes}--\eqref{eq:threehigh}.
Then by \eqref{eq:roughpairproduct} and \eqref{eq:roughpairweights}, \eqref{eq:threehigh} gives the product of the inverse squares at the two smaller frequencies in each case. The six possible products are exactly the six terms on the right-hand side of \eqref{eq:threehightarget}.
Applying \eqref{eq:quarticexpdecay} to the two terms separately and adding proves \eqref{eq:threehightarget}.
\end{proof}

\begin{proof}[Proof of Proposition~\ref{prop:quartic}: estimate of $I_5+I_7$]
We first estimate $I_5+I_7$ in \eqref{eq:pair}.
All time norms are taken over $(0,T_1)$. By \eqref{eq:weightedamplitudes},
\eqref{eq:controlenergy}, and \eqref{eq:minimumcontrol},
\[
\begin{aligned}
 &\sum_{s=\pm1}\Bigl(
 \|F^s\|_{L^\infty_tL^2}^2
 +\|\sqrt gF^s\|_{L^2_tL^2}^2
 +\|\sqrt\rho|D|F^s\|_{L^2_tL^2}^2\Bigr)
 \le C\mathcal E(T_1).
\end{aligned}
\]
We shall use the following inequalities.
Since $N-9>1$, we have
\begin{equation}
 \|\Lambda_t^{9-N}(k,\eta)F^s(k,\eta)\|_1
 \le C_N\|F^s\|_2.
 \label{eq:actualsobolevL1}
\end{equation}
The same inequality holds with $F^s$ replaced by $\sqrt gF^s$ or $\sqrt\rho|D|F^s$.
For every fixed $b\ge0$, the substitution $\tau=\rho^{1/3}t$ also gives
\begin{equation}
\begin{aligned}
 \|e^{-\rho^{1/3}t/16}\av{\rho^{1/3}t}^{b}\|_{L^1_t}
 &\le C_b\rho^{-1/3},\\
 \|e^{-\rho^{1/3}t/16}\av{\rho^{1/3}t}^{b}\|_{L^2_t}
 &\le C_b\rho^{-1/6},\\
 \|\av t e^{-\rho^{1/3}t/16}\av{\rho^{1/3}t}^{b}\|_{L^2_t}
 &\le C_b\rho^{-1/2},\\
 \|\av t^2e^{-\rho^{1/3}t/16}\av{\rho^{1/3}t}^{b}\|_{L^\infty_t}
 &\le C_b\rho^{-2/3}.
\end{aligned}
\label{eq:quartictimebounds}
\end{equation}
By $\av t\le\rho^{-1/3}\av{\rho^{1/3}t}$ and $\tau=\rho^{1/3}t$, it also holds
\begin{align}
 &\bigl\|e^{-\rho^{1/3}t/16}
   [\av t^2+\rho^{-1/3}\av t\av{\rho^{1/3}t}^4]\bigr\|_{L^\infty_t}
       \le C\rho^{-2/3},\nonumber\\
 &\bigl\|e^{-\rho^{1/3}t/16}
   [\av{\rho^{1/3}t}^2+\av t+\rho^{-1/3}\av{\rho^{1/3}t}^4]\bigr\|_{L^1_t}
       \le C\rho^{-2/3}.
\label{eq:quarticscalarbounds}
\end{align}
Set $k=j+\ell+p$ and $\eta=\alpha+\xi+\gamma$.
Lemmas~\ref{lem:lowpair}, \ref{lem:meanpair}, \ref{lem:quarticinterface}, \ref{lem:otherpairs}, and \ref{lem:highpartition} give pointwise estimates for 
\begin{align*}
 e^{\rho^{1/3}t\mathbf1_{k\ne0}\mathbf1_{j+p\ne0}/8}
 (h(k,\eta)-h(j+p,\alpha+\gamma))
 a(k,\eta;\ell,\xi)a(j+p,\alpha+\gamma;p,\gamma)\\
 -e^{\rho^{1/3}t\mathbf1_{j+\ell\ne0}\mathbf1_{j\ne0}/8}
 (h(j+\ell,\alpha+\xi)-h(j,\alpha))
 a(k,\eta;p,\gamma)a(j+\ell,\alpha+\xi;\ell,\xi).
\end{align*}
Substituting \eqref{eq:weightedamplitudes} and applying Young's inequality, we obtain
\[
 \int_0^{T_1}|I_5(t)+I_7(t)|\dd t\le\sum_{r=1}^7K_r,
\]
where $K_1,\ldots,K_7$ are the upper bounds defined below. We estimate these terms separately.

\noindent\emph{Case 1. Region 1.}
We distinguish between $p=0$ and $p\ne0$. By \eqref{eq:fexplicit}, \eqref{eq:weightedamplitudes}, Lemmas~\ref{lem:lowpair}, \ref{lem:meanpair}, \ref{lem:quarticinterface}, and Young's inequality, we obtain
\begin{align*}
K_1 ={} &\frac{C_N}{|\beta|}\sum_s\int_0^{T_1} e^{-\rho^{1/3}t/16}
 \Bigl[
 \big[\av t^2+\rho^{-1/3}\av t\av{\rho^{1/3}t}^4\big]
       \|F^s\|_2^2\|\Lambda_t^{9-N}\sqrt gF^{-s}_{\ne}\|_1^2\\
 &+\av t^2\|\sqrt gF^s\|_2^2
                          \|\Lambda_t^{9-N}F^{-s}\|_1^2
 +\big[\av{\rho^{1/3}t}^2+\av t+\rho^{-1/3}\av{\rho^{1/3}t}^4\big]
       \|F^s\|_2^2\|\Lambda_t^{9-N}F^{-s}_{\ne}\|_1
                    \|\Lambda_t^{9-N}F^{-s}_0\|_1\\
 &+\|F^s\|_2\||D|F^s\|_2
       \|\Lambda_t^{9-N}\sqrt gF^{-s}_{\ne}\|_1
       \|\Lambda_t^{9-N}F^{-s}_{\ne}\|_1\Bigr]\dd t.
\end{align*}
Thus $\rho\le1$, \eqref{eq:actualsobolevL1}, \eqref{eq:quartictimebounds}, and \eqref{eq:quarticscalarbounds} give
\begin{align*}
K_1 \le{} &\frac{C_N}{|\beta|}\sum_s
 \Bigl[\rho^{-2/3}
   \bigl[
   \|F^s\|_{L^\infty_tL^2}^2\|\sqrt gF^{-s}_{\ne}\|_{L^2_tL^2}^2
   +\|\sqrt gF^s\|_{L^2_tL^2}^2\|F^{-s}\|_{L^\infty_tL^2}^2\nonumber\\*
 &+
   \|F^s\|_{L^\infty_tL^2}^2
   \|F^{-s}_{\ne}\|_{L^\infty_tL^2}
   \|F^{-s}_0\|_{L^\infty_tL^2}\bigr]\nonumber\\
 &+\rho^{-1/2}
   \|F^s\|_{L^\infty_tL^2}\|F^{-s}_{\ne}\|_{L^\infty_tL^2}
   \|\sqrt\rho|D|F^s\|_{L^2_tL^2}
   \|\sqrt gF^{-s}_{\ne}\|_{L^2_tL^2}\Bigr]\nonumber\\
 \le{} &\frac{C_N}{|\beta|}\rho^{-2/3}\mathcal E(T_1)^2.
\end{align*}

\noindent\emph{Case 2. Region 2.}
Applying the bound for region 2 in \eqref{eq:otherpairstarget}, with the factors at $(k,\eta)$ and $(j,\alpha)$ in $L^1$ and the other two in $L^2$, we obtain
\begin{align*}
K_2 ={} &\frac{C_N}{|\beta|}\sum_s\int_0^{T_1} e^{-\rho^{1/3}t/16}
 \|\Lambda_t^{9-N}F^s\|_1^2
 \bigl[\|F^{-s}\|_2^2+\av t\|F^{-s}\|_2\|\sqrt gF^{-s}\|_2
                    +\av t^2\|\sqrt gF^{-s}\|_2^2\bigr]\dd t\nonumber\\
 \le{} &\frac{C_N}{|\beta|}\sum_s\|F^s\|_{L^\infty_tL^2}^2
 \Bigl[
 \|e^{-\rho^{1/3}t/16}\|_{L^1_t}\|F^{-s}\|_{L^\infty_tL^2}^2\nonumber\\*
 &+\|\av t e^{-\rho^{1/3}t/16}\|_{L^2_t}
         \|F^{-s}\|_{L^\infty_tL^2}\|\sqrt gF^{-s}\|_{L^2_tL^2}
 +\|\av t^2e^{-\rho^{1/3}t/16}\|_{L^\infty_t}
         \|\sqrt gF^{-s}\|_{L^2_tL^2}^2\Bigr]\nonumber\\
 \le{} &\frac{C_N}{|\beta|}
       (\rho^{-1/3}+\rho^{-1/2}+\rho^{-2/3})\mathcal E(T_1)^2.
\end{align*}

\noindent\emph{Case 3. Region 3.}
Applying the bound for region 3 in \eqref{eq:otherpairstarget}, with the factors at $(j,\alpha)$ and $(p,\gamma)$ in $L^1$ and the other two in $L^2$, we obtain
\begin{align*}
K_3 ={} &\frac{C_N}{|\beta|}\sum_s\int_0^{T_1} e^{-\rho^{1/3}t/16}
 \|\Lambda_t^{9-N}F^s\|_1\|F^s\|_2
 \Bigl[\|\Lambda_t^{9-N}F^{-s}\|_1
       \big[\|F^{-s}\|_2+\av t\|\sqrt gF^{-s}\|_2\big]\nonumber\\
 &+\|\Lambda_t^{9-N}\sqrt gF^{-s}\|_1
       \big[\||D|F^{-s}\|_2+\av t^2\|\sqrt gF^{-s}\|_2\big]\Bigr]\dd t\nonumber\\
 \le{} &\frac{C_N}{|\beta|}\sum_s\|F^s\|_{L^\infty_tL^2}^2
 \Bigl[
 \|e^{-\rho^{1/3}t/16}\|_{L^1_t}\|F^{-s}\|_{L^\infty_tL^2}^2
 +\|\av t e^{-\rho^{1/3}t/16}\|_{L^2_t}
       \|F^{-s}\|_{L^\infty_tL^2}\|\sqrt gF^{-s}\|_{L^2_tL^2}\nonumber\\
 &+\rho^{-1/2}
       \|\sqrt\rho|D|F^{-s}\|_{L^2_tL^2}
       \|\sqrt gF^{-s}\|_{L^2_tL^2}
 +\|\av t^2e^{-\rho^{1/3}t/16}\|_{L^\infty_t}
       \|\sqrt gF^{-s}\|_{L^2_tL^2}^2\Bigr]\nonumber\\
 \le{} &\frac{C_N}{|\beta|}
       (\rho^{-1/3}+\rho^{-1/2}+\rho^{-2/3})\mathcal E(T_1)^2.
\end{align*}

\noindent\emph{Case 4. Region 4.}
We apply the bound for region 4 in \eqref{eq:otherpairstarget}, with the factors at $(k,\eta)$ and $(p,\gamma)$ in $L^1$ and the other two in $L^2$. Arguing as in Case 3, we obtain
\begin{align*}
K_4 ={} &\frac{C_N}{|\beta|}\sum_s\int_0^{T_1} e^{-\rho^{1/3}t/16}
 \|\Lambda_t^{9-N}F^s\|_1\|F^s\|_2
 \Bigl[\|\Lambda_t^{9-N}F^{-s}\|_1
       \big[\|F^{-s}\|_2+\av t\|\sqrt gF^{-s}\|_2\big]\nonumber\\
 &+\|\Lambda_t^{9-N}\sqrt gF^{-s}\|_1
      \big[\||D|F^{-s}\|_2
                              +\av t^2\|\sqrt gF^{-s}\|_2\big]\Bigr]\dd t\nonumber\\
 \le{} &\frac{C_N}{|\beta|}
       (\rho^{-1/3}+\rho^{-1/2}+\rho^{-2/3})\mathcal E(T_1)^2.
\end{align*}

\noindent\emph{Case 5. Region 5.}
Applying the bound for region 5 in \eqref{eq:otherpairstarget}, with the factors at $(j,\alpha)$ and $(\ell,\xi)$ in $L^1$ and the other two in $L^2$, we obtain
\begin{align*}
K_5 ={} &\frac{C_N}{|\beta|}\sum_s\int_0^{T_1} e^{-\rho^{1/3}t/16}
 \|\Lambda_t^{9-N}F^s\|_1\|\Lambda_t^{9-N}F^{-s}\|_1\nonumber\\
 &\qquad\times\Bigl[
 \av{\rho^{1/3}t}^2\|F^s\|_2\|F^{-s}\|_2
 +\av t\|F^s\|_2\|\sqrt gF^{-s}\|_2
 +\av t^2\|\sqrt gF^s\|_2\|\sqrt gF^{-s}\|_2\Bigr]\dd t\nonumber\\
 \le{} &\frac{C_N}{|\beta|}\sum_s
 \|F^s\|_{L^\infty_tL^2}\|F^{-s}\|_{L^\infty_tL^2}
 \Bigl[
 \|e^{-\rho^{1/3}t/16}\av{\rho^{1/3}t}^2\|_{L^1_t}
       \|F^s\|_{L^\infty_tL^2}\|F^{-s}\|_{L^\infty_tL^2}\nonumber\\
 &+\|\av t e^{-\rho^{1/3}t/16}\|_{L^2_t}
       \|F^s\|_{L^\infty_tL^2}\|\sqrt gF^{-s}\|_{L^2_tL^2}
 +\|\av t^2e^{-\rho^{1/3}t/16}\|_{L^\infty_t}
       \|\sqrt gF^s\|_{L^2_tL^2}\|\sqrt gF^{-s}\|_{L^2_tL^2}
 \Bigr]\nonumber\\
 \le{} &\frac{C_N}{|\beta|}
       (\rho^{-1/3}+\rho^{-1/2}+\rho^{-2/3})\mathcal E(T_1)^2.
\end{align*}

\noindent\emph{Case 6. Region 6.}
We apply the bound for region 6 in \eqref{eq:otherpairstarget}, with the factors at $(k,\eta)$ and $(\ell,\xi)$ in $L^1$ and the other two in $L^2$. Arguing as in Case 5, we obtain
\begin{align*}
K_6 ={} &\frac{C_N}{|\beta|}\sum_s\int_0^{T_1} e^{-\rho^{1/3}t/16}
 \|\Lambda_t^{9-N}F^s\|_1\|\Lambda_t^{9-N}F^{-s}\|_1\nonumber\\*
 &\qquad\times\Bigl[
 \av{\rho^{1/3}t}^2\|F^s\|_2\|F^{-s}\|_2
 +\av t\|F^s\|_2\|\sqrt gF^{-s}\|_2
 +\av t^2\|\sqrt gF^s\|_2\|\sqrt gF^{-s}\|_2\Bigr]\dd t\nonumber\\
 \le{} &\frac{C_N}{|\beta|}
       (\rho^{-1/3}+\rho^{-1/2}+\rho^{-2/3})\mathcal E(T_1)^2.
\end{align*}

\noindent\emph{Case 7. Region 7.}
By Lemma~\ref{lem:highpartition} and Young's inequality, we obtain
\begin{align*}
K_7 ={} &\frac{C_N}{|\beta|}\sum_s\int_0^{T_1} e^{-\rho^{1/3}t/16}
 \Bigl[
 \|F^s\|_2^2\|\Lambda_t^{-2}F^{-s}\|_1^2
 +\|\Lambda_t^{-2}F^s\|_1^2\|F^{-s}\|_2^2\\*
 &\qquad+\|\Lambda_t^{-2}F^s\|_1\|F^s\|_2
          \|\Lambda_t^{-2}F^{-s}\|_1\|F^{-s}\|_2
 \Bigr]\dd t\\
 \le{} &\frac{C_N}{|\beta|}
 \|e^{-\rho^{1/3}t/16}\|_{L^1_t}
 \sum_s\|F^s\|_{L^\infty_tL^2}^2\|F^{-s}\|_{L^\infty_tL^2}^2\\
 \le{} &\frac{C_N}{|\beta|}\rho^{-1/3}\mathcal E(T_1)^2.
\end{align*}
Adding the above estimates and using $0<\rho\le1$, we obtain
\begin{equation}
 \int_0^{T_1}|I_5(t)+I_7(t)|\dd t
 \le\sum_{r=1}^7K_r
 \le\frac{C_N}{|\beta|}\rho^{-2/3}\mathcal E(T_1)^2,
 \label{eq:pairedquarticintegral}
\end{equation}
which proves the desired conclusion.
\end{proof}

\subsection{Estimates for $I_6$}
\label{sec:quarticsecond}

For \eqref{eq:adv}, $k=p+q+m$, $\eta=\gamma+\delta+\zeta$,
and $p+q\ne0$. Thus
\[
 (k,\eta+kt)=(p,\gamma+pt)+(q,\delta+qt)+(m,\zeta+mt).
\]
We again divide the admissible frequencies into \emph{seven
disjoint regions}. Regions 1--6 are defined, in order, by
\begin{equation}
\begin{gathered}
8\big(\sqrt{p^2+(\gamma+pt)^2}+\sqrt{q^2+(\delta+qt)^2}\bigr)\le\min\big\{\sqrt{k^2+(\eta+kt)^2},\sqrt{m^2+(\zeta+mt)^2}\big\},\\
8\bigl(\sqrt{k^2+(\eta+kt)^2}+\sqrt{m^2+(\zeta+mt)^2}\bigr)\le\min\bigl\{\sqrt{p^2+(\gamma+pt)^2},\sqrt{q^2+(\delta+qt)^2}\bigr\},\\
8\bigl(\sqrt{m^2+(\zeta+mt)^2}+\sqrt{q^2+(\delta+qt)^2}\bigr)\le\min\bigl\{\sqrt{k^2+(\eta+kt)^2},\sqrt{p^2+(\gamma+pt)^2}\bigr\},\\
8\bigl(\sqrt{k^2+(\eta+kt)^2}+\sqrt{q^2+(\delta+qt)^2}\bigr)\le\min\bigl\{\sqrt{m^2+(\zeta+mt)^2},\sqrt{p^2+(\gamma+pt)^2}\bigr\},\\
8\bigl(\sqrt{m^2+(\zeta+mt)^2}+\sqrt{p^2+(\gamma+pt)^2}\bigr)\le\min\bigl\{\sqrt{k^2+(\eta+kt)^2},\sqrt{q^2+(\delta+qt)^2}\bigr\},\\
8\bigl(\sqrt{k^2+(\eta+kt)^2}+\sqrt{p^2+(\gamma+pt)^2}\bigr)\le\min\bigl\{\sqrt{m^2+(\zeta+mt)^2},\sqrt{q^2+(\delta+qt)^2}\bigr\}.
\end{gathered}
\label{eq:advpartition}
\end{equation}
Region 7 is the complement of their union. It consists exactly of the frequencies for which all six inequalities in \eqref{eq:advpartition} fail. Since $p+q\ne0$, at least one of $p,q$ is nonzero, and the frequency relation implies that at least two of the four lengths are positive. The disjointness argument in Subsection~\ref{sec:quarticfirst} applies without change. 
The exponential weights satisfy
\begin{equation}
 \frac{e^{\rho^{1/3}t\mathbf1_{k\ne0}\mathbf1_{m\ne0}/8}}
 {e^{\rho^{1/3}t(\mathbf1_{k\ne0}+\mathbf1_{p\ne0}
                    +\mathbf1_{q\ne0}+\mathbf1_{m\ne0})/16}}
 \le e^{-\rho^{1/3}t/16}.
\label{eq:advectorquarticexpdecay}
\end{equation}
If $k,m\ne0$, then $p+q\ne0$ provides a third nonzero horizontal frequency. Otherwise, the relations $p+q\ne0$ and $k=p+q+m$ still give at least two nonzero horizontal frequencies.
The following two lemmas give the kernel estimates on Regions 1--6 and Region 7, respectively.

\begin{lemma}\label{lem:advectorpairs}
Let $k=p+q+m$, $\eta=\gamma+\delta+\zeta$, and $p+q\ne0$.
On the six regions in \eqref{eq:advpartition}, the weights and kernels defined in \eqref{eq:weight} and \eqref{eq:kernelgeometric} satisfy
\begin{align}
&\frac{e^{\rho^{1/3}t\mathbf1_{k\ne0}\mathbf1_{m\ne0}/8}}{e^{\rho^{1/3}t(\mathbf1_{k\ne0}+\mathbf1_{p\ne0}
+\mathbf1_{q\ne0}+\mathbf1_{m\ne0})/16}}\frac{|h(k,\eta)-h(m,\zeta)|}
 {|p+q|
\Lambda_t^N(k,\eta)\Lambda_t^N(p,\gamma)
\Lambda_t^N(q,\delta)\Lambda_t^N(m,\zeta)}\nonumber\\
&\times
 |a(k,\eta;p+q,\gamma+\delta)a(p+q,\gamma+\delta;p,\gamma)|\nonumber\\
\le{} &C_Ne^{-\rho^{1/3}t/16}
\begin{cases}
\begin{aligned}
 &\Lambda_t^{9-N}(p,\gamma)\Lambda_t^{9-N}(q,\delta)
 \big[\av{\rho^{1/3}t}
                  +\av t(\sqrt{g(k,\eta)}+\sqrt{g(m,\zeta)})\big],
\end{aligned}&\text{region }1,\\
 \Lambda_t^{9-N}(k,\eta)\Lambda_t^{9-N}(m,\zeta)
                         \big[1+\av t\sqrt{g(p,\gamma)}\big],
 &\text{region }2,\\
\begin{aligned}
 &\Lambda_t^{9-N}(m,\zeta)\Lambda_t^{9-N}(q,\delta)\\
 &\quad\times\big[1+\av t(\sqrt{g(k,\eta)}+\sqrt{g(p,\gamma)})
                         +\av t^2\sqrt{g(k,\eta)g(p,\gamma)}\big],
\end{aligned}&\text{region }3,\\
\begin{aligned}
 &\Lambda_t^{9-N}(k,\eta)\Lambda_t^{9-N}(q,\delta)\\
 &\quad\times\big[1+\av t(\sqrt{g(m,\zeta)}+\sqrt{g(p,\gamma)})
                         +\av t^2\sqrt{g(m,\zeta)g(p,\gamma)}\big],
\end{aligned}&\text{region }4,\\
\begin{aligned}
 &\Lambda_t^{9-N}(m,\zeta)\Lambda_t^{9-N}(p,\gamma)\\
 &\quad\times\big[1+\av t\sqrt{g(k,\eta)}
      +\sqrt{q^2+\delta^2}\sqrt{g(p,\gamma)}
                         +\av t^2\sqrt{g(k,\eta)g(p,\gamma)}\big],
\end{aligned}&\text{region }5,\\
\begin{aligned}
 &\Lambda_t^{9-N}(k,\eta)\Lambda_t^{9-N}(p,\gamma)\\
 &\quad\times\big[1+\av t\sqrt{g(m,\zeta)}
      +\sqrt{q^2+\delta^2}\sqrt{g(p,\gamma)}
                         +\av t^2\sqrt{g(m,\zeta)g(p,\gamma)}\big],
\end{aligned}&\text{region }6.
\end{cases}
\label{eq:advectorpairstarget}
\end{align}
The region numbers refer to \eqref{eq:advpartition}.
\end{lemma}

\begin{proof}
Throughout the proof, $\ell=p+q\ne0$ and $\xi=\gamma+\delta$.
By \eqref{eq:endpointsym} and \eqref{eq:elementarysymbols},
\begin{align}
 |a(k,\eta;\ell,\xi)|
 &\le C\Lambda_t(m,\zeta)(1+\av t\kappa(\ell,\xi)),\nonumber\\
 |a(k,\eta;\ell,\xi)|
 &=|a(m,\zeta;-\ell,-\xi)|
 \le C\Lambda_t(k,\eta)(1+\av t\kappa(\ell,\xi)).
\label{eq:advouterelementary}
\end{align}
Applying the same formulas gives
\begin{align}
 |a(\ell,\xi;p,\gamma)|
 &\le C\Lambda_t(q,\delta)(1+\av t\kappa(p,\gamma)),\nonumber\\
 |a(\ell,\xi;p,\gamma)|
 &=|a(q,\delta;-p,-\gamma)|
 \le C\Lambda_t(\ell,\xi)(1+\av t\kappa(p,\gamma)).
\label{eq:advinnerelementary}
\end{align}
For $p=0$, these inequalities follow instead from
\eqref{eq:meanA}, since $|q|=|\ell|$ and $\kappa(0,\gamma)=0$.
 \eqref{eq:A} also gives
\begin{align}
 |a(\ell,\xi;p,\gamma)|
 &\le\frac{|p\delta-q\gamma|}{\sqrt{p^2+\gamma^2}}
 \le\Lambda_t(p,\gamma)\Lambda_t(q,\delta)\kappa(p,\gamma)
 \quad(p\ne0),\nonumber\\
 |a(\ell,\xi;0,\gamma)|
 &\le|q|\le\Lambda_t(q,\delta).
\label{eq:advinnerdeterminant}
\end{align}
\noindent\emph{Case 1.}
Here $\ell=p+q$, $\xi=\gamma+\delta$. Since $\kappa\le1$ and $\Lambda_t\ge1$, \eqref{eq:advinnerdeterminant} bounds $|a(\ell,\xi;p,\gamma)|$ by $\Lambda_t(p,\gamma)\Lambda_t(q,\delta)$ for both $p\ne0$ and $p=0$.
By \eqref{eq:tkappa}, we have
\[
\av t\kappa(\ell,\xi)\le C\Lambda_t(\ell,\xi)
 \le C\Lambda_t(p,\gamma)\Lambda_t(q,\delta).
\]
The first condition in \eqref{eq:advpartition} gives
\[
 8\sqrt{\ell^2+(\xi+\ell t)^2}
 \le8\big[\sqrt{p^2+(\gamma+pt)^2}+\sqrt{q^2+(\delta+qt)^2}\big]
 \le\sqrt{m^2+(\zeta+mt)^2}.
\]
Applying \eqref{eq:FD1} at $(m,\zeta)$ and $(k,\eta)$ and using \eqref{eq:lowA}, we obtain
\begin{equation}
\begin{aligned}
&\frac{|h(k,\eta)-h(m,\zeta)|
 |a(k,\eta;\ell,\xi)a(\ell,\xi;p,\gamma)|}
 {|\ell|\Lambda_t^N(k,\eta)\Lambda_t^N(m,\zeta)}\\
\le{}& C_N\Lambda_t^{2}(\ell,\xi)\Lambda_t(p,\gamma)\Lambda_t(q,\delta)\big[\av t\av{\rho^{1/3}t}+\av t^2\sqrt{g(k,\eta)+g(m,\zeta)}\big] \kappa(\ell,\xi)\\
\le{}& C\Lambda_t^{3}(\ell,\xi)\Lambda_t(p,\gamma)\Lambda_t(q,\delta)\big[\av{\rho^{1/3}t}+\av t\sqrt{g(k,\eta)+g(m,\zeta)}\big]\\
\le{}& C\Lambda_t^{4}(p,\gamma)\Lambda_t^{4}(q,\delta)
       \big[\av{\rho^{1/3}t}+\av t\sqrt{g(k,\eta)+g(m,\zeta)}\big].
\label{eq:advKMpolynomial}
\end{aligned}
\end{equation}
Dividing \eqref{eq:advKMpolynomial} by $\Lambda_t^N(p,\gamma)\Lambda_t^N(q,\delta)$ and applying \eqref{eq:advectorquarticexpdecay}, we bound the left-hand side of \eqref{eq:advectorpairstarget} in region 1 by 
\[
 C_Ne^{-\rho^{1/3}t/16}\Lambda_t^{4-N}(p,\gamma)\Lambda_t^{4-N}(q,\delta)
 \big[\av{\rho^{1/3}t}+\av t(\sqrt{g(k,\eta)}+\sqrt{g(m,\zeta)})\big].
\]
Since $\Lambda_t^{4-N}\le\Lambda_t^{9-N}$, this proves region 1 of \eqref{eq:advectorpairstarget}.

\noindent\emph{Case 2.}
Here $\ell=k-m$, $\xi=\eta-\zeta$ and $\ell=p+q\neq0$. Equations \eqref{eq:A} and \eqref{eq:advinnerelementary} give
\[
\begin{aligned}
&|a(k,\eta;\ell,\xi)| \le C\Lambda_t(k,\eta)\Lambda_t(m,\zeta),\\
&|a(\ell,\xi;p,\gamma)|\le C\Lambda_t(\ell,\xi)(1+\av t\kappa(p,\gamma)).
\end{aligned}
\]
Notice that $\Lambda_t(\ell,\xi)\le2\Lambda_t(k,\eta)\Lambda_t(m,\zeta)$ due to $\ell=k-m$ and $\xi=\eta-\zeta$.
Multiplication of the two preceding bounds therefore gives
\begin{equation}
 |a(k,\eta;\ell,\xi)a(\ell,\xi;p,\gamma)|
 \le C\Lambda_t^2(k,\eta)\Lambda_t^2(m,\zeta)
             (1+\av t\kappa(p,\gamma)).
\label{eq:advPQproduct}
\end{equation}
The second condition of \eqref{eq:advpartition} and $k-m=p+q$ give
\[
\begin{aligned}
 \big|\sqrt{p^2+(\gamma+pt)^2}-\sqrt{q^2+(\delta+qt)^2}\big|
 &\le\sqrt{k^2+(\eta+kt)^2}+\sqrt{m^2+(\zeta+mt)^2}\\
 &\le\frac18\min\big\{\sqrt{p^2+(\gamma+pt)^2},\sqrt{q^2+(\delta+qt)^2}\big\}.
\end{aligned}
\]
Hence $\Lambda_t(p,\gamma)$ and $\Lambda_t(q,\delta)$ are comparable and dominate the other two weights.
By \eqref{eq:weight} and \eqref{eq:thetabounds},
\begin{equation}
 \frac{|h(k,\eta)-h(m,\zeta)|}
 {\Lambda_t^N(p,\gamma)\Lambda_t^N(q,\delta)}
 \le\frac{\Lambda_t^{2N}(k,\eta)+\Lambda_t^{2N}(m,\zeta)}
 {\Lambda_t^N(p,\gamma)\Lambda_t^N(q,\delta)}
 \le C_N.
\label{eq:advPQenergy}
\end{equation}
Combining \eqref{eq:advPQproduct}, \eqref{eq:advPQenergy}, $|\ell|\ge1$, \eqref{eq:gprimitive}, and \eqref{eq:advectorquarticexpdecay}, the left-hand side of \eqref{eq:advectorpairstarget} in region 2 is at most 
\[
 C_Ne^{-\rho^{1/3}t/16}\Lambda_t^{2-N}(k,\eta)\Lambda_t^{2-N}(m,\zeta)\big[1+\av t\sqrt{g(p,\gamma)}\big].
\]
This implies region 2 of \eqref{eq:advectorpairstarget}.

\noindent\emph{Case 3.}
Here $\ell=p+q$, $\xi=\gamma+\delta$. The third condition of \eqref{eq:advpartition} implies
\begin{align*}
\big|\sqrt{k^2+(\eta+kt)^2}-\sqrt{p^2+(\gamma+pt)^2}\big|
 &\le\sqrt{q^2+(\delta+qt)^2}+\sqrt{m^2+(\zeta+mt)^2}\\
 &\le\frac18\min\big\{\sqrt{k^2+(\eta+kt)^2},\sqrt{p^2+(\gamma+pt)^2}\big\}.
\end{align*}
Thus \eqref{eq:weight} and
\eqref{eq:thetabounds} give
\begin{equation}
 \frac{|h(k,\eta)-h(m,\zeta)|}
 {\Lambda_t^N(k,\eta)\Lambda_t^N(p,\gamma)}
 \le\frac{\Lambda_t^{2N}(k,\eta)+\Lambda_t^{2N}(m,\zeta)}
 {\Lambda_t^N(k,\eta)\Lambda_t^N(p,\gamma)}
 \le C_N.
\label{eq:advKPenergy}
\end{equation}
Multiplying the first inequalities of \eqref{eq:advouterelementary}
and \eqref{eq:advinnerelementary}, then using
\eqref{eq:advKPenergy} and $|\ell|\ge1$, yields
\begin{equation}
\begin{aligned}
 \frac{|h(k,\eta)-h(m,\zeta)|
 |a(k,\eta;\ell,\xi)a(\ell,\xi;p,\gamma)|}
 {|\ell|\Lambda_t^N(k,\eta)\Lambda_t^N(p,\gamma)}
 \le C_N\Lambda_t(m,\zeta)\Lambda_t(q,\delta)
 (1+\av t\kappa(\ell,\xi))(1+\av t\kappa(p,\gamma)).
\end{aligned}
\label{eq:advKP}
\end{equation}
For $p=0$, \eqref{eq:meanA} gives the same inequality with $\kappa(p,\gamma)=0$. Since $\ell=k-m$ and $\xi=\eta-\zeta$, \eqref{eq:gprimitive} and \eqref{eq:shiftkappa} give
\[
 \kappa(\ell,\xi)\le C\Lambda_t^5(m,\zeta)\sqrt{g(k,\eta)},
 \qquad
 \kappa(p,\gamma)\le\sqrt{g(p,\gamma)}.
\]
Substituting these two inequalities in \eqref{eq:advKP},
and expanding the product, gives
\begin{equation}
\begin{aligned}
&C\Lambda_t(m,\zeta)\Lambda_t(q,\delta)(1+\av t\kappa(\ell,\xi))(1+\av t\kappa(p,\gamma))\\
\le{} &C\Lambda_t^{6}(m,\zeta)\Lambda_t(q,\delta) \big[1+\av t\sqrt{g(k,\eta)}+\av t\sqrt{g(p,\gamma)}
+\av t^2\sqrt{g(k,\eta)g(p,\gamma)}\big].
\label{eq:advKPpolynomial}
\end{aligned}
\end{equation}
Dividing \eqref{eq:advKPpolynomial} by $\Lambda_t^N(m,\zeta)\Lambda_t^N(q,\delta)$ and applying \eqref{eq:advectorquarticexpdecay}, we bound the left-hand side of \eqref{eq:advectorpairstarget} in region 3 by 
\[
 C_Ne^{-\rho^{1/3}t/16}\Lambda_t^{6-N}(m,\zeta)\Lambda_t^{1-N}(q,\delta)\big[1+\av t\sqrt{g(k,\eta)}+\av t\sqrt{g(p,\gamma)}+\av t^2\sqrt{g(k,\eta)g(p,\gamma)}\big].
\]
This implies region 3 of \eqref{eq:advectorpairstarget}.

\noindent\emph{Case 4.}
Here \(\ell=k-m=p+q\), \(\xi=\eta-\zeta=\gamma+\delta\).
By the second bound of \eqref{eq:advouterelementary} and the first bound of \eqref{eq:advinnerelementary},
\begin{align*}
|a(k,\eta;\ell,\xi)a(\ell,\xi;p,\gamma)|
 &=|a(m,\zeta;-\ell,-\xi)a(\ell,\xi;p,\gamma)|\\
 &\le C\Lambda_t(k,\eta)\Lambda_t(q,\delta)
              (1+\av t\kappa(\ell,\xi))(1+\av t\kappa(p,\gamma)).
\end{align*}
The fourth condition of \eqref{eq:advpartition} and
$p+m=k-q$ imply
\begin{align*}
 \big|\sqrt{m^2+(\zeta+mt)^2}-\sqrt{p^2+(\gamma+pt)^2}\big|
 &\le\sqrt{k^2+(\eta+kt)^2}+\sqrt{q^2+(\delta+qt)^2}\\
 &\le \frac18\min\big\{\sqrt{m^2+(\zeta+mt)^2},\sqrt{p^2+(\gamma+pt)^2}\big\}.
\end{align*}
Thus $\Lambda_t(m,\zeta)$ and $\Lambda_t(p,\gamma)$ are comparable and dominate the other two weights. Applying \eqref{eq:weight}, \eqref{eq:thetabounds}, and \eqref{eq:shiftkappa} with $\ell=-m+k$, $\xi=-\zeta+\eta$ and the evenness of $g$ from \eqref{eq:g} gives
\begin{align*}
 \frac{|h(k,\eta)-h(m,\zeta)|}
 {\Lambda_t^N(m,\zeta)\Lambda_t^N(p,\gamma)}&\le C_N,\\
 \kappa(\ell,\xi)\kappa(p,\gamma)
 &\le C\Lambda_t^5(k,\eta)\sqrt{g(m,\zeta)g(p,\gamma)}.
\end{align*}
By \eqref{eq:shiftkappa} and the evenness of $g$ from \eqref{eq:g}, $\kappa(\ell,\xi)\le C\Lambda_t^5(k,\eta)\sqrt{g(m,\zeta)}$. Hence by $\ell=p+q\neq0$, we get
\begin{equation}
\begin{aligned}
&\frac{|h(k,\eta)-h(m,\zeta)|
 |a(k,\eta;\ell,\xi)a(\ell,\xi;p,\gamma)|}
 {|\ell|\Lambda_t^N(m,\zeta)\Lambda_t^N(p,\gamma)}\\
\le{} &C\Lambda_t^{6}(k,\eta)\Lambda_t(q,\delta) \big[1+\av t\sqrt{g(m,\zeta)}+\av t\sqrt{g(p,\gamma)}+\av t^2\sqrt{g(m,\zeta)g(p,\gamma)}\big].
\label{eq:advMPpolynomial}
\end{aligned}
\end{equation}
Dividing \eqref{eq:advMPpolynomial} by $\Lambda_t^N(k,\eta)\Lambda_t^N(q,\delta)$ and applying \eqref{eq:advectorquarticexpdecay}, we bound the left-hand side of \eqref{eq:advectorpairstarget} in region 4 by
\[
 C_Ne^{-\rho^{1/3}t/16}\Lambda_t^{6-N}(k,\eta)\Lambda_t^{1-N}(q,\delta)\big[1+\av t\sqrt{g(m,\zeta)}+\av t\sqrt{g(p,\gamma)}+\av t^2\sqrt{g(m,\zeta)g(p,\gamma)}\big].
\]
This implies region 4 of \eqref{eq:advectorpairstarget}.

\noindent\emph{Case 5.}
Here $\ell=p+q$, $\xi=\gamma+\delta$. For $p\ne0$, \eqref{eq:A}, the fifth condition of \eqref{eq:advpartition} and \eqref{eq:advinnerdeterminant} give
\begin{align*}
|a(k,\eta;\ell,\xi)|
 &\le C\frac{\Lambda_t(m,\zeta)\Lambda_t(\ell,\xi)}{\sqrt{\ell^2+\xi^2}}\le C\frac{\Lambda_t(m,\zeta)\Lambda_t(q,\delta)}{\sqrt{\ell^2+\xi^2}},\\
|a(\ell,\xi;p,\gamma)|
 &\le C\Lambda_t(p,\gamma)\Lambda_t(q,\delta)\kappa(p,\gamma),\\
\sqrt{\ell^2+\xi^2}
 &\le\sqrt{q^2+\delta^2}+\sqrt{p^2+\gamma^2}
 \le\sqrt{q^2+\delta^2}+C\av t\Lambda_t(p,\gamma).
\end{align*}
The fifth condition of \eqref{eq:advpartition} gives $\sqrt{p^2+(\gamma+pt)^2}\le\frac18\sqrt{q^2+(\delta+qt)^2}$.
Thus $\Lambda_t(q,\delta)\le C\Lambda_t(\ell,\xi)$.
Using $\Lambda_t(\ell,\xi)^2\le C(\ell^2+\xi^2+t^2\ell^2)$ and $|\ell|\ge1$, we obtain
\begin{equation}
\begin{aligned}
 \frac{\Lambda_t^2(q,\delta)\kappa(p,\gamma)}
 {|\ell|\sqrt{\ell^2+\xi^2}}
 &\le C\bigg(\frac{\sqrt{\ell^2+\xi^2}}{|\ell|}
                         +t^2\kappa(\ell,\xi)\bigg)\kappa(p,\gamma)\\
 &\le C\bigl[\sqrt{q^2+\delta^2}\kappa(p,\gamma)
       +\av t\Lambda_t(p,\gamma)\kappa(p,\gamma)
       +\av t^2\kappa(\ell,\xi)\kappa(p,\gamma)\bigr].
\end{aligned}
\label{eq:advKQ}
\end{equation}
For the last term, \eqref{eq:shiftkappa} with $\ell=k-m$ and $\xi=\eta-\zeta$ and \eqref{eq:gprimitive} give
\[
 \kappa(\ell,\xi)\kappa(p,\gamma)
 \le C\Lambda_t^5(m,\zeta)\sqrt{g(k,\eta)g(p,\gamma)}.
\]
When $p=0$, \eqref{eq:meanA} gives
$|a(\ell,\xi;0,\gamma)|\le|q|=|\ell|$,
which cancels the denominator $|\ell|$.
By \eqref{eq:tkappa},
\[
 \av t\Lambda_t(p,\gamma)\kappa(p,\gamma)
 \le C\Lambda_t^2(p,\gamma).
\]
Multiplying \eqref{eq:advKQ} by the factor $C\Lambda_t(m,\zeta)\Lambda_t(p,\gamma)$ from the two kernel bounds, and using \eqref{eq:shiftkappa}, \eqref{eq:tkappa}, and \eqref{eq:gprimitive}, gives, for $p\ne0$,
\begin{equation}
\begin{aligned}
\frac{|a(k,\eta;\ell,\xi)a(\ell,\xi;p,\gamma)|}{|\ell|}
\le{} &C\Lambda_t(m,\zeta)\Lambda_t(p,\gamma)\sqrt{q^2+\delta^2}\sqrt{g(p,\gamma)}
+C\Lambda_t(m,\zeta)\Lambda_t^{3}(p,\gamma)\\
&+C\Lambda_t^{6}(m,\zeta)\Lambda_t(p,\gamma)\av t^2\sqrt{g(k,\eta)g(p,\gamma)}.
\label{eq:advKQpolynomial}
\end{aligned}
\end{equation}
For $p=0$, cancellation of $|\ell|$ gives instead
\[
 C\Lambda_t(m,\zeta)(1+\av t\kappa(\ell,\xi))
       \le C\Lambda_t^{6}(m,\zeta)\big(1+\av t\sqrt{g(k,\eta)}\big).
\]
Due to $|\sqrt{k^2+(\eta+kt)^2}-\sqrt{q^2+(\delta+qt)^2}|\le\sqrt{p^2+(\gamma+pt)^2}+\sqrt{m^2+(\zeta+mt)^2}$, $\Lambda_t(k,\eta)$ and $\Lambda_t(q,\delta)$ are comparable and dominate the other two weights due to region 5 of \eqref{eq:advpartition}. Equations \eqref{eq:weight} and \eqref{eq:thetabounds} then give
\begin{equation}
 \frac{|h(k,\eta)-h(m,\zeta)|}
 {\Lambda_t^N(k,\eta)\Lambda_t^N(q,\delta)}
 \le\frac{\Lambda_t^{2N}(k,\eta)+\Lambda_t^{2N}(m,\zeta)}
 {\Lambda_t^N(k,\eta)\Lambda_t^N(q,\delta)}
 \le C_N.
\label{eq:advKQenergy}
\end{equation}
Multiplying \eqref{eq:advKQpolynomial} and the $p=0$ bound by \eqref{eq:advKQenergy}, dividing by $\Lambda_t^N(m,\zeta)\Lambda_t^N(p,\gamma)$, and applying \eqref{eq:advectorquarticexpdecay}, we bound the left-hand side of \eqref{eq:advectorpairstarget} in region 5 by
\begin{align*}
C_Ne^{-\rho^{1/3}t/16}\Lambda_t^{6-N}(m,\zeta)\Lambda_t^{3-N}(p,\gamma)
\big[1+\av t\sqrt{g(k,\eta)}
 +\sqrt{q^2+\delta^2}\sqrt{g(p,\gamma)}
 +\av t^2\sqrt{g(k,\eta)g(p,\gamma)}\big].
\end{align*}
This is region 5 of \eqref{eq:advectorpairstarget}.

\noindent\emph{Case 6.}
For $p\ne0$, \eqref{eq:endpointsym}, \eqref{eq:A}, the sixth condition of \eqref{eq:advpartition} and \eqref{eq:advinnerdeterminant} give
\[
\begin{aligned}
&|a(k,\eta;\ell,\xi)| =|a(m,\zeta;-\ell,-\xi)|
\le C\frac{\Lambda_t(k,\eta) \Lambda_t(\ell,\xi)}{\sqrt{\ell^2+\xi^2}}\le C\frac{\Lambda_t(k,\eta) \Lambda_t(q,\delta)}{\sqrt{\ell^2+\xi^2}},\\
&|a(\ell,\xi;p,\gamma)|\le C\Lambda_t(p,\gamma) \Lambda_t(q,\delta)\kappa(p,\gamma).
\end{aligned}
\]
Since $\ell=-m+k=p+q$, $\xi=-\zeta+\eta=\gamma+\delta$, \eqref{eq:shiftkappa} and \eqref{eq:gprimitive} give
\begin{align*}
 \sqrt{\ell^2+\xi^2}
 &\le\sqrt{q^2+\delta^2}+C\av t\Lambda_t(p,\gamma),\\
 \kappa(\ell,\xi)\kappa(p,\gamma)
 &\le C\Lambda_t^5(k,\eta)\sqrt{g(m,\zeta)g(p,\gamma)}.
\end{align*}
The sixth condition of \eqref{eq:advpartition} gives $\sqrt{p^2+(\gamma+pt)^2}\le\frac18\sqrt{q^2+(\delta+qt)^2}$.
Thus $\Lambda_t(q,\delta)\le C\Lambda_t(\ell,\xi)$. Multiplying \eqref{eq:advKQ} by $C\Lambda_t(k,\eta)\Lambda_t(p,\gamma)$, substituting these two inequalities, and applying \eqref{eq:tkappa} and \eqref{eq:gprimitive}, we obtain
\begin{equation}
\begin{aligned}
\frac{|a(k,\eta;\ell,\xi)a(\ell,\xi;p,\gamma)|}{|\ell|}\le{} &C\Lambda_t(k,\eta)\Lambda_t(p,\gamma)\sqrt{q^2+\delta^2}\sqrt{g(p,\gamma)}
+C\Lambda_t(k,\eta)\Lambda_t^{3}(p,\gamma)\\
&+C\Lambda_t^{6}(k,\eta)\Lambda_t(p,\gamma)\av t^2\sqrt{g(m,\zeta)g(p,\gamma)}.
\label{eq:advMQpolynomial}
\end{aligned}
\end{equation}
For $p=0$, \eqref{eq:meanA}, the second inequality of \eqref{eq:advouterelementary}, and \eqref{eq:shiftkappa} give
\[
 \frac{|a(k,\eta;\ell,\xi)a(\ell,\xi;0,\gamma)|}{|\ell|}
 \le C\Lambda_t(k,\eta)(1+\av t\kappa(\ell,\xi))
 \le C\Lambda_t^6(k,\eta)\big(1+\av t\sqrt{g(m,\zeta)}\big).
\]
The sixth condition in \eqref{eq:advpartition} gives
\begin{align*}
 \big|\sqrt{m^2+(\zeta+mt)^2}-\sqrt{q^2+(\delta+qt)^2}\big|
 &\le\sqrt{k^2+(\eta+kt)^2}+\sqrt{p^2+(\gamma+pt)^2}\\
 &\le\frac18\min\bigl\{\sqrt{m^2+(\zeta+mt)^2},\sqrt{q^2+(\delta+qt)^2}\bigr\}.
\end{align*}
Thus \eqref{eq:weight} and \eqref{eq:thetabounds} imply
\begin{equation}
 \frac{|h(k,\eta)-h(m,\zeta)|}
 {\Lambda_t^N(m,\zeta)\Lambda_t^N(q,\delta)}
 \le\frac{\Lambda_t^{2N}(k,\eta)+\Lambda_t^{2N}(m,\zeta)}
 {\Lambda_t^N(m,\zeta)\Lambda_t^N(q,\delta)}
 \le C_N.
\label{eq:advMQenergy}
\end{equation}
Combining \eqref{eq:advMQpolynomial}, the $p=0$ bound, and \eqref{eq:advMQenergy}, then dividing by $\Lambda_t^N(k,\eta)\Lambda_t^N(p,\gamma)$ and using \eqref{eq:advectorquarticexpdecay}, we find that the left-hand side of \eqref{eq:advectorpairstarget} in region 6 is at most
\begin{align*}
C_Ne^{-\rho^{1/3}t/16}\Lambda_t^{6-N}(k,\eta)\Lambda_t^{3-N}(p,\gamma)
\big[1+\av t\sqrt{g(m,\zeta)}
 +\sqrt{q^2+\delta^2}\sqrt{g(p,\gamma)}
 +\av t^2\sqrt{g(m,\zeta)g(p,\gamma)}\big].
\end{align*}
This proves region 6 of \eqref{eq:advectorpairstarget}.
\end{proof}

\begin{lemma}\label{lem:advectorhighpartition}
Let $k=p+q+m$, $\eta=\gamma+\delta+\zeta$, and $p+q\ne0$.
If all six inequalities in \eqref{eq:advpartition} fail, then  the weights and kernels defined in \eqref{eq:weight} and \eqref{eq:kernelgeometric} satisfy
\begin{align}
&\frac{e^{\rho^{1/3}t\mathbf1_{k\ne0}\mathbf1_{m\ne0}/8}}{e^{\rho^{1/3}t(\mathbf1_{k\ne0}+\mathbf1_{p\ne0}
+\mathbf1_{q\ne0}+\mathbf1_{m\ne0})/16}}\frac{|h(k,\eta)-h(m,\zeta)|}
 {|p+q|\Lambda_t^N(k,\eta)\Lambda_t^N(p,\gamma)\Lambda_t^N(q,\delta)\Lambda_t^N(m,\zeta)}\nonumber\\
&\times
 |a(k,\eta;p+q,\gamma+\delta)a(p+q,\gamma+\delta;p,\gamma)|\nonumber\\
\le{} &C_Ne^{-\rho^{1/3}t/16}\Bigl[
\Lambda_t^{-2}(k,\eta)\Lambda_t^{-2}(m,\zeta)
 +\Lambda_t^{-2}(k,\eta)\Lambda_t^{-2}(p,\gamma)
 +\Lambda_t^{-2}(k,\eta)\Lambda_t^{-2}(q,\delta)\nonumber\\
&+\Lambda_t^{-2}(m,\zeta)\Lambda_t^{-2}(p,\gamma)
 +\Lambda_t^{-2}(m,\zeta)\Lambda_t^{-2}(q,\delta)
 +\Lambda_t^{-2}(p,\gamma)\Lambda_t^{-2}(q,\delta)\Bigr].
\label{eq:threehighadvector}
\end{align}
\end{lemma}
\begin{proof}
Since $p+q\ne0$, \eqref{eq:kernelgeometric} gives
\begin{align}
|a(k,\eta;p+q,\gamma+\delta)|
 \le\frac{\Lambda_t(p+q,\gamma+\delta)\Lambda_t(m,\zeta)}
             {\sqrt{(p+q)^2+(\gamma+\delta)^2}}
 \le[\Lambda_t(p,\gamma)+\Lambda_t(q,\delta)]\Lambda_t(m,\zeta).
\label{eq:roughadvouter}
\end{align}
For $p\ne0$, the same formula yields
\[
 |a(p+q,\gamma+\delta;p,\gamma)|
 \le\frac{\Lambda_t(p,\gamma)\Lambda_t(q,\delta)}
             {\sqrt{p^2+\gamma^2}}
 \le\Lambda_t(p,\gamma)\Lambda_t(q,\delta).
\]
For $p=0$, \eqref{eq:meanA} instead gives
$|a(q,\gamma+\delta;0,\gamma)|\le|q|\le\Lambda_t(q,\delta)$.
These estimates and \eqref{eq:roughadvouter} imply
\begin{align}
|a(k,\eta;p+q,\gamma+\delta)a(p+q,\gamma+\delta;p,\gamma)|
\le C\max\{\Lambda_t(k,\eta),\Lambda_t(m,\zeta),
                 \Lambda_t(p,\gamma),\Lambda_t(q,\delta)\}^{4}.
\label{eq:roughadvproduct}
\end{align}
By \eqref{eq:weight} and \eqref{eq:thetabounds},
\begin{equation}
 |h(k,\eta)-h(m,\zeta)|
 \le\Lambda_t^{2N}(k,\eta)+\Lambda_t^{2N}(m,\zeta).
\label{eq:roughadvweight}
\end{equation}
Consider first
\[
 \sqrt{k^2+(\eta+kt)^2}\ge\sqrt{m^2+(\zeta+mt)^2}
 \ge\sqrt{p^2+(\gamma+pt)^2}\ge\sqrt{q^2+(\delta+qt)^2}.
\]
Since $k=p+q+m$ and $\eta=\gamma+\delta+\zeta$, and the first
inequality of \eqref{eq:advpartition} fails, we have
\begin{align}
 \sqrt{k^2+(\eta+kt)^2}
 &\le3\sqrt{m^2+(\zeta+mt)^2},\nonumber\\
 \sqrt{m^2+(\zeta+mt)^2}
 &<8\big[\sqrt{p^2+(\gamma+pt)^2}
                  +\sqrt{q^2+(\delta+qt)^2}\big]
 \le16\sqrt{p^2+(\gamma+pt)^2}.
\label{eq:advthreehighsizes}
\end{align}
Using \eqref{eq:roughadvproduct}, \eqref{eq:roughadvweight},
$|p+q|\ge1$, and then \eqref{eq:advthreehighsizes}, we obtain
\begin{align}
&\frac{|h(k,\eta)-h(m,\zeta)|
 |a(k,\eta;p+q,\gamma+\delta)a(p+q,\gamma+\delta;p,\gamma)|}
 {|p+q|\Lambda_t^N(k,\eta)\Lambda_t^N(m,\zeta)
          \Lambda_t^N(p,\gamma)\Lambda_t^N(q,\delta)}\nonumber\\
\le{}&\frac{C_N\Lambda_t^{2N+4}(k,\eta)}
 {\Lambda_t^N(k,\eta)\Lambda_t^N(m,\zeta)
          \Lambda_t^N(p,\gamma)\Lambda_t^N(q,\delta)}\nonumber\\
\le{}& C_N\Lambda_t^{4-N}(p,\gamma)\Lambda_t^{-N}(q,\delta)
 \le C_N\Lambda_t^{-2}(p,\gamma)\Lambda_t^{-2}(q,\delta).
\label{eq:advthreehighnormalized}
\end{align}
For each other ordering, the frequency relation again bounds
the largest length by three times the second largest.
Failure of the corresponding inequality in \eqref{eq:advpartition}
bounds the second largest by sixteen times the third largest.
The calculation \eqref{eq:advthreehighnormalized}, using the
uniform bounds \eqref{eq:roughadvproduct} and \eqref{eq:roughadvweight},
therefore gives the product of inverse squares at the two
smaller frequencies. The six choices give the six products in
\eqref{eq:threehighadvector}. Multiplying by the exponential
bound \eqref{eq:advectorquarticexpdecay} proves that inequality.
\end{proof}

\begin{proof}[Proof of Proposition~\ref{prop:quartic}: estimate of $I_6$]
We next estimate $I_6$ in \eqref{eq:adv}, using \eqref{eq:controlenergy}, \eqref{eq:minimumcontrol}, \eqref{eq:actualsobolevL1}, and \eqref{eq:quartictimebounds}. All time norms are taken over $(0,T_1)$.

Set $k=p+q+m$, $\eta=\gamma+\delta+\zeta$, and $p+q\ne0$.
Lemmas~\ref{lem:advectorpairs} and \ref{lem:advectorhighpartition}
give pointwise estimates for the kernel
\begin{align*}
 e^{\rho^{1/3}t\mathbf1_{k\ne0}\mathbf1_{m\ne0}/8}
 (h(k,\eta)-h(m,\zeta))
 a(k,\eta;p+q,\gamma+\delta)a(p+q,\gamma+\delta;p,\gamma).
\end{align*}
Substituting \eqref{eq:weightedamplitudes} and applying Young's inequality, we obtain
\[
 \int_0^{T_1}|I_6(t)|\dd t\le\sum_{r=1}^7L_r.
\]
Here $L_1,\ldots,L_7$ denote the upper bounds for the second integral. We estimate them as follows.

\noindent\emph{Case 1. Region 1.}
Applying the first bound in \eqref{eq:advectorpairstarget}, with the factors at $(p,\gamma)$ and $(q,\delta)$ in $L^1$ and the other two in $L^2$, we obtain
\begin{align*}
L_1 &=\frac{C_N}{|\beta|}\sum_s\int_0^{T_1} e^{-\rho^{1/3}t/16}
 \|\Lambda_t^{9-N}F^s\|_1\|\Lambda_t^{9-N}F^{-s}\|_1
 \big[\av{\rho^{1/3}t}\|F^s\|_2^2
              +\av t\|\sqrt gF^s\|_2\|F^s\|_2\big]\dd t\nonumber\\
 &\le\frac{C_N}{|\beta|}\sum_s
 \|F^s\|_{L^\infty_tL^2}^2\|F^{-s}\|_{L^\infty_tL^2}
 \Bigl[
 \|e^{-\rho^{1/3}t/16}\av{\rho^{1/3}t}\|_{L^1_t}
       \|F^s\|_{L^\infty_tL^2}
 +\|\av t e^{-\rho^{1/3}t/16}\|_{L^2_t}
       \|\sqrt gF^s\|_{L^2_tL^2}\Bigr]\nonumber\\
 &\le\frac{C_N}{|\beta|}
       (\rho^{-1/3}+\rho^{-1/2})\mathcal E(T_1)^2.
\end{align*}

\noindent\emph{Case 2. Region 2.}
Applying the second bound in \eqref{eq:advectorpairstarget}, with the factors at $(k,\eta)$ and $(m,\zeta)$ in $L^1$ and the other two in $L^2$, we obtain
\begin{align*}
L_2 &=\frac{C_N}{|\beta|}\sum_s\int_0^{T_1} e^{-\rho^{1/3}t/16}
 \|\Lambda_t^{9-N}F^s\|_1^2\|F^{-s}\|_2
       \big[\|F^s\|_2+\av t\|\sqrt gF^s\|_2\big]\dd t\nonumber\\
 &\le\frac{C_N}{|\beta|}\sum_s
 \|F^s\|_{L^\infty_tL^2}^2\|F^{-s}\|_{L^\infty_tL^2}
 \Bigl[
 \|e^{-\rho^{1/3}t/16}\|_{L^1_t}\|F^s\|_{L^\infty_tL^2}
 +\|\av t e^{-\rho^{1/3}t/16}\|_{L^2_t}
       \|\sqrt gF^s\|_{L^2_tL^2}\Bigr]\nonumber\\
 &\le\frac{C_N}{|\beta|}
       (\rho^{-1/3}+\rho^{-1/2})\mathcal E(T_1)^2.
\end{align*}

\noindent\emph{Case 3. Region 3.}
Applying the third bound in \eqref{eq:advectorpairstarget}, with the factors at $(m,\zeta)$ and $(q,\delta)$ in $L^1$ and the other two in $L^2$, we obtain
\begin{align*}
L_3 ={} &\frac{C_N}{|\beta|}\sum_s\int_0^{T_1} e^{-\rho^{1/3}t/16}
 \|\Lambda_t^{9-N}F^s\|_1\|\Lambda_t^{9-N}F^{-s}\|_1
 \bigl[\|F^s\|_2^2+\av t\|\sqrt gF^s\|_2\|F^s\|_2
                       +\av t^2\|\sqrt gF^s\|_2^2\bigr]\dd t\nonumber\\
 \le{} &\frac{C_N}{|\beta|}\sum_s
 \|F^s\|_{L^\infty_tL^2}\|F^{-s}\|_{L^\infty_tL^2}
 \Bigl[
 \|e^{-\rho^{1/3}t/16}\|_{L^1_t}\|F^s\|_{L^\infty_tL^2}^2\nonumber\\
 &\qquad+\|\av t e^{-\rho^{1/3}t/16}\|_{L^2_t}
       \|\sqrt gF^s\|_{L^2_tL^2}\|F^s\|_{L^\infty_tL^2}
 +\|\av t^2e^{-\rho^{1/3}t/16}\|_{L^\infty_t}
       \|\sqrt gF^s\|_{L^2_tL^2}^2\Bigr]\nonumber\\
 \le{} &\frac{C_N}{|\beta|}
       (\rho^{-1/3}+\rho^{-1/2}+\rho^{-2/3})\mathcal E(T_1)^2.
\end{align*}

\noindent\emph{Case 4. Region 4.}
We apply the fourth bound in \eqref{eq:advectorpairstarget}, with the factors at $(k,\eta)$ and $(q,\delta)$ in $L^1$ and the other two in $L^2$. Arguing as in Case 3, we obtain
\begin{align*}
L_4 ={} &\frac{C_N}{|\beta|}\sum_s\int_0^{T_1} e^{-\rho^{1/3}t/16}
 \|\Lambda_t^{9-N}F^s\|_1\|\Lambda_t^{9-N}F^{-s}\|_1
 \bigl[\|F^s\|_2^2+\av t\|\sqrt gF^s\|_2\|F^s\|_2
                       +\av t^2\|\sqrt gF^s\|_2^2\bigr]\dd t\nonumber\\
 \le{} &\frac{C_N}{|\beta|}
       (\rho^{-1/3}+\rho^{-1/2}+\rho^{-2/3})\mathcal E(T_1)^2.
\end{align*}

\noindent\emph{Case 5. Region 5.}
Applying the fifth bound in \eqref{eq:advectorpairstarget}, with the factors at $(m,\zeta)$ and $(p,\gamma)$ in $L^1$ and the other two in $L^2$, we obtain
\begin{align*}
 L_5 ={} &\frac{C_N}{|\beta|}\sum_s\int_0^{T_1} e^{-\rho^{1/3}t/16}
 \|\Lambda_t^{9-N}F^s\|_1
 \Bigl[\|\Lambda_t^{9-N}F^s\|_1\|F^{-s}\|_2
       \big[\|F^s\|_2+\av t\|\sqrt gF^s\|_2\big]\nonumber\\
 &\qquad+\|\Lambda_t^{9-N}\sqrt gF^s\|_1
       \big[\|F^s\|_2\||D|F^{-s}\|_2
                      +\av t^2\|\sqrt gF^s\|_2\|F^{-s}\|_2\big]\Bigr]\dd t\nonumber\\
 \le{} &\frac{C_N}{|\beta|}\sum_s
 \Bigl[
 \|e^{-\rho^{1/3}t/16}\|_{L^1_t}
       \|F^s\|_{L^\infty_tL^2}^3\|F^{-s}\|_{L^\infty_tL^2}\nonumber\\
 &\qquad+\|\av t e^{-\rho^{1/3}t/16}\|_{L^2_t}
       \|\sqrt gF^s\|_{L^2_tL^2}
       \|F^s\|_{L^\infty_tL^2}^2\|F^{-s}\|_{L^\infty_tL^2}\nonumber\\
 &\qquad+\rho^{-1/2}
       \|F^s\|_{L^\infty_tL^2}^2
       \|\sqrt\rho|D|F^{-s}\|_{L^2_tL^2}
       \|\sqrt gF^s\|_{L^2_tL^2}\nonumber\\
 &\qquad+\|\av t^2e^{-\rho^{1/3}t/16}\|_{L^\infty_t}
       \|\sqrt gF^s\|_{L^2_tL^2}^2
       \|F^s\|_{L^\infty_tL^2}\|F^{-s}\|_{L^\infty_tL^2}\Bigr]\nonumber\\
 \le{} &\frac{C_N}{|\beta|}
       (\rho^{-1/3}+\rho^{-1/2}+\rho^{-2/3})\mathcal E(T_1)^2.
\end{align*}

\noindent\emph{Case 6. Region 6.}
We apply the sixth bound in \eqref{eq:advectorpairstarget}, with the factors at $(k,\eta)$ and $(p,\gamma)$ in $L^1$ and the other two in $L^2$. Arguing as in Case 5, we obtain
\begin{align*}
L_6 ={} &\frac{C_N}{|\beta|}\sum_s\int_0^{T_1} e^{-\rho^{1/3}t/16}
 \|\Lambda_t^{9-N}F^s\|_1
 \Bigl[\|\Lambda_t^{9-N}F^s\|_1\|F^{-s}\|_2
       \big[\|F^s\|_2+\av t\|\sqrt gF^s\|_2\big]\nonumber\\
 &\qquad+\|\Lambda_t^{9-N}\sqrt gF^s\|_1
       \big[\|F^s\|_2\||D|F^{-s}\|_2
                      +\av t^2\|\sqrt gF^s\|_2\|F^{-s}\|_2\big]\Bigr]\dd t\nonumber\\
 \le{} &\frac{C_N}{|\beta|}
       (\rho^{-1/3}+\rho^{-1/2}+\rho^{-2/3})\mathcal E(T_1)^2.
\end{align*}

\noindent\emph{Case 7. Region 7.}
By Lemma~\ref{lem:advectorhighpartition} and Young's inequality, we obtain
\begin{align*}
L_7 &=\frac{C_N}{|\beta|}\sum_s\int_0^{T_1} e^{-\rho^{1/3}t/16}
 \Bigl[
 \|F^s\|_2^2\|\Lambda_t^{-2}F^s\|_1\|\Lambda_t^{-2}F^{-s}\|_1
 +\|F^s\|_2\|F^{-s}\|_2\|\Lambda_t^{-2}F^s\|_1^2
 \Bigr]\dd t\\
 &\le\frac{C_N}{|\beta|}
 \|e^{-\rho^{1/3}t/16}\|_{L^1_t}
 \sum_s\|F^s\|_{L^\infty_tL^2}^3\|F^{-s}\|_{L^\infty_tL^2}\\
 &\le\frac{C_N}{|\beta|}\rho^{-1/3}\mathcal E(T_1)^2.
\end{align*}
Adding the above estimates and using $0<\rho\le1$, we obtain
\begin{equation}
 \int_0^{T_1}|I_6(t)|\dd t
 \le\sum_{r=1}^7L_r
 \le\frac{C_N}{|\beta|}\rho^{-2/3}\mathcal E(T_1)^2.
 \label{eq:advectorquarticintegral}
\end{equation}
Adding \eqref{eq:pairedquarticintegral} and
\eqref{eq:advectorquarticintegral} proves \eqref{eq:quarticbound}.
\end{proof}

\section{Proof of the energy estimate}
\label{sec:bootstrap}

In this section, we prove Proposition~\ref{prop:closedenergy} and Corollary~\ref{cor:decay} separately. 

\subsection{Proof of Proposition~\ref{prop:closedenergy}}

\begin{lemma}
\label{lem:allremainders}
With the notation of Section~\ref{sec:equations}, let $(u,b)$ be a smooth solution of the MHD equations \eqref{eq:perturbation} on $[0,T_1]$. Then, for every $t\in[0,T_1]$, we have
\begin{align}
 |\mathcal B(t)|&\le\frac{C_N}{|\beta|}\rho^{-1/3}
                              \mathcal E(T_1)^{1/2}E(t),
                         \label{eq:boundaryestimate}\\
 |\mathcal B(0)|&\le\frac{C_N}{|\beta|}\varepsilon^3,
 \qquad E(0)\le4\varepsilon^2.
                         \label{eq:initialenergy}
\end{align}
Moreover, it holds that
\begin{equation}\label{eq:allremainders}
\begin{split}
 &E(t)-\mathcal B(t)
   +\int_0^t\Big[2\nu\|\nabla\mathcal A\widetilde u\|_2^2
                         +2\mu\|\nabla\mathcal Ab\|_2^2
 +\frac{2(\mu-\nu)}\beta
         \langle\nabla\mathcal A\widetilde u^1,
                      \nabla\mathcal A\partial_x^{-1}b^2_{\ne}\rangle
        -4\langle\mathcal A\widetilde u,
                       \mathcal A\partial_x\nabla\Delta^{-1}u^2\rangle\\
 &\qquad+2\sum_k\int_\R(7g+(\partial_t-k\partial_\eta)V)\mathcal A^2
                (|\widehat{\widetilde u}|^2+|\widehat b|^2)\frac{\dd\eta}{2\pi}
        -\frac{\rho^{1/3}}8
            \big(\|\mathcal A\widetilde u_{\ne}\|_2^2+\|\mathcal Ab_{\ne}\|_2^2\big)
        \Big]\dd\tau\\
 \le{} &E(0)-\mathcal B(0)
        +C_N\rho^{-1/3}\mathcal E(T_1)^{3/2}
        +C_N\rho^{-2/3}\mathcal E(T_1)^2.
\end{split}
\end{equation}
\end{lemma}

\begin{proof}
By \eqref{eq:Bbound}, \eqref{eq:quadraticenergy} and \eqref{eq:energyequivalence},
\begin{align*}
 |\mathcal B(t)|
 &\le\frac{C_N}{|\beta|}\av t e^{-\rho^{1/3}t/16}
 \Bigl(\sum_s\|F^s(t)\|_2^2\Bigr)^{1/2}
                         \sum_s\|F^s(t)\|_2^2\\
 &\le\frac{C_N}{|\beta|}\rho^{-1/3}\mathcal E(T_1)^{1/2}E(t).
\end{align*}
This is \eqref{eq:boundaryestimate}. At $t=0$, $\theta=1$ and $\mathcal A=\Lambda_0^N$. Hence
\begin{align*}
 E(0)=\|\Lambda_0^N\widetilde u(0)\|_2^2
                      +\|\Lambda_0^Nb(0)\|_2^2
 \le4\bigl(\|u(0)\|_{H^N}^2+\|b(0)\|_{H^N}^2\bigr)
 \le4\varepsilon^2,\qquad
 |\mathcal B(0)|\le\frac{C_N}{|\beta|}\varepsilon^3,
\end{align*}
by \eqref{eq:physicalenergycomparison} and \eqref{eq:Binitial}.

Integrating the identity in Proposition~\ref{prop:identity} from $0$
to $t\le T_1$ and applying \eqref{eq:Nbound}, \eqref{eq:T0bound},
\eqref{eq:TJbound}, \eqref{eq:R3bound}, \eqref{eq:mixeddiffusionbound},
and \eqref{eq:quarticbound}, we obtain
\begin{align*}
 &E(t)-\mathcal B(t)
   +\int_0^t\Big[2\nu\|\nabla\mathcal A\widetilde u\|_2^2
                         +2\mu\|\nabla\mathcal Ab\|_2^2
 +\frac{2(\mu-\nu)}\beta
         \langle\nabla\mathcal A\widetilde u^1,
                      \nabla\mathcal A\partial_x^{-1}b^2_{\ne}\rangle
        -4\langle\mathcal A\widetilde u,
                       \mathcal A\partial_x\nabla\Delta^{-1}u^2\rangle\\
 &\quad+2\sum_k\int_\R(7g+(\partial_t-k\partial_\eta)V)\mathcal A^2
                (|\widehat{\widetilde u}|^2+|\widehat b|^2)\frac{\dd\eta}{2\pi}
        -\frac{\rho^{1/3}}8
            (\|\mathcal A\widetilde u_{\ne}\|_2^2+\|\mathcal Ab_{\ne}\|_2^2)
        \Big]\dd\tau\\
\le{} & E(0)-\mathcal B(0)+C_N\Big(1+\frac1{|\beta|}
       +\frac{\nu+\mu}{|\beta|\sqrt{\nu\mu}}\Big)
                         \rho^{-1/3}\mathcal E(T_1)^{3/2}
 +\frac{C_N}{|\beta|}\rho^{-2/3}\mathcal E(T_1)^2\\
\le{} & E(0)-\mathcal B(0)+C_N\rho^{-1/3}\mathcal E(T_1)^{3/2}
                   +C_N\rho^{-2/3}\mathcal E(T_1)^2,
\end{align*}
where we also used \eqref{eq:coefficientratios} to conclude $|\beta|^{-1}\le3/2$ and $(\nu+\mu)/(|\beta|\sqrt{\nu\mu})\le3$.
The proof is thus complete.
\end{proof}

\begin{proof}[Proof of Proposition~\ref{prop:closedenergy}]
\noindent\emph{Step 1. The diffusion terms.}
Due to $\widetilde u=u+\beta^{-1}\partial_x^{-1}b_{\ne}^2e_1$, we have
\begin{align*}
 2\nu\|\nabla\mathcal A\widetilde u\|_2^2
 &=2\nu\|\nabla\mathcal Au\|_2^2
 +\frac{4\nu}{\beta}\langle\nabla\mathcal Au^1,
                          \nabla\mathcal A\partial_x^{-1}b^2_{\ne}\rangle
 +\frac{2\nu}{\beta^2}\|\nabla\mathcal A\partial_x^{-1}b^2_{\ne}\|_2^2,\\
 \frac{2(\mu-\nu)}\beta
 \langle\nabla\mathcal A\widetilde u^1,
              \nabla\mathcal A\partial_x^{-1}b^2_{\ne}\rangle
 &=\frac{2(\mu-\nu)}\beta
 \langle\nabla\mathcal Au^1,\nabla\mathcal A\partial_x^{-1}b^2_{\ne}\rangle
 +\frac{2(\mu-\nu)}{\beta^2}
                        \|\nabla\mathcal A\partial_x^{-1}b^2_{\ne}\|_2^2.
\end{align*}
Incompressibility then gives
\[
 \|\nabla\mathcal A\partial_x^{-1}b_{\ne}^2\|_2^2
 =\sum_{k\ne0}\int_\R
      \mathcal A^2(k,\eta)\frac{k^2+\eta^2}{k^2}
                          |\widehat b^2(k,\eta)|^2
                          \frac{\dd\eta}{2\pi}
 =\|\mathcal Ab_{\ne}\|_2^2.
\]
Hence we get
\begin{equation}
\begin{aligned}
 &2\nu\|\nabla\mathcal A\widetilde u\|_2^2
      +2\mu\|\nabla\mathcal Ab\|_2^2+\frac{2(\mu-\nu)}{\beta}
     \left\langle\nabla\mathcal A\widetilde u^1,
           \nabla\mathcal A\partial_x^{-1}b_{\ne}^2\right\rangle
                                                        \\
 ={}&2\nu\|\nabla\mathcal Au\|_2^2+2\mu\|\nabla\mathcal Ab\|_2^2
       +\frac{2\mu}{\beta^2}\|\mathcal Ab_{\ne}\|_2^2+\frac{2(\nu+\mu)}{\beta}
     \left\langle\nabla\mathcal Au^1,
           \nabla\mathcal A\partial_x^{-1}b_{\ne}^2\right\rangle.
 \label{eq:combinedquadraticidentity}
\end{aligned}
\end{equation}
For $k\ne0$, incompressibility gives
\begin{align}
 |\widehat u^1(k,\eta)|
       =\frac{|\eta|}{\sqrt{k^2+\eta^2}}|\widehat u(k,\eta)|,
 \qquad
 |\widehat b^2(k,\eta)|
       =\frac{|k|}{\sqrt{k^2+\eta^2}}|\widehat b(k,\eta)|.
\end{align}
Thus \eqref{eq:strongfield} and $2|\eta|\le k^2+\eta^2 $ give
\begin{equation}
\begin{aligned}
 \Big|\frac{2(\nu+\mu)}{\beta}
     \left\langle\nabla\mathcal Au^1,
           \nabla\mathcal A\partial_x^{-1}b_{\ne}^2\right\rangle\Big|
 &\le\frac{2(\nu+\mu)}{|\beta|}
       \sum_{k\ne0}\int_\R \mathcal A^2(k,\eta)|\eta|
               |\widehat u(k,\eta)||\widehat b(k,\eta)|
                                      \frac{\dd\eta}{2\pi}\\
 &\le\frac{\nu+\mu}{|\beta|\sqrt{\nu\mu}}
       \sum_{k\ne0}\int_\R \mathcal A^2(k,\eta)|\eta|
       \bigl(\nu|\widehat u(k,\eta)|^2+\mu|\widehat b(k,\eta)|^2\bigr) \frac{\dd\eta}{2\pi}\\
 &\le\frac{\nu+\mu}{2|\beta|\sqrt{\nu\mu}}
       \bigl(\nu\|\nabla\mathcal Au\|_2^2
                        +\mu\|\nabla\mathcal Ab\|_2^2\bigr)\\
 &\le\frac32
       \bigl(\nu\|\nabla\mathcal Au\|_2^2
                        +\mu\|\nabla\mathcal Ab\|_2^2\bigr).
\end{aligned}
\end{equation}
Combining the preceding estimate with \eqref{eq:combinedquadraticidentity} and using $\frac{2\mu}{\beta^2}\|\mathcal Ab_{\ne}\|_2^2\ge0$, we obtain
\begin{equation}
\begin{aligned}
2\nu\|\nabla\mathcal A\widetilde u\|_2^2
      +2\mu\|\nabla\mathcal Ab\|_2^2+\frac{2(\mu-\nu)}{\beta}
     \left\langle\nabla\mathcal A\widetilde u^1,
           \nabla\mathcal A\partial_x^{-1}b_{\ne}^2\right\rangle
\ge\frac12
       \bigl(\nu\|\nabla\mathcal Au\|_2^2
                        +\mu\|\nabla\mathcal Ab\|_2^2\bigr).
\end{aligned}
\end{equation}
For $k=0$, the correction and both mixed terms vanish, and the physical diffusion is unchanged.

\noindent\emph{Step 2. The linear pressure.}
Due to \eqref{eq:pressurebound}, we have
\[
\big|\widehat{\partial_x\nabla\Delta^{-1}u^2}(k,\eta)\big|
       \le\kappa^2(k,\eta)|\widehat u(k,\eta)|.
\]
This symbol vanishes at $k=0$. Since $\kappa^2\le g$ and
$|\widehat u|\le|\widehat{\widetilde u}|+|\beta|^{-1}|\widehat b|$,
\begin{equation}
\begin{aligned}
 \big|4\left\langle\mathcal A\widetilde u,
              \mathcal A\partial_x\nabla\Delta^{-1}u^2\right\rangle\big|
 &\le4\sum_k\int_\R
       \kappa^2(k,\eta)\mathcal A^2(k,\eta)
       |\widehat{\widetilde u}(k,\eta)||\widehat u(k,\eta)|
                                      \frac{\dd\eta}{2\pi}\\
 &\le(4+2|\beta|^{-1})\sum_k\int_\R
       g(k,\eta)\mathcal A^2(k,\eta)
       \bigl(|\widehat{\widetilde u}(k,\eta)|^2
                            +|\widehat b(k,\eta)|^2\bigr)
                                      \frac{\dd\eta}{2\pi}\\
 &\le7\sum_k\int_\R
       g(k,\eta)\mathcal A^2(k,\eta)
       \bigl(|\widehat{\widetilde u}(k,\eta)|^2
                            +|\widehat b(k,\eta)|^2\bigr)
                                      \frac{\dd\eta}{2\pi}.
\end{aligned}
\end{equation}
Here we used $|\beta|^{-1}\le3/2$ from
\eqref{eq:coefficientratios}. Subtraction from the $14g$ term
in \eqref{eq:allremainders} leaves $7g$.

\noindent\emph{Step 3. The nonzero horizontal modes.}
$\nu,\mu\ge\rho$ and
\eqref{eq:physicalenergycomparison} give
\begin{equation*}
\begin{aligned}
 \frac14(k^2+\eta^2)
       \bigl(\nu|\widehat u(k,\eta)|^2+\mu|\widehat b(k,\eta)|^2\bigr)
  \ge\frac{\rho}{16}(k^2+\eta^2)
       \bigl(|\widehat{\widetilde u}(k,\eta)|^2
                             +|\widehat b(k,\eta)|^2\bigr).
\end{aligned}
\end{equation*}
For $k\ne0$, using $k^2\ge1$ and \eqref{eq:physicalV}, we have
\[
\begin{aligned}
\frac{\rho(k^2+\eta^2)}{16\rho^{1/3}}
 +\frac{2(\partial_t-k\partial_\eta)V(k,\eta)}{\rho^{1/3}}-\frac12
\ge\frac{\rho^{2/3}\eta^2}{16k^2}
 +\frac2{1+\rho^{2/3}\eta^2/k^2}-\frac12
=\frac{(\rho^{2/3}\eta^2/k^2-7/2)^2+47/4}
 {16(1+\rho^{2/3}\eta^2/k^2)}\ge0,
\end{aligned}
\]
from which we infer
\begin{equation}
 \frac{\rho}{16}(k^2+\eta^2)+2(\partial_t-k\partial_{\eta})V(k,\eta)
                         -\frac18\rho^{1/3}
       \ge\frac38\rho^{1/3}.
\end{equation}
For $k=0$, both $(\partial_t-k\partial_\eta)V$ and $\rho^{1/3}\mathbf1_{k\ne0}$ vanish. Combining Steps 1 and 2 with the preceding inequality, we obtain the following lower bound for the integrand in \eqref{eq:allremainders}:
\begin{equation}
\begin{aligned}
 \frac14\bigl(\nu\|\nabla\mathcal Au\|_2^2
                       +\mu\|\nabla\mathcal Ab\|_2^2\bigr)
 +\sum_k\int_\R
       \Big(7g(k,\eta)+\frac38\rho^{1/3}1_{k\ne0}\Big)
       \mathcal A^2(k,\eta)
       \bigl(|\widehat{\widetilde u}(k,\eta)|^2
                            +|\widehat b(k,\eta)|^2\bigr)
                                      \frac{\dd\eta}{2\pi}.
\end{aligned}
\end{equation}
Substitution in \eqref{eq:allremainders} of Lemma~\ref{lem:allremainders} proves
\eqref{eq:correctedenergyestimate}.

\noindent\emph{Step 4. Bounds for $\mathcal B(t)$ and $\mathcal E(T_1)$.}
We first estimate $\mathcal B(t)$. By
\eqref{eq:boundaryestimate}, \eqref{eq:energyequivalence}, and
$|\beta|^{-1}\le3/2$,
\begin{equation}
\begin{aligned}
 \sup_{0\le t\le T_1}|\mathcal B(t)|
 \le\frac{C_N}{|\beta|}\rho^{-1/3}\mathcal E(T_1)^{1/2}
                  \sup_{0\le t\le T_1}E(t)
 \le C_N\rho^{-1/3}\mathcal E(T_1)^{3/2}.
\end{aligned}
\label{eq:boundarysupremum}
\end{equation}
Together with \eqref{eq:initialenergy}, this proves \eqref{eq:correctionboundarycontrol}.

At the initial time, \eqref{eq:initialenergy} gives
\[
 E(0)-\mathcal B(0)
 \le E(0)+|\mathcal B(0)|
 \le4\varepsilon^2+\frac{C_N}{|\beta|}\varepsilon^3
 \le C_N\varepsilon^2,
\]
where $\varepsilon\le1$ and \eqref{eq:coefficientratios} were used.
Moving $\mathcal B(t)$ to the right-hand side of \eqref{eq:correctedenergyestimate}
and using the initial estimate and \eqref{eq:boundarysupremum}, we obtain,
for every $0\le t\le T_1$,
\begin{equation}
\begin{aligned}
 &E(t)
 +\frac14\int_0^t
       \bigl(\nu\|\nabla\mathcal A(s)u(s)\|_2^2
                     +\mu\|\nabla\mathcal A(s)b(s)\|_2^2\bigr)\dd s
                                                        \\
 &+\int_0^t\sum_k\int_\R
       \Big(7g(s,k,\eta)+\frac38\rho^{1/3}1_{k\ne0}\Big)
       \mathcal A^2(s,k,\eta)
       \bigl(|\widehat{\widetilde u}(s,k,\eta)|^2
                            +|\widehat b(s,k,\eta)|^2\bigr)
                         \frac{\dd\eta}{2\pi}\dd s\\
 \le{} &C_N\big[
       \varepsilon^2+\rho^{-1/3}\mathcal E(T_1)^{3/2}
                         +\rho^{-2/3}\mathcal E(T_1)^2\big].
\end{aligned}
\end{equation}
By \eqref{eq:thetabounds} and \eqref{eq:energyequivalence},
\begin{align*}
&\sup_{t\le T_1}\sum_s\|F^s(t)\|_2^2
 \le C\sup_{t\le T_1}E(t),\\
&\int_0^{T_1}\sum_s\bigl(\|\sqrt gF^s\|_2^2
                    +\rho^{1/3}\|F^s_{\ne}\|_2^2\bigr)\dd t
\le C\int_0^{T_1}\sum_k\int_\R
       \Big(7g+\frac38\rho^{1/3}1_{k\ne0}\Big)\mathcal A^2
       (|\widehat{\widetilde u}|^2+|\widehat b|^2)
                          \frac{\dd\eta}{2\pi}\dd t,\\
&2\int_0^{T_1}\sum_k\int_\R e^{\rho^{1/3}t\mathbf 1_{k\ne0}/8}\Lambda_t^{2N}
 (k^2+\eta^2)(\nu|\widehat u|^2+\mu|\widehat b|^2)
                          \frac{\dd\eta}{2\pi}\dd t
\le C\int_0^{T_1}\big(\nu\|\nabla\mathcal Au\|_2^2
                   +\mu\|\nabla\mathcal Ab\|_2^2\big)\dd t.
\end{align*}
Combining the above estimates, we obtain
\[
 \mathcal E(T_1)\le C_0\bigl[\varepsilon^2
       +\rho^{-1/3}\mathcal E(T_1)^{3/2}
       +\rho^{-2/3}\mathcal E(T_1)^2\bigr].
\]
This proves \eqref{eq:closedenergy}, where $C_0\ge4$ depends only on $N$ and is independent of $T_1,\nu,\mu,\beta$ under \eqref{eq:strongfield}.
\end{proof}

\subsection{Proof of Corollary~\ref{cor:decay}}

\begin{proof}
The enhanced-dissipation estimate \eqref{eq:fielddecay} follows directly from \eqref{eq:mainbound}.
Let $w=u$ or $b$. By incompressibility, $w^2=w_{\ne}^2$.
For $k\ne0$, the divergence constraint
$k\widehat w^1+\eta\widehat w^2=0$ yields
\[
 |\widehat w^2(t,k,\eta)|
 =\frac{|k|}{\sqrt{k^2+\eta^2}}|\widehat w(t,k,\eta)|
 =\kappa(k,\eta)|\widehat w(t,k,\eta)|.
\]
By \eqref{eq:tkappa}, we have $(1+t)\kappa(k,\eta)\le C\Lambda_t(k,\eta)$ for $k\ne0$.
Plancherel's identity and \eqref{eq:fielddecay} therefore give
\begin{align*}
 \|\Lambda_t^{N-1}(u^2,b^2)(t)\|_{L^2}
 \le C(1+t)^{-1}
       \|\Lambda_t^N(u_{\ne},b_{\ne})(t)\|_{L^2}
 \le C_N\varepsilon(1+t)^{-1}e^{-\rho^{1/3}t/16},
\end{align*}
which proves \eqref{eq:verticaldecay}.
\end{proof}

\section*{Acknowledgements}
J. Wu was supported by the National Science Foundation of USA (DMS 2104682 and DMS 2309748). T. Yu was supported by the National Natural Science Foundation of China (Grant No. 12501297) and the Natural Science Foundation of Liaoning Province in China (Grant No. 2026-BS-0217) and the Fundamental Research Funds for the Central Universities (Grant No. 3132026351). N. Zhu was supported by the National Natural Science Foundation of China (Grants No. 12301285, No. 12171010) and project ZR2023QA002 supported by Shandong Provincial Natural Science Foundation.

\subsection*{Data availability}
No data was used for the research described in the article.

\subsection*{Disclosure statement}
The authors declare that they have no conflict of interest.

\end{document}